# ON THE BORDERLINE OF FIELDS AND HYPERFIELDS, PART II – ENUMERATION AND CLASSIFICATION OF THE HYPERFIELDS OF ORDER 7

Christos G. Massouros and Gerasimos G. Massouros

Abstract: The quotient hyperfield is a landmark on the borderline of fields and hyperfields. In this paper, which is the second part of our previously published paper, all the hyperfields of order 7 are constructed, enumerated and presented, in the course of which an important family of 7-element canonical hypergroups is revealed. The study of these hyperfields proved the existence of both quotient and non-quotient ones among them. Their construction became feasible because it is based on a new definition of the hyperfield with fewer axioms, which is introduced in this paper following our proof that the axiom of reversibility can derive from the other axioms of the hyperfield. Hence, the processing power needed for a computer to test whether a structure is a hyperfield or not is much less. This paper also proves properties and contains examples of skew hyperfields, strongly canonical hyperfields/hyperrings and superiorly canonical hyperfields/hyperrings that wrap up and complete its previously published first part's conclusions and results.

## 1. Introduction

This paper is the succession of our previous paper [1] and as such, in its introduction to the topics discussed, it refers to the detailed introduction of [1]. As described in it, hypercompositional Algebra began its existence in Mathematics with the hypergroup which was introduced by F. Marty in 1934 [2] and advanced with the hyperfield which was introduced in 1956 by M. Krasner [3].

The basic concept in hypercompositional Algebra is the hypercomposition. A hypercomposition or a hyperoperation over a non-empty set is a mapping from the cartesian product $E \times E$ into the power set $P(E)$ of $E$. A *hypergroup* is a non-empty set $E$ enriched with a hypercomposition "·" which satisfies the following two axioms:

(i) the axiom of *associativity*:

$$a \cdot (b \cdot c) = (a \cdot b) \cdot c \, , \text{ for all } a,b,c \in E$$

(ii) the axiom of *reproductivity*:

$$a \cdot E = E \cdot a = E \, , \text{ for all } a \in E$$

Papers [4, 5] present in detail that the group is defined with exactly the same axioms as above. Namely, a group is a non-empty set $E$ which is enriched with a composition (i.e. a mapping from the cartesian product $E \times E$ into the set $E$) that satisfies the axioms (i) and (ii).

If "$\cdot$" is an internal composition on a set $E$ and $A$, $B$ are subsets of $E$, then $A \cdot B$ signifies the set $\left\{ a \cdot b \mid (a,b) \in A \times B \right\}$ while if "$\cdot$" is a hypercomposition then $A \cdot B$ is the union $\bigcup\limits_{(a,b) \in A \times B} a \cdot b$. $A \cdot b$ and $a \cdot B$ have the same meaning as $A \cdot \{b\}$ and $\{a\} \cdot B$ respectively. In general, the singleton $\{a\}$ is identified with its member $a$.

**Theorem 1.** *Ref. [1] If either $A = \varnothing$ or $B = \varnothing$, then $AB = \varnothing$ and vice versa.*

**Theorem 2.** *Ref. [5,6] The result of the hypercomposition of any two elements in a hypergroup is always non-void.*

**Definition 1.** *Ref. [3,7]* A *hyperfield* is an algebraic structure $(H, +, \cdot)$ where $H$ is a non-empty set, "$\cdot$" is an internal composition on $H$, and "$+$" is a hypercomposition on $H$ which satisfies the axioms:

> ***I. Multiplicative axiom***
>
> $H = H^* \cup \{0\}$ where $(H^*, \cdot)$ is a multiplicative group and $0$ is a bilaterally absorbing element of $H$, i.e., $0 \cdot a = a \cdot 0 = 0$, for all $a \in H$.

> ***II. Additive axioms***
>
> i. associativity:
> $a + (b + c) = (a + b) + c$, for all $a, b, c \in H$
>
> ii. commutativity:
> $a + b = b + a$, for all $a, b \in H$
>
> iii. for every $a \in H$ there exists one and only one $a' \in H$ such that $0 \in a + a'$. $a'$ is written $-a$ and called the opposite of $a$; moreover, instead of $a + (-b)$ we write $a - b$
>
> iv. reversibility:
> if $a \in b + c$ then $c \in a - b$

> ***III. Distributive axiom***
>
> $a \cdot (b + c) = a \cdot b + a \cdot c$, $(b + c) \cdot a = b \cdot a + c \cdot a$, for all $a, b, c \in H$

If the multiplicative axiom ***I*** is replaced by the axiom:

> ***I′.*** $H = H^* \cup \{0\}$ where $(H^*, \cdot)$ is a multiplicative semigroup and $0$ is a bilaterally absorbing element of $H$,

then a more general structure is obtained which is called *hyperring* [7].

**Theorem 3.** *Ref. [4,8] A non-empty set $H$ enriched with the additive axioms II is a hypergroup.*

**Proof.** The associativity is valid. Next, due to Theorem 2, for every $a$, $b$ in $H$ there exists $x$ in $H$ such that $x \in b - a$ or, because of the reversibility, $b \in a + x$. Therefore, $H \subseteq a + H$. Moreover $a + H \subseteq H$. Consequently $a + H = H$. Hence the reproductivity is valid as well. Thus, $H$ is a hypergroup. □

The special hypergroup of the hyperfield was named *canonical hypergroup* by J. Mittas who studied it in depth and presented his research results through a multitude of papers e.g. [8-10].

**Theorem 4.** *Ref. [11] In a canonical hypergroup $H$ the equality $a + 0 = a$, for all $a$ in $H$, is valid.*

**Proof.** $0 \in a - a$, therefore $a \in a + 0$. Next, if $x \in a + 0$ for all $x \in H$, then $0 \in a - x$, and hence $x = a$.

**Theorem 5.** *In a canonical hypergroup $H$ the axiom of reversibility is equivalent to the validity of the equality*

$$-(a + b) = -a - b$$

*for all a,b in H.*

**Proof.** Suppose that the equality $-(a + b) = -a - b$ is valid for all $a,b$ in $H$. If $a \in b + c$ then equivalently we have that $0 \in a - (b + c)$ or $0 \in a + (-b - c)$ or $0 \in (a - b) - c$. Therefore $c \in a - b$, hence the reversibility holds. Vice versa now. Suppose that the reversibility holds, and $w \in -(a + b)$. Then $0 \in w + (a + b)$. Applying the associativity and the reversibility we have:

$$0 \in w + (a + b) \Rightarrow 0 \in (w + a) + b \Rightarrow -b \in w + a \Rightarrow w \in -a - b.$$

So, $-(a + b) \subseteq -a - b$. Next suppose that $w \in -a - b$. Then we have the sequence of implications:

$$w \in -a - b \Rightarrow -a \in w + b \Rightarrow 0 \in (w + b) + a \Rightarrow 0 \in w + (b + a) \Rightarrow w \in -(a + b).$$

Thus $-a - b \subseteq -(a + b)$. Therefore $-(a + b) = -a - b$. □

In the hyperfields, due to distributivity, the equality $-(a + b) = -a - b$ is always valid. Therefore Theorem 5 reduces the axioms of the hyperfield, since in the additive axioms the reversibility can be deduced from the rest axioms. Hence, we have the Definition:

**Definition 2.** A *hyperfield* is an algebraic structure $(H, +, \cdot)$ where $H$ is a non-empty set, "$\cdot$" is an internal composition on $H$, and "+" is a hypercomposition on $H$ which satisfies the axioms:

 *I. Multiplicative axiom*

  $H = H^* \cup \{0\}$ where $(H^*, \cdot)$ is a multiplicative group and $0$ is a bilaterally absorbing element of $H$, i.e., $0 \cdot a = a \cdot 0 = 0$, for all $a \in H$.

 *II. Additive axioms*

  i. associativity:
   $a + (b + c) = (a + b) + c$, for all $a,b,c \in H$

  ii. commutativity:
   $a + b = b + a$, for all $a,b \in H$

  iii. for every $a \in H$ there exists one and only one $a' \in H$ such that $0 \in a + a'$. $a'$ is written $-a$ and called the opposite of $a$; moreover, instead of $a + (-b)$ we write $a - b$

 *III. Distributive axiom*

  $a \cdot (b + c) = a \cdot b + a \cdot c$,  $(b + c) \cdot a = b \cdot a + c \cdot a$, for all $a,b,c \in H$

The elimination of one redundant axiom (reversibility) makes the verification that a structure is a hyperfield more concise with proofs that do not include extra or unnecessary steps, which in turn facilitates the development of computer packages capable of performing it with reduced computational resources. So, based on the above definition, it was feasible to enumerate the hyperfields of order 7 and present the 277 non-isomorphic hyperfields in the Appendix at the end. The notation used is $HF_7^k$, $1 \leq k \leq 277$ and $HF_7^1$ is the field $\mathbb{Z}_7$.

*Special notation*: In the following pages, in addition to the typical algebraic notations, we are using Krasner's notation for the complement and the difference [12]. So, $A \cdot\cdot B$ denotes the set of elements that are in the set $A$, but not in the set $B$. If $K$ is a field or a hyperfield then, $K^*$ denotes the set $K \cdot\cdot \{0\}$.

## 2. The skew hyperfield

The term division rings, also called skew fields, is used to define nontrivial rings in which every nonzero element has a multiplicative inverse. Historically, division rings were sometimes referred to as fields, while fields were called "commutative fields". Krasner defined the hyperfield considering that its multiplicative group is not necessarily commutative [3,7]. But, in the relevant bibliography the term hyperfield is used, almost from the very first time, to indicate that its multiplicative group is commutative. Consequently, it is necessary to introduce the term "skew hyperfield" in order to refer to hypercompositional structures that are analogous to skew fields. The constructions of hyperfields given in [11,13] enrich a multiplicative group, which need not be commutative, with a hypercomposition. Thus, we have the following three classes of skew hyperfields:

**Theorem 6**. *Let $G$ be a multiplicative group and $0$ a bilaterally absorbing element, i.e. $x \cdot 0 = 0 \cdot x = 0$, for all $x \in G \cup \{0\}$. Then $G \cup \{0\}$ equipped with the hypercompositions*:

$$(1) \quad \left|
\begin{array}{l}
x + y = \{x, y\}, \qquad \text{for all } x, y \in G, \text{ with } x \neq y \\
x + 0 = 0 + x = x, \quad \text{for all } x \in G \cup \{0\} \\
x + x = G \cup \{0\}, \quad \text{for all } x \in G
\end{array}
\right.$$

$$(2) \quad \left|
\begin{array}{l}
x + y = \{x, y\}, \qquad\qquad \text{for all } x, y \in G, \text{ with } x \neq y \\
x + 0 = 0 + x = x, \qquad\quad \text{for all } x \in G \cup \{0\} \\
x + x = \left[ G \cup \{0\} \right] \cdot\cdot \{x\}, \quad \text{for all } x \in G
\end{array}
\right.$$

*and*

$$(3) \quad \left|
\begin{array}{l}
x + y = G \cdot\cdot \{0, x, y\}, \quad \text{for all } x, y \in G, \text{ with } x \neq y \\
x + 0 = 0 + x = x, \qquad \text{for all } x \in G \cup \{0\} \\
x + x = \{0, x\}, \qquad\quad \text{for all } x \in G
\end{array}
\right.$$

*creates three non-isomorphic skew hyperfields. For the sake of convenience, these hyperfields will be denoted by $SHF_1(G)$, $SHF_2(G)$ and $SHF_3(G)$ respectively.*

The above constructions reveal that in the case of hyperfields, the equivalent of Wedderburn's little theorem [14] does not apply i.e. although every finite skew field is commutative, there exist non-commutative finite hyperfields. Indeed, when the constructions of the above theorem are applied to finite non-commutative multiplicative groups, they produce non-commutative finite hyperfields. For instance, the dihedral group $D_3$ is the non-Abelian group having the smallest group order. Therefore, $SHF_1(D_3)$, $SHF_2(D_3)$ and $SHF_3(D_3)$ are three skew hyperfields having the smallest order.

Moreover, it can be noted that both Construction I from [13] and Proposition 1 from [15] do not necessitate that the multiplicative group is abelian. Consequently, Theorem 7 is valid:

**Theorem 7**. *Let G be a non-unitary multiplicative group and let $(H^*, \cdot)$ be its direct product with the multiplicative group $\{-1, 1\}$. Consider the set $H = H^* \cup \{0\}$, where 0 is a bilaterally absorbing element in H, ie. $0w = w0 = 0$, for all $w \in H$. Then H, equipped with the hypercompositions:*

$$(1) \begin{vmatrix} (x,i) \,\hat{+}\, (y,j) = \{(x,i),(y,j)\}, & \text{if } (y,j) \neq (x,-i) \\[6pt] (x,i) \,\hat{+}\, (x,-i) = H, & \text{for all } (x,i) \in H^* \\[6pt] (x,i) \,\hat{+}\, 0 = 0 \,\hat{+}\, (x,i) = (x,i) \ \text{ and } \ 0\,\hat{+}\,0 = 0 & \text{for all } (x,i) \in H^* \end{vmatrix}$$

$$(2) \begin{vmatrix} (x,i) \,\hat{+}\, (w,j) = \{(x,i),(w,j),(x,-i),(w,-j)\}, & \text{if } (w,j) \neq (x,i),\ (x,-i) \\[6pt] (x,i) \,\hat{+}\, (x,i) = H \cdots \{(x,i),\ (x,-i),0\} & \text{for all } (x,i) \in H^* \\[6pt] (x,i) \,\hat{+}\, (x,-i) = H \cdots \{(x,i),\ (x,-i)\} & \text{for all } (x,i) \in H^* \\[6pt] (x,i) \,\hat{+}\, 0 = 0 \,\hat{+}\, (x,i) = (x,i) \ \text{ and } \ 0\,\hat{+}\,0 = 0 & \text{for all } (x,i) \in H^* \end{vmatrix}$$

*creates two non-isomorphic skew hyperfields. These hyperfields will be denoted by $SHF_4(G)$ and $SHF_5(G)$ respectively.*

Furthermore, the Construction II of [13] can be applied to a skew field or a skew hyperfield and therefore the following theorems hold:

**Theorem 8**. *Let $(K,+,\cdot)$ be a skew field. If we define the hypercomposition $\dot{+}$ on F as follows:*

$$x \dot{+} y = \{x, y, x+y\}, \qquad \text{if } y \neq -x,\ x,y \neq 0$$
$$x \dot{+} (-x) = K, \qquad \text{for all } x \in K^*$$
$$x \dot{+} 0 = 0 \dot{+} x = x, \qquad \text{for all } x \in K$$

*then $(K, \dot{+}, \cdot)$ is a skew hyperfield.*

**Theorem 9**. *Let $(H,+,\cdot)$ be a skew hyperfield. If we define a new hypercomposition «$\dot{+}$» on H as follows:*

$$x \dotplus y = \{x, y\} \cup x+y, \qquad \text{for all } x,y \in H^*, \text{ with } y \neq -x$$

$$x \dotplus (-x) = H, \qquad\qquad \text{for all } x \in H^*$$

$$x \dotplus 0 = 0 \dotplus x = x, \qquad\quad \text{for all } x \in H$$

then, $(H, \dotplus, \cdot)$ is a skew hyperfield.

As with the hyperfields, the skew hyperfields of Theorems 8 and 9 are constructed over skew fields and hyperfields using an extensive enlargement of their composition or hypercomposition ([5], example 8) and they will be called *augmented skew hyperfields* in accordance with the terminology established in [1]. In the following, the term *hyperfield* is used to indicate that its multiplicative group is commutative.

## 3. The quotient hyperfield/hyperring and the non-quotient hyperfields/hyperrings

M. Krasner constructed the quotient hyperfield and quotient hyperring, which are based on a field and a ring, respectively [7]. This construction is presented in detail in [1]. However, certain aspects of its content are repeated here for the self-sufficiency of the present paper. If $F$ is a field and $G$ a subgroup of $F$'s multiplicative group $F^*$, then the multiplicative classes modulo $G$ in $F$ form a partition of $F$. Krasner observed that the product of two such classes, considered as subsets of $F$, is also a class modulo $G$, while their sum is a union of such classes and subsequently, he proved that the set $F/G$ of the classes of this partition becomes a hyperfield if the multiplication and the addition are defined as follows:

$$xG \cdot yG = xyG$$

$$xG \dotplus yG = \left\{ (xp + yq)G \ \middle| \ p, q \in G \right\}$$

for all $xG, yG \in F/G$.

Moreover, Krasner showed that if $R$ is a ring and $G$ is a normal subgroup of its multiplicative group, then the above construction gives a hyperring [7]. Ch. Massouros in [15] generalized this construction using not normal multiplicative subgroups, since he proved that in rings there exist multiplicative subgroups $G$ which satisfy the property $xG \cdot yG = xyG$, even when they are not normal.

Crucial to the independence of the theory of hyperfields/hyperrings from the corresponding theory of fields/rings was the existence of hyperfields/hyperrings which do not belong to the class of quotient hyperfields/hyperrings. The existence of such hyperfields/hyperrings was proved by Ch. Massouros in [11,15,16]. The basic theorems that prove the existence of non-quotient hyperfields/hyperrings are the following:

**Theorem 10**. *Ref. [15]  Let $\Theta$ be a multiplicative group which has more than two elements and let $(K^*, \cdot)$ be its direct product with the multiplicative group $\{-1,1\}$. Consider the set $K = K^* \cup \{0\}$, where $0$ is a bilaterally absorbing element in $K$, i.e. $0w = w0 = 0$, for all $w \in K$. The following hypercomposition is introduced on $K$:*

$$(x,i) \dotplus (y,j) = \{(x,i),\ (y,j),\ (x,-i),\ (y,-j)\},\ \ if\ (y,j) \neq (x,i),\ (x,-i)$$

$$(x,i) \dotplus (x,i) = K \smallsetminus \{(x,i),\ (x,-i),\ 0\}$$

$$(x,i) \dotplus (x,-i) = K \smallsetminus \{(x,i),\ (x,-i)\}$$

$$(x,i) \dotplus 0 = 0 \dotplus (x,i) = (x,i)\ \ and\ \ 0 \dotplus 0 = 0$$

*Then, the triplet $K(\Theta) = (K, +, \cdot)$ is a hyperfield which does not belong to the class of quotient hyperfields when $\Theta$ is a periodic group.*

**Corollary 1.** *The hyperfield $HF_7^{81}$ of the Appendix is not a quotient hyperfield.*

**Table 1.** The canonical hypergroup of the non-quotient hyperfield $HF_7^{81}$

| $HF_7^{81}$ | 0 | 1 | a | b | c | d | e |
|---|---|---|---|---|---|---|---|
| 0 | 0 | 1 | a | b | c | d | e |
| 1 | 1 | {a,b,d,e} | {1,a,c,d} | {1,b,c,e} | {0,a,b,d,e} | {1,a,c,d} | {1,b,c,e} |
| a | a | {1,a,c,d} | {1,b,c,e} | {a,b,d,e} | {1,a,c,d} | {0,1,b,c,e} | {a,b,d,e} |
| b | b | {1,b,c,e} | {a,b,d,e} | {1,a,c,d} | {1,b,c,e} | {a,b,d,e} | {0,1,a,c,d} |
| c | c | {0,a,b,d,e} | {1,a,c,d} | {1,b,c,e} | {a,b,d,e} | {1,a,c,d} | {1,b,c,e} |
| d | d | {1,a,c,d} | {0,1,b,c,e} | {a,b,d,e} | {1,a,c,d} | {1,b,c,e} | {a,b,d,e} |
| e | e | {1,b,c,e} | {a,b,d,e} | {0,1,a,c,d} | {1,b,c,e} | {a,b,d,e} | {1,a,c,d} |

**Theorem 11.** *Ref. [11, 16] Let $\Theta$ be a multiplicative group which has more than two elements and let 0 be a multiplicatively bilaterally absorbing element. If a hypercomposition $+$ is defined on $H = \Theta \cup \{0\}$ as follows:*

$$x + y = \{x, y\}, \qquad \text{for all } x, y \in \Theta, \text{ with } y \neq x$$

$$x + x = H \cdots \{x\}, \qquad \text{for all } x \in \Theta$$

$$x + 0 = 0 + x = x, \qquad \text{for all } x \in H$$

*then, the triplet $H(\Theta) = (\Theta \cup \{0\}, +, \cdot)$ is a hyperfield which is not isomorphic to a quotient hyperfield when $\Theta$ is a periodic group.*

**Corollary 2.** *The hyperfield $HF_7^{258}$ of the Appendix is not a quotient hyperfield.*

**Table 2.** The canonical hypergroup of the non-quotient hyperfield $HF_7^{258}$

| $HF_7^{258}$ | 0 | 1 | a | b | c | d | e |
|---|---|---|---|---|---|---|---|
| 0 | 0 | 1 | a | b | c | d | e |
| 1 | 1 | {0,a,b,c,d,e} | {1,a} | {1,b} | {1,c} | {1,d} | {1,e} |
| a | a | {1,a} | {0,1,b,c,d,e} | {a,b} | {a,c} | {a,d} | {a,e} |
| b | b | {1,b} | {a,b} | {0,1,a,c,d,e} | {b,c} | {b,d} | {b,e} |
| c | c | {1,c} | {a,c} | {b,c} | {0,1,a,b,d,e} | {c,d} | {c,e} |
| d | d | {1,d} | {a,d} | {b,d} | {c,d} | {0,1,a,b,c,e} | {d,e} |
| e | e | {1,e} | {a,e} | {b,e} | {c,e} | {d,e} | {0,1,a,b,c,d} |

**Theorem 12**. *Ref. [15] The direct sum S of the hyperrings $S_i$, $i \in I$, is not isomorphic to a sub-hyperring of a quotient hyperring if at least one of the $S_i$ is a non-quotient hyperfield.*

The existence of some more classes of non-quotient hyperfields/hyperrings was subsequently made by Nakassis [17].

**Theorem 13.** *Ref. [17] Let $(T, \cdot)$ be a multiplicative group of order m, with m > 3. Also let $H = T \cup \{0\}$, where 0 is a multiplicatively absorbing element. If H is equipped with the hypercomposition:*

$$x + y = H \cdots \{0, x, y\} \quad \text{for all } x, y \in T, \text{ with } y \neq x$$

$$x + x = \{0, x\}, \quad \quad \text{for all } x \in T$$

$$x + 0 = 0 + x = x, \quad \text{for all } x \in H$$

*then, $H(T) = (T \cup \{0\}, +, \cdot)$ is a hyperfield which is a non-quotient hyperfield if the cardinality of T is properly chosen.*

The next two Propositions are crucial for proving that hyperfields of the type outlined in Theorem 13 are not quotient hyperfields.

**Proposition 1.** *In a quotient hyperfield $F/Q$ the cardinality of the sum of any two elements cannot exceed the cardinality of Q.*

**Proof.** Suppose that $xQ, yQ$ are two arbitrary elements of $F/Q$. Then,

$$xQ + yQ = \left\{ (x + yq)Q \mid q \in Q \right\}$$

and the function $f : Q \to xQ + yQ$ with $f(q) = (x + yq)Q$ is a surjection. □

**Proposition 2.** *If in a quotient hyperfield $F/Q$ the differences $xQ - xQ$, $xQ \in F/Q$ have only 0 in common, then the cardinality of the sum of any two non-opposite elements as well as that of any two non-equal elements is equal to the cardinality of Q.*

**Proof.** Let $xQ, yQ$ be two non-opposite and non-equal elements in $F/Q$. Then

$$xQ + yQ = \left\{ (x + yq)Q \mid q \in Q \right\} .$$

Next, if $(x + yq)Q = (x + yp)Q$, then

$$x + yq = (x + yp)r \Leftrightarrow x - xr = yq - ypr \Rightarrow (xQ - xQ) \cap (yQ - yQ) \neq \varnothing$$

But, from the validity of the equality $(xQ - xQ) \cap (yQ - yQ) = \{0\}$, it follows that $x - xr = 0$. Therefore $r = 1$ and consequently $yq - yp = 0$ or equivalently $q = p$. Hence

$$card(xQ + yQ) = cardQ$$

and so the Proposition. □

We will subsequently present some non-quotient hyperfields from the list of 277 seven-element hyperfields detailed in the Appendix, utilizing these Propositions.

**Proposition 3.**

*i.    The hyperfield $HF_7^3$ of the Appendix is not a quotient hyperfield.*

**Table 3.** The canonical hypergroup of the non-quotient hyperfield $HF_7^3$

| $HF_7^3$ | 0 | 1 | a | b | c | d | e |
|----------|---|---|---|---|---|---|---|
| 0 | 0 | 1 | a | b | c | d | e |
| 1 | 1 | {1,a,b,d,e} | {b,c,d,e} | {a,c,d,e} | {0,1,c} | {a,b,c,e} | {a,b,c,d} |
| a | a | {b,c,d,e} | {1,a,b,c,e} | {1,c,d,e} | {1,b,d,e} | {0,a,d} | {1,b,c,d} |
| b | b | {a,c,d,e} | {1,c,d,e} | {1,a,b,c,d} | {1,a,d,e} | {1,a,c,e} | {0,b,e} |
| c | c | {0,1,c} | {1,b,d,e} | {1,a,d,e} | {a,b,c,d,e} | {1,a,b,e} | {1,a,b,d} |
| d | d | {a,b,c,e} | {0,a,d} | {1,a,c,e} | {1,a,b,e} | {1,b,c,d,e} | {1,a,b,c} |
| e | e | {a,b,c,d} | {1,b,c,d} | {0,b,e} | {1,a,b,d} | {1,a,b,c} | {1,a,c,d,e} |

*ii.    The hyperfield $HF_7^4$ of the Appendix is not a quotient hyperfield.*

**Table 4.** The canonical hypergroup of the non-quotient hyperfield $HF_7^4$

| $HF_7^4$ | 0 | 1 | a | b | c | d | e |
|----------|---|---|---|---|---|---|---|
| 0 | 0 | 1 | a | b | c | d | e |
| 1 | 1 | {1,a,b,c,d,e} | {b,c,d,e} | {a,c,d,e} | {0,1,c} | {a,b,c,e} | {a,b,c,d} |
| a | a | {b,c,d,e} | {1,a,b,c,d,e} | {1,c,d,e} | {1,b,d,e} | {0,a,d} | {1,b,c,d} |
| b | b | {a,c,d,e} | {1,c,d,e} | {1,a,b,c,d,e} | {1,a,d,e} | {1,a,c,e} | {0,b,e} |
| c | c | {0,1,c} | {1,b,d,e} | {1,a,d,e} | {1,a,b,c,d,e} | {1,a,b,e} | {1,a,b,d} |
| d | d | {a,b,c,e} | {0,a,d} | {1,a,c,e} | {1,a,b,e} | {1,a,b,c,d,e} | {1,a,b,c} |
| e | e | {a,b,c,d} | {1,b,c,d} | {0,b,e} | {1,a,b,d} | {1,a,b,c} | {1,a,b,c,d,e} |

**Proof.** (*i*) The opposite of 1 is *c*, the opposite of *a* is *d*, the opposite of *b* is *e* and the differences 1-*c*, *a*-*d* and *b*-*e* have only 0 in common. Therefore, if $HF_7^3$ were isomorphic to a quotient hyperfield *F/Q*, then, according to Proposition 2, the cardinality of *Q* would be equal to the cardinality of the sum of any two non-opposite and non-equal elements, which is 4. Moreover, according to Proposition 1, the cardinality of the sum of any two elements cannot exceed the cardinality of *Q*. However, the cardinality of the sum *x*+*x* is equal to 5 for any $x \in HF_7^3$. Hence $HF_7^3$ is not a quotient hyperfield. Similar is the proof of (***ii***). □

**Proposition 4.**

*i.   The hyperfield $HF_7^5$ of the Appendix is not a quotient hyperfield.*

**Table 5.** The canonical hypergroup of the non-quotient hyperfield $HF_7^5$

| $HF_7^5$ | 0 | 1 | a | b | c | d | e |
|---|---|---|---|---|---|---|---|
| 0 | 0 | 1 | a | b | c | d | e |
| 1 | 1 | {a,b,c,e} | {1,b,d,e} | {a,b,c,d,e} | {0,a,d} | {1,a,b,c,e} | {a,c,d,e} |
| a | a | {1,b,d,e} | {1,b,c,d} | {1,a,c,e} | {1,b,c,d,e} | {0,b,e} | {1,a,b,c,d} |
| b | b | {a,b,c,d,e} | {1,a,c,e} | {a,c,d,e} | {1,a,b,d} | {1,a,c,d,e} | {0,1,c} |
| c | c | {0,a,d} | {1,b,c,d,e} | {1,a,b,d} | {1,b,d,e} | {a,b,c,e} | {1,a,b,d,e} |
| d | d | {1,a,b,c,e} | {0,b,e} | {1,a,c,d,e} | {a,b,c,e} | {1,a,c,e} | {1,b,c,d} |
| e | e | {a,c,d,e} | {1,a,b,c,d} | {0,1,c} | {1,a,b,d,e} | {1,b,c,d} | {1,a,b,d} |

*ii.   The hyperfield $HF_7^6$ of the Appendix is not a quotient hyperfield.*

**Table 6.** The canonical hypergroup of the non-quotient hyperfield $HF_7^6$

| $HF_7^6$ | 0 | 1 | a | b | c | d | e |
|---|---|---|---|---|---|---|---|
| 0 | 0 | 1 | a | b | c | d | e |
| 1 | 1 | {a,b,c,d,e} | {1,b,c,d,e} | {b,c,d,e} | {0,a,d} | {1,a,b,c} | {a,b,c,d,e} |
| a | a | {1,b,c,d,e} | {1,b,c,d,e} | {1,a,c,d,e} | {1,c,d,e} | {0,b,e} | {a,b,c,d} |
| b | b | {b,c,d,e} | {1,a,c,d,e} | {1,a,c,d,e} | {1,a,b,d,e} | {1,a,d,e} | {0,1,c} |
| c | c | {0,a,d} | {1,c,d,e} | {1,a,b,d,e} | {1,a,b,d,e} | {1,a,b,c,e} | {1,a,b,e} |
| d | d | {1,a,b,c} | {0,b,e} | {1,a,d,e} | {1,a,b,c,e} | {1,a,b,c,e} | {1,a,b,c,d} |
| e | e | {a,b,c,d,e} | {a,b,c,d} | {0,1,c} | {1,a,b,e} | {1,a,b,c,d} | {1,a,b,c,d} |

*iii.   The hyperfield $HF_7^{11}$ of the Appendix is not a quotient hyperfield.*

**Table 7.** The canonical hypergroup of the non-quotient hyperfield $HF_7^{11}$

| $HF_7^{11}$ | 0 | 1 | a | b | c | d | e |
|---|---|---|---|---|---|---|---|
| 0 | 0 | 1 | a | b | c | d | e |
| 1 | 1 | {a,b,c,d} | {a,b,c,d,e} | {1,a,d,e} | {0,b,e} | {b,c,d,e} | {1,a,b,c,d} |
| a | a | {a,b,c,d,e} | {b,c,d,e} | {1,b,c,d,e} | {1,a,b,e} | {0,1,c} | {1,c,d,e} |
| b | b | {1,a,d,e} | {1,b,c,d,e} | {1,c,d,e} | {1,a,c,d,e} | {1,a,b,c} | {0,a,d} |
| c | c | {0,b,e} | {1,a,b,e} | {1,a,c,d,e} | {1,a,d,e} | {1,a,b,d,e} | {a,b,c,d} |
| d | d | {b,c,d,e} | {0,1,c} | {1,a,b,c} | {1,a,b,d,e} | {1,a,b,e} | {1,a,b,c,e} |
| e | e | {1,a,b,c,d} | {1,c,d,e} | {0,a,d} | {a,b,c,d} | {1,a,b,c,e} | {1,a,b,c} |

**Proof.** (***i***) The opposite of 1 is $c$, the opposite of $a$ is $d$, the opposite of $b$ is $e$ and the differences $1\text{-}c$, $a\text{-}d$ and $b\text{-}e$ have only 0 in common. Therefore, if $HF_7^5$ were isomorphic to a quotient hyperfield $F/Q$, then, according to Proposition 2, the cardinality of $Q$ would be equal to the cardinality of the sum of any two non-opposite and non-equal elements. Therefore, the cardinalities of the sum of any two non-opposite and non-equal elements would have been equal. But, card(1+$a$) = 4, while card($b$+$c$) = 5. Consequently $HF_7^5$ is not a quotient hyperfield. Similar is the proof of (***ii***) and (***iii***). □

**Remark 1.** In the light of Propositions 1 and 2 it is preferable to keep the first proof that $NQHF_5^2$ in [1] is not a quotient hyperfield when its multiplicative group is the Vierergruppe and it is classified according to Theorem 14 of [1], when its multiplicative group is cyclic.

### 4. Strongly canonical and superiorly canonical hyperfields/hyperrings

As in the case of fields, a valuation theory was developed for the hyperfields as well. The concept of the valuation of hyperfields was introduced by M. Krasner [3]. [18-21] contain a recent study on valuations in hyperfields, while in [22] A. Linzi gives an extremely detailed and in-depth presentation of the valuation theory of hyperfields. Among other things in his paper, A. Linzi clarifies several points of J. Mittas' earlier research on valued hyperfields [23-27]. J. Mittas proved [23-27] that a necessary and sufficient condition for a canonical hypergroup, and consequently for a hyperfield, to be valuated or hypervaluated is the validity of certain additional properties of a purely algebraic type, i.e., properties that can be expressed without the intervention of the valuation or the hypervaluation, respectively. This led him to the definition of the following two special canonical hypergroups:

(a) The *strongly canonical hypergroup* which is a canonical hypergroup that also satisfies the axioms:

$S_1$: If $x \in x+y$, then $x+y=x$, for all $x,y \in H$

$S_2$: If $(x+y) \cap (z+w) \neq \varnothing$, then either $x+y \subseteq z+w$ or $z+w \subseteq x+y$

(b) The *superiorly canonical hypergroup* which is a strongly canonical hypergroup that also satisfies the axioms:

$S_3$: If $z,w \in x-y$ and $x \neq y$, then $z-z=w-w$

$S_4$: If $x \in z-z$ and $y \notin z-z$ then $x-x \subseteq y-y$

Thus, depending on the type of their additive hypergroup, strongly canonical and superiorly canonical hyperfields are defined respectively. Mittas proved the following theorem (e.g. [29], Theorem 4.2]:

**Theorem 14.** *A hyperfield can be valuated if and only if its additive hypergroup is superiorly canonical.*

A detailed and thorough proof of this theorem is given by A. Linzi in [22] (Theorem 4.22).

**Proposition 5.** *In a strongly canonical hyperfield if* $x \neq y$ *, then*:

$$(x-x) \cap (y-x) = \varnothing \quad and \quad (y-y) \cap (y-x) = \varnothing.$$

**Proof.** Suppose that $x \neq y$ and $(x-x) \cap (y-x) \neq \varnothing$. Let $w \in (x-x) \cap (y-x)$. Then, $w \in y-x$ implies that $y \in x+w$. Moreover $w \in x-x$ implies that $x \in x+w$. Therefore $x = x+w$. Thus $x = y$, which contradicts our assumption. □

**Proposition 6.** *In a strongly canonical hyperfield K it holds that*:

*i.* $x + (x-x) = x$, *for all* $x \in K$

*ii.* $x \notin x-x$ *and equivalently* $x \notin x+x$, *for all* $x \in K \cdot\cdot \{0\}$

**Proof.** (i) Suppose that $w \in x-x$. Then $x \in x+w$ and therefore $x = x+w$. So

$$x + (x-x) = \bigcup_{w \in x-x} x+w = x.$$

(ii) Suppose that $x$ is a non-zero element such that $x \in x-x$. Then $0 \in x+(x-x)$. But, according to (i), $x = x+(x-x)$. Consequently $x = 0$, which is absurd. □

**Lemma 1.** *Ref. [28] In a canonical hypergroup H it holds that*:

$$(x+y) \cap (z+w) \neq \varnothing \Rightarrow (x-z) \cap (w-y) \neq \varnothing$$

*for all* $x, y, z, w \in H$.

**Theorem 15.** *In a strongly canonical hyperfield K, the sets* $a+x$, $x \in K$ *form a partition of K, for any fixed point* $a \in K$.

**Proof.** Let $x \neq y$ and suppose that $(a+x) \cap (a+y) \neq \varnothing$. Then, because of lemma 1, $(x-y) \cap (a-a) \neq \varnothing$. If $w \in (x-y) \cap (a-a)$, then:

$$w \in x-y \Rightarrow x \in w+y \quad and \quad w \in a-a \Rightarrow a \in w+a \Rightarrow a = w+a.$$

Hence $\qquad\qquad a+x \subseteq a+(w+y) = (a+w)+y = a+y.$

Moreover, $\qquad\quad w \in x-y \Rightarrow y \in x-w \quad and \quad w \in a-a \Rightarrow a \in a-w \Rightarrow a = a-w.$

So $\qquad\qquad\quad a+y \subseteq a+(x-w) = (a-w)+x = a+x.$

Therefore, $a+x = a+y$.

**Example 1.** Let H be a set which is totally ordered and symmetric around a center, denoted by $0 \in H$. Then H, equipped with the hypercomposition:

$$x+y = \begin{cases} y, & if \ |x| < |y| \\ x+x = x \\ x-x = \left[-|x|, \ |x|\right] \end{cases}$$

is a canonical hypergroup. Now, if $H \cdot\cdot \{0\}$ is a multiplicative group and 0 is bilaterally absorbing with regards to the multiplication, then H becomes a hyperfield. Moreover, if the hypercomposition on H is defined as follows:

$$x + y = \begin{cases} y, & \text{if } |x| < |y| \\ x + x = -x \\ x - x = \left(-|x|, \ |x|\right) \end{cases}$$

then H is a strongly canonical hyperfield.

**Example 2**. Let $(E, \cdot)$ be a totally ordered multiplicative semigroup, having a minimum element 0, which is bilaterally absorbing with regards to the multiplication. The following hypercomposition is defined on E:

$$x \,\widehat{+}\, y = \begin{cases} \max\{x, y\} & \text{if } x \neq y \\ \{z \in E \,|\, z \leq x\} & \text{if } x = y \end{cases}$$

Then $(E, \widehat{+}, \cdot)$ is a hyperring. If $E \cdot \cdot \{0\}$ is a multiplicative group, then $(E, \widehat{+}, \cdot)$ is a skew hyperfield while if it is an abelian group, $(E, \widehat{+}, \cdot)$ is a hyperfield. This hyperfield was introduced by J. Mittas (see [10] page 86 and [30] page 370) and nowadays it is called tropical hyperfield (see e.g. [31-35]). Moreover, if the hypercomposition on E is defined as follows:

$$x \,\widecheck{+}\, y = \begin{cases} \max\{x, y\} & \text{if } x \neq y \\ \{z \in E \,|\, z < x\} & \text{if } x = y \end{cases}$$

then $(E, \widecheck{+}, \cdot)$ is a strongly canonical hyperring. If $E \cdot \cdot \{0\}$ is a multiplicative group, then $(E, \widecheck{+}, \cdot)$ is a strongly canonical skew hyperfield while if it is an abelian group, $(E, \widecheck{+}, \cdot)$ is a strongly canonical hyperfield.

## Corrigendum on [1].

Theorem 12 in [1] states that the above hyperring $(E, \widecheck{+}, \cdot)$ is a non-quotient hyperring. However, this is not true because in the proof there is a slight confusion. Specifically, while it is initially shown that the assumption that $E$ is isomorphic to a quotient hyperfield $R/Q$ implies that « $4Q$ *is a new class different from $Q$ and* $2Q$ » just six lines below it is mistakenly considered that $4Q=Q$, which led to the incorrect reductio ad absurdum, that $7 \in Q$ and $7 \notin Q$. This was observed by David Hobby and Jaiung Jun, who subsequently examined whether the family of these hyperfields contains quotient hyperfields and constructed an interesting example that answers this question in the affirmative [36].

## 5. Construction of the hyperfields of order 7

This section presents the algorithmic method which is used for the construction of the hyperfields that have 7 elements. The detailed list of these hyperfields is given in the Appendix at the end. There are 277 hyperfields of order 7, and they are denoted with the notation $HF_7^k$, $1 \leq k \leq 277$. $HF_7^1$ is the field $\mathbb{Z}_7$. The multiplicative group $HF_7^*$ of these hyperfields consists of six elements. As it is known there are two groups with six elements, the cyclic group $\mathbb{Z}_6$ and the dihedral group $D_3$. Since $\mathbb{Z}_6$ is cyclic, it is abelian, while $D_3$ is not abelian. Consequently, the multiplicative group of hyperfields of order 7 is the cyclic group $\mathbb{Z}_6$. Of course, $D_3$ can be used to form skew hyperfields, as it is indicated in Section 2.

Since the multiplicative group $HF_7^*$ of hyperfields of order 7 is the cyclic group $\mathbb{Z}_6$ we have that $HF_7^* = \left\{1, a, a^2, a^3, a^4, a^5\right\}$, where $a$ is a generator of $HF_7^*$. Note that the generators of the multiplicative group of $\mathbb{Z}_7$ are $a = 3$ and $a = 5$.

**Theorem 16.** *For the hyperfields of order 7 the following apply:*

$$\textit{i.} \quad 0 \notin a+1, \quad \textit{ii.} \quad 0 \notin a^2+1, \quad \textit{iii.} \quad 0 \notin a^4+1, \quad \textit{iv.} \quad 0 \notin a^5+1$$

*where $a$ is the generator of the multiplicative group.*

**Proof.** i. Suppose that $0 \in a+1$, then $a = -1$ and therefore $a^2 = 1$. Hence $a$ is not a generator of $HF_7^*$, which is a contradiction.

ii. If $0 \in a^2+1$, then $a^2 = -1$. Hence $a^4 = 1$ and therefore $a$ is not a generator of $HF_7^*$, which contradicts the assumption for $a$.

iii. If $0 \in a^4+1$, then $a^4 = -1$. Also, $a^8 = a^6 a^2 = a^2$. So, $a^2 = 1$ and therefore $a$ is not a generator of $HF_7^*$, which is absurd.

iv. If $0 \in a^5+1$, then $a^5 = -1$. From $a^{10} = a^6 a^4 = a^4$ it derives that $a^4 = 1$ and therefore $a$ is not a generator of $HF_7^*$, absurd. □

**Corollary 1.** *For the hyperfields of order 7, it holds that either $0 \in 1+1$ or $0 \in a^3+1$, where $a$ is the generator of the multiplicative group.*

**Corollary 2.** *For the hyperfields of order 7 it holds that:*

*i.* $\quad 1+a \subseteq \left\{1, a, a^2, a^3, a^4, a^5\right\}$

*ii.* $\quad 1+a^2 \subseteq \left\{1, a, a^2, a^3, a^4, a^5\right\}$

*where $a$ is the generator of the multiplicative group.*

**Proposition 7.** *The addition in $HF_7^k$ is completely determined by the sums:*

$$1+1, \quad 1+a, \quad 1+a^2, \quad 1+a^3$$

*where $a$ is the generator of the multiplicative group.*

**Proof.** Distributivity implies the following:

| | | | |
|---|---|---|---|
| i. | $1+a^4 = a^4(a^2+1)$ | vii. $a+a^5 = a(a^4+1)$ | xiii. $a^3+a^4 = a^3(a+1)$ |
| ii. | $1+a^5 = a^5(a+1)$ | viii. $a^2+a^2 = a^2(1+1)$ | xiv. $a^3+a^5 = a^3(a^2+1)$ |
| iii. | $a+a = a(1+1)$ | ix. $a^2+a^3 = a^2(a+1)$ | xv. $a^4+a^4 = a^4(1+1)$ |
| iv. | $a+a^2 = a(a+1)$ | x. $a^2+a^4 = a^2(a^2+1)$ | xvi. $a^4+a^5 = a^4(a+1)$ |
| v. | $a+a^3 = a(a^2+1)$ | xi. $a^2+a^5 = a^2(a^3+1)$ | xvii. $a^5+a^5 = a^5(1+1)$ |
| vi. | $a+a^4 = a(a^3+1)$ | xii. $a^3+a^3 = a^3(1+1)$ | |

Using the notation $b, c, d, e$ for the elements $a^2, a^3, a^4, a^5$ respectively, the following Cayley table summarizes the above results:

**Table 8.** The addition in $HF_7^k$

|   | 0 | 1 | $a$ | $b$ | $c$ | $d$ | $e$ |
|---|---|---|-----|-----|-----|-----|-----|
| 0 | 0 | 1 | $a$ | $b$ | $c$ | $d$ | $e$ |
| 1 | 1 | $1+1 \subseteq$ $\{0\}\cup\{1,a,b,c,d,e\}$ | $1+a \subseteq$ $\{1,a,b,c,d,e\}$ | $1+b \subseteq$ $\{1,a,b,c,d,e\}$ | $1+c \subseteq$ $\{0\}\cup\{1,a,b,c,d,e\}$ | $d(b+1)$ | $e(a+1)$ |
| $a$ | $a$ | $1+a \subseteq$ $\{1,a,b,c,d,e\}$ | $a(1+1)$ | $a(a+1)$ | $a(b+1)$ | $a(c+1)$ | $a(d+1)$ |
| $b$ | $b$ | $1+b \subseteq$ $\{1,a,b,c,d,e\}$ | $a(a+1)$ | $b(1+1)$ | $b(a+1)$ | $b(b+1)$ | $b(c+1)$ |
| $c$ | $c$ | $1+c \subseteq$ $\{0\}\cup\{1,a,b,c,d,e\}$ | $a(b+1)$ | $b(a+1)$ | $c(1+1)$ | $c(a+1)$ | $c(b+1)$ |
| $d$ | $d$ | $d(b+1)$ | $a(c+1)$ | $b(b+1)$ | $c(a+1)$ | $d(1+1)$ | $d(a+1)$ |
| $e$ | $e$ | $e(a+1)$ | $a(d+1)$ | $b(c+1)$ | $c(b+1)$ | $d(a+1)$ | $e(1+1)$ |

# 6. Classification of the hyperfields of order 7

Some of the 277 seven-element hyperfields which are listed in the Appendix are classified in this section. Krasner's question regarding the existence of non-quotient hyperfields [7] led directly to the question and subsequently to a study of the conditions under which a field $F$ can be expressed as the difference $G - G$, where $G$ is a multiplicative subgroup of the field [15, 16, 37, 38]. Thus, if a hyperfield of order 7 is isomorphic to a quotient hyperfield F/G, then G is a multiplicative subgroup of $F$ having index 6 and, in this case, the following Theorem ([1], Theorem 15) has been proven to hold:

**Theorem 17.** *Ref. [1]* *If F is a finite field and G is a subgroup of its multiplicative group of index 6 and order $m$, then $G \cdot G = F$, if and only if:*

$-1 \notin G$ *and* $m \geq 11$,

$-1 \in G$, *charF* $= 11$ *and* $m \geq 20$,

$-1 \in G$, *charF* $= 13$ *and* $m \geq 28$,

$-1 \in G$, *charF* $\neq 11,13$ *and* $m \geq 30$.

Moreover, it is known that if $G$ is a subgroup of finite index in the multiplicative group of an infinite field $F$, then the equality $G - G = F$ is valid [1, 39, 40]. Consequently, the quotient hyperfields $F/G$ for which $G - G \neq F/G$ holds are fully and explicitly determined by the aforementioned Theorem 17. Therefore, the following Theorem is valid:

**Theorem 18.** *The quotient hyperfields of order 7 that satisfy the condition $G - G \neq F/G$ are as follows:*

**i.**   $Z_7$, $Z_{19}/G$, $Z_{31}/G$, $Z_{43}/G$, $GF\left[5^2\right]/G$ and $GF\left[7^2\right]/G$, when the multiplicative subgroup $G$ does not contain $-1$.

**ii.**   $Z_{13}/G$, $Z_{37}/G$, $Z_{61}/G$, $Z_{73}/G$, $Z_{97}/G$ and $Z_{109}/G$, when the multiplicative subgroup $G$ contains $-1$.

**Remark 1.**  As shown below $Z_{157}/G$ is isomorphic to $Z_{97}/G$ and for this reason it is not included in case (ii) of the above Theorem.

The following subsections focus on the examination of the isomorphisms of the hyperfields presented in Theorem 18 in relation to the seven-element hyperfields listed in the Appendix.  The objective is to establish their correspondence and to identify which hyperfields from this list are not categorized as quotient hyperfields as well as which ones are classified as quotient hyperfields.

### 6.1. The quotient hyperfields F/G of order 7 that are constructed from the prime fields and for which $G - G \neq F/G$ and $-1 \notin G$ hold.

The field $Z_7$, which may be regarded as a trivial hyperfield, could be considered as the initial member of this category. The remaining members, according to Theorem 18, are $Z_{19}/G$, $Z_{31}/G$ and $Z_{43}/G$.

**Proposition 8.**

**i.**   The multiplicative subgroup of index 6 of the field $Z_{19}$ is

$$G = \left\{1, 7, 11\right\}$$

and the hyperfield $Z_{19}/G$ is isomorphic to $HF_7^9$.

**Table 9.** The canonical hypergroup of the hyperfield $Z_{19}/G$

| $HF_7^9$ | 0 | 1 | a | b | c | d | e |
|---|---|---|---|---|---|---|---|
| 0 | 0 | 1 | a | b | c | d | e |
| 1 | 1 | {a,c} | {a,b,e} | {1,d,e} | {0,b,e} | {b,c,d} | {1,a,d} |
| a | a | {a,b,e} | {b,d} | {1,b,c} | {1,a,e} | {0,1,c} | {c,d,e} |
| b | b | {1,d,e} | {1,b,c} | {c,e} | {a,c,d} | {1,a,b} | {0,a,d} |
| c | c | {0,b,e} | {1,a,e} | {a,c,d} | {1,d} | {b,d,e} | {a,b,c} |
| d | d | {b,c,d} | {0,1,c} | {1,a,b} | {b,d,e} | {a,e} | {1,c,e} |
| e | e | {1,a,d} | {c,d,e} | {0,a,d} | {a,b,c} | {1,c,e} | {1,b} |

**ii.**   The multiplicative subgroup of index 6 of the field $Z_{31}$ is

$$G = \left\{1, 2, 4, 8, 16\right\}$$

*and the hyperfield* $Z_{31}/G$ *is isomorphic to* $HF_7^{13}$ .

**Table 10.** The canonical hypergroup of the hyperfield $Z_{31}/G$

| $HF_7^{13}$ | 0 | 1 | a | b | c | d | e |
|---|---|---|---|---|---|---|---|
| 0 | 0 | 1 | a | b | c | d | e |
| 1 | 1 | {1,a,b} | {1,b,d,e} | {a,b,d,e} | {0,1,a,c,d} | {1,b,c,e} | {a,c,d,e} |
| a | a | {1,b,d,e} | {a,b,c} | {1,a,c,e} | {1,b,c,e} | {0,a,b,d,e} | {1,a,c,d} |
| b | b | {a,b,d,e} | {1,a,c,e} | {b,c,d} | {1,a,b,d} | {1,a,c,d} | {0,1,b,c,e} |
| c | c | {0,1,a,c,d} | {1,b,c,e} | {1,a,b,d} | {c,d,e} | {a,b,c,e} | {a,b,d,e} |
| d | d | {1,b,c,e} | {0,a,b,d,e} | {1,a,c,d} | {a,b,c,e} | {1,d,e} | {1,b,c,d} |
| e | e | {a,c,d,e} | {1,a,c,d} | {0,1,b,c,e} | {a,b,d,e} | {1,b,c,d} | {1,a,e} |

***iii.*** *The multiplicative subgroup of index 6 of the field* $Z_{43}$ *is*

$$G = \{1, 4, 11, 16, 21, 35, 41\}$$

*and the hyperfield* $Z_{43}/G$ *is isomorphic to* $HF_7^{61}$ .

**Table 11.** The canonical hypergroup of the hyperfield $Z_{43}/G$

| $HF_7^{61}$ | 0 | 1 | a | b | c | d | e |
|---|---|---|---|---|---|---|---|
| 0 | 0 | 1 | a | b | c | d | e |
| 1 | 1 | {a,b,c} | {1,a,b,d,e} | {1,a,b,d,e} | {0,a,b,d,e} | {1,b,c,d,e} | {1,a,c,d,e} |
| a | a | {1,a,b,d,e} | {b,c,d} | {1,a,b,c,e} | {1,a,b,c,e} | {0,1,b,c,e} | {1,a,c,d,e} |
| b | b | {1,a,b,d,e} | {1,a,b,c,e} | {c,d,e} | {1,a,b,c,d} | {1,a,b,c,d} | {0,1,a,c,d} |
| c | c | {0,a,b,d,e} | {1,a,b,c,e} | {1,a,b,c,d} | {1,d,e} | {a,b,c,d,e} | {a,b,c,d,e} |
| d | d | {1,b,c,d,e} | {0,1,b,c,e} | {1,a,b,c,d} | {a,b,c,d,e} | {1,a,e} | {1,b,c,d,e} |
| e | e | {1,a,c,d,e} | {1,a,c,d,e} | {0,1,a,c,d} | {a,b,c,d,e} | {1,b,c,d,e} | {1,a,b} |

**6.2.** *The quotient hyperfields F/G of order 7 that are constructed from the prime fields and for which* $G - G \neq F/G$ *and* $-1 \in G$ *hold*.

The members of this category, according to Theorem 18, are $Z_{13}/G$, $Z_{37}/G$, $Z_{61}/G$, $Z_{73}/G$, $Z_{97}/G$ and $Z_{109}/G$.

**Proposition 9.**

***i.*** *The multiplicative subgroup of index 6 of the field* $Z_{13}$ *is*

$$G = \{1,12\}$$

*and the hyperfield $Z_{13}/G$ is isomorphic to $HF_7^{143}$.*

**Table 12.** The canonical hypergroup of the hyperfield $Z_{13}/G$

| $HF_7^{143}$ | 0 | 1 | a | b | c | d | e |
|---|---|---|---|---|---|---|---|
| 0 | 0 | 1 | a | b | c | d | e |
| 1 | 1 | {0,a} | {1,d} | {c,d} | {b,e} | {a,b} | {c,e} |
| a | a | {1,d} | {0,b} | {a,e} | {d,e} | {1,c} | {b,c} |
| b | b | {c,d} | {a,e} | {0,c} | {1,b} | {1,e} | {a,d} |
| c | c | {b,e} | {d,e} | {1,b} | {0,d} | {a,c} | {1,a} |
| d | d | {a,b} | {1,c} | {1,e} | {a,c} | {0,e} | {b,d} |
| e | e | {c,e} | {b,c} | {a,d} | {1,a} | {b,d} | {0,1} |

***ii.*** *The multiplicative subgroup of index 6 of the field $Z_{37}$ is*

$$G = \{1,\ 10,\ 11,\ 26,\ 27,\ 36\}$$

*and the hyperfield $Z_{37}/G$ is isomorphic to $HF_7^{160}$.*

**Table 13.** The canonical hypergroup of the hyperfield $Z_{37}/G$

| $HF_7^{160}$ | 0 | 1 | a | b | c | d | e |
|---|---|---|---|---|---|---|---|
| 0 | 0 | 1 | a | b | c | d | e |
| 1 | 1 | {0,1,a,d} | {1,b,c,d,e} | {a,b,c,e} | {a,b,d,e} | {1,a,c,e} | {a,b,c,d,e} |
| a | a | {1,b,c,d,e} | {0,a,b,e} | {1,a,c,d,e} | {1,b,c,d} | {1,b,c,e} | {1,a,b,d} |
| b | b | {a,b,c,e} | {1,a,c,d,e} | {0,1,b,c} | {1,a,b,d,e} | {a,c,d,e} | {1,a,c,d} |
| c | c | {a,b,d,e} | {1,b,c,d} | {1,a,b,d,e} | {0,a,c,d} | {1,a,b,c,e} | {1,b,d,e} |
| d | d | {1,a,c,e} | {1,b,c,e} | {a,c,d,e} | {1,a,b,c,e} | {0,b,d,e} | {1,a,b,c,d} |
| e | e | {a,b,c,d,e} | {1,a,b,d} | {1,a,c,d} | {1,b,d,e} | {1,a,b,c,d} | {0,1,c,e} |

***iii.*** *The multiplicative subgroup of index 6 of the field $Z_{61}$ is*

$$G = \{1,\ 3,\ 9,\ 20,\ 27,\ 34,\ 41,\ 52,\ 58,\ 60\}$$

*and the hyperfield $Z_{61}/G$ is isomorphic to $HF_7^{234}$.*

**Table 14.** The canonical hypergroup of the hyperfield $Z_{61}/G$

| $HF_7^{234}$ | 0 | 1 | a | b | c | d | e |
|---|---|---|---|---|---|---|---|
| 0 | 0 | 1 | a | b | c | d | e |
| 1 | 1 | {0,a,b,c,e} | {1,a,b,c,d,e} | {1,a,c,d,e} | {1,a,b,c,d,e} | {a,b,c,d,e} | {1,a,b,c,d,e} |
| a | a | {1,a,b,c,d,e} | {0,1,b,c,d} | {1,a,b,c,d,e} | {1,a,b,d,e} | {1,a,b,c,d,e} | {1,b,c,d,e} |
| b | b | {1,a,c,d,e} | {1,a,b,c,d,e} | {0,a,c,d,e} | {1,a,b,c,d,e} | {1,a,b,c,e} | {1,a,b,c,d,e} |
| c | c | {1,a,b,c,d,e} | {1,a,b,d,e} | {1,a,b,c,d,e} | {0,1,b,d,e} | {1,a,b,c,d,e} | {1,a,b,c,d} |
| d | d | {a,b,c,d,e} | {1,a,b,c,d,e} | {1,a,b,c,e} | {1,a,b,c,d,e} | {0,1,a,c,e} | {1,a,b,c,d,e} |
| e | e | {1,a,b,c,d,e} | {1,b,c,d,e} | {1,a,b,c,d,e} | {1,a,b,c,d} | {1,a,b,c,d,e} | {0,1,a,b,d} |

***iv.*** *The multiplicative subgroup of index 6 of the field $Z_{73}$ is*

$$G = \{1, 3, 8, 9, 24, 27, 46, 49, 64, 65, 70, 72\}$$

*and the hyperfield $Z_{73}/G$ is isomorphic to $HF_7^{245}$.*

**Table 15.** The canonical hypergroup of the hyperfield $Z_{73}/G$

| $HF_7^{245}$ | 0 | 1 | a | b | c | d | e |
|---|---|---|---|---|---|---|---|
| 0 | 0 | 1 | a | b | c | d | e |
| 1 | 1 | {0,1,a,b,c,d} | {1,b,c,d,e} | {1,a,b,c,e} | {1,a,b,c,d,e} | {1,a,c,d,e} | {a,b,c,d,e} |
| a | a | {1,b,c,d,e} | {0,a,b,c,d,e} | {1,a,c,d,e} | {1,a,b,c,d} | {1,a,b,c,d,e} | {1,a,b,d,e} |
| b | b | {1,a,b,c,e} | {1,a,c,d,e} | {0,1,b,c,d,e} | {1,a,b,d,e} | {a,b,c,d,e} | {1,a,b,c,d,e} |
| c | c | {1,a,b,c,d,e} | {1,a,b,c,d} | {1,a,b,d,e} | {0,1,a,c,d,e} | {1,a,b,c,e} | {1,b,c,d,e} |
| d | d | {1,a,c,d,e} | {1,a,b,c,d,e} | {a,b,c,d,e} | {1,a,b,c,e} | {0,1,a,b,d,e} | {1,a,b,c,d} |
| e | e | {a,b,c,d,e} | {1,a,b,d,e} | {1,a,b,c,d,e} | {1,b,c,d,e} | {1,a,b,c,d} | {0,1,a,b,c,e} |

***v.*** *The multiplicative subgroups of index 6 of the fields $Z_{97}$, $Z_{157}$ are respectively:*

$$G = \{1, 8, 12, 18, 22, 27, 33, 47, 50, 64, 70, 75, 79, 85, 89, 96\}$$
$$G = \begin{Bmatrix} 1, 4, 14, 16, 27, 39, 46, 49, 56, 58, 64, 67, 75, 82, \\ 90, 93, 99, 101, 108, 111, 118, 130, 141, 143, 153, 156 \end{Bmatrix}$$

*and the hyperfields $Z_{97}/G$, $Z_{157}/G$ are isomorphic to $HF_7^{267}$.*

**Table 16.** The canonical hypergroup of the hyperfields $Z_{97}/G$ and $Z_{157}/G$

| $HF_7^{267}$ | 0 | 1 | a | b | c | d | e |
|---|---|---|---|---|---|---|---|
| 0 | 0 | 1 | a | b | c | d | e |
| 1 | 1 | {0,a,b,c,d,e} | {1,a,b,c,d,e} | {1,a,b,c,d,e} | {1,a,b,c,d,e} | {1,a,b,c,d,e} | {1,a,b,c,d,e} |
| a | a | {1,a,b,c,d,e} | {0,1,b,c,d,e} | {1,a,b,c,d,e} | {1,a,b,c,d,e} | {1,a,b,c,d,e} | {1,a,b,c,d,e} |
| b | b | {1,a,b,c,d,e} | {1,a,b,c,d,e} | {0,1,a,c,d,e} | {1,a,b,c,d,e} | {1,a,b,c,d,e} | {1,a,b,c,d,e} |
| c | c | {1,a,b,c,d,e} | {1,a,b,c,d,e} | {1,a,b,c,d,e} | {0,1,a,b,d,e} | {1,a,b,c,d,e} | {1,a,b,c,d,e} |
| d | d | {1,a,b,c,d,e} | {1,a,b,c,d,e} | {1,a,b,c,d,e} | {1,a,b,c,d,e} | {0,1,a,b,c,e} | {1,a,b,c,d,e} |
| e | e | {1,a,b,c,d,e} | {1,a,b,c,d,e} | {1,a,b,c,d,e} | {1,a,b,c,d,e} | {1,a,b,c,d,e} | {0,1,a,b,c,d} |

***vi.*** *The multiplicative subgroup of index 6 of the field* $Z_{109}$ *is*

$$G = \{1, 4, 16, 27, 34, 38, 43, 45, 46, 63, 64, 66, 71, 75, 82, 93, 105, 108\}$$

*and the hyperfield* $Z_{109}/G$ *is isomorphic to* $HF_7^{246}$.

**Table 17.** The canonical hypergroup of the hyperfield $Z_{109}/G$

| $HF_7^{246}$ | 0 | 1 | a | b | c | d | e |
|---|---|---|---|---|---|---|---|
| 0 | 0 | 1 | a | b | c | d | e |
| 1 | 1 | {0,1,a,b,c,d} | {1,b,c,d,e} | {1,a,b,c,d,e} | {1,a,b,c,d,e} | {1,a,b,c,d,e} | {a,b,c,d,e} |
| a | a | {1,b,c,d,e} | {0,a,b,c,d,e} | {1,a,c,d,e} | {1,a,b,c,d,e} | {1,a,b,c,d,e} | {1,a,b,c,d,e} |
| b | b | {1,a,b,c,d,e} | {1,a,c,d,e} | {0,1,b,c,d,e} | {1,a,b,d,e} | {1,a,b,c,d,e} | {1,a,b,c,d,e} |
| c | c | {1,a,b,c,d,e} | {1,a,b,c,d,e} | {1,a,b,d,e} | {0,1,a,c,d,e} | {1,a,b,c,e} | {1,a,b,c,d,e} |
| d | d | {1,a,b,c,d,e} | {1,a,b,c,d,e} | {1,a,b,c,d,e} | {1,a,b,c,e} | {0,1,a,b,d,e} | {1,a,b,c,d} |
| e | e | {a,b,c,d,e} | {1,a,b,c,d,e} | {1,a,b,c,d,e} | {1,a,b,c,d,e} | {1,a,b,c,d} | {0,1,a,b,c,e} |

### 6.3. The quotient hyperfields F/G of order 7 that are constructed from the finite fields GF[$p^n$], n>1 and for which $G - G \neq F/G$ holds.

According to Theorem 18, this category has two members, which are the hyperfields $GF\left[5^2\right]/G$ and $GF\left[7^2\right]/G$.

The field $GF[5^2]$ consists of all the linear polynomials with coefficients in the field of residues modulo 5. In $GF[5^2]$ the polynomial $x^2+3x+4$ is irreducible. Thus, in the multiplication of the polynomials we set $x^2 = -3x-4 = 2x+1$ and then they are combined according to the ordinary rules, working modulo 5. The multiplicative subgroup of index 6 in $GF[5^2]$ is G = {1,2,3,4} and its cosets are:

$$G, xG, (x+1)G, (2x+1)G, (3x+1)G, (4x+1)G$$

The results of the hypercomposition in the above hyperfield lead to the following Proposition:

**Proposition 10.** *The hyperfield*

$$GF\left[5^2\right]\Big/ G = \left\{0, G, xG, (x+1)G, (2x+1)G, (3x+1)G, (4x+1)G\right\}$$

*is isomorphic to* $HF_7^{142}$.

**Table 18.** The canonical hypergroup of the hyperfield $GF[5^2]/G$

| $HF_7^{142}$ | 0 | 1 | a | b | c | d | e |
|---|---|---|---|---|---|---|---|
| 0 | 0 | 1 | a | b | c | d | e |
| 1 | 1 | {0,1} | {b,c,d,e} | {a,c,d,e} | {a,b,d,e} | {a,b,c,e} | {a,b,c,d} |
| a | a | {b,c,d,e} | {0,a} | {1,c,d,e} | {1,b,d,e} | {1,b,c,e} | {1,b,c,d} |
| b | b | {a,c,d,e} | {1,c,d,e} | {0,b} | {1,a,d,e} | {1,a,c,e} | {1,a,c,d} |
| c | c | {a,b,d,e} | {1,b,d,e} | {1,a,d,e} | {0,c} | {1,a,b,e} | {1,a,b,d} |
| d | d | {a,b,c,e} | {1,b,c,e} | {1,a,c,e} | {1,a,b,e} | {0,d} | {1,a,b,c} |
| e | e | {a,b,c,d} | {1,b,c,d} | {1,a,c,d} | {1,a,b,d} | {1,a,b,c} | {0,e} |

In the field $GF[7^2]$ of all linear polynomials with coefficients from $\mathbb{Z}_7$, the addition and the multiplication of the polynomials are defined in the usual way, by replacing $x^2$ with 6, since $x^2+1$ is the irreducible polynomial of degree 2. The multiplicative subgroup of index 6 in the field $GF[7^2]$ is

$$G = \{1, 6, x, 6x, 2x+2, 2x+5, 5x+2, 5x+5\}$$

and its cosets are:

$$G, \ (x+2)G, \ (x+2)^2G = (4x+3)G, \ (x+2)^3G = (4x+2)G, \ (x+2)^4G = (3x)G, \ (x+2)^5G = (6x+4)G$$

The results of the hypercomposition in the hyperfield $GF[7^2]/G$ lead to the following Proposition:

**Proposition 11.** *The hyperfield*

$$GF\left[7^2\right]\Big/ G = \left\{0, G, (x+2)G, (4x+3)G, (4x+2)G, (3x)G, (6x+4)G, \right\}$$

*is isomorphic to* $HF_7^{225}$.

**Table 19.** The canonical hypergroup of the hyperfield $GF[7^2]/G$

| $HF_7^{225}$ | 0 | 1 | a | b | c | d | e |
|---|---|---|---|---|---|---|---|
| 0 | 0 | 1 | a | b | c | d | e |
| 1 | 1 | {0,a,b,c,d} | {1,b,d,e} | {1,a,b,c,d} | {1,b,c,e} | {1,a,b,d,e} | {a,c,d,e} |
| a | a | {1,b,d,e} | {0,b,c,d,e} | {1,a,c,e} | {a,b,c,d,e} | {1,a,c,d} | {1,a,b,c,e} |
| b | b | {1,a,b,c,d} | {1,a,c,e} | {0,1,c,d,e} | {1,a,b,d} | {1,b,c,d,e} | {a,b,d,e} |
| c | c | {1,b,c,e} | {a,b,c,d,e} | {1,a,b,d} | {0,1,a,d,e} | {a,b,c,e} | {1,a,c,d,e} |
| d | d | {1,a,b,d,e} | {1,a,c,d} | {1,b,c,d,e} | {a,b,c,e} | {0,1,a,b,e} | {1,b,c,d} |
| e | e | {a,c,d,e} | {1,a,b,c,e} | {a,b,d,e} | {1,a,c,d,e} | {1,b,c,d} | {0,1,a,b,c} |

### 6.4. Quotient hyperfields F/G of order 7 for which G – G = F/G. Augmented hyperfields.

This subsection presents some seven-element hyperfields that derive as quotients of finite fields, wherein the sum of two opposite elements yields the whole hyperfield.

**Proposition 12.** *The multiplicative subgroup of index 6 of the fields* $Z_{67}$, $Z_{79}$ *and* $Z_{139}$ *are respectively:*

$$G = \{1, 9, 14, 15, 22, 24, 25, 40, 59, 62, 64\}$$
$$G = \{1, 8, 10, 21, 24, 25, 38, 46, 52, 62, 64, 65, 67\}$$
$$G = \begin{cases} 1, 6, 34, 36, 44, 45, 52, 55, 57, 63, 64, 65, \\ 77, 79, 80, 91, 100, 106, 112, 116, 125, 129, 131 \end{cases}$$

*and the hyperfields* $Z_{67}/G$, $Z_{79}/G$ *and* $Z_{139}/G$ *are isomorphic to* $HF_7^{137}$.

**Table 20.** The canonical hypergroup of the hyperfields $Z_{67}/G$, $Z_{79}/G$ *and* $Z_{139}/G$

| $HF_7^{137}$ | 0 | 1 | a | b | c | d | e |
|---|---|---|---|---|---|---|---|
| 0 | 0 | 1 | a | b | c | d | e |
| 1 | 1 | {1,a,b,d,e} | {1,a,b,c,d,e} | {1,a,b,c,d,e} | {0,1,a,b,c,d,e} | {1,a,b,c,d,e} | {1,a,b,c,d,e} |
| a | a | {1,a,b,c,d,e} | {1,a,b,c,e} | {1,a,b,c,d,e} | {1,a,b,c,d,e} | {0,1,a,b,c,d,e} | {1,a,b,c,d,e} |
| b | b | {1,a,b,c,d,e} | {1,a,b,c,d,e} | {1,a,b,c,d} | {1,a,b,c,d,e} | {1,a,b,c,d,e} | {0,1,a,b,c,d,e} |
| c | c | {0,1,a,b,c,d,e} | {1,a,b,c,d,e} | {1,a,b,c,d,e} | {a,b,c,d,e} | {1,a,b,c,d,e} | {1,a,b,c,d,e} |
| d | d | {1,a,b,c,d,e} | {0,1,a,b,c,d,e} | {1,a,b,c,d,e} | {1,a,b,c,d,e} | {1,b,c,d,e} | {1,a,b,c,d,e} |
| e | e | {1,a,b,c,d,e} | {1,a,b,c,d,e} | {0,1,a,b,c,d,e} | {1,a,b,c,d,e} | {1,a,b,c,d,e} | {1,a,c,d,e} |

**Proposition 13.**

*i.* *The multiplicative subgroups of index 6 of the fields* $Z_{103}$, $Z_{127}$, $Z_{151}$, $Z_{163}$ *are respectively:*

$$G = \{1, 8, 9, 13, 14, 23, 30, 34, 61, 64, 66, 72, 76, 79, 81, 93, 100\}$$

$$G = \begin{Bmatrix} 1, 2, 4, 8, 16, 19, 25, 32, 38, 47, 50, \\ 61, 64, 73, 76, 87, 94, 100, 107, 117, 122 \end{Bmatrix}$$

$$G = \begin{Bmatrix} 1, 8, 9, 19, 20, 29, 44, 50, 59, 64, 68, 72, 78, 81, \\ 84, 86, 91, 94, 98, 110, 123 \quad, 124, 125, 127, 148 \end{Bmatrix}$$

$$G = \begin{Bmatrix} 1, 6, 21, 22, 25, 36, 38, 40, 53, 58, 61, 64, 65, 77, 85, \\ 104, 115, 126, 132, 133, 135, 136, 140, 146, 150, 155, 158 \end{Bmatrix}$$

*and the hyperfields $Z_{103}/G$, $Z_{127}/G$, $Z_{151}/G$, $Z_{163}/G$ are isomorphic to $HF_7^{141}$.*

**Table 21.** The canonical hypergroup of the hyperfields $Z_{103}/G$, $Z_{127}/G$, $Z_{151}/G$ and $Z_{163}/G$

| $HF_7^{141}$ | 0 | 1 | a | b | c | d | e |
|---|---|---|---|---|---|---|---|
| 0 | 0 | 1 | a | b | c | d | e |
| 1 | 1 | {1,a,b,c,d,e} | {1,a,b,c,d,e} | {1,a,b,c,d,e} | {0,1,a,b,c,d,e} | {1,a,b,c,d,e} | {1,a,b,c,d,e} |
| a | a | {1,a,b,c,d,e} | {1,a,b,c,d,e} | {1,a,b,c,d,e} | {1,a,b,c,d,e} | {0,1,a,b,c,d,e} | {1,a,b,c,d,e} |
| b | b | {1,a,b,c,d,e} | {1,a,b,c,d,e} | {1,a,b,c,d,e} | {1,a,b,c,d,e} | {1,a,b,c,d,e} | {0,1,a,b,c,d,e} |
| c | c | {0,1,a,b,c,d,e} | {1,a,b,c,d,e} | {1,a,b,c,d,e} | {1,a,b,c,d,e} | {1,a,b,c,d,e} | {1,a,b,c,d,e} |
| d | d | {1,a,b,c,d,e} | {0,1,a,b,c,d,e} | {1,a,b,c,d,e} | {1,a,b,c,d,e} | {1,a,b,c,d,e} | {1,a,b,c,d,e} |
| e | e | {1,a,b,c,d,e} | {1,a,b,c,d,e} | {0,1,a,b,c,d,e} | {1,a,b,c,d,e} | {1,a,b,c,d,e} | {1,a,b,c,d,e} |

***ii.*** *The multiplicative subgroup of index 6 of the field $Z_{181}$ is*

$$G = \begin{Bmatrix} 1, 5, 25, 27, 29, 36, 42, 46, 48, 49, 56, 59, 64, 67, 82, 99, 114, \\ 117, 122, 125, 132, 133, 135, 139, 145, 152, 154, 156, 176, 180 \end{Bmatrix}$$

*and the hyperfield $Z_{181}/G$ is isomorphic to $HF_7^{277}$.*

**Table 22.** The canonical hypergroup of the hyperfield $Z_{181}/G$

| $HF_7^{277}$ | 0 | 1 | a | b | c | d | e |
|---|---|---|---|---|---|---|---|
| 0 | 0 | 1 | a | b | c | d | e |
| 1 | 1 | {0,1,a,b,c,d,e} | {1,a,b,c,d,e} | {1,a,b,c,d,e} | {1,a,b,c,d,e} | {1,a,b,c,d,e} | {1,a,b,c,d,e} |
| a | a | {1,a,b,c,d,e} | {0,1,a,b,c,d,e} | {1,a,b,c,d,e} | {1,a,b,c,d,e} | {1,a,b,c,d,e} | {1,a,b,c,d,e} |
| b | b | {1,a,b,c,d,e} | {1,a,b,c,d,e} | {0,1,a,b,c,d,e} | {1,a,b,c,d,e} | {1,a,b,c,d,e} | {1,a,b,c,d,e} |
| c | c | {1,a,b,c,d,e} | {1,a,b,c,d,e} | {1,a,b,c,d,e} | {0,1,a,b,c,d,e} | {1,a,b,c,d,e} | {1,a,b,c,d,e} |
| d | d | {1,a,b,c,d,e} | {1,a,b,c,d,e} | {1,a,b,c,d,e} | {1,a,b,c,d,e} | {0,1,a,b,c,d,e} | {1,a,b,c,d,e} |
| e | e | {1,a,b,c,d,e} | {1,a,b,c,d,e} | {1,a,b,c,d,e} | {1,a,b,c,d,e} | {1,a,b,c,d,e} | {0,1,a,b,c,d,e} |

The field $GF[11^2]$ consists of all the linear polynomials with coefficients from the field of residues modulo 11. In $GF[11^2]$ the polynomial $x^2 + 1$ is irreducible. Thus, the polynomials are combined according to the ordinary rules, working modulo 11, by setting $x^2 = -1 = 10$. The multiplicative subgroup of index 6 in $GF[11^2]$ is

$$G = \{1, 2, 3, 4, 5, 6, 7, 8, 9, 10, x, 2x, 3x, 4x, 5x, 6x, 7x, 8x, 9x, 10x\}$$

and its cosets are:

$$G, \ (3x+1)G, \ (3x+1)^2 G = (6x+3)G, \ (3x+1)^3 G = (4x+7)G,$$

$$(3x+1)^4 G = (3x+6)G, \ (3x+1)^5 G = (10x+8)G.$$

The outcomes of the hypercomposition within the hyperfield $GF\left[11^2\right] / G$ give rise to the next Proposition:

**Proposition 14.** *The hyperfield* $GF\left[11^2\right] / G$ *is isomorphic to* $HF_7^{277}$ .

In the field $GF[13^2]$ the irreducible polynomial is $x^2 + 2$. Thus, in the multiplication of the polynomials we set $x^2 = -2 = 11$. The multiplicative subgroup of index 6 in $GF[13^2]$ is

$$G = \begin{cases} 1, \ 5, \ 8, \ 12, \\ 5x+1, \ 8x+1, \ 2x+2, \ 11x+2, \ 3x+3, \ 10x+3, \ 5x+4, \ 8x+4, \ x+5, \ 12x+5, \\ x+6, \ 12x+6, \ x+7, \ 12x+7, \ x+8, \ 12x+8, \ 5x+9, \ 8x+9, \\ 3x+10, \ 10x+10, \ 2x+11, \ 11x+11, \ 5x+12, \ 8x+12 \end{cases}$$

which has 28 elements and its cosets are:

$$G, \ xG, \ x^2 G = 11G, \ x^3 G = 11xG, \ x^4 G = 4G, \ x^5 G = 11xG.$$

Therefore, the Proposition holds:

**Proposition 15.** *The hyperfield* $GF\left[13^2\right] / G$ *is isomorphic to* $HF_7^{277}$ .

A hyperfield of the type of $HF_7^{277}$ or $HF_7^{141}$ is called *total hyperfield*.

**Example 3**. The multiplicative subgroup of index 6 of the field $Z_{193}$ is:

$$G = \begin{cases} 1, \ 3, \ 8, \ 9, \ 14, \ 23, \ 24, \ 27, \ 42, \ 43, \ 50, \ 64, \ 67, \ 69, \ 72, \ 81, \ 112, \ 121, \ 124, \\ 126, \ 129, \ 143, \ 150, \ 151, \ 166, \ 169, \ 170, \ 179, \ 184, \ 185, \ 190, \ 192 \end{cases}$$

Its cosets are:

$$G, \ 5G, \ 5^2 G, \ 5^3 G, \ 5^4 G, \ 5^5 G$$

and the hyperfield $\mathbb{Z}_{193} / G$ is total with self-opposite elements, that is isomorphic to $HF_7^{277}$ .

**Example 4**. The multiplicative subgroup of index 6 of the field $Z_{199}$ is:

$$G = \begin{cases} 1, \ 5, \ 8, \ 18, \ 25, \ 28, \ 40, \ 52, \ 61, \ 62, \ 63, \ 64, \ 90, \ 92, \ 98, \ 103, \ 106, \ 111, \ 114, \\ 116, \ 117, \ 121, \ 123, \ 125, \ 132, \ 139, \ 140, \ 144, \ 157, \ 172, \ 182, \ 187, \ 188 \end{cases}$$

Its cosets are:

$$G, \ 3G, \ 3^2G, \ 3^3G, \ 3^4G, \ 3^5G$$

and the hyperfield $\mathbb{Z}_{199} / G$ is total with no self-opposite elements, that is isomorphic to $HF_7^{141}$

The combination of Propositions 13, 14 and 15 along with Examples 3 and 4, gives rise to the conjecture:

**Conjecture**: *If the order of the multiplicative subgroup of a finite field exceeds a certain number, then the generated quotient hyperfield is total.*

The next Theorems refer to the augmented hyperfields of the above quotient hyperfields. Bear in mind that if $(H,+,\cdot)$ is a field or a hyperfield, then its augmented hyperfield is the hyperfield in which the sum of any two non-zero and non-opposite elements is extended (augmented) to include the two addends (extensive enlargement of the hypercomposition [5, example 8]). Proposition 1 in [13] (see also Proposition 2 in [1]) shows that, in this case, the sum of two opposite elements is equal to the entire set $H$. Thus, the augmented hyperfield of a hyperfield is endowed with the following hypercomposition «$\dotplus$»:

$x \dotplus y = \{x, y\} \cup x+y,$      for all $x,y \in H^*$, with $y \neq -x$,

$x \dotplus (-x) = H,$      for all $x \in H^*$,

$x \dotplus 0 = 0 \dotplus x = x,$      for all $x \in H$.

The augmented hyperfield of a hyperfield $H$ is denoted by $[H]$ and as proved in [13] (Proposition 4) and in [1] (Theorems 4 and 5) if $H$ is a field or a quotient hyperfield, then $[H]$ is also a quotient hyperfield.

**Theorem 19.** *The augmented hyperfield of the field $\mathbb{Z}_7$ is isomorphic to $HF_7^2$ .*

**Table 23.** The canonical hypergroup of the augmented hyperfield of the field $Z_7$

| $[\mathbb{Z}_7]$ $\cong$ $HF_7^2$ | 0 | 1 | a | b | c | d | e |
|---|---|---|---|---|---|---|---|
| 0 | 0 | 1 | a | b | c | d | e |
| 1 | 1 | {1,b} | {1,a,d} | {1,a,b} | {0,1,a,b,c,d,e} | {1,d,e} | {1,c,e} |
| a | a | {1,a,d} | {a,c} | {a,b,e} | {a,b,c} | {0,1,a,b,c,d,e} | {1,a,e} |
| b | b | {1,a,b} | {a,b,e} | {b,d} | {1,b,c} | {b,c,d} | {0,1,a,b,c,d,e} |
| c | c | {0,1,a,b,c,d,e} | {a,b,c} | {1,b,c} | {c,e} | {a,c,d} | {c,d,e} |
| d | d | {1,d,e} | {0,1,a,b,c,d,e} | {b,c,d} | {a,c,d} | {1,d} | {b,d,e} |
| e | e | {1,c,e} | {1,a,e} | {0,1,a,b,c,d,e} | {c,d,e} | {b,d,e} | {a,e} |

**Theorem 20.** *If G is the multiplicative subgroup of index 6 in the following fields, then:*

***i.*** *the augmented hyperfield of the hyperfield $\mathbb{Z}_{13} / G$ is isomorphic to $HF_7^{271}$ , i.e.*

$$\left[\mathbb{Z}_{13} / G\right] \cong \left[HF_7^{143}\right] \cong HF_7^{271}$$

**ii.** the augmented hyperfield of the hyperfield $\mathbb{Z}_{19} / G$ is isomorphic to $HF_7^{96}$, i.e.

$$\left[\mathbb{Z}_{19} / G\right] \cong \left[HF_7^9\right] \cong HF_7^{96}$$

**iii.** the augmented hyperfield of the hyperfield $\mathbb{Z}_{31} / G$ is isomorphic to $HF_7^{95}$, i.e.

$$\left[\mathbb{Z}_{31} / G\right] \cong \left[HF_7^{13}\right] \cong HF_7^{95}$$

**iv.** the augmented hyperfield of the hyperfield $\mathbb{Z}_{37} / G$ is isomorphic to $HF_7^{276}$, i.e.

$$\left[\mathbb{Z}_{37} / G\right] \cong \left[HF_7^{160}\right] \cong HF_7^{276}$$

**v.** the augmented hyperfield of the hyperfield $\mathbb{Z}_{43} / G$ is isomorphic to $HF_7^{101}$, i.e.

$$\left[\mathbb{Z}_{43} / G\right] \cong \left[HF_7^{61}\right] \cong HF_7^{110}$$

**vi.** the augmented hyperfield of the hyperfield $\mathbb{Z}_{61} / G$ is isomorphic to $HF_7^{277}$, i.e.

$$\left[\mathbb{Z}_{61} / G\right] \cong \left[HF_7^{234}\right] \cong HF_7^{277}$$

**vii.** for the hyperfields $\mathbb{Z}_{67} / G$, $\mathbb{Z}_{79} / G$ and $\mathbb{Z}_{139} / G$ it holds that:

$$\left[\mathbb{Z}_{67} / G\right] \cong \left[\mathbb{Z}_{79} / G\right] \cong \left[\mathbb{Z}_{139} / G\right] \cong \left[HF_7^{137}\right] \equiv HF_7^{137}$$

**viii.** the augmented hyperfield of the hyperfield $\mathbb{Z}_{73} / G$ is isomorphic to $HF_7^{276}$, i.e.

$$\left[\mathbb{Z}_{73} / G\right] \cong \left[HF_7^{245}\right] \cong HF_7^{276}$$

**ix.** the augmented hyperfield of:

    (a)   the isomorphic hyperfields $\mathbb{Z}_{97} / G$ and $\mathbb{Z}_{157} / G$,

    (b)   the hyperfield $\mathbb{Z}_{109} / G$ and

    (c)   the hyperfield $\mathbb{Z}_{181} / G$

is isomorphic to $HF_7^{277}$, that is

$$\left[\mathbb{Z}_{97} / G\right] \cong \left[\mathbb{Z}_{157} / G\right] \cong \left[HF_7^{267}\right] \cong HF_7^{277},$$

$$\left[\mathbb{Z}_{109} / G\right] \cong \left[HF_7^{246}\right] \cong HF_7^{277}$$

and
$$\left[\mathbb{Z}_{181} / G\right] \cong \left[HF_7^{277}\right] \equiv HF_7^{277}.$$

**x.** for the hyperfields $\mathbb{Z}_{103} / G$, $\mathbb{Z}_{127} / G$, $\mathbb{Z}_{151} / G$, $\mathbb{Z}_{163} / G$ it holds that:

$$\left[\mathbb{Z}_{103} / G\right] \cong \left[\mathbb{Z}_{127} / G\right] \cong \left[\mathbb{Z}_{151} / G\right] \cong \left[\mathbb{Z}_{163} / G\right] \cong \left[HF_7^{141}\right] \equiv HF_7^{141}$$

***xi.*** *the augmented hyperfield of the hyperfield* $GF\left[5^2\right]/G$ *is isomorphic to* $HF_7^{277}$ *, i.e.*

$$\left[GF\left[5^2\right]/G\right] \cong \left[HF_7^{142}\right] \equiv HF_7^{277}$$

***xii.*** *the augmented hyperfield of the hyperfield* $GF\left[7^2\right]/G$ *is isomorphic to* $HF_7^{275}$ *, i.e.*

$$\left[GF\left[7^2\right]/G\right] \cong \left[HF_7^{225}\right] \equiv HF_7^{275}$$

The proofs for cases (i) through (xi) of the aforementioned Theorem are straightforward; however, the proof for the last case is included in the Appendix along with the necessary analysis.

### 6.5. Non quotient hyperfields of order 7.

The analysis and study of the various cases, which is carried out through Propositions 8, 9, 10 and 11 leads to the following Theorem, which classifies all the hyperfields of the sections A1-A3 and B1-B5 of the Appendix:

**Theorem 21.**

***i.*** *Section A1 consists of the field $Z_7$ and its augmented hyperfield.*

***ii.*** *All hyperfields in sections A2i and A2ii are non-quotient hyperfields.*

***iii.*** *All hyperfields in section A3i are non-quotient hyperfields with the exception of $HF_7^{13}$ , which is a quotient hyperfield.*

***iv.*** *All hyperfields in section A3ii are non-quotient hyperfields with the exception of $HF_7^{61}$ , which is a quotient hyperfield.*

***v.*** *Section B1 consists of quotient hyperfields.*

***vi.*** *Section B2 consists of non-quotient hyperfields.*

***vii.*** *All hyperfields in section B3i are non-quotient hyperfields with the exception of $HF_7^{160}$ , which is a quotient hyperfield.*

***viii.*** *All hyperfields in sections B3ii and B4i are non-quotient hyperfields.*

***ix.*** *All hyperfields in section B4ii are non-quotient hyperfields with the exception of $HF_7^{225}$ and $HF_7^{234}$ , which are quotient hyperfields.*

***x.*** *All hyperfields in section B5i are non-quotient hyperfields with the exception of $HF_7^{245}$ and $HF_7^{246}$ , which are quotient hyperfields.*

***xi.*** *All hyperfields in section B5ii are non-quotient hyperfields with the exception of $HF_7^{267}$ , which is a quotient hyperfield.*

# Appendix

The enumeration of hypercompositional structures has been a research area since the 1980s. The initial publication addressing this subject was authored by M. De Salvo and D. Freni [41]. Subsequently, R. Migliorato focused on the enumeration of 3-element hypergroups [42], with his findings being validated by later studies [43, 44], which progressively refined the computational techniques employed in the enumeration process. Since then, this topic has been the subject of numerous papers (e.g. [45-64]). The hyperfields of order less than 6 are enumerated in [1, 52-55], from which [1] presents their detailed classification with the necessary extensions and corrections of the previous results and summarizes all of them in its Table 26. Considering the field as a special case of the hyperfield, there exist two hyperfields of order 2, five hyperfields of order 3 and 27 hyperfields of order 5. Interestingly, all the hyperfields of order 2 and 3 are quotient hyperfields [1].

This Appendix presents the hyperfields of order 7, which are produced as described in the above section 5 and with the use of the Mathematica [65] packages that are developed in [44, 56-59]. A collateral direct consequence is that the corresponding family of canonical 7-element hypergroups is also revealed. The symbols **0**, **1**, **a**, **b**, **c**, **d** and **e** are used to denote these hyperfields' elements, where **a** is the generator of their multiplicative subgroup and **b**, **c**, **d** and **e** represent the elements $a^2$, $a^3$, $a^4$ and $a^5$ respectively.

The table of the order 7 hyperfields' multiplicative subgroup is the following one:

**Table 24.** The multiplicative subgroup of the hyperfields of order 7

| $\left(HF_7^*, \cdot\right)$ | 1 | a | b=α² | c=α³ | d=α⁴ | e=α⁵ |
|---|---|---|---|---|---|---|
| 1 | 1 | a | b | c | d | e |
| a | a | b | c | d | e | 1 |
| b=α² | b | c | d | e | 1 | a |
| c=α³ | c | d | e | 1 | a | b |
| d=α⁴ | d | e | 1 | a | b | c |
| e=α⁵ | e | 1 | a | b | c | d |

Note that the isomorphic for each hyperfield in the list below can be obtained by interchanging **a** with **e** and **b** with **d**. For example, the augmented hyperfield of $HF_7^{225}$, referenced in Theorem 20.xii, has the following canonical additive hypergroup:

**Table 25.** The canonical hypergroup of the augmented hyperfield $\left[ HF_7^{225} \right]$ of $HF_7^{225}$

| $\left[ HF_7^{225} \right]$ | 0 | 1 | a | b | c | d | e |
|---|---|---|---|---|---|---|---|
| 0 | 0 | 1 | a | b | c | d | e |
| 1 | 1 | {0,1,a,b,c,d,e} | {1,a,b,d,e} | {1,a,b,c,d} | {1,b,c,e} | {1,a,b,d,e} | {1,a,c,d,e} |
| a | a | {1,a,b,d,e} | {0,1,a,b,c,d,e} | {1,a,b,c,e} | {a,b,c,d,e} | {1,a,c,d} | {1,a,b,c,e} |
| b | b | {1,a,b,c,d} | {1,a,b,c,e} | {0,1,a,b,c,d,e} | {1,a,b,c,d} | {1,b,c,d,e} | {a,b,d,e} |
| c | c | {1,b,c,e} | {a,b,c,d,e} | {1,a,b,c,d} | {0,1,a,b,c,d,e} | {a,b,c,d,e} | {1,a,c,d,e} |
| d | d | {1,a,b,d,e} | {1,a,c,d} | {1,b,c,d,e} | {a,b,c,d,e} | {0,1,a,b,c,d,e} | {1,b,c,d,e} |
| e | e | {1,a,c,d,e} | {1,a,b,c,e} | {a,b,d,e} | {1,a,c,d,e} | {1,b,c,d,e} | {0,1,a,b,c,d,e} |

This hyperfield is not among the ones in the list below. However, the list contains its isomorphic hyperfield $HF_7^{275}$, which is obtained from $\left[ HF_7^{225} \right]$ by interchanging **a** with **e** and **b** with **d** in the above Cayley table.

The categories in which the hyperfields of order 7 are classified in the following list are:

A. the hyperfields of order 7 which do not have self-opposite elements. For these $0 \in 1+c$.

B. the hyperfields of order 7 with self-opposite elements. For these $0 \in 1+1$.

The above categories were divided into subcategories according to the number of elements included in the difference $x$-$x$. Observe that when card($x$-$x$)=7, then the hypercomposition is closed, i.e. $x,y \in x+y$ (see Proposition 1 [13], Proposition 2 [1]). The following tables summarize the number of hyperfields of order 7 in each category and subcategory.

**Table 26.**

| A. The class of hyperfields with no self-opposite elements | Number of elements: 141 |
|---|---|
| Cardinality of x-x | Number of hyperfields |
| 1 | 1 (the field $Z_7$) |
| 2 | 0 |
| 3 | 9 |
| 4 | 0 |
| 5 | 77 |
| 6 | 0 |
| 7 | 54 |

**Table 27.**

| B. The class of hyperfields with self-opposite elements | Number of elements: 136 |
|---|---|
| Cardinality of x-x | Number of hyperfields |
| 1 | 0 |
| 2 | 2 |
| 3 | 9 |
| 4 | 40 |
| 5 | 46 |
| 6 | 29 |
| 7 | 10 |

## A. The hyperfields which do not have self-opposite elements

## A1. The field $\mathbb{Z}_7$ and its augmented hyperfield.

| $HF_7^1 = \mathbb{Z}_7$ | 0 | 1 | a | b | c | d | e |
|---|---|---|---|---|---|---|---|
| 0 | 0 | 1 | a | b | c | d | e |
| 1 | 1 | b | d | a | 0 | e | c |
| a | a | d | c | e | b | 0 | 1 |
| b | b | a | e | d | 1 | c | 0 |
| c | c | 0 | b | 1 | e | a | d |
| d | d | e | 0 | c | a | 1 | b |
| e | e | c | 1 | 0 | d | b | a |

| $HF_7^2$ $= [\mathbb{Z}_7]$ | 0 | 1 | a | b | c | d | e |
|---|---|---|---|---|---|---|---|
| 0 | 0 | 1 | a | b | c | d | e |
| 1 | 1 | {1,b} | {1,a,d} | {1,a,b} | {0,1,a,b,c,d,e} | {1,d,e} | {1,c,e} |
| a | a | {1,a,d} | {a,c} | {a,b,e} | {a,b,c} | {0,1,a,b,c,d,e} | {1,a,e} |
| b | b | {1,a,b} | {a,b,e} | {b,d} | {1,b,c} | {b,c,d} | {0,1,a,b,c,d,e} |
| c | c | {0,1,a,b,c,d,e} | {a,b,c} | {1,b,c} | {c,e} | {a,c,d} | {c,d,e} |
| d | d | {1,d,e} | {0,1,a,b,c,d,e} | {b,c,d} | {a,c,d} | {1,d} | {b,d,e} |
| e | e | {1,c,e} | {1,a,e} | {0,1,a,b,c,d,e} | {c,d,e} | {b,d,e} | {a,e} |

**A2i. Hyperfields for which card(x-x) = 3 and x, -x ∈ x-x, for every non-zero element x.**

| $HF_7^3$ | 0 | 1 | a | b | c | d | e |
|---|---|---|---|---|---|---|---|
| 0 | 0 | 1 | a | b | c | d | e |
| 1 | 1 | {1,a,b,d,e} | {b,c,d,e} | {a,c,d,e} | {0,1,c} | {a,b,c,e} | {a,b,c,d} |
| a | a | {b,c,d,e} | {1,a,b,c,e} | {1,c,d,e} | {1,b,d,e} | {0,a,d} | {1,b,c,d} |
| b | b | {a,c,d,e} | {1,c,d,e} | {1,a,b,c,d} | {1,a,d,e} | {1,a,c,e} | {0,b,e} |
| c | c | {0,1,c} | {1,b,d,e} | {1,a,d,e} | {a,b,c,d,e} | {1,a,b,e} | {1,a,b,d} |
| d | d | {a,b,c,e} | {0,a,d} | {1,a,c,e} | {1,a,b,e} | {1,b,c,d,e} | {1,a,b,c} |
| e | e | {a,b,c,d} | {1,b,c,d} | {0,b,e} | {1,a,b,d} | {1,a,b,c} | {1,a,c,d,e} |

| $HF_7^4$ | 0 | 1 | a | b | c | d | e |
|---|---|---|---|---|---|---|---|
| 0 | 0 | 1 | a | b | c | d | e |
| 1 | 1 | {1,a,b,c,d,e} | {b,c,d,e} | {a,c,d,e} | {0,1,c} | {a,b,c,e} | {a,b,c,d} |
| a | a | {b,c,d,e} | {1,a,b,c,d,e} | {1,c,d,e} | {1,b,d,e} | {0,a,d} | {1,b,c,d} |
| b | b | {a,c,d,e} | {1,c,d,e} | {1,a,b,c,d,e} | {1,a,d,e} | {1,a,c,e} | {0,b,e} |
| c | c | {0,1,c} | {1,b,d,e} | {1,a,d,e} | {1,a,b,c,d,e} | {1,a,b,e} | {1,a,b,d} |
| d | d | {a,b,c,e} | {0,a,d} | {1,a,c,e} | {1,a,b,e} | {1,a,b,c,d,e} | {1,a,b,c} |
| e | e | {a,b,c,d} | {1,b,c,d} | {0,b,e} | {1,a,b,d} | {1,a,b,c} | {1,a,b,c,d,e} |

**A2ii. Hyperfields for which card(x-x) = 3 and x, -x ∉ x-x, for every non-zero element x.**

| $HF_7^5$ | 0 | 1 | a | b | c | d | e |
|---|---|---|---|---|---|---|---|
| 0 | 0 | 1 | a | b | c | d | e |
| 1 | 1 | {a,b,c,e} | {1,b,d,e} | {a,b,c,d,e} | {0,a,d} | {1,a,b,c,e} | {a,c,d,e} |
| a | a | {1,b,d,e} | {1,b,c,d} | {1,a,c,e} | {1,b,c,d,e} | {0,b,e} | {1,a,b,c,d} |
| b | b | {a,b,c,d,e} | {1,a,c,e} | {a,c,d,e} | {1,a,b,d} | {1,a,c,d,e} | {0,1,c} |
| c | c | {0,a,d} | {1,b,c,d,e} | {1,a,b,d} | {1,b,d,e} | {a,b,c,e} | {1,a,b,d,e} |
| d | d | {1,a,b,c,e} | {0,b,e} | {1,a,c,d,e} | {a,b,c,e} | {1,a,c,e} | {1,b,c,d} |
| e | e | {a,c,d,e} | {1,a,b,c,d} | {0,1,c} | {1,a,b,d,e} | {1,b,c,d} | {1,a,b,d} |

| $HF_7^6$ | 0 | 1 | a | b | c | d | e |
|---|---|---|---|---|---|---|---|
| 0 | 0 | 1 | a | b | c | d | e |
| 1 | 1 | {a,b,c,d,e} | {1,b,c,d,e} | {b,c,d,e} | {0,a,d} | {1,a,b,c} | {a,b,c,d,e} |
| a | a | {1,b,c,d,e} | {1,b,c,d,e} | {1,a,c,d,e} | {1,c,d,e} | {0,b,e} | {a,b,c,d} |
| b | b | {b,c,d,e} | {1,a,c,d,e} | {1,a,c,d,e} | {1,a,b,d,e} | {1,a,d,e} | {0,1,c} |
| c | c | {0,a,d} | {1,c,d,e} | {1,a,b,d,e} | {1,a,b,d,e} | {1,a,b,c,e} | {1,a,b,e} |
| d | d | {1,a,b,c} | {0,b,e} | {1,a,d,e} | {1,a,b,c,e} | {1,a,b,c,e} | {1,a,b,c,d} |
| e | e | {a,b,c,d,e} | {a,b,c,d} | {0,1,c} | {1,a,b,e} | {1,a,b,c,d} | {1,a,b,c,d} |

| $HF_7^7$ | 0 | 1 | a | b | c | d | e |
|---|---|---|---|---|---|---|---|
| 0 | 0 | 1 | a | b | c | d | e |
| 1 | 1 | {a,b,c,d,e} | {1,b,c,d,e} | {a,b,c,d,e} | {0,a,d} | {1,a,b,c,e} | {a,b,c,d,e} |
| a | a | {1,b,c,d,e} | {1,b,c,d,e} | {1,a,c,d,e} | {1,b,c,d,e} | {0,b,e} | {1,a,b,c,d} |
| b | b | {a,b,c,d,e} | {1,a,c,d,e} | {1,a,c,d,e} | {1,a,b,d,e} | {1,a,c,d,e} | {0,1,c} |
| c | c | {0,a,d} | {1,b,c,d,e} | {1,a,b,d,e} | {1,a,b,d,e} | {1,a,b,c,e} | {1,a,b,d,e} |
| d | d | {1,a,b,c,e} | {0,b,e} | {1,a,c,d,e} | {1,a,b,c,e} | {1,a,b,c,e} | {1,a,b,c,d} |
| e | e | {a,b,c,d,e} | {1,a,b,c,d} | {0,1,c} | {1,a,b,d,e} | {1,a,b,c,d} | {1,a,b,c,d} |

| $HF_7^8$ | 0 | 1 | a | b | c | d | e |
|---|---|---|---|---|---|---|---|
| 0 | 0 | 1 | a | b | c | d | e |
| 1 | 1 | {a,b,c,e} | {1,b,d,e} | {b,c,d,e} | {0,a,d} | {1,a,b,c} | {a,c,d,e} |
| a | a | {1,b,d,e} | {1,b,c,d} | {1,a,c,e} | {1,c,d,e} | {0,b,e} | {a,b,c,d} |
| b | b | {b,c,d,e} | {1,a,c,e} | {a,c,d,e} | {1,a,b,d} | {1,a,d,e} | {0,1,c} |
| c | c | {0,a,d} | {1,c,d,e} | {1,a,b,d} | {1,b,d,e} | {a,b,c,e} | {1,a,b,e} |
| d | d | {1,a,b,c} | {0,b,e} | {1,a,d,e} | {a,b,c,e} | {1,a,c,e} | {1,b,c,d} |
| e | e | {a,c,d,e} | {a,b,c,d} | {0,1,c} | {1,a,b,e} | {1,b,c,d} | {1,a,b,d} |

| $HF_7^9$ | 0 | 1 | a | b | c | d | e |
|---|---|---|---|---|---|---|---|
| 0 | 0 | 1 | a | b | c | d | e |
| 1 | 1 | {a,c} | {a,b,e} | {1,d,e} | {0,b,e} | {b,c,d} | {1,a,d} |
| a | a | {a,b,e} | {b,d} | {1,b,c} | {1,a,e} | {0,1,c} | {c,d,e} |
| b | b | {1,d,e} | {1,b,c} | {c,e} | {a,c,d} | {1,a,b} | {0,a,d} |
| c | c | {0,b,e} | {1,a,e} | {a,c,d} | {1,d} | {b,d,e} | {a,b,c} |
| d | d | {b,c,d} | {0,1,c} | {1,a,b} | {b,d,e} | {a,e} | {1,c,e} |
| e | e | {1,a,d} | {c,d,e} | {0,a,d} | {a,b,c} | {1,c,e} | {1,b} |

| $HF_7^{10}$ | 0 | 1 | a | b | c | d | e |
|---|---|---|---|---|---|---|---|
| 0 | 0 | 1 | a | b | c | d | e |
| 1 | 1 | {a,c,d} | {a,b,c,e} | {1,a,d,e} | {0,b,e} | {b,c,d,e} | {1,a,b,d} |
| a | a | {a,b,c,e} | {b,d,e} | {1,b,c,d} | {1,a,b,e} | {0,1,c} | {1,c,d,e} |
| b | b | {1,a,d,e} | {1,b,c,d} | {1,c,e} | {a,c,d,e} | {1,a,b,c} | {0,a,d} |
| c | c | {0,b,e} | {1,a,b,e} | {a,c,d,e} | {1,a,d} | {1,b,d,e} | {a,b,c,d} |
| d | d | {b,c,d,e} | {0,1,c} | {1,a,b,c} | {1,b,d,e} | {a,b,e} | {1,a,c,e} |
| e | e | {1,a,b,d} | {1,c,d,e} | {0,a,d} | {a,b,c,d} | {1,a,c,e} | {1,b,c} |

| $HF_7^{11}$ | 0 | 1 | a | b | c | d | e |
|---|---|---|---|---|---|---|---|
| 0 | 0 | 1 | a | b | c | d | e |
| 1 | 1 | {a,b,c,d} | {a,b,c,d,e} | {1,a,d,e} | {0,b,e} | {b,c,d,e} | {1,a,b,c,d} |
| a | a | {a,b,c,d,e} | {b,c,d,e} | {1,b,c,d,e} | {1,a,b,e} | {0,1,c} | {1,c,d,e} |
| b | b | {1,a,d,e} | {1,b,c,d,e} | {1,c,d,e} | {1,a,c,d,e} | {1,a,b,c} | {0,a,d} |
| c | c | {0,b,e} | {1,a,b,e} | {1,a,c,d,e} | {1,a,d,e} | {1,a,b,d,e} | {a,b,c,d} |
| d | d | {b,c,d,e} | {0,1,c} | {1,a,b,c} | {1,a,b,d,e} | {1,a,b,e} | {1,a,b,c,e} |
| e | e | {1,a,b,c,d} | {1,c,d,e} | {0,a,d} | {a,b,c,d} | {1,a,b,c,e} | {1,a,b,c} |

**A3i. Hyperfields for which card(x-x) = 5 and  x, -x ∈ x-x, for every non-zero element x.**

| $HF_7^{12}$ | 0 | 1 | a | b | c | d | e |
|---|---|---|---|---|---|---|---|
| 0 | 0 | 1 | a | b | c | d | e |
| 1 | 1 | {1,a,b} | {1,b,d,e} | {b,d,e} | {0,1,a,c,d} | {1,b,c} | {a,c,d,e} |
| a | a | {1,b,d,e} | {a,b,c} | {1,a,c,e} | {1,c,e} | {0,a,b,d,e} | {a,c,d} |
| b | b | {b,d,e} | {1,a,c,e} | {b,c,d} | {1,a,b,d} | {1,a,d} | {0,1,b,c,e} |
| c | c | {0,1,a,c,d} | {1,c,e} | {1,a,b,d} | {c,d,e} | {a,b,c,e} | {a,b,e} |
| d | d | {1,b,c} | {0,a,b,d,e} | {1,a,d} | {a,b,c,e} | {1,d,e} | {1,b,c,d} |
| e | e | {a,c,d,e} | {a,c,d} | {0,1,b,c,e} | {a,b,e} | {1,b,c,d} | {1,a,e} |

| $HF_7^{13}$ | 0 | 1 | a | b | c | d | e |
|---|---|---|---|---|---|---|---|
| 0 | 0 | 1 | a | b | c | d | e |
| 1 | 1 | {1,a,b} | {1,b,d,e} | {a,b,d,e} | {0,1,a,c,d} | {1,b,c,e} | {a,c,d,e} |
| a | a | {1,b,d,e} | {a,b,c} | {1,a,c,e} | {1,b,c,e} | {0,a,b,d,e} | {1,a,c,d} |
| b | b | {a,b,d,e} | {1,a,c,e} | {b,c,d} | {1,a,b,d} | {1,a,c,d} | {0,1,b,c,e} |
| c | c | {0,1,a,c,d} | {1,b,c,e} | {1,a,b,d} | {c,d,e} | {a,b,c,e} | {a,b,d,e} |
| d | d | {1,b,c,e} | {0,a,b,d,e} | {1,a,c,d} | {a,b,c,e} | {1,d,e} | {1,b,c,d} |
| e | e | {a,c,d,e} | {1,a,c,d} | {0,1,b,c,e} | {a,b,d,e} | {1,b,c,d} | {1,a,e} |

| $HF_7^{14}$ | 0 | 1 | a | b | c | d | e |
|---|---|---|---|---|---|---|---|
| 0 | 0 | 1 | a | b | c | d | e |
| 1 | 1 | {1,a,c} | {1,b,e} | {a,b,d,e} | {0,1,a,c,d} | {1,b,c,e} | {a,d,e} |
| a | a | {1,b,e} | {a,b,d} | {1,a,c} | {1,b,c,e} | {0,a,b,d,e} | {1,a,c,d} |
| b | b | {a,b,d,e} | {1,a,c} | {b,c,e} | {a,b,d} | {1,a,c,d} | {0,1,b,c,e} |
| c | c | {0,1,a,c,d} | {1,b,c,e} | {a,b,d} | {1,c,d} | {b,c,e} | {a,b,d,e} |
| d | d | {1,b,c,e} | {0,a,b,d,e} | {1,a,c,d} | {b,c,e} | {a,d,e} | {1,c,d} |
| e | e | {a,d,e} | {1,a,c,d} | {0,1,b,c,e} | {a,b,d,e} | {1,c,d} | {1,b,e} |

| $HF_7^{15}$ | 0 | 1 | a | b | c | d | e |
|---|---|---|---|---|---|---|---|
| 0 | 0 | 1 | a | b | c | d | e |
| 1 | 1 | {1,a,d} | {1,b,c,e} | {b,d,e} | {0,1,a,c,d} | {1,b,c} | {a,b,d,e} |
| a | a | {1,b,c,e} | {a,b,e} | {1,a,c,d} | {1,c,e} | {0,a,b,d,e} | {a,c,d} |
| b | b | {b,d,e} | {1,a,c,d} | {1,b,c} | {a,b,d,e} | {1,a,d} | {0,1,b,c,e} |
| c | c | {0,1,a,c,d} | {1,c,e} | {a,b,d,e} | {a,c,d} | {1,b,c,e} | {a,b,e} |
| d | d | {1,b,c} | {0,a,b,d,e} | {1,a,d} | {1,b,c,e} | {b,d,e} | {1,a,c,d} |
| e | e | {a,b,d,e} | {a,c,d} | {0,1,b,c,e} | {a,b,e} | {1,a,c,d} | {1,c,e} |

| $HF_7^{16}$ | 0 | 1 | a | b | c | d | e |
|---|---|---|---|---|---|---|---|
| 0 | 0 | 1 | a | b | c | d | e |
| 1 | 1 | {1,a,d} | {1,b,c,e} | {a,b,d,e} | {0,1,a,c,d} | {1,b,c,e} | {a,b,d,e} |
| a | a | {1,b,c,e} | {a,b,e} | {1,a,c,d} | {1,b,c,e} | {0,a,b,d,e} | {1,a,c,d} |
| b | b | {a,b,d,e} | {1,a,c,d} | {1,b,c} | {a,b,d,e} | {1,a,c,d} | {0,1,b,c,e} |
| c | c | {0,1,a,c,d} | {1,b,c,e} | {a,b,d,e} | {a,c,d} | {1,b,c,e} | {a,b,d,e} |
| d | d | {1,b,c,e} | {0,a,b,d,e} | {1,a,c,d} | {1,b,c,e} | {b,d,e} | {1,a,c,d} |
| e | e | {a,b,d,e} | {1,a,c,d} | {0,1,b,c,e} | {a,b,d,e} | {1,a,c,d} | {1,c,e} |

| $HF_7^{17}$ | 0 | 1 | a | b | c | d | e |
|---|---|---|---|---|---|---|---|
| 0 | 0 | 1 | a | b | c | d | e |
| 1 | 1 | {1,a,e} | {1,b,e} | {b,c,d,e} | {0,1,a,c,d} | {1,a,b,c} | {a,d,e} |
| a | a | {1,b,e} | {1,a,b} | {1,a,c} | {1,c,d,e} | {0,a,b,d,e} | {a,b,c,d} |
| b | b | {b,c,d,e} | {1,a,c} | {a,b,c} | {a,b,d} | {1,a,d,e} | {0,1,b,c,e} |
| c | c | {0,1,a,c,d} | {1,c,d,e} | {a,b,d} | {b,c,d} | {b,c,e} | {1,a,b,e} |
| d | d | {1,a,b,c} | {0,a,b,d,e} | {1,a,d,e} | {b,c,e} | {c,d,e} | {1,c,d} |
| e | e | {a,d,e} | {a,b,c,d} | {0,1,b,c,e} | {1,a,b,e} | {1,c,d} | {1,d,e} |

| $HF_7^{18}$ | 0 | 1 | a | b | c | d | e |
|---|---|---|---|---|---|---|---|
| 0 | 0 | 1 | a | b | c | d | e |
| 1 | 1 | {1,a,e} | {1,b,e} | {a,b,c,d,e} | {0,1,a,c,d} | {1,a,b,c,e} | {a,d,e} |
| a | a | {1,b,e} | {1,a,b} | {1,a,c} | {1,b,c,d,e} | {0,a,b,d,e} | {1,a,b,c,d} |
| b | b | {a,b,c,d,e} | {1,a,c} | {a,b,c} | {a,b,d} | {1,a,c,d,e} | {0,1,b,c,e} |
| c | c | {0,1,a,c,d} | {1,b,c,d,e} | {a,b,d} | {b,c,d} | {b,c,e} | {1,a,b,d,e} |
| d | d | {1,a,b,c,e} | {0,a,b,d,e} | {1,a,c,d,e} | {b,c,e} | {c,d,e} | {1,c,d} |
| e | e | {a,d,e} | {1,a,b,c,d} | {0,1,b,c,e} | {1,a,b,d,e} | {1,c,d} | {1,d,e} |

| $HF_7^{19}$ | 0 | 1 | a | b | c | d | e |
|---|---|---|---|---|---|---|---|
| 0 | 0 | 1 | a | b | c | d | e |
| 1 | 1 | {1,b,c} | {1,b,d,e} | {a,b,d} | {0,1,a,c,d} | {1,b,e} | {a,c,d,e} |
| a | a | {1,b,d,e} | {a,c,d} | {1,a,c,e} | {b,c,e} | {0,a,b,d,e} | {1,a,c} |
| b | b | {a,b,d} | {1,a,c,e} | {b,d,e} | {1,a,b,d} | {1,c,d} | {0,1,b,c,e} |
| c | c | {0,1,a,c,d} | {b,c,e} | {1,a,b,d} | {1,c,e} | {a,b,c,e} | {a,d,e} |
| d | d | {1,b,e} | {0,a,b,d,e} | {1,c,d} | {a,b,c,e} | {1,a,d} | {1,b,c,d} |
| e | e | {a,c,d,e} | {1,a,c} | {0,1,b,c,e} | {a,d,e} | {1,b,c,d} | {a,b,e} |

| $HF_7^{20}$ | 0 | 1 | a | b | c | d | e |
|---|---|---|---|---|---|---|---|
| 0 | 0 | 1 | a | b | c | d | e |
| 1 | 1 | {1,b,d} | {1,b,c,d,e} | {b,d} | {0,1,a,c,d} | {1,b} | {a,b,c,d,e} |
| a | a | {1,b,c,d,e} | {a,c,e} | {1,a,c,d,e} | {c,e} | {0,a,b,d,e} | {a,c} |
| b | b | {b,d} | {1,a,c,d,e} | {1,b,d} | {1,a,b,d,e} | {1,d} | {0,1,b,c,e} |
| c | c | {0,1,a,c,d} | {c,e} | {1,a,b,d,e} | {a,c,e} | {1,a,b,c,e} | {a,e} |
| d | d | {1,b} | {0,a,b,d,e} | {1,d} | {1,a,b,c,e} | {1,b,d} | {1,a,b,c,d} |
| e | e | {a,b,c,d,e} | {a,c} | {0,1,b,c,e} | {a,e} | {1,a,b,c,d} | {a,c,e} |

| $HF_7^{21}$ | 0 | 1 | a | b | c | d | e |
|---|---|---|---|---|---|---|---|
| 0 | 0 | 1 | a | b | c | d | e |
| 1 | 1 | {1,a,b,c} | {1,b,d,e} | {b,d,e} | {0,1,a,c,d} | {1,b,c} | {a,c,d,e} |
| a | a | {1,b,d,e} | {a,b,c,d} | {1,a,c,e} | {1,c,e} | {0,a,b,d,e} | {a,c,d} |
| b | b | {b,d,e} | {1,a,c,e} | {b,c,d,e} | {1,a,b,d} | {1,a,d} | {0,1,b,c,e} |
| c | c | {0,1,a,c,d} | {1,c,e} | {1,a,b,d} | {1,c,d,e} | {a,b,c,e} | {a,b,e} |
| d | d | {1,b,c} | {0,a,b,d,e} | {1,a,d} | {a,b,c,e} | {1,a,d,e} | {1,b,c,d} |
| e | e | {a,c,d,e} | {a,c,d} | {0,1,b,c,e} | {a,b,e} | {1,b,c,d} | {1,a,b,e} |

| $HF_7^{22}$ | 0 | 1 | a | b | c | d | e |
|---|---|---|---|---|---|---|---|
| 0 | 0 | 1 | a | b | c | d | e |
| 1 | 1 | {1,a,b,c} | {1,b,d,e} | {a,b,d,e} | {0,1,a,c,d} | {1,b,c,e} | {a,c,d,e} |
| a | a | {1,b,d,e} | {a,b,c,d} | {1,a,c,e} | {1,b,c,e} | {0,a,b,d,e} | {1,a,c,d} |
| b | b | {a,b,d,e} | {1,a,c,e} | {b,c,d,e} | {1,a,b,d} | {1,a,c,d} | {0,1,b,c,e} |
| c | c | {0,1,a,c,d} | {1,b,c,e} | {1,a,b,d} | {1,c,d,e} | {a,b,c,e} | {a,b,d,e} |
| d | d | {1,b,c,e} | {0,a,b,d,e} | {1,a,c,d} | {a,b,c,e} | {1,a,d,e} | {1,b,c,d} |
| e | e | {a,c,d,e} | {1,a,c,d} | {0,1,b,c,e} | {a,b,d,e} | {1,b,c,d} | {1,a,b,e} |

| $HF_7^{23}$ | 0 | 1 | a | b | c | d | e |
|---|---|---|---|---|---|---|---|
| 0 | 0 | 1 | a | b | c | d | e |
| 1 | 1 | {1,a,b,d} | {1,b,c,d,e} | {b,d,e} | {0,1,a,c,d} | {1,b,c} | {a,b,c,d,e} |
| a | a | {1,b,c,d,e} | {a,b,c,e} | {1,a,c,d,e} | {1,c,e} | {0,a,b,d,e} | {a,c,d} |
| b | b | {b,d,e} | {1,a,c,d,e} | {1,b,c,d} | {1,a,b,d,e} | {1,a,d} | {0,1,b,c,e} |
| c | c | {0,1,a,c,d} | {1,c,e} | {1,a,b,d,e} | {a,c,d,e} | {1,a,b,c,e} | {a,b,e} |
| d | d | {1,b,c} | {0,a,b,d,e} | {1,a,d} | {1,a,b,c,e} | {1,b,d,e} | {1,a,b,c,d} |
| e | e | {a,b,c,d,e} | {a,c,d} | {0,1,b,c,e} | {a,b,e} | {1,a,b,c,d} | {1,a,c,e} |

| $HF_7^{24}$ | 0 | 1 | a | b | c | d | e |
|---|---|---|---|---|---|---|---|
| 0 | 0 | 1 | a | b | c | d | e |
| 1 | 1 | {1,a,b,d} | {1,b,c,d,e} | {a,b,d,e} | {0,1,a,c,d} | {1,b,c,e} | {a,b,c,d,e} |
| a | a | {1,b,c,d,e} | {a,b,c,e} | {1,a,c,d,e} | {1,b,c,e} | {0,a,b,d,e} | {1,a,c,d} |
| b | b | {a,b,d,e} | {1,a,c,d,e} | {1,b,c,d} | {1,a,b,d,e} | {1,a,c,d} | {0,1,b,c,e} |
| c | c | {0,1,a,c,d} | {1,b,c,e} | {1,a,b,d,e} | {a,c,d,e} | {1,a,b,c,e} | {a,b,d,e} |
| d | d | {1,b,c,e} | {0,a,b,d,e} | {1,a,c,d} | {1,a,b,c,e} | {1,b,d,e} | {1,a,b,c,d} |
| e | e | {a,b,c,d,e} | {1,a,c,d} | {0,1,b,c,e} | {a,b,d,e} | {1,a,b,c,d} | {1,a,c,e} |

| $HF_7^{25}$ | 0 | 1 | a | b | c | d | e |
|---|---|---|---|---|---|---|---|
| 0 | 0 | 1 | a | b | c | d | e |
| 1 | 1 | {1,a,b,e} | {1,b,d,e} | {b,c,d,e} | {0,1,a,c,d} | {1,a,b,c} | {a,c,d,e} |
| a | a | {1,b,d,e} | {1,a,b,c} | {1,a,c,e} | {1,c,d,e} | {0,a,b,d,e} | {a,b,c,d} |
| b | b | {b,c,d,e} | {1,a,c,e} | {a,b,c,d} | {1,a,b,d} | {1,a,d,e} | {0,1,b,c,e} |
| c | c | {0,1,a,c,d} | {1,c,d,e} | {1,a,b,d} | {b,c,d,e} | {a,b,c,e} | {1,a,b,e} |
| d | d | {1,a,b,c} | {0,a,b,d,e} | {1,a,d,e} | {a,b,c,e} | {1,c,d,e} | {1,b,c,d} |
| e | e | {a,c,d,e} | {a,b,c,d} | {0,1,b,c,e} | {1,a,b,e} | {1,b,c,d} | {1,a,d,e} |

| $HF_7^{26}$ | 0 | 1 | a | b | c | d | e |
|---|---|---|---|---|---|---|---|
| 0 | 0 | 1 | a | b | c | d | e |
| 1 | 1 | {1,a,b,e} | {1,b,d,e} | {a,b,c,d,e} | {0,1,a,c,d} | {1,a,b,c,e} | {a,c,d,e} |
| a | a | {1,b,d,e} | {1,a,b,c} | {1,a,c,e} | {1,b,c,d,e} | {0,a,b,d,e} | {1,a,b,c,d} |
| b | b | {a,b,c,d,e} | {1,a,c,e} | {a,b,c,d} | {1,a,b,d} | {1,a,c,d,e} | {0,1,b,c,e} |
| c | c | {0,1,a,c,d} | {1,b,c,d,e} | {1,a,b,d} | {b,c,d,e} | {a,b,c,e} | {1,a,b,d,e} |
| d | d | {1,a,b,c,e} | {0,a,b,d,e} | {1,a,c,d,e} | {a,b,c,e} | {1,c,d,e} | {1,b,c,d} |
| e | e | {a,c,d,e} | {1,a,b,c,d} | {0,1,b,c,e} | {1,a,b,d,e} | {1,b,c,d} | {1,a,d,e} |

| $HF_7^{27}$ | 0 | 1 | a | b | c | d | e |
|---|---|---|---|---|---|---|---|
| 0 | 0 | 1 | a | b | c | d | e |
| 1 | 1 | {1,a,c,d} | {1,b,c,e} | {b,d,e} | {0,1,a,c,d} | {1,b,c} | {a,b,d,e} |
| a | a | {1,b,c,e} | {a,b,d,e} | {1,a,c,d} | {1,c,e} | {0,a,b,d,e} | {a,c,d} |
| b | b | {b,d,e} | {1,a,c,d} | {1,b,c,e} | {a,b,d,e} | {1,a,d} | {0,1,b,c,e} |
| c | c | {0,1,a,c,d} | {1,c,e} | {a,b,d,e} | {1,a,c,d} | {1,b,c,e} | {a,b,e} |
| d | d | {1,b,c} | {0,a,b,d,e} | {1,a,d} | {1,b,c,e} | {a,b,d,e} | {1,a,c,d} |
| e | e | {a,b,d,e} | {a,c,d} | {0,1,b,c,e} | {a,b,e} | {1,a,c,d} | {1,b,c,e} |

| $HF_7^{28}$ | 0 | 1 | a | b | c | d | e |
|---|---|---|---|---|---|---|---|
| 0 | 0 | 1 | a | b | c | d | e |
| 1 | 1 | {1,a,c,d} | {1,b,c,e} | {a,b,d,e} | {0,1,a,c,d} | {1,b,c,e} | {a,b,d,e} |
| a | a | {1,b,c,e} | {a,b,d,e} | {1,a,c,d} | {1,b,c,e} | {0,a,b,d,e} | {1,a,c,d} |
| b | b | {a,b,d,e} | {1,a,c,d} | {1,b,c,e} | {a,b,d,e} | {1,a,c,d} | {0,1,b,c,e} |
| c | c | {0,1,a,c,d} | {1,b,c,e} | {a,b,d,e} | {1,a,c,d} | {1,b,c,e} | {a,b,d,e} |
| d | d | {1,b,c,e} | {0,a,b,d,e} | {1,a,c,d} | {1,b,c,e} | {a,b,d,e} | {1,a,c,d} |
| e | e | {a,b,d,e} | {1,a,c,d} | {0,1,b,c,e} | {a,b,d,e} | {1,a,c,d} | {1,b,c,e} |

| $HF_7^{29}$ | 0 | 1 | a | b | c | d | e |
|---|---|---|---|---|---|---|---|
| 0 | 0 | 1 | a | b | c | d | e |
| 1 | 1 | {1,a,c,e} | {1,b,e} | {b,c,d,e} | {0,1,a,c,d} | {1,a,b,c} | {a,d,e} |
| a | a | {1,b,e} | {1,a,b,d} | {1,a,c} | {1,c,d,e} | {0,a,b,d,e} | {a,b,c,d} |
| b | b | {b,c,d,e} | {1,a,c} | {a,b,c,e} | {a,b,d} | {1,a,d,e} | {0,1,b,c,e} |
| c | c | {0,1,a,c,d} | {1,c,d,e} | {a,b,d} | {1,b,c,d} | {b,c,e} | {1,a,b,e} |
| d | d | {1,a,b,c} | {0,a,b,d,e} | {1,a,d,e} | {b,c,e} | {a,c,d,e} | {1,c,d} |
| e | e | {a,d,e} | {a,b,c,d} | {0,1,b,c,e} | {1,a,b,e} | {1,c,d} | {1,b,d,e} |

| $HF_7^{30}$ | 0 | 1 | a | b | c | d | e |
|---|---|---|---|---|---|---|---|
| 0 | 0 | 1 | a | b | c | d | e |
| 1 | 1 | {1,a,c,e} | {1,b,e} | {a,b,c,d,e} | {0,1,a,c,d} | {1,a,b,c,e} | {a,d,e} |
| a | a | {1,b,e} | {1,a,b,d} | {1,a,c} | {1,b,c,d,e} | {0,a,b,d,e} | {1,a,b,c,d} |
| b | b | {a,b,c,d,e} | {1,a,c} | {a,b,c,e} | {a,b,d} | {1,a,c,d,e} | {0,1,b,c,e} |
| c | c | {0,1,a,c,d} | {1,b,c,d,e} | {a,b,d} | {1,b,c,d} | {b,c,e} | {1,a,b,d,e} |
| d | d | {1,a,b,c,e} | {0,a,b,d,e} | {1,a,c,d,e} | {b,c,e} | {a,c,d,e} | {1,c,d} |
| e | e | {a,d,e} | {1,a,b,c,d} | {0,1,b,c,e} | {1,a,b,d,e} | {1,c,d} | {1,b,d,e} |

| $HF_7^{31}$ | 0 | 1 | a | b | c | d | e |
|---|---|---|---|---|---|---|---|
| 0 | 0 | 1 | a | b | c | d | e |
| 1 | 1 | {1,b,c,d} | {1,b,c,d,e} | {a,b,d} | {0,1,a,c,d} | {1,b,e} | {a,b,c,d,e} |
| a | a | {1,b,c,d,e} | {a,c,d,e} | {1,a,c,d,e} | {b,c,e} | {0,a,b,d,e} | {1,a,c} |
| b | b | {a,b,d} | {1,a,c,d,e} | {1,b,d,e} | {1,a,b,d,e} | {1,c,d} | {0,1,b,c,e} |
| c | c | {0,1,a,c,d} | {b,c,e} | {1,a,b,d,e} | {1,a,c,e} | {1,a,b,c,e} | {a,d,e} |
| d | d | {1,b,e} | {0,a,b,d,e} | {1,c,d} | {1,a,b,c,e} | {1,a,b,d} | {1,a,b,c,d} |
| e | e | {a,b,c,d,e} | {1,a,c} | {0,1,b,c,e} | {a,d,e} | {1,a,b,c,d} | {a,b,c,e} |

| $HF_7^{32}$ | 0 | 1 | a | b | c | d | e |
|---|---|---|---|---|---|---|---|
| 0 | 0 | 1 | a | b | c | d | e |
| 1 | 1 | {1,a,b,c,d} | {1,b,c,d,e} | {b,d,e} | {0,1,a,c,d} | {1,b,c} | {a,b,c,d,e} |
| a | a | {1,b,c,d,e} | {a,b,c,d,e} | {1,a,c,d,e} | {1,c,e} | {0,a,b,d,e} | {a,c,d} |
| b | b | {b,d,e} | {1,a,c,d,e} | {1,b,c,d,e} | {1,a,b,d,e} | {1,a,d} | {0,1,b,c,e} |
| c | c | {0,1,a,c,d} | {1,c,e} | {1,a,b,d,e} | {1,a,c,d,e} | {1,a,b,c,e} | {a,b,e} |
| d | d | {1,b,c} | {0,a,b,d,e} | {1,a,d} | {1,a,b,c,e} | {1,a,b,d,e} | {1,a,b,c,d} |
| e | e | {a,b,c,d,e} | {a,c,d} | {0,1,b,c,e} | {a,b,e} | {1,a,b,c,d} | {1,a,b,c,e} |

| $HF_7^{33}$ | 0 | 1 | a | b | c | d | e |
|---|---|---|---|---|---|---|---|
| 0 | 0 | 1 | a | b | c | d | e |
| 1 | 1 | {1,a,b,c,d} | {1,b,c,d,e} | {a,b,d,e} | {0,1,a,c,d} | {1,b,c,e} | {a,b,c,d,e} |
| a | a | {1,b,c,d,e} | {a,b,c,d,e} | {1,a,c,d,e} | {1,b,c,e} | {0,a,b,d,e} | {1,a,c,d} |
| b | b | {a,b,d,e} | {1,a,c,d,e} | {1,b,c,d,e} | {1,a,b,d,e} | {1,a,c,d} | {0,1,b,c,e} |
| c | c | {0,1,a,c,d} | {1,b,c,e} | {1,a,b,d,e} | {1,a,c,d,e} | {1,a,b,c,e} | {a,b,d,e} |
| d | d | {1,b,c,e} | {0,a,b,d,e} | {1,a,c,d} | {1,a,b,c,e} | {1,a,b,d,e} | {1,a,b,c,d} |
| e | e | {a,b,c,d,e} | {1,a,c,d} | {0,1,b,c,e} | {a,b,d,e} | {1,a,b,c,d} | {1,a,b,c,e} |

| $HF_7^{34}$ | 0 | 1 | a | b | c | d | e |
|---|---|---|---|---|---|---|---|
| 0 | 0 | 1 | a | b | c | d | e |
| 1 | 1 | {1,a,b,c,e} | {1,b,d,e} | {b,c,d,e} | {0,1,a,c,d} | {1,a,b,c} | {a,c,d,e} |
| a | a | {1,b,d,e} | {1,a,b,c,d} | {1,a,c,e} | {1,c,d,e} | {0,a,b,d,e} | {a,b,c,d} |
| b | b | {b,c,d,e} | {1,a,c,e} | {a,b,c,d,e} | {1,a,b,d} | {1,a,d,e} | {0,1,b,c,e} |
| c | c | {0,1,a,c,d} | {1,c,d,e} | {1,a,b,d} | {1,b,c,d,e} | {a,b,c,e} | {1,a,b,e} |
| d | d | {1,a,b,c} | {0,a,b,d,e} | {1,a,d,e} | {a,b,c,e} | {1,a,c,d,e} | {1,b,c,d} |
| e | e | {a,c,d,e} | {a,b,c,d} | {0,1,b,c,e} | {1,a,b,e} | {1,b,c,d} | {1,a,b,d,e} |

| $HF_7^{35}$ | 0 | 1 | a | b | c | d | e |
|---|---|---|---|---|---|---|---|
| 0 | 0 | 1 | a | b | c | d | e |
| 1 | 1 | {1,a,b,c,e} | {1,b,d,e} | {a,b,c,d,e} | {0,1,a,c,d} | {1,a,b,c,e} | {a,c,d,e} |
| a | a | {1,b,d,e} | {1,a,b,c,d} | {1,a,c,e} | {1,b,c,d,e} | {0,a,b,d,e} | {1,a,b,c,d} |
| b | b | {a,b,c,d,e} | {1,a,c,e} | {a,b,c,d,e} | {1,a,b,d} | {1,a,c,d,e} | {0,1,b,c,e} |
| c | c | {0,1,a,c,d} | {1,b,c,d,e} | {1,a,b,d} | {1,b,c,d,e} | {a,b,c,e} | {1,a,b,d,e} |
| d | d | {1,a,b,c,e} | {0,a,b,d,e} | {1,a,c,d,e} | {a,b,c,e} | {1,a,c,d,e} | {1,b,c,d} |
| e | e | {a,c,d,e} | {1,a,b,c,d} | {0,1,b,c,e} | {1,a,b,d,e} | {1,b,c,d} | {1,a,b,d,e} |

| $HF_7^{36}$ | 0 | 1 | a | b | c | d | e |
|---|---|---|---|---|---|---|---|
| 0 | 0 | 1 | a | b | c | d | e |
| 1 | 1 | {1,a,b,d,e} | {1,c,d} | {b,c,e} | {0,1,a,c,d} | {1,a,c} | {b,c,e} |
| a | a | {1,c,d} | {1,a,b,c,e} | {a,d,e} | {1,c,d} | {0,a,b,d,e} | {a,b,d} |
| b | b | {b,c,e} | {a,d,e} | {1,a,b,c,d} | {1,b,e} | {a,d,e} | {0,1,b,c,e} |
| c | c | {0,1,a,c,d} | {1,c,d} | {1,b,e} | {a,b,c,d,e} | {1,a,c} | {1,b,e} |
| d | d | {1,a,c} | {0,a,b,d,e} | {a,d,e} | {1,a,c} | {1,b,c,d,e} | {a,b,d} |
| e | e | {b,c,e} | {a,b,d} | {0,1,b,c,e} | {1,b,e} | {a,b,d} | {1,a,c,d,e} |

| $HF_7^{37}$ | 0 | 1 | a | b | c | d | e |
|---|---|---|---|---|---|---|---|
| 0 | 0 | 1 | a | b | c | d | e |
| 1 | 1 | {1,a,b,d,e} | {1,c,d} | {a,b,c,e} | {0,1,a,c,d} | {1,a,c,e} | {b,c,e} |
| a | a | {1,c,d} | {1,a,b,c,e} | {a,d,e} | {1,b,c,d} | {0,a,b,d,e} | {1,a,b,d} |
| b | b | {a,b,c,e} | {a,d,e} | {1,a,b,c,d} | {1,b,e} | {a,c,d,e} | {0,1,b,c,e} |
| c | c | {0,1,a,c,d} | {1,b,c,d} | {1,b,e} | {a,b,c,d,e} | {1,a,c} | {1,b,d,e} |
| d | d | {1,a,c,e} | {0,a,b,d,e} | {a,c,d,e} | {1,a,c} | {1,b,c,d,e} | {a,b,d} |
| e | e | {b,c,e} | {1,a,b,d} | {0,1,b,c,e} | {1,b,d,e} | {a,b,d} | {1,a,c,d,e} |

| $HF_7^{38}$ | 0 | 1 | a | b | c | d | e |
|---|---|---|---|---|---|---|---|
| 0 | 0 | 1 | a | b | c | d | e |
| 1 | 1 | {1,a,b,d,e} | {1,b,c,d,e} | {b,c,d,e} | {0,1,a,c,d} | {1,a,b,c} | {a,b,c,d,e} |
| a | a | {1,b,c,d,e} | {1,a,b,c,e} | {1,a,c,d,e} | {1,c,d,e} | {0,a,b,d,e} | {a,b,c,d} |
| b | b | {b,c,d,e} | {1,a,c,d,e} | {1,a,b,c,d} | {1,a,b,d,e} | {1,a,d,e} | {0,1,b,c,e} |
| c | c | {0,1,a,c,d} | {1,c,d,e} | {1,a,b,d,e} | {a,b,c,d,e} | {1,a,b,c,e} | {1,a,b,e} |
| d | d | {1,a,b,c} | {0,a,b,d,e} | {1,a,d,e} | {1,a,b,c,e} | {1,b,c,d,e} | {1,a,b,c,d} |
| e | e | {a,b,c,d,e} | {a,b,c,d} | {0,1,b,c,e} | {1,a,b,e} | {1,a,b,c,d} | {1,a,c,d,e} |

| $HF_7^{39}$ | 0 | 1 | a | b | c | d | e |
|---|---|---|---|---|---|---|---|
| 0 | 0 | 1 | a | b | c | d | e |
| 1 | 1 | {1,a,b,d,e} | {1,b,c,d,e} | {a,b,c,d,e} | {0,1,a,c,d} | {1,a,b,c,e} | {a,b,c,d,e} |
| a | a | {1,b,c,d,e} | {1,a,b,c,e} | {1,a,c,d,e} | {1,b,c,d,e} | {0,a,b,d,e} | {1,a,b,c,d} |
| b | b | {a,b,c,d,e} | {1,a,c,d,e} | {1,a,b,c,d} | {1,a,b,d,e} | {1,a,c,d,e} | {0,1,b,c,e} |
| c | c | {0,1,a,c,d} | {1,b,c,d,e} | {1,a,b,d,e} | {a,b,c,d,e} | {1,a,b,c,e} | {1,a,b,d,e} |
| d | d | {1,a,b,c,e} | {0,a,b,d,e} | {1,a,c,d,e} | {1,a,b,c,e} | {1,b,c,d,e} | {1,a,b,c,d} |
| e | e | {a,b,c,d,e} | {1,a,b,c,d} | {0,1,b,c,e} | {1,a,b,d,e} | {1,a,b,c,d} | {1,a,c,d,e} |

| $HF_7^{40}$ | 0 | 1 | a | b | c | d | e |
|---|---|---|---|---|---|---|---|
| 0 | 0 | 1 | a | b | c | d | e |
| 1 | 1 | {1,a,b,c,d,e} | {1,c,d} | {b,c,e} | {0,1,a,c,d} | {1,a,c} | {b,c,e} |
| a | a | {1,c,d} | {1,a,b,c,d,e} | {a,d,e} | {1,c,d} | {0,a,b,d,e} | {a,b,d} |
| b | b | {b,c,e} | {a,d,e} | {1,a,b,c,d,e} | {1,b,e} | {a,d,e} | {0,1,b,c,e} |
| c | c | {0,1,a,c,d} | {1,c,d} | {1,b,e} | {1,a,b,c,d,e} | {1,a,c} | {1,b,e} |
| d | d | {1,a,c} | {0,a,b,d,e} | {a,d,e} | {1,a,c} | {1,a,b,c,d,e} | {a,b,d} |
| e | e | {b,c,e} | {a,b,d} | {0,1,b,c,e} | {1,b,e} | {a,b,d} | {1,a,b,c,d,e} |

| $HF_7^{41}$ | 0 | 1 | a | b | c | d | e |
|---|---|---|---|---|---|---|---|
| 0 | 0 | 1 | a | b | c | d | e |
| 1 | 1 | {1,a,b,c,d,e} | {1,c,d} | {a,b,c,e} | {0,1,a,c,d} | {1,a,c,e} | {b,c,e} |
| a | a | {1,c,d} | {1,a,b,c,d,e} | {a,d,e} | {1,b,c,d} | {0,a,b,d,e} | {1,a,b,d} |
| b | b | {a,b,c,e} | {a,d,e} | {1,a,b,c,d,e} | {1,b,e} | {a,c,d,e} | {0,1,b,c,e} |
| c | c | {0,1,a,c,d} | {1,b,c,d} | {1,b,e} | {1,a,b,c,d,e} | {1,a,c} | {1,b,d,e} |
| d | d | {1,a,c,e} | {0,a,b,d,e} | {a,c,d,e} | {1,a,c} | {1,a,b,c,d,e} | {a,b,d} |
| e | e | {b,c,e} | {1,a,b,d} | {0,1,b,c,e} | {1,b,d,e} | {a,b,d} | {1,a,b,c,d,e} |

| $HF_7^{42}$ | 0 | 1 | a | b | c | d | e |
|---|---|---|---|---|---|---|---|
| 0 | 0 | 1 | a | b | c | d | e |
| 1 | 1 | {1,a,b,c,d,e} | {1,b,c,d,e} | {b,c,d,e} | {0,1,a,c,d} | {1,a,b,c} | {a,b,c,d,e} |
| a | a | {1,b,c,d,e} | {1,a,b,c,d,e} | {1,a,c,d,e} | {1,c,d,e} | {0,a,b,d,e} | {a,b,c,d} |
| b | b | {b,c,d,e} | {1,a,c,d,e} | {1,a,b,c,d,e} | {1,a,b,d,e} | {1,a,d,e} | {0,1,b,c,e} |
| c | c | {0,1,a,c,d} | {1,c,d,e} | {1,a,b,d,e} | {1,a,b,c,d,e} | {1,a,b,c,e} | {1,a,b,e} |
| d | d | {1,a,b,c} | {0,a,b,d,e} | {1,a,d,e} | {1,a,b,c,e} | {1,a,b,c,d,e} | {1,a,b,c,d} |
| e | e | {a,b,c,d,e} | {a,b,c,d} | {0,1,b,c,e} | {1,a,b,e} | {1,a,b,c,d} | {1,a,b,c,d,e} |

| $HF_7^{43}$ | 0 | 1 | a | b | c | d | e |
|---|---|---|---|---|---|---|---|
| 0 | 0 | 1 | a | b | c | d | e |
| 1 | 1 | {1,a,b,c,d,e} | {1,b,c,d,e} | {a,b,c,d,e} | {0,1,a,c,d} | {1,a,b,c,e} | {a,b,c,d,e} |
| a | a | {1,b,c,d,e} | {1,a,b,c,d,e} | {1,a,c,d,e} | {1,b,c,d,e} | {0,a,b,d,e} | {1,a,b,c,d} |
| b | b | {a,b,c,d,e} | {1,a,c,d,e} | {1,a,b,c,d,e} | {1,a,b,d,e} | {1,a,c,d,e} | {0,1,b,c,e} |
| c | c | {0,1,a,c,d} | {1,b,c,d,e} | {1,a,b,d,e} | {1,a,b,c,d,e} | {1,a,b,c,e} | {1,a,b,d,e} |
| d | d | {1,a,b,c,e} | {0,a,b,d,e} | {1,a,c,d,e} | {1,a,b,c,e} | {1,a,b,c,d,e} | {1,a,b,c,d} |
| e | e | {a,b,c,d,e} | {1,a,b,c,d} | {0,1,b,c,e} | {1,a,b,d,e} | {1,a,b,c,d} | {1,a,b,c,d,e} |

| $HF_7^{44}$ | 0 | 1 | a | b | c | d | e |
|---|---|---|---|---|---|---|---|
| 0 | 0 | 1 | a | b | c | d | e |
| 1 | 1 | {1,b} | {a,b,d,e} | {1,d} | {0,1,b,c,e} | {b,d} | {1,a,c,d} |
| a | a | {a,b,d,e} | {a,c} | {1,b,c,e} | {a,e} | {0,1,a,c,d} | {c,e} |
| b | b | {1,d} | {1,b,c,e} | {b,d} | {1,a,c,d} | {1,b} | {0,a,b,d,e} |
| c | c | {0,1,b,c,e} | {a,e} | {1,a,c,d} | {c,e} | {a,b,d,e} | {a,c} |
| d | d | {b,d} | {0,1,a,c,d} | {1,b} | {a,b,d,e} | {1,d} | {1,b,c,e} |
| e | e | {1,a,c,d} | {c,e} | {0,a,b,d,e} | {a,c} | {1,b,c,e} | {a,e} |

| $HF_7^{45}$ | 0 | 1 | a | b | c | d | e |
|---|---|---|---|---|---|---|---|
| 0 | 0 | 1 | a | b | c | d | e |
| 1 | 1 | {1,a,b} | {a,b,d,e} | {1,a,d,e} | {0,1,b,c,e} | {b,c,d,e} | {1,a,c,d} |
| a | a | {a,b,d,e} | {a,b,c} | {1,b,c,e} | {1,a,b,e} | {0,1,a,c,d} | {1,c,d,e} |
| b | b | {1,a,d,e} | {1,b,c,e} | {b,c,d} | {1,a,c,d} | {1,a,b,c} | {0,a,b,d,e} |
| c | c | {0,1,b,c,e} | {1,a,b,e} | {1,a,c,d} | {c,d,e} | {a,b,d,e} | {a,b,c,d} |
| d | d | {b,c,d,e} | {0,1,a,c,d} | {1,a,b,c} | {a,b,d,e} | {1,d,e} | {1,b,c,e} |
| e | e | {1,a,c,d} | {1,c,d,e} | {0,a,b,d,e} | {a,b,c,d} | {1,b,c,e} | {1,a,e} |

| $HF_7^{46}$ | 0 | 1 | a | b | c | d | e |
|---|---|---|---|---|---|---|---|
| 0 | 0 | 1 | a | b | c | d | e |
| 1 | 1 | {1,a,c} | {a,b,e} | {1,a,d,e} | {0,1,b,c,e} | {b,c,d,e} | {1,a,d} |
| a | a | {a,b,e} | {a,b,d} | {1,b,c} | {1,a,b,e} | {0,1,a,c,d} | {1,c,d,e} |
| b | b | {1,a,d,e} | {1,b,c} | {b,c,e} | {a,c,d} | {1,a,b,c} | {0,a,b,d,e} |
| c | c | {0,1,b,c,e} | {1,a,b,e} | {a,c,d} | {1,c,d} | {b,d,e} | {a,b,c,d} |
| d | d | {b,c,d,e} | {0,1,a,c,d} | {1,a,b,c} | {b,d,e} | {a,d,e} | {1,c,e} |
| e | e | {1,a,d} | {1,c,d,e} | {0,a,b,d,e} | {a,b,c,d} | {1,c,e} | {1,b,e} |

| $HF_7^{47}$ | 0 | 1 | a | b | c | d | e |
|---|---|---|---|---|---|---|---|
| 0 | 0 | 1 | a | b | c | d | e |
| 1 | 1 | {1,a,d} | {a,b,c,e} | {1,a,d,e} | {0,1,b,c,e} | {b,c,d,e} | {1,a,b,d} |
| a | a | {a,b,c,e} | {a,b,e} | {1,b,c,d} | {1,a,b,e} | {0,1,a,c,d} | {1,c,d,e} |
| b | b | {1,a,d,e} | {1,b,c,d} | {1,b,c} | {a,c,d,e} | {1,a,b,c} | {0,a,b,d,e} |
| c | c | {0,1,b,c,e} | {1,a,b,e} | {a,c,d,e} | {a,c,d} | {1,b,d,e} | {a,b,c,d} |
| d | d | {b,c,d,e} | {0,1,a,c,d} | {1,a,b,c} | {1,b,d,e} | {b,d,e} | {1,a,c,e} |
| e | e | {1,a,b,d} | {1,c,d,e} | {0,a,b,d,e} | {a,b,c,d} | {1,a,c,e} | {1,c,e} |

| $HF_7^{48}$ | 0 | 1 | a | b | c | d | e |
|---|---|---|---|---|---|---|---|
| 0 | 0 | 1 | a | b | c | d | e |
| 1 | 1 | {1,b,c} | {a,b,d,e} | {1,a,d} | {0,1,b,c,e} | {b,d,e} | {1,a,c,d} |
| a | a | {a,b,d,e} | {a,c,d} | {1,b,c,e} | {a,b,e} | {0,1,a,c,d} | {1,c,e} |
| b | b | {1,a,d} | {1,b,c,e} | {b,d,e} | {1,a,c,d} | {1,b,c} | {0,a,b,d,e} |
| c | c | {0,1,b,c,e} | {a,b,e} | {1,a,c,d} | {1,c,e} | {a,b,d,e} | {a,c,d} |
| d | d | {b,d,e} | {0,1,a,c,d} | {1,b,c} | {a,b,d,e} | {1,a,d} | {1,b,c,e} |
| e | e | {1,a,c,d} | {1,c,e} | {0,a,b,d,e} | {a,c,d} | {1,b,c,e} | {a,b,e} |

| $HF_7^{49}$ | 0 | 1 | a | b | c | d | e |
|---|---|---|---|---|---|---|---|
| 0 | 0 | 1 | a | b | c | d | e |
| 1 | 1 | {1,a,b,c} | {a,b,d,e} | {1,a,d,e} | {0,1,b,c,e} | {b,c,d,e} | {1,a,c,d} |
| a | a | {a,b,d,e} | {a,b,c,d} | {1,b,c,e} | {1,a,b,e} | {0,1,a,c,d} | {1,c,d,e} |
| b | b | {1,a,d,e} | {1,b,c,e} | {b,c,d,e} | {1,a,c,d} | {1,a,b,c} | {0,a,b,d,e} |
| c | c | {0,1,b,c,e} | {1,a,b,e} | {1,a,c,d} | {1,c,d,e} | {a,b,d,e} | {a,b,c,d} |
| d | d | {b,c,d,e} | {0,1,a,c,d} | {1,a,b,c} | {a,b,d,e} | {1,a,d,e} | {1,b,c,e} |
| e | e | {1,a,c,d} | {1,c,d,e} | {0,a,b,d,e} | {a,b,c,d} | {1,b,c,e} | {1,a,b,e} |

| $HF_7^{50}$ | 0 | 1 | a | b | c | d | e |
|---|---|---|---|---|---|---|---|
| 0 | 0 | 1 | a | b | c | d | e |
| 1 | 1 | {1,a,b,d} | {a,b,c,d,e} | {1,a,d,e} | {0,1,b,c,e} | {b,c,d,e} | {1,a,b,c,d} |
| a | a | {a,b,c,d,e} | {a,b,c,e} | {1,b,c,d,e} | {1,a,b,e} | {0,1,a,c,d} | {1,c,d,e} |
| b | b | {1,a,d,e} | {1,b,c,d,e} | {1,b,c,d} | {1,a,c,d,e} | {1,a,b,c} | {0,a,b,d,e} |
| c | c | {0,1,b,c,e} | {1,a,b,e} | {1,a,c,d,e} | {a,c,d,e} | {1,a,b,d,e} | {a,b,c,d} |
| d | d | {b,c,d,e} | {0,1,a,c,d} | {1,a,b,c} | {1,a,b,d,e} | {1,b,d,e} | {1,a,b,c,e} |
| e | e | {1,a,b,c,d} | {1,c,d,e} | {0,a,b,d,e} | {a,b,c,d} | {1,a,b,c,e} | {1,a,c,e} |

| $HF_7^{51}$ | 0 | 1 | a | b | c | d | e |
|---|---|---|---|---|---|---|---|
| 0 | 0 | 1 | a | b | c | d | e |
| 1 | 1 | {1,a,b,e} | {a,b,d,e} | {1,c,d,e} | {0,1,b,c,e} | {a,b,c,d} | {1,a,c,d} |
| a | a | {a,b,d,e} | {1,a,b,c} | {1,b,c,e} | {1,a,d,e} | {0,1,a,c,d} | {b,c,d,e} |
| b | b | {1,c,d,e} | {1,b,c,e} | {a,b,c,d} | {1,a,c,d} | {1,a,b,e} | {0,a,b,d,e} |
| c | c | {0,1,b,c,e} | {1,a,d,e} | {1,a,c,d} | {b,c,d,e} | {a,b,d,e} | {1,a,b,c} |
| d | d | {a,b,c,d} | {0,1,a,c,d} | {1,a,b,e} | {a,b,d,e} | {1,c,d,e} | {1,b,c,e} |
| e | e | {1,a,c,d} | {b,c,d,e} | {0,a,b,d,e} | {1,a,b,c} | {1,b,c,e} | {1,a,d,e} |

| $HF_7^{52}$ | 0 | 1 | a | b | c | d | e |
|---|---|---|---|---|---|---|---|
| 0 | 0 | 1 | a | b | c | d | e |
| 1 | 1 | {1,a,b,e} | {a,b,d,e} | {1,a,c,d,e} | {0,1,b,c,e} | {a,b,c,d,e} | {1,a,c,d} |
| a | a | {a,b,d,e} | {1,a,b,c} | {1,b,c,e} | {1,a,b,d,e} | {0,1,a,c,d} | {1,b,c,d,e} |
| b | b | {1,a,c,d,e} | {1,b,c,e} | {a,b,c,d} | {1,a,c,d} | {1,a,b,c,e} | {0,a,b,d,e} |
| c | c | {0,1,b,c,e} | {1,a,b,d,e} | {1,a,c,d} | {b,c,d,e} | {a,b,d,e} | {1,a,b,c,d} |
| d | d | {a,b,c,d,e} | {0,1,a,c,d} | {1,a,b,c,e} | {a,b,d,e} | {1,c,d,e} | {1,b,c,e} |
| e | e | {1,a,c,d} | {1,b,c,d,e} | {0,a,b,d,e} | {1,a,b,c,d} | {1,b,c,e} | {1,a,d,e} |

| $HF_7^{53}$ | 0 | 1 | a | b | c | d | e |
|---|---|---|---|---|---|---|---|
| 0 | 0 | 1 | a | b | c | d | e |
| 1 | 1 | {1,a,c,d} | {a,b,c,e} | {1,a,d,e} | {0,1,b,c,e} | {b,c,d,e} | {1,a,b,d} |
| a | a | {a,b,c,e} | {a,b,d,e} | {1,b,c,d} | {1,a,b,e} | {0,1,a,c,d} | {1,c,d,e} |
| b | b | {1,a,d,e} | {1,b,c,d} | {1,b,c,e} | {a,c,d,e} | {1,a,b,c} | {0,a,b,d,e} |
| c | c | {0,1,b,c,e} | {1,a,b,e} | {a,c,d,e} | {1,a,c,d} | {1,b,d,e} | {a,b,c,d} |
| d | d | {b,c,d,e} | {0,1,a,c,d} | {1,a,b,c} | {1,b,d,e} | {a,b,d,e} | {1,a,c,e} |
| e | e | {1,a,b,d} | {1,c,d,e} | {0,a,b,d,e} | {a,b,c,d} | {1,a,c,e} | {1,b,c,e} |

| $HF_7^{54}$ | 0 | 1 | a | b | c | d | e |
|---|---|---|---|---|---|---|---|
| 0 | 0 | 1 | a | b | c | d | e |
| 1 | 1 | {1,a,b,c,d} | {a,b,c,d,e} | {1,a,d,e} | {0,1,b,c,e} | {b,c,d,e} | {1,a,b,c,d} |
| a | a | {a,b,c,d,e} | {a,b,c,d,e} | {1,b,c,d,e} | {1,a,b,e} | {0,1,a,c,d} | {1,c,d,e} |
| b | b | {1,a,d,e} | {1,b,c,d,e} | {1,b,c,d,e} | {1,a,c,d,e} | {1,a,b,c} | {0,a,b,d,e} |
| c | c | {0,1,b,c,e} | {1,a,b,e} | {1,a,c,d,e} | {1,a,c,d,e} | {1,a,b,d,e} | {a,b,c,d} |
| d | d | {b,c,d,e} | {0,1,a,c,d} | {1,a,b,c} | {1,a,b,d,e} | {1,a,b,d,e} | {1,a,b,c,e} |
| e | e | {1,a,b,c,d} | {1,c,d,e} | {0,a,b,d,e} | {a,b,c,d} | {1,a,b,c,e} | {1,a,b,c,e} |

| $HF_7^{55}$ | 0 | 1 | a | b | c | d | e |
|---|---|---|---|---|---|---|---|
| 0 | 0 | 1 | a | b | c | d | e |
| 1 | 1 | {1,a,b,c,e} | {a,b,d,e} | {1,c,d,e} | {0,1,b,c,e} | {a,b,c,d} | {1,a,c,d} |
| a | a | {a,b,d,e} | {1,a,b,c,d} | {1,b,c,e} | {1,a,d,e} | {0,1,a,c,d} | {b,c,d,e} |
| b | b | {1,c,d,e} | {1,b,c,e} | {a,b,c,d,e} | {1,a,c,d} | {1,a,b,e} | {0,a,b,d,e} |
| c | c | {0,1,b,c,e} | {1,a,d,e} | {1,a,c,d} | {1,b,c,d,e} | {a,b,d,e} | {1,a,b,c} |
| d | d | {a,b,c,d} | {0,1,a,c,d} | {1,a,b,e} | {a,b,d,e} | {1,a,c,d,e} | {1,b,c,e} |
| e | e | {1,a,c,d} | {b,c,d,e} | {0,a,b,d,e} | {1,a,b,c} | {1,b,c,e} | {1,a,b,d,e} |

| $HF_7^{56}$ | 0 | 1 | a | b | c | d | e |
|---|---|---|---|---|---|---|---|
| 0 | 0 | 1 | a | b | c | d | e |
| 1 | 1 | {1,a,b,c,e} | {a,b,d,e} | {1,a,c,d,e} | {0,1,b,c,e} | {a,b,c,d,e} | {1,a,c,d} |
| a | a | {a,b,d,e} | {1,a,b,c,d} | {1,b,c,e} | {1,a,b,d,e} | {0,1,a,c,d} | {1,b,c,d,e} |
| b | b | {1,a,c,d,e} | {1,b,c,e} | {a,b,c,d,e} | {1,a,c,d} | {1,a,b,c,e} | {0,a,b,d,e} |
| c | c | {0,1,b,c,e} | {1,a,b,d,e} | {1,a,c,d} | {1,b,c,d,e} | {a,b,d,e} | {1,a,b,c,d} |
| d | d | {a,b,c,d,e} | {0,1,a,c,d} | {1,a,b,c,e} | {a,b,d,e} | {1,a,c,d,e} | {1,b,c,e} |
| e | e | {1,a,c,d} | {1,b,c,d,e} | {0,a,b,d,e} | {1,a,b,c,d} | {1,b,c,e} | {1,a,b,d,e} |

**A3ii. Hyperfields for which card(x-x) = 5 and  x, -x ∉ x-x, for every non-zero element x.**

| $HF_7^{57}$ | 0 | 1 | a | b | c | d | e |
|---|---|---|---|---|---|---|---|
| 0 | 0 | 1 | a | b | c | d | e |
| 1 | 1 | {b,d} | {1,a,c,d} | {1,b} | {0,a,b,d,e} | {1,d} | {1,b,c,e} |
| a | a | {1,a,c,d} | {c,e} | {a,b,d,e} | {a,c} | {0,1,b,c,e} | {a,e} |
| b | b | {1,b} | {a,b,d,e} | {1,d} | {1,b,c,e} | {b,d} | {0,1,a,c,d} |
| c | c | {0,a,b,d,e} | {a,c} | {1,b,c,e} | {a,e} | {1,a,c,d} | {c,e} |
| d | d | {1,d} | {0,1,b,c,e} | {b,d} | {1,a,c,d} | {1,b} | {a,b,d,e} |
| e | e | {1,b,c,e} | {a,e} | {0,1,a,c,d} | {c,e} | {a,b,d,e} | {a,c} |

| $HF_7^{58}$ | 0 | 1 | a | b | c | d | e |
|---|---|---|---|---|---|---|---|
| 0 | 0 | 1 | a | b | c | d | e |
| 1 | 1 | {b,d} | {1,a,b,c,d,e} | {1,b,d} | {0,a,b,d,e} | {1,b,d} | {1,a,b,c,d,e} |
| a | a | {1,a,b,c,d,e} | {c,e} | {1,a,b,c,d,e} | {a,c,e} | {0,1,b,c,e} | {a,c,e} |
| b | b | {1,b,d} | {1,a,b,c,d,e} | {1,d} | {1,a,b,c,d,e} | {1,b,d} | {0,1,a,c,d} |
| c | c | {0,a,b,d,e} | {a,c,e} | {1,a,b,c,d,e} | {a,e} | {1,a,b,c,d,e} | {a,c,e} |
| d | d | {1,b,d} | {0,1,b,c,e} | {1,b,d} | {1,a,b,c,d,e} | {1,b} | {1,a,b,c,d,e} |
| e | e | {1,a,b,c,d,e} | {a,c,e} | {0,1,a,c,d} | {a,c,e} | {1,a,b,c,d,e} | {a,c} |

| $HF_7^{59}$ | 0 | 1 | a | b | c | d | e |
|---|---|---|---|---|---|---|---|
| 0 | 0 | 1 | a | b | c | d | e |
| 1 | 1 | {a,b,c} | {1,a,d} | {1,a,b,e} | {0,a,b,d,e} | {1,c,d,e} | {1,c,e} |
| a | a | {1,a,d} | {b,c,d} | {a,b,e} | {1,a,b,c} | {0,1,b,c,e} | {1,a,d,e} |
| b | b | {1,a,b,e} | {a,b,e} | {c,d,e} | {1,b,c} | {a,b,c,d} | {0,1,a,c,d} |
| c | c | {0,a,b,d,e} | {1,a,b,c} | {1,b,c} | {1,d,e} | {a,c,d} | {b,c,d,e} |
| d | d | {1,c,d,e} | {0,1,b,c,e} | {a,b,c,d} | {a,c,d} | {1,a,e} | {b,d,e} |
| e | e | {1,c,e} | {1,a,d,e} | {0,1,a,c,d} | {b,c,d,e} | {b,d,e} | {1,a,b} |

| $HF_7^{60}$ | 0 | 1 | a | b | c | d | e |
|---|---|---|---|---|---|---|---|
| 0 | 0 | 1 | a | b | c | d | e |
| 1 | 1 | {a,b,c} | {1,a,b,d,e} | {1,b,d,e} | {0,a,b,d,e} | {1,b,c,d} | {1,a,c,d,e} |
| a | a | {1,a,b,d,e} | {b,c,d} | {1,a,b,c,e} | {1,a,c,e} | {0,1,b,c,e} | {a,c,d,e} |
| b | b | {1,b,d,e} | {1,a,b,c,e} | {c,d,e} | {1,a,b,c,d} | {1,a,b,d} | {0,1,a,c,d} |
| c | c | {0,a,b,d,e} | {1,a,c,e} | {1,a,b,c,d} | {1,d,e} | {a,b,c,d,e} | {a,b,c,e} |
| d | d | {1,b,c,d} | {0,1,b,c,e} | {1,a,b,d} | {a,b,c,d,e} | {1,a,e} | {1,b,c,d,e} |
| e | e | {1,a,c,d,e} | {a,c,d,e} | {0,1,a,c,d} | {a,b,c,e} | {1,b,c,d,e} | {1,a,b} |

| $HF_7^{61}$ | 0 | 1 | a | b | c | d | e |
|---|---|---|---|---|---|---|---|
| 0 | 0 | 1 | a | b | c | d | e |
| 1 | 1 | {a,b,c} | {1,a,b,d,e} | {1,a,b,d,e} | {0,a,b,d,e} | {1,b,c,d,e} | {1,a,c,d,e} |
| a | a | {1,a,b,d,e} | {b,c,d} | {1,a,b,c,e} | {1,a,b,c,e} | {0,1,b,c,e} | {1,a,c,d,e} |
| b | b | {1,a,b,d,e} | {1,a,b,c,e} | {c,d,e} | {1,a,b,c,d} | {1,a,b,c,d} | {0,1,a,c,d} |
| c | c | {0,a,b,d,e} | {1,a,b,c,e} | {1,a,b,c,d} | {1,d,e} | {a,b,c,d,e} | {a,b,c,d,e} |
| d | d | {1,b,c,d,e} | {0,1,b,c,e} | {1,a,b,c,d} | {a,b,c,d,e} | {1,a,e} | {1,b,c,d,e} |
| e | e | {1,a,c,d,e} | {1,a,c,d,e} | {0,1,a,c,d} | {a,b,c,d,e} | {1,b,c,d,e} | {1,a,b} |

| $HF_7^{62}$ | 0 | 1 | a | b | c | d | e |
|---|---|---|---|---|---|---|---|
| 0 | 0 | 1 | a | b | c | d | e |
| 1 | 1 | {a,b,d} | {1,a,c,d} | {1,b,e} | {0,a,b,d,e} | {1,c,d} | {1,b,c,e} |
| a | a | {1,a,c,d} | {b,c,e} | {a,b,d,e} | {1,a,c} | {0,1,b,c,e} | {a,d,e} |
| b | b | {1,b,e} | {a,b,d,e} | {1,c,d} | {1,b,c,e} | {a,b,d} | {0,1,a,c,d} |
| c | c | {0,a,b,d,e} | {1,a,c} | {1,b,c,e} | {a,d,e} | {1,a,c,d} | {b,c,e} |
| d | d | {1,c,d} | {0,1,b,c,e} | {a,b,d} | {1,a,c,d} | {1,b,e} | {a,b,d,e} |
| e | e | {1,b,c,e} | {a,d,e} | {0,1,a,c,d} | {b,c,e} | {a,b,d,e} | {1,a,c} |

| $HF_7^{63}$ | 0 | 1 | a | b | c | d | e |
|---|---|---|---|---|---|---|---|
| 0 | 0 | 1 | a | b | c | d | e |
| 1 | 1 | {a,b,d} | {1,a,c,d} | {1,a,b,e} | {0,a,b,d,e} | {1,c,d,e} | {1,b,c,e} |
| a | a | {1,a,c,d} | {b,c,e} | {a,b,d,e} | {1,a,b,c} | {0,1,b,c,e} | {1,a,d,e} |
| b | b | {1,a,b,e} | {a,b,d,e} | {1,c,d} | {1,b,c,e} | {a,b,c,d} | {0,1,a,c,d} |
| c | c | {0,a,b,d,e} | {1,a,b,c} | {1,b,c,e} | {a,d,e} | {1,a,c,d} | {b,c,d,e} |
| d | d | {1,c,d,e} | {0,1,b,c,e} | {a,b,c,d} | {1,a,c,d} | {1,b,e} | {a,b,d,e} |
| e | e | {1,b,c,e} | {1,a,d,e} | {0,1,a,c,d} | {b,c,d,e} | {a,b,d,e} | {1,a,c} |

| $HF_7^{64}$ | 0 | 1 | a | b | c | d | e |
|---|---|---|---|---|---|---|---|
| 0 | 0 | 1 | a | b | c | d | e |
| 1 | 1 | {a,b,d} | {1,a,b,c,d,e} | {1,b,d,e} | {0,a,b,d,e} | {1,b,c,d} | {1,a,b,c,d,e} |
| a | a | {1,a,b,c,d,e} | {b,c,e} | {1,a,b,c,d,e} | {1,a,c,e} | {0,1,b,c,e} | {a,c,d,e} |
| b | b | {1,b,d,e} | {1,a,b,c,d,e} | {1,c,d} | {1,a,b,c,d,e} | {1,a,b,d} | {0,1,a,c,d} |
| c | c | {0,a,b,d,e} | {1,a,c,e} | {1,a,b,c,d,e} | {a,d,e} | {1,a,b,c,d,e} | {a,b,c,e} |
| d | d | {1,b,c,d} | {0,1,b,c,e} | {1,a,b,d} | {1,a,b,c,d,e} | {1,b,e} | {1,a,b,c,d,e} |
| e | e | {1,a,b,c,d,e} | {a,c,d,e} | {0,1,a,c,d} | {a,b,c,e} | {1,a,b,c,d,e} | {1,a,c} |

| $HF_7^{65}$ | 0 | 1 | a | b | c | d | e |
|---|---|---|---|---|---|---|---|
| 0 | 0 | 1 | a | b | c | d | e |
| 1 | 1 | {a,b,d} | {1,a,b,c,d,e} | {1,a,b,d,e} | {0,a,b,d,e} | {1,b,c,d,e} | {1,a,b,c,d,e} |
| a | a | {1,a,b,c,d,e} | {b,c,e} | {1,a,b,c,d,e} | {1,a,b,c,e} | {0,1,b,c,e} | {1,a,c,d,e} |
| b | b | {1,a,b,d,e} | {1,a,b,c,d,e} | {1,c,d} | {1,a,b,c,d,e} | {1,a,b,c,d} | {0,1,a,c,d} |
| c | c | {0,a,b,d,e} | {1,a,b,c,e} | {1,a,b,c,d,e} | {a,d,e} | {1,a,b,c,d,e} | {a,b,c,d,e} |
| d | d | {1,b,c,d,e} | {0,1,b,c,e} | {1,a,b,c,d} | {1,a,b,c,d,e} | {1,b,e} | {1,a,b,c,d,e} |
| e | e | {1,a,b,c,d,e} | {1,a,c,d,e} | {0,1,a,c,d} | {a,b,c,d,e} | {1,a,b,c,d,e} | {1,a,c} |

| $HF_7^{66}$ | 0 | 1 | a | b | c | d | e |
|---|---|---|---|---|---|---|---|
| 0 | 0 | 1 | a | b | c | d | e |
| 1 | 1 | {a,b,e} | {1,a,d} | {1,b,c,e} | {0,a,b,d,e} | {1,a,c,d} | {1,c,e} |
| a | a | {1,a,d} | {1,b,c} | {a,b,e} | {1,a,c,d} | {0,1,b,c,e} | {a,b,d,e} |
| b | b | {1,b,c,e} | {a,b,e} | {a,c,d} | {1,b,c} | {a,b,d,e} | {0,1,a,c,d} |
| c | c | {0,a,b,d,e} | {1,a,c,d} | {1,b,c} | {b,d,e} | {a,c,d} | {1,b,c,e} |
| d | d | {1,a,c,d} | {0,1,b,c,e} | {a,b,d,e} | {a,c,d} | {1,c,e} | {b,d,e} |
| e | e | {1,c,e} | {a,b,d,e} | {0,1,a,c,d} | {1,b,c,e} | {b,d,e} | {1,a,d} |

| $HF_7^{67}$ | 0 | 1 | a | b | c | d | e |
|---|---|---|---|---|---|---|---|
| 0 | 0 | 1 | a | b | c | d | e |
| 1 | 1 | {a,b,e} | {1,a,d} | {1,a,b,c,e} | {0,a,b,d,e} | {1,a,c,d,e} | {1,c,e} |
| a | a | {1,a,d} | {1,b,c} | {a,b,e} | {1,a,b,c,d} | {0,1,b,c,e} | {1,a,b,d,e} |
| b | b | {1,a,b,c,e} | {a,b,e} | {a,c,d} | {1,b,c} | {a,b,c,d,e} | {0,1,a,c,d} |
| c | c | {0,a,b,d,e} | {1,a,b,c,d} | {1,b,c} | {b,d,e} | {a,c,d} | {1,b,c,d,e} |
| d | d | {1,a,c,d,e} | {0,1,b,c,e} | {a,b,c,d,e} | {a,c,d} | {1,c,e} | {b,d,e} |
| e | e | {1,c,e} | {1,a,b,d,e} | {0,1,a,c,d} | {1,b,c,d,e} | {b,d,e} | {1,a,d} |

| $HF_7^{68}$ | 0 | 1 | a | b | c | d | e |
|---|---|---|---|---|---|---|---|
| 0 | 0 | 1 | a | b | c | d | e |
| 1 | 1 | {a,b,e} | {1,a,b,d,e} | {1,b,c,d,e} | {0,a,b,d,e} | {1,a,b,c,d} | {1,a,c,d,e} |
| a | a | {1,a,b,d,e} | {1,b,c} | {1,a,b,c,e} | {1,a,c,d,e} | {0,1,b,c,e} | {a,b,c,d,e} |
| b | b | {1,b,c,d,e} | {1,a,b,c,e} | {a,c,d} | {1,a,b,c,d} | {1,a,b,d,e} | {0,1,a,c,d} |
| c | c | {0,a,b,d,e} | {1,a,c,d,e} | {1,a,b,c,d} | {b,d,e} | {a,b,c,d,e} | {1,a,b,c,e} |
| d | d | {1,a,b,c,d} | {0,1,b,c,e} | {1,a,b,d,e} | {a,b,c,d,e} | {1,c,e} | {1,b,c,d,e} |
| e | e | {1,a,c,d,e} | {a,b,c,d,e} | {0,1,a,c,d} | {1,a,b,c,e} | {1,b,c,d,e} | {1,a,d} |

| $HF_7^{69}$ | 0 | 1 | a | b | c | d | e |
|---|---|---|---|---|---|---|---|
| 0 | 0 | 1 | a | b | c | d | e |
| 1 | 1 | {a,b,e} | {1,a,b,d,e} | {1,a,b,c,d,e} | {0,a,b,d,e} | {1,a,b,c,d,e} | {1,a,c,d,e} |
| a | a | {1,a,b,d,e} | {1,b,c} | {1,a,b,c,e} | {1,a,b,c,d,e} | {0,1,b,c,e} | {1,a,b,c,d,e} |
| b | b | {1,a,b,c,d,e} | {1,a,b,c,e} | {a,c,d} | {1,a,b,c,d} | {1,a,b,c,d,e} | {0,1,a,c,d} |
| c | c | {0,a,b,d,e} | {1,a,b,c,d,e} | {1,a,b,c,d} | {b,d,e} | {a,b,c,d,e} | {1,a,b,c,d,e} |
| d | d | {1,a,b,c,d,e} | {0,1,b,c,e} | {1,a,b,c,d,e} | {a,b,c,d,e} | {1,c,e} | {1,b,c,d,e} |
| e | e | {1,a,c,d,e} | {1,a,b,c,d,e} | {0,1,a,c,d} | {1,a,b,c,d,e} | {1,b,c,d,e} | {1,a,d} |

| $HF_7^{70}$ | 0 | 1 | a | b | c | d | e |
|---|---|---|---|---|---|---|---|
| 0 | 0 | 1 | a | b | c | d | e |
| 1 | 1 | {a,c,d} | {1,a,b,c,e} | {1,b,d,e} | {0,a,b,d,e} | {1,b,c,d} | {1,a,b,d,e} |
| a | a | {1,a,b,c,e} | {b,d,e} | {1,a,b,c,d} | {1,a,c,e} | {0,1,b,c,e} | {a,c,d,e} |
| b | b | {1,b,d,e} | {1,a,b,c,d} | {1,c,e} | {a,b,c,d,e} | {1,a,b,d} | {0,1,a,c,d} |
| c | c | {0,a,b,d,e} | {1,a,c,e} | {a,b,c,d,e} | {1,a,d} | {1,b,c,d,e} | {a,b,c,e} |
| d | d | {1,b,c,d} | {0,1,b,c,e} | {1,a,b,d} | {1,b,c,d,e} | {a,b,e} | {1,a,c,d,e} |
| e | e | {1,a,b,d,e} | {a,c,d,e} | {0,1,a,c,d} | {a,b,c,e} | {1,a,c,d,e} | {1,b,c} |

| $HF_7^{71}$ | 0 | 1 | a | b | c | d | e |
|---|---|---|---|---|---|---|---|
| 0 | 0 | 1 | a | b | c | d | e |
| 1 | 1 | {a,c,d} | {1,a,b,c,e} | {1,a,b,d,e} | {0,a,b,d,e} | {1,b,c,d,e} | {1,a,b,d,e} |
| a | a | {1,a,b,c,e} | {b,d,e} | {1,a,b,c,d} | {1,a,b,c,e} | {0,1,b,c,e} | {1,a,c,d,e} |
| b | b | {1,a,b,d,e} | {1,a,b,c,d} | {1,c,e} | {a,b,c,d,e} | {1,a,b,c,d} | {0,1,a,c,d} |
| c | c | {0,a,b,d,e} | {1,a,b,c,e} | {a,b,c,d,e} | {1,a,d} | {1,b,c,d,e} | {a,b,c,d,e} |
| d | d | {1,b,c,d,e} | {0,1,b,c,e} | {1,a,b,c,d} | {1,b,c,d,e} | {a,b,e} | {1,a,c,d,e} |
| e | e | {1,a,b,d,e} | {1,a,c,d,e} | {0,1,a,c,d} | {a,b,c,d,e} | {1,a,c,d,e} | {1,b,c} |

| $HF_7^{72}$ | 0 | 1 | a | b | c | d | e |
|---|---|---|---|---|---|---|---|
| 0 | 0 | 1 | a | b | c | d | e |
| 1 | 1 | {b,c,d} | {1,a,b,c,d,e} | {1,a,b,d} | {0,a,b,d,e} | {1,b,d,e} | {1,a,b,c,d,e} |
| a | a | {1,a,b,c,d,e} | {c,d,e} | {1,a,b,c,d,e} | {a,b,c,e} | {0,1,b,c,e} | {1,a,c,e} |
| b | b | {1,a,b,d} | {1,a,b,c,d,e} | {1,d,e} | {1,a,b,c,d,e} | {1,b,c,d} | {0,1,a,c,d} |
| c | c | {0,a,b,d,e} | {a,b,c,e} | {1,a,b,c,d,e} | {1,a,e} | {1,a,b,c,d,e} | {a,c,d,e} |
| d | d | {1,b,d,e} | {0,1,b,c,e} | {1,b,c,d} | {1,a,b,c,d,e} | {1,a,b} | {1,a,b,c,d,e} |
| e | e | {1,a,b,c,d,e} | {1,a,c,e} | {0,1,a,c,d} | {a,c,d,e} | {1,a,b,c,d,e} | {a,b,c} |

| $HF_7^{73}$ | 0 | 1 | a | b | c | d | e |
|---|---|---|---|---|---|---|---|
| 0 | 0 | 1 | a | b | c | d | e |
| 1 | 1 | {a,b,c,d} | {1,a,c,d} | {1,b,e} | {0,a,b,d,e} | {1,c,d} | {1,b,c,e} |
| a | a | {1,a,c,d} | {b,c,d,e} | {a,b,d,e} | {1,a,c} | {0,1,b,c,e} | {a,d,e} |
| b | b | {1,b,e} | {a,b,d,e} | {1,c,d,e} | {1,b,c,e} | {a,b,d} | {0,1,a,c,d} |
| c | c | {0,a,b,d,e} | {1,a,c} | {1,b,c,e} | {1,a,d,e} | {1,a,c,d} | {b,c,e} |
| d | d | {1,c,d} | {0,1,b,c,e} | {a,b,d} | {1,a,c,d} | {1,a,b,e} | {a,b,d,e} |
| e | e | {1,b,c,e} | {a,d,e} | {0,1,a,c,d} | {b,c,e} | {a,b,d,e} | {1,a,b,c} |

| $HF_7^{74}$ | 0 | 1 | a | b | c | d | e |
|---|---|---|---|---|---|---|---|
| 0 | 0 | 1 | a | b | c | d | e |
| 1 | 1 | {a,b,c,d} | {1,a,c,d} | {1,a,b,e} | {0,a,b,d,e} | {1,c,d,e} | {1,b,c,e} |
| a | a | {1,a,c,d} | {b,c,d,e} | {a,b,d,e} | {1,a,b,c} | {0,1,b,c,e} | {1,a,d,e} |
| b | b | {1,a,b,e} | {a,b,d,e} | {1,c,d,e} | {1,b,c,e} | {a,b,c,d} | {0,1,a,c,d} |
| c | c | {0,a,b,d,e} | {1,a,b,c} | {1,b,c,e} | {1,a,d,e} | {1,a,c,d} | {b,c,d,e} |
| d | d | {1,c,d,e} | {0,1,b,c,e} | {a,b,c,d} | {1,a,c,d} | {1,a,b,e} | {a,b,d,e} |
| e | e | {1,b,c,e} | {1,a,d,e} | {0,1,a,c,d} | {b,c,d,e} | {a,b,d,e} | {1,a,b,c} |

| $HF_7^{75}$ | 0 | 1 | a | b | c | d | e |
|---|---|---|---|---|---|---|---|
| 0 | 0 | 1 | a | b | c | d | e |
| 1 | 1 | {a,b,c,d} | {1,a,b,c,d,e} | {1,b,d,e} | {0,a,b,d,e} | {1,b,c,d} | {1,a,b,c,d,e} |
| a | a | {1,a,b,c,d,e} | {b,c,d,e} | {1,a,b,c,d,e} | {1,a,c,e} | {0,1,b,c,e} | {a,c,d,e} |
| b | b | {1,b,d,e} | {1,a,b,c,d,e} | {1,c,d,e} | {1,a,b,c,d,e} | {1,a,b,d} | {0,1,a,c,d} |
| c | c | {0,a,b,d,e} | {1,a,c,e} | {1,a,b,c,d,e} | {1,a,d,e} | {1,a,b,c,d,e} | {a,b,c,e} |
| d | d | {1,b,c,d} | {0,1,b,c,e} | {1,a,b,d} | {1,a,b,c,d,e} | {1,a,b,e} | {1,a,b,c,d,e} |
| e | e | {1,a,b,c,d,e} | {a,c,d,e} | {0,1,a,c,d} | {a,b,c,e} | {1,a,b,c,d,e} | {1,a,b,c} |

| $HF_7^{76}$ | 0 | 1 | a | b | c | d | e |
|---|---|---|---|---|---|---|---|
| 0 | 0 | 1 | a | b | c | d | e |
| 1 | 1 | {a,b,c,d} | {1,a,b,c,d,e} | {1,a,b,d,e} | {0,a,b,d,e} | {1,b,c,d,e} | {1,a,b,c,d,e} |
| a | a | {1,a,b,c,d,e} | {b,c,d,e} | {1,a,b,c,d,e} | {1,a,b,c,e} | {0,1,b,c,e} | {1,a,c,d,e} |
| b | b | {1,a,b,d,e} | {1,a,b,c,d,e} | {1,c,d,e} | {1,a,b,c,d,e} | {1,a,b,c,d} | {0,1,a,c,d} |
| c | c | {0,a,b,d,e} | {1,a,b,c,e} | {1,a,b,c,d,e} | {1,a,d,e} | {1,a,b,c,d,e} | {a,b,c,d,e} |
| d | d | {1,b,c,d,e} | {0,1,b,c,e} | {1,a,b,c,d} | {1,a,b,c,d,e} | {1,a,b,e} | {1,a,b,c,d,e} |
| e | e | {1,a,b,c,d,e} | {1,a,c,d,e} | {0,1,a,c,d} | {a,b,c,d,e} | {1,a,b,c,d,e} | {1,a,b,c} |

| $HF_7^{77}$ | 0 | 1 | a | b | c | d | e |
|---|---|---|---|---|---|---|---|
| 0 | 0 | 1 | a | b | c | d | e |
| 1 | 1 | {a,b,c,e} | {1,a,d} | {1,b,c,e} | {0,a,b,d,e} | {1,a,c,d} | {1,c,e} |
| a | a | {1,a,d} | {1,b,c,d} | {a,b,e} | {1,a,c,d} | {0,1,b,c,e} | {a,b,d,e} |
| b | b | {1,b,c,e} | {a,b,e} | {a,c,d,e} | {1,b,c} | {a,b,d,e} | {0,1,a,c,d} |
| c | c | {0,a,b,d,e} | {1,a,c,d} | {1,b,c} | {1,b,d,e} | {a,c,d} | {1,b,c,e} |
| d | d | {1,a,c,d} | {0,1,b,c,e} | {a,b,d,e} | {a,c,d} | {1,a,c,e} | {b,d,e} |
| e | e | {1,c,e} | {a,b,d,e} | {0,1,a,c,d} | {1,b,c,e} | {b,d,e} | {1,a,b,d} |

| $HF_7^{78}$ | 0 | 1 | a | b | c | d | e |
|---|---|---|---|---|---|---|---|
| 0 | 0 | 1 | a | b | c | d | e |
| 1 | 1 | {a,b,c,e} | {1,a,d} | {1,a,b,c,e} | {0,a,b,d,e} | {1,a,c,d,e} | {1,c,e} |
| a | a | {1,a,d} | {1,b,c,d} | {a,b,e} | {1,a,b,c,d} | {0,1,b,c,e} | {1,a,b,d,e} |
| b | b | {1,a,b,c,e} | {a,b,e} | {a,c,d,e} | {1,b,c} | {a,b,c,d,e} | {0,1,a,c,d} |
| c | c | {0,a,b,d,e} | {1,a,b,c,d} | {1,b,c} | {1,b,d,e} | {a,c,d} | {1,b,c,d,e} |
| d | d | {1,a,c,d,e} | {0,1,b,c,e} | {a,b,c,d,e} | {a,c,d} | {1,a,c,e} | {b,d,e} |
| e | e | {1,c,e} | {1,a,b,d,e} | {0,1,a,c,d} | {1,b,c,d,e} | {b,d,e} | {1,a,b,d} |

| $HF_7^{79}$ | 0 | 1 | a | b | c | d | e |
|---|---|---|---|---|---|---|---|
| 0 | 0 | 1 | a | b | c | d | e |
| 1 | 1 | {a,b,c,e} | {1,a,b,d,e} | {1,b,c,d,e} | {0,a,b,d,e} | {1,a,b,c,d} | {1,a,c,d,e} |
| a | a | {1,a,b,d,e} | {1,b,c,d} | {1,a,b,c,e} | {1,a,c,d,e} | {0,1,b,c,e} | {a,b,c,d,e} |
| b | b | {1,b,c,d,e} | {1,a,b,c,e} | {a,c,d,e} | {1,a,b,c,d} | {1,a,b,d,e} | {0,1,a,c,d} |
| c | c | {0,a,b,d,e} | {1,a,c,d,e} | {1,a,b,c,d} | {1,b,d,e} | {a,b,c,d,e} | {1,a,b,c,e} |
| d | d | {1,a,b,c,d} | {0,1,b,c,e} | {1,a,b,d,e} | {a,b,c,d,e} | {1,a,c,e} | {1,b,c,d,e} |
| e | e | {1,a,c,d,e} | {a,b,c,d,e} | {0,1,a,c,d} | {1,a,b,c,e} | {1,b,c,d,e} | {1,a,b,d} |

| $HF_7^{80}$ | 0 | 1 | a | b | c | d | e |
|---|---|---|---|---|---|---|---|
| 0 | 0 | 1 | a | b | c | d | e |
| 1 | 1 | {a,b,c,e} | {1,a,b,d,e} | {1,a,b,c,d,e} | {0,a,b,d,e} | {1,a,b,c,d,e} | {1,a,c,d,e} |
| a | a | {1,a,b,d,e} | {1,b,c,d} | {1,a,b,c,e} | {1,a,b,c,d,e} | {0,1,b,c,e} | {1,a,b,c,d,e} |
| b | b | {1,a,b,c,d,e} | {1,a,b,c,e} | {a,c,d,e} | {1,a,b,c,d} | {1,a,b,c,d,e} | {0,1,a,c,d} |
| c | c | {0,a,b,d,e} | {1,a,b,c,d,e} | {1,a,b,c,d} | {1,b,d,e} | {a,b,c,d,e} | {1,a,b,c,d,e} |
| d | d | {1,a,b,c,d,e} | {0,1,b,c,e} | {1,a,b,c,d,e} | {a,b,c,d,e} | {1,a,c,e} | {1,b,c,d,e} |
| e | e | {1,a,c,d,e} | {1,a,b,c,d,e} | {0,1,a,c,d} | {1,a,b,c,d,e} | {1,b,c,d,e} | {1,a,b,d} |

| $HF_7^{81}$ | 0 | 1 | a | b | c | d | e |
|---|---|---|---|---|---|---|---|
| 0 | 0 | 1 | a | b | c | d | e |
| 1 | 1 | {a,b,d,e} | {1,a,c,d} | {1,b,c,e} | {0,a,b,d,e} | {1,a,c,d} | {1,b,c,e} |
| a | a | {1,a,c,d} | {1,b,c,e} | {a,b,d,e} | {1,a,c,d} | {0,1,b,c,e} | {a,b,d,e} |
| b | b | {1,b,c,e} | {a,b,d,e} | {1,a,c,d} | {1,b,c,e} | {a,b,d,e} | {0,1,a,c,d} |
| c | c | {0,a,b,d,e} | {1,a,c,d} | {1,b,c,e} | {a,b,d,e} | {1,a,c,d} | {1,b,c,e} |
| d | d | {1,a,c,d} | {0,1,b,c,e} | {a,b,d,e} | {1,a,c,d} | {1,b,c,e} | {a,b,d,e} |
| e | e | {1,b,c,e} | {a,b,d,e} | {0,1,a,c,d} | {1,b,c,e} | {a,b,d,e} | {1,a,c,d} |

| $HF_7^{82}$ | 0 | 1 | a | b | c | d | e |
|---|---|---|---|---|---|---|---|
| 0 | 0 | 1 | a | b | c | d | e |
| 1 | 1 | {a,b,d,e} | {1,a,c,d} | {1,a,b,c,e} | {0,a,b,d,e} | {1,a,c,d,e} | {1,b,c,e} |
| a | a | {1,a,c,d} | {1,b,c,e} | {a,b,d,e} | {1,a,b,c,d} | {0,1,b,c,e} | {1,a,b,d,e} |
| b | b | {1,a,b,c,e} | {a,b,d,e} | {1,a,c,d} | {1,b,c,e} | {a,b,c,d,e} | {0,1,a,c,d} |
| c | c | {0,a,b,d,e} | {1,a,b,c,d} | {1,b,c,e} | {a,b,d,e} | {1,a,c,d} | {1,b,c,d,e} |
| d | d | {1,a,c,d,e} | {0,1,b,c,e} | {a,b,c,d,e} | {1,a,c,d} | {1,b,c,e} | {a,b,d,e} |
| e | e | {1,b,c,e} | {1,a,b,d,e} | {0,1,a,c,d} | {1,b,c,d,e} | {a,b,d,e} | {1,a,c,d} |

| $HF_7^{83}$ | 0 | 1 | a | b | c | d | e |
|---|---|---|---|---|---|---|---|
| 0 | 0 | 1 | a | b | c | d | e |
| 1 | 1 | {a,b,d,e} | {1,a,b,c,d,e} | {1,b,c,d,e} | {0,a,b,d,e} | {1,a,b,c,d} | {1,a,b,c,d,e} |
| a | a | {1,a,b,c,d,e} | {1,b,c,e} | {1,a,b,c,d,e} | {1,a,c,d,e} | {0,1,b,c,e} | {a,b,c,d,e} |
| b | b | {1,b,c,d,e} | {1,a,b,c,d,e} | {1,a,c,d} | {1,a,b,c,d,e} | {1,a,b,d,e} | {0,1,a,c,d} |
| c | c | {0,a,b,d,e} | {1,a,c,d,e} | {1,a,b,c,d,e} | {a,b,d,e} | {1,a,b,c,d,e} | {1,a,b,c,e} |
| d | d | {1,a,b,c,d} | {0,1,b,c,e} | {1,a,b,d,e} | {1,a,b,c,d,e} | {1,b,c,e} | {1,a,b,c,d,e} |
| e | e | {1,a,b,c,d,e} | {a,b,c,d,e} | {0,1,a,c,d} | {1,a,b,c,e} | {1,a,b,c,d,e} | {1,a,c,d} |

| $HF_7^{84}$ | 0 | 1 | a | b | c | d | e |
|---|---|---|---|---|---|---|---|
| 0 | 0 | 1 | a | b | c | d | e |
| 1 | 1 | {a,b,d,e} | {1,a,b,c,d,e} | {1,a,b,c,d,e} | {0,a,b,d,e} | {1,a,b,c,d,e} | {1,a,b,c,d,e} |
| a | a | {1,a,b,c,d,e} | {1,b,c,e} | {1,a,b,c,d,e} | {1,a,b,c,d,e} | {0,1,b,c,e} | {1,a,b,c,d,e} |
| b | b | {1,a,b,c,d,e} | {1,a,b,c,d,e} | {1,a,c,d} | {1,a,b,c,d,e} | {1,a,b,c,d,e} | {0,1,a,c,d} |
| c | c | {0,a,b,d,e} | {1,a,b,c,d,e} | {1,a,b,c,d,e} | {a,b,d,e} | {1,a,b,c,d,e} | {1,a,b,c,d,e} |
| d | d | {1,a,b,c,d,e} | {0,1,b,c,e} | {1,a,b,c,d,e} | {1,a,b,c,d,e} | {1,b,c,e} | {1,a,b,c,d,e} |
| e | e | {1,a,b,c,d,e} | {1,a,b,c,d,e} | {0,1,a,c,d} | {1,a,b,c,d,e} | {1,a,b,c,d,e} | {1,a,c,d} |

| $HF_7^{85}$ | 0 | 1 | a | b | c | d | e |
|---|---|---|---|---|---|---|---|
| 0 | 0 | 1 | a | b | c | d | e |
| 1 | 1 | {a,b,c,d,e} | {1,a,c,d} | {1,b,c,e} | {0,a,b,d,e} | {1,a,c,d} | {1,b,c,e} |
| a | a | {1,a,c,d} | {1,b,c,d,e} | {a,b,d,e} | {1,a,c,d} | {0,1,b,c,e} | {a,b,d,e} |
| b | b | {1,b,c,e} | {a,b,d,e} | {1,a,c,d,e} | {1,b,c,e} | {a,b,d,e} | {0,1,a,c,d} |
| c | c | {0,a,b,d,e} | {1,a,c,d} | {1,b,c,e} | {1,a,b,d,e} | {1,a,c,d} | {1,b,c,e} |
| d | d | {1,a,c,d} | {0,1,b,c,e} | {a,b,d,e} | {1,a,c,d} | {1,a,b,c,e} | {a,b,d,e} |
| e | e | {1,b,c,e} | {a,b,d,e} | {0,1,a,c,d} | {1,b,c,e} | {a,b,d,e} | {1,a,b,c,d} |

| $HF_7^{86}$ | 0 | 1 | a | b | c | d | e |
|---|---|---|---|---|---|---|---|
| 0 | 0 | 1 | a | b | c | d | e |
| 1 | 1 | {a,b,c,d,e} | {1,a,c,d} | {1,a,b,c,e} | {0,a,b,d,e} | {1,a,c,d,e} | {1,b,c,e} |
| a | a | {1,a,c,d} | {1,b,c,d,e} | {a,b,d,e} | {1,a,b,c,d} | {0,1,b,c,e} | {1,a,b,d,e} |
| b | b | {1,a,b,c,e} | {a,b,d,e} | {1,a,c,d,e} | {1,b,c,e} | {a,b,c,d,e} | {0,1,a,c,d} |
| c | c | {0,a,b,d,e} | {1,a,b,c,d} | {1,b,c,e} | {1,a,b,d,e} | {1,a,c,d} | {1,b,c,d,e} |
| d | d | {1,a,c,d,e} | {0,1,b,c,e} | {a,b,c,d,e} | {1,a,c,d} | {1,a,b,c,e} | {a,b,d,e} |
| e | e | {1,b,c,e} | {1,a,b,d,e} | {0,1,a,c,d} | {1,b,c,d,e} | {a,b,d,e} | {1,a,b,c,d} |

| $HF_7^{87}$ | 0 | 1 | a | b | c | d | e |
|---|---|---|---|---|---|---|---|
| 0 | 0 | 1 | a | b | c | d | e |
| 1 | 1 | {a,b,c,d,e} | {1,a,b,c,d,e} | {1,b,c,d,e} | {0,a,b,d,e} | {1,a,b,c,d} | {1,a,b,c,d,e} |
| a | a | {1,a,b,c,d,e} | {1,b,c,d,e} | {1,a,b,c,d,e} | {1,a,c,d,e} | {0,1,b,c,e} | {a,b,c,d,e} |
| b | b | {1,b,c,d,e} | {1,a,b,c,d,e} | {1,a,c,d,e} | {1,a,b,c,d,e} | {1,a,b,d,e} | {0,1,a,c,d} |
| c | c | {0,a,b,d,e} | {1,a,c,d,e} | {1,a,b,c,d,e} | {1,a,b,d,e} | {1,a,b,c,d,e} | {1,a,b,c,e} |
| d | d | {1,a,b,c,d} | {0,1,b,c,e} | {1,a,b,d,e} | {1,a,b,c,d,e} | {1,a,b,c,e} | {1,a,b,c,d,e} |
| e | e | {1,a,b,c,d,e} | {a,b,c,d,e} | {0,1,a,c,d} | {1,a,b,c,e} | {1,a,b,c,d,e} | {1,a,b,c,d} |

| $HF_7^{88}$ | 0 | 1 | a | b | c | d | e |
|---|---|---|---|---|---|---|---|
| 0 | 0 | 1 | a | b | c | d | e |
| 1 | 1 | {a,b,c,d,e} | {1,a,b,c,d,e} | {1,a,b,c,d,e} | {0,a,b,d,e} | {1,a,b,c,d,e} | {1,a,b,c,d,e} |
| a | a | {1,a,b,c,d,e} | {1,b,c,d,e} | {1,a,b,c,d,e} | {1,a,b,c,d,e} | {0,1,b,c,e} | {1,a,b,c,d,e} |
| b | b | {1,a,b,c,d,e} | {1,a,b,c,d,e} | {1,a,c,d,e} | {1,a,b,c,d,e} | {1,a,b,c,d,e} | {0,1,a,c,d} |
| c | c | {0,a,b,d,e} | {1,a,b,c,d,e} | {1,a,b,c,d,e} | {1,a,b,d,e} | {1,a,b,c,d,e} | {1,a,b,c,d,e} |
| d | d | {1,a,b,c,d,e} | {0,1,b,c,e} | {1,a,b,c,d,e} | {1,a,b,c,d,e} | {1,a,b,c,e} | {1,a,b,c,d,e} |
| e | e | {1,a,b,c,d,e} | {1,a,b,c,d,e} | {0,1,a,c,d} | {1,a,b,c,d,e} | {1,a,b,c,d,e} | {1,a,b,c,d} |

## 4. Hyperfields for which card(x-x) = 7, for every non-zero element x.

| $HF_7^2$ $= [\mathbb{Z}_7]$ | 0 | 1 | a | b | c | d | e |
|---|---|---|---|---|---|---|---|
| 0 | 0 | 1 | a | b | c | d | e |
| 1 | 1 | {1,b} | {1,a,d} | {1,a,b} | {0,1,a,b,c,d,e} | {1,d,e} | {1,c,e} |
| a | a | {1,a,d} | {a,c} | {a,b,e} | {a,b,c} | {0,1,a,b,c,d,e} | {1,a,e} |
| b | b | {1,a,b} | {a,b,e} | {b,d} | {1,b,c} | {b,c,d} | {0,1,a,b,c,d,e} |
| c | c | {0,1,a,b,c,d,e} | {a,b,c} | {1,b,c} | {c,e} | {a,c,d} | {c,d,e} |
| d | d | {1,d,e} | {0,1,a,b,c,d,e} | {b,c,d} | {a,c,d} | {1,d} | {b,d,e} |
| e | e | {1,c,e} | {1,a,e} | {0,1,a,b,c,d,e} | {c,d,e} | {b,d,e} | {a,e} |

| $HF_7^{89}$ | 0 | 1 | a | b | c | d | e |
|---|---|---|---|---|---|---|---|
| 0 | 0 | 1 | a | b | c | d | e |
| 1 | 1 | 1 | {1,a} | {1,b} | {0,1,a,b,c,d,e} | {1,d} | {1,e} |
| a | a | {1,a} | a | {a,b} | {a,c} | {0,1,a,b,c,d,e} | {a,e} |
| b | b | {1,b} | {a,b} | b | {b,c} | {b,d} | {0,1,a,b,c,d,e} |
| c | c | {0,1,a,b,c,d,e} | {a,c} | {b,c} | c | {c,d} | {c,e} |
| d | d | {1,d} | {0,1,a,b,c,d,e} | {b,d} | {c,d} | d | {d,e} |
| e | e | {1,e} | {a,e} | {0,1,a,b,c,d,e} | {c,e} | {d,e} | e |

| $HF_7^{90}$ | 0 | 1 | a | b | c | d | e |
|---|---|---|---|---|---|---|---|
| 0 | 0 | 1 | a | b | c | d | e |
| 1 | 1 | {1,c} | {1,a} | {1,b} | {0,1,a,b,c,d,e} | {1,d} | {1,e} |
| a | a | {1,a} | {a,d} | {a,b} | {a,c} | {0,1,a,b,c,d,e} | {a,e} |
| b | b | {1,b} | {a,b} | {b,e} | {b,c} | {b,d} | {0,1,a,b,c,d,e} |
| c | c | {0,1,a,b,c,d,e} | {a,c} | {b,c} | {1,c} | {c,d} | {c,e} |
| d | d | {1,d} | {0,1,a,b,c,d,e} | {b,d} | {c,d} | {a,d} | {d,e} |
| e | e | {1,e} | {a,e} | {0,1,a,b,c,d,e} | {c,e} | {d,e} | {b,e} |

| $HF_7^{91}$ | 0 | 1 | a | b | c | d | e |
|---|---|---|---|---|---|---|---|
| 0 | 0 | 1 | a | b | c | d | e |
| 1 | 1 | 1 | {1,a} | {1,a,b} | {0,1,a,b,c,d,e} | {1,d,e} | {1,e} |
| a | a | {1,a} | a | {a,b} | {a,b,c} | {0,1,a,b,c,d,e} | {1,a,e} |
| b | b | {1,a,b} | {a,b} | b | {b,c} | {b,c,d} | {0,1,a,b,c,d,e} |
| c | c | {0,1,a,b,c,d,e} | {a,b,c} | {b,c} | c | {c,d} | {c,d,e} |
| d | d | {1,d,e} | {0,1,a,b,c,d,e} | {b,c,d} | {c,d} | d | {d,e} |
| e | e | {1,e} | {1,a,e} | {0,1,a,b,c,d,e} | {c,d,e} | {d,e} | e |

| $HF_7^{92}$ | 0 | 1 | a | b | c | d | e |
|---|---|---|---|---|---|---|---|
| 0 | 0 | 1 | a | b | c | d | e |
| 1 | 1 | {1,b} | {1,a,b,d,e} | {1,b,d} | {0,1,a,b,c,d,e} | {1,b,d} | {1,a,c,d,e} |
| a | a | {1,a,b,d,e} | {a,c} | {1,a,b,c,e} | {a,c,e} | {0,1,a,b,c,d,e} | {a,c,e} |
| b | b | {1,b,d} | {1,a,b,c,e} | {b,d} | {1,a,b,c,d} | {1,b,d} | {0,1,a,b,c,d,e} |
| c | c | {0,1,a,b,c,d,e} | {a,c,e} | {1,a,b,c,d} | {c,e} | {a,b,c,d,e} | {a,c,e} |
| d | d | {1,b,d} | {0,1,a,b,c,d,e} | {1,b,d} | {a,b,c,d,e} | {1,d} | {1,b,c,d,e} |
| e | e | {1,a,c,d,e} | {a,c,e} | {0,1,a,b,c,d,e} | {a,c,e} | {1,b,c,d,e} | {a,e} |

| $HF_7^{93}$ | 0 | 1 | a | b | c | d | e |
|---|---|---|---|---|---|---|---|
| 0 | 0 | 1 | a | b | c | d | e |
| 1 | 1 | {1,a,b} | {1,a,d} | {1,a,b,e} | {0,1,a,b,c,d,e} | {1,c,d,e} | {1,c,e} |
| a | a | {1,a,d} | {a,b,c} | {a,b,e} | {1,a,b,c} | {0,1,a,b,c,d,e} | {1,a,d,e} |
| b | b | {1,a,b,e} | {a,b,e} | {b,c,d} | {1,b,c} | {a,b,c,d} | {0,1,a,b,c,d,e} |
| c | c | {0,1,a,b,c,d,e} | {1,a,b,c} | {1,b,c} | {c,d,e} | {a,c,d} | {b,c,d,e} |
| d | d | {1,c,d,e} | {0,1,a,b,c,d,e} | {a,b,c,d} | {a,c,d} | {1,d,e} | {b,d,e} |
| e | e | {1,c,e} | {1,a,d,e} | {0,1,a,b,c,d,e} | {b,c,d,e} | {b,d,e} | {1,a,e} |

| $HF_7^{94}$ | 0 | 1 | a | b | c | d | e |
|---|---|---|---|---|---|---|---|
| 0 | 0 | 1 | a | b | c | d | e |
| 1 | 1 | {1,a,b} | {1,a,b,d,e} | {1,b,d,e} | {0,1,a,b,c,d,e} | {1,b,c,d} | {1,a,c,d,e} |
| a | a | {1,a,b,d,e} | {a,b,c} | {1,a,b,c,e} | {1,a,c,e} | {0,1,a,b,c,d,e} | {a,c,d,e} |
| b | b | {1,b,d,e} | {1,a,b,c,e} | {b,c,d} | {1,a,b,c,d} | {1,a,b,d} | {0,1,a,b,c,d,e} |
| c | c | {0,1,a,b,c,d,e} | {1,a,c,e} | {1,a,b,c,d} | {c,d,e} | {a,b,c,d,e} | {a,b,c,e} |
| d | d | {1,b,c,d} | {0,1,a,b,c,d,e} | {1,a,b,d} | {a,b,c,d,e} | {1,d,e} | {1,b,c,d,e} |
| e | e | {1,a,c,d,e} | {a,c,d,e} | {0,1,a,b,c,d,e} | {a,b,c,e} | {1,b,c,d,e} | {1,a,e} |

| $HF_7^{95}$ | 0 | 1 | a | b | c | d | e |
|---|---|---|---|---|---|---|---|
| 0 | 0 | 1 | a | b | c | d | e |
| 1 | 1 | {1,a,b} | {1,a,b,d,e} | {1,a,b,d,e} | {0,1,a,b,c,d,e} | {1,b,c,d,e} | {1,a,c,d,e} |
| a | a | {1,a,b,d,e} | {a,b,c} | {1,a,b,c,e} | {1,a,b,c,e} | {0,1,a,b,c,d,e} | {1,a,c,d,e} |
| b | b | {1,a,b,d,e} | {1,a,b,c,e} | {b,c,d} | {1,a,b,c,d} | {1,a,b,c,d} | {0,1,a,b,c,d,e} |
| c | c | {0,1,a,b,c,d,e} | {1,a,b,c,e} | {1,a,b,c,d} | {c,d,e} | {a,b,c,d,e} | {a,b,c,d,e} |
| d | d | {1,b,c,d,e} | {0,1,a,b,c,d,e} | {1,a,b,c,d} | {a,b,c,d,e} | {1,d,e} | {1,b,c,d,e} |
| e | e | {1,a,c,d,e} | {1,a,c,d,e} | {0,1,a,b,c,d,e} | {a,b,c,d,e} | {1,b,c,d,e} | {1,a,e} |

| $HF_7^{96}$ | 0 | 1 | a | b | c | d | e |
|---|---|---|---|---|---|---|---|
| 0 | 0 | 1 | a | b | c | d | e |
| 1 | 1 | {1,a,c} | {1,a,b,e} | {1,b,d,e} | {0,1,a,b,c,d,e} | {1,b,c,d} | {1,a,d,e} |
| a | a | {1,a,b,e} | {a,b,d} | {1,a,b,c} | {1,a,c,e} | {0,1,a,b,c,d,e} | {a,c,d,e} |
| b | b | {1,b,d,e} | {1,a,b,c} | {b,c,e} | {a,b,c,d} | {1,a,b,d} | {0,1,a,b,c,d,e} |
| c | c | {0,1,a,b,c,d,e} | {1,a,c,e} | {a,b,c,d} | {1,c,d} | {b,c,d,e} | {a,b,c,e} |
| d | d | {1,b,c,d} | {0,1,a,b,c,d,e} | {1,a,b,d} | {b,c,d,e} | {a,d,e} | {1,c,d,e} |
| e | e | {1,a,d,e} | {a,c,d,e} | {0,1,a,b,c,d,e} | {a,b,c,e} | {1,c,d,e} | {1,b,e} |

| $HF_7^{97}$ | 0 | 1 | a | b | c | d | e |
|---|---|---|---|---|---|---|---|
| 0 | 0 | 1 | a | b | c | d | e |
| 1 | 1 | {1,a,c} | {1,a,b,e} | {1,a,b,d,e} | {0,1,a,b,c,d,e} | {1,b,c,d,e} | {1,a,d,e} |
| a | a | {1,a,b,e} | {a,b,d} | {1,a,b,c} | {1,a,b,c,e} | {0,1,a,b,c,d,e} | {1,a,c,d,e} |
| b | b | {1,a,b,d,e} | {1,a,b,c} | {b,c,e} | {a,b,c,d} | {1,a,b,c,d} | {0,1,a,b,c,d,e} |
| c | c | {0,1,a,b,c,d,e} | {1,a,b,c,e} | {a,b,c,d} | {1,c,d} | {b,c,d,e} | {a,b,c,d,e} |
| d | d | {1,b,c,d,e} | {0,1,a,b,c,d,e} | {1,a,b,c,d} | {b,c,d,e} | {a,d,e} | {1,c,d,e} |
| e | e | {1,a,d,e} | {1,a,c,d,e} | {0,1,a,b,c,d,e} | {a,b,c,d,e} | {1,c,d,e} | {1,b,e} |

| $HF_7^{98}$ | 0 | 1 | a | b | c | d | e |
|---|---|---|---|---|---|---|---|
| 0 | 0 | 1 | a | b | c | d | e |
| 1 | 1 | {1,a,d} | {1,a,c} | {1,b,e} | {0,1,a,b,c,d,e} | {1,c,d} | {1,b,e} |
| a | a | {1,a,c} | {a,b,e} | {a,b,d} | {1,a,c} | {0,1,a,b,c,d,e} | {a,d,e} |
| b | b | {1,b,e} | {a,b,d} | {1,b,c} | {b,c,e} | {a,b,d} | {0,1,a,b,c,d,e} |
| c | c | {0,1,a,b,c,d,e} | {1,a,c} | {b,c,e} | {a,c,d} | {1,c,d} | {b,c,e} |
| d | d | {1,c,d} | {0,1,a,b,c,d,e} | {a,b,d} | {1,c,d} | {b,d,e} | {a,d,e} |
| e | e | {1,b,e} | {a,d,e} | {0,1,a,b,c,d,e} | {b,c,e} | {a,d,e} | {1,c,e} |

| $HF_7^{99}$ | 0 | 1 | a | b | c | d | e |
|---|---|---|---|---|---|---|---|
| 0 | 0 | 1 | a | b | c | d | e |
| 1 | 1 | {1,a,d} | {1,a,c} | {1,a,b,e} | {0,1,a,b,c,d,e} | {1,c,d,e} | {1,b,e} |
| a | a | {1,a,c} | {a,b,e} | {a,b,d} | {1,a,b,c} | {0,1,a,b,c,d,e} | {1,a,d,e} |
| b | b | {1,a,b,e} | {a,b,d} | {1,b,c} | {b,c,e} | {a,b,c,d} | {0,1,a,b,c,d,e} |
| c | c | {0,1,a,b,c,d,e} | {1,a,b,c} | {b,c,e} | {a,c,d} | {1,c,d} | {b,c,d,e} |
| d | d | {1,c,d,e} | {0,1,a,b,c,d,e} | {a,b,c,d} | {1,c,d} | {b,d,e} | {a,d,e} |
| e | e | {1,b,e} | {1,a,d,e} | {0,1,a,b,c,d,e} | {b,c,d,e} | {a,d,e} | {1,c,e} |

| $HF_7^{100}$ | 0 | 1 | a | b | c | d | e |
|---|---|---|---|---|---|---|---|
| 0 | 0 | 1 | a | b | c | d | e |
| 1 | 1 | {1,a,d} | {1,a,b,c,e} | {1,b,d,e} | {0,1,a,b,c,d,e} | {1,b,c,d} | {1,a,b,d,e} |
| a | a | {1,a,b,c,e} | {a,b,e} | {1,a,b,c,d} | {1,a,c,e} | {0,1,a,b,c,d,e} | {a,c,d,e} |
| b | b | {1,b,d,e} | {1,a,b,c,d} | {1,b,c} | {a,b,c,d,e} | {1,a,b,d} | {0,1,a,b,c,d,e} |
| c | c | {0,1,a,b,c,d,e} | {1,a,c,e} | {a,b,c,d,e} | {a,c,d} | {1,b,c,d,e} | {a,b,c,e} |
| d | d | {1,b,c,d} | {0,1,a,b,c,d,e} | {1,a,b,d} | {1,b,c,d,e} | {b,d,e} | {1,a,c,d,e} |
| e | e | {1,a,b,d,e} | {a,c,d,e} | {0,1,a,b,c,d,e} | {a,b,c,e} | {1,a,c,d,e} | {1,c,e} |

| $HF_7^{101}$ | 0 | 1 | a | b | c | d | e |
|---|---|---|---|---|---|---|---|
| 0 | 0 | 1 | a | b | c | d | e |
| 1 | 1 | {1,a,d} | {1,a,b,c,e} | {1,a,b,d,e} | {0,1,a,b,c,d,e} | {1,b,c,d,e} | {1,a,b,d,e} |
| a | a | {1,a,b,c,e} | {a,b,e} | {1,a,b,c,d} | {1,a,b,c,e} | {0,1,a,b,c,d,e} | {1,a,c,d,e} |
| b | b | {1,a,b,d,e} | {1,a,b,c,d} | {1,b,c} | {a,b,c,d,e} | {1,a,b,c,d} | {0,1,a,b,c,d,e} |
| c | c | {0,1,a,b,c,d,e} | {1,a,b,c,e} | {a,b,c,d,e} | {a,c,d} | {1,b,c,d,e} | {a,b,c,d,e} |
| d | d | {1,b,c,d,e} | {0,1,a,b,c,d,e} | {1,a,b,c,d} | {1,b,c,d,e} | {b,d,e} | {1,a,c,d,e} |
| e | e | {1,a,b,d,e} | {1,a,c,d,e} | {0,1,a,b,c,d,e} | {a,b,c,d,e} | {1,a,c,d,e} | {1,c,e} |

| $HF_7^{102}$ | 0 | 1 | a | b | c | d | e |
|---|---|---|---|---|---|---|---|
| 0 | 0 | 1 | a | b | c | d | e |
| 1 | 1 | {1,a,e} | {1,a,b,e} | {1,b,c,d,e} | {0,1,a,b,c,d,e} | {1,a,b,c,d} | {1,a,d,e} |
| a | a | {1,a,b,e} | {1,a,b} | {1,a,b,c} | {1,a,c,d,e} | {0,1,a,b,c,d,e} | {a,b,c,d,e} |
| b | b | {1,b,c,d,e} | {1,a,b,c} | {a,b,c} | {a,b,c,d} | {1,a,b,d,e} | {0,1,a,b,c,d,e} |
| c | c | {0,1,a,b,c,d,e} | {1,a,c,d,e} | {a,b,c,d} | {b,c,d} | {b,c,d,e} | {1,a,b,c,e} |
| d | d | {1,a,b,c,d} | {0,1,a,b,c,d,e} | {1,a,b,d,e} | {b,c,d,e} | {c,d,e} | {1,c,d,e} |
| e | e | {1,a,d,e} | {a,b,c,d,e} | {0,1,a,b,c,d,e} | {1,a,b,c,e} | {1,c,d,e} | {1,d,e} |

| $HF_7^{103}$ | 0 | 1 | a | b | c | d | e |
|---|---|---|---|---|---|---|---|
| 0 | 0 | 1 | a | b | c | d | e |
| 1 | 1 | {1,a,e} | {1,a,b,e} | {1,a,b,c,d,e} | {0,1,a,b,c,d,e} | {1,a,b,c,d,e} | {1,a,d,e} |
| a | a | {1,a,b,e} | {1,a,b} | {1,a,b,c} | {1,a,b,c,d,e} | {0,1,a,b,c,d,e} | {1,a,b,c,d,e} |
| b | b | {1,a,b,c,d,e} | {1,a,b,c} | {a,b,c} | {a,b,c,d} | {1,a,b,c,d,e} | {0,1,a,b,c,d,e} |
| c | c | {0,1,a,b,c,d,e} | {1,a,b,c,d,e} | {a,b,c,d} | {b,c,d} | {b,c,d,e} | {1,a,b,c,d,e} |
| d | d | {1,a,b,c,d,e} | {0,1,a,b,c,d,e} | {1,a,b,c,d,e} | {b,c,d,e} | {c,d,e} | {1,c,d,e} |
| e | e | {1,a,d,e} | {1,a,b,c,d,e} | {0,1,a,b,c,d,e} | {1,a,b,c,d,e} | {1,c,d,e} | {1,d,e} |

| $HF_7^{104}$ | 0 | 1 | a | b | c | d | e |
|---|---|---|---|---|---|---|---|
| 0 | 0 | 1 | a | b | c | d | e |
| 1 | 1 | {1,b,c} | {1,a,b,d,e} | {1,a,b,d} | {0,1,a,b,c,d,e} | {1,b,d,e} | {1,a,c,d,e} |
| a | a | {1,a,b,d,e} | {a,c,d} | {1,a,b,c,e} | {a,b,c,e} | {0,1,a,b,c,d,e} | {1,a,c,e} |
| b | b | {1,a,b,d} | {1,a,b,c,e} | {b,d,e} | {1,a,b,c,d} | {1,b,c,d} | {0,1,a,b,c,d,e} |
| c | c | {0,1,a,b,c,d,e} | {a,b,c,e} | {1,a,b,c,d} | {1,c,e} | {a,b,c,d,e} | {a,c,d,e} |
| d | d | {1,b,d,e} | {0,1,a,b,c,d,e} | {1,b,c,d} | {a,b,c,d,e} | {1,a,d} | {1,b,c,d,e} |
| e | e | {1,a,c,d,e} | {1,a,c,e} | {0,1,a,b,c,d,e} | {a,c,d,e} | {1,b,c,d,e} | {a,b,e} |

| $HF_7^{105}$ | 0 | 1 | a | b | c | d | e |
|---|---|---|---|---|---|---|---|
| 0 | 0 | 1 | a | b | c | d | e |
| 1 | 1 | {1,b,d} | {1,a,c,d} | {1,b} | {0,1,a,b,c,d,e} | {1,d} | {1,b,c,e} |
| a | a | {1,a,c,d} | {a,c,e} | {a,b,d,e} | {a,c} | {0,1,a,b,c,d,e} | {a,e} |
| b | b | {1,b} | {a,b,d,e} | {1,b,d} | {1,b,c,e} | {b,d} | {0,1,a,b,c,d,e} |
| c | c | {0,1,a,b,c,d,e} | {a,c} | {1,b,c,e} | {a,c,e} | {1,a,c,d} | {c,e} |
| d | d | {1,d} | {0,1,a,b,c,d,e} | {b,d} | {1,a,c,d} | {1,b,d} | {a,b,d,e} |
| e | e | {1,b,c,e} | {a,e} | {0,1,a,b,c,d,e} | {c,e} | {a,b,d,e} | {a,c,e} |

| $HF_7^{106}$ | 0 | 1 | a | b | c | d | e |
|---|---|---|---|---|---|---|---|
| 0 | 0 | 1 | a | b | c | d | e |
| 1 | 1 | {1,b,d} | {1,a,c,d} | {1,a,b} | {0,1,a,b,c,d,e} | {1,d,e} | {1,b,c,e} |
| a | a | {1,a,c,d} | {a,c,e} | {a,b,d,e} | {a,b,c} | {0,1,a,b,c,d,e} | {1,a,e} |
| b | b | {1,a,b} | {a,b,d,e} | {1,b,d} | {1,b,c,e} | {b,c,d} | {0,1,a,b,c,d,e} |
| c | c | {0,1,a,b,c,d,e} | {a,b,c} | {1,b,c,e} | {a,c,e} | {1,a,c,d} | {c,d,e} |
| d | d | {1,d,e} | {0,1,a,b,c,d,e} | {b,c,d} | {1,a,c,d} | {1,b,d} | {a,b,d,e} |
| e | e | {1,b,c,e} | {1,a,e} | {0,1,a,b,c,d,e} | {c,d,e} | {a,b,d,e} | {a,c,e} |

| $HF_7^{107}$ | 0 | 1 | a | b | c | d | e |
|---|---|---|---|---|---|---|---|
| 0 | 0 | 1 | a | b | c | d | e |
| 1 | 1 | {1,b,d} | {1,a,b,c,d,e} | {1,b,d} | {0,1,a,b,c,d,e} | {1,b,d} | {1,a,b,c,d,e} |
| a | a | {1,a,b,c,d,e} | {a,c,e} | {1,a,b,c,d,e} | {a,c,e} | {0,1,a,b,c,d,e} | {a,c,e} |
| b | b | {1,b,d} | {1,a,b,c,d,e} | {1,b,d} | {1,a,b,c,d,e} | {1,b,d} | {0,1,a,b,c,d,e} |
| c | c | {0,1,a,b,c,d,e} | {a,c,e} | {1,a,b,c,d,e} | {a,c,e} | {1,a,b,c,d,e} | {a,c,e} |
| d | d | {1,b,d} | {0,1,a,b,c,d,e} | {1,b,d} | {1,a,b,c,d,e} | {1,b,d} | {1,a,b,c,d,e} |
| e | e | {1,a,b,c,d,e} | {a,c,e} | {0,1,a,b,c,d,e} | {a,c,e} | {1,a,b,c,d,e} | {a,c,e} |

| $HF_7^{108}$ | 0 | 1 | a | b | c | d | e |
|---|---|---|---|---|---|---|---|
| 0 | 0 | 1 | a | b | c | d | e |
| 1 | 1 | {1,a,b,c} | {1,a,d} | {1,a,b,e} | {0,1,a,b,c,d,e} | {1,c,d,e} | {1,c,e} |
| a | a | {1,a,d} | {a,b,c,d} | {a,b,e} | {1,a,b,c} | {0,1,a,b,c,d,e} | {1,a,d,e} |
| b | b | {1,a,b,e} | {a,b,e} | {b,c,d,e} | {1,b,c} | {a,b,c,d} | {0,1,a,b,c,d,e} |
| c | c | {0,1,a,b,c,d,e} | {1,a,b,c} | {1,b,c} | {1,c,d,e} | {a,c,d} | {b,c,d,e} |
| d | d | {1,c,d,e} | {0,1,a,b,c,d,e} | {a,b,c,d} | {a,c,d} | {1,a,d,e} | {b,d,e} |
| e | e | {1,c,e} | {1,a,d,e} | {0,1,a,b,c,d,e} | {b,c,d,e} | {b,d,e} | {1,a,b,e} |

| $HF_7^{109}$ | 0 | 1 | a | b | c | d | e |
|---|---|---|---|---|---|---|---|
| 0 | 0 | 1 | a | b | c | d | e |
| 1 | 1 | {1,a,b,c} | {1,a,b,d,e} | {1,b,d,e} | {0,1,a,b,c,d,e} | {1,b,c,d} | {1,a,c,d,e} |
| a | a | {1,a,b,d,e} | {a,b,c,d} | {1,a,b,c,e} | {1,a,c,e} | {0,1,a,b,c,d,e} | {a,c,d,e} |
| b | b | {1,b,d,e} | {1,a,b,c,e} | {b,c,d,e} | {1,a,b,c,d} | {1,a,b,d} | {0,1,a,b,c,d,e} |
| c | c | {0,1,a,b,c,d,e} | {1,a,c,e} | {1,a,b,c,d} | {1,c,d,e} | {a,b,c,d,e} | {a,b,c,e} |
| d | d | {1,b,c,d} | {0,1,a,b,c,d,e} | {1,a,b,d} | {a,b,c,d,e} | {1,a,d,e} | {1,b,c,d,e} |
| e | e | {1,a,c,d,e} | {a,c,d,e} | {0,1,a,b,c,d,e} | {a,b,c,e} | {1,b,c,d,e} | {1,a,b,e} |

| $HF_7^{110}$ | 0 | 1 | a | b | c | d | e |
|---|---|---|---|---|---|---|---|
| 0 | 0 | 1 | a | b | c | d | e |
| 1 | 1 | {1,a,b,c} | {1,a,b,d,e} | {1,a,b,d,e} | {0,1,a,b,c,d,e} | {1,b,c,d,e} | {1,a,c,d,e} |
| a | a | {1,a,b,d,e} | {a,b,c,d} | {1,a,b,c,e} | {1,a,b,c,e} | {0,1,a,b,c,d,e} | {1,a,c,d,e} |
| b | b | {1,a,b,d,e} | {1,a,b,c,e} | {b,c,d,e} | {1,a,b,c,d} | {1,a,b,c,d} | {0,1,a,b,c,d,e} |
| c | c | {0,1,a,b,c,d,e} | {1,a,b,c,e} | {1,a,b,c,d} | {1,c,d,e} | {a,b,c,d,e} | {a,b,c,d,e} |
| d | d | {1,b,c,d,e} | {0,1,a,b,c,d,e} | {1,a,b,c,d} | {a,b,c,d,e} | {1,a,d,e} | {1,b,c,d,e} |
| e | e | {1,a,c,d,e} | {1,a,c,d,e} | {0,1,a,b,c,d,e} | {a,b,c,d,e} | {1,b,c,d,e} | {1,a,b,e} |

| $HF_7^{111}$ | 0 | 1 | a | b | c | d | e |
|---|---|---|---|---|---|---|---|
| 0 | 0 | 1 | a | b | c | d | e |
| 1 | 1 | {1,a,b,d} | {1,a,c,d} | {1,b,e} | {0,1,a,b,c,d,e} | {1,c,d} | {1,b,c,e} |
| a | a | {1,a,c,d} | {a,b,c,e} | {a,b,d,e} | {1,a,c} | {0,1,a,b,c,d,e} | {a,d,e} |
| b | b | {1,b,e} | {a,b,d,e} | {1,b,c,d} | {1,b,c,e} | {a,b,d} | {0,1,a,b,c,d,e} |
| c | c | {0,1,a,b,c,d,e} | {1,a,c} | {1,b,c,e} | {a,c,d,e} | {1,a,c,d} | {b,c,e} |
| d | d | {1,c,d} | {0,1,a,b,c,d,e} | {a,b,d} | {1,a,c,d} | {1,b,d,e} | {a,b,d,e} |
| e | e | {1,b,c,e} | {a,d,e} | {0,1,a,b,c,d,e} | {b,c,e} | {a,b,d,e} | {1,a,c,e} |

| $HF_7^{112}$ | 0 | 1 | a | b | c | d | e |
|---|---|---|---|---|---|---|---|
| 0 | 0 | 1 | a | b | c | d | e |
| 1 | 1 | {1,a,b,d} | {1,a,c,d} | {1,a,b,e} | {0,1,a,b,c,d,e} | {1,c,d,e} | {1,b,c,e} |
| a | a | {1,a,c,d} | {a,b,c,e} | {a,b,d,e} | {1,a,b,c} | {0,1,a,b,c,d,e} | {1,a,d,e} |
| b | b | {1,a,b,e} | {a,b,d,e} | {1,b,c,d} | {1,b,c,e} | {a,b,c,d} | {0,1,a,b,c,d,e} |
| c | c | {0,1,a,b,c,d,e} | {1,a,b,c} | {1,b,c,e} | {a,c,d,e} | {1,a,c,d} | {b,c,d,e} |
| d | d | {1,c,d,e} | {0,1,a,b,c,d,e} | {a,b,c,d} | {1,a,c,d} | {1,b,d,e} | {a,b,d,e} |
| e | e | {1,b,c,e} | {1,a,d,e} | {0,1,a,b,c,d,e} | {b,c,d,e} | {a,b,d,e} | {1,a,c,e} |

| $HF_7^{113}$ | 0 | 1 | a | b | c | d | e |
|---|---|---|---|---|---|---|---|
| 0 | 0 | 1 | a | b | c | d | e |
| 1 | 1 | {1,a,b,d} | {1,a,b,c,d,e} | {1,b,d,e} | {0,1,a,b,c,d,e} | {1,b,c,d} | {1,a,b,c,d,e} |
| a | a | {1,a,b,c,d,e} | {a,b,c,e} | {1,a,b,c,d,e} | {1,a,c,e} | {0,1,a,b,c,d,e} | {a,c,d,e} |
| b | b | {1,b,d,e} | {1,a,b,c,d,e} | {1,b,c,d} | {1,a,b,c,d,e} | {1,a,b,d} | {0,1,a,b,c,d,e} |
| c | c | {0,1,a,b,c,d,e} | {1,a,c,e} | {1,a,b,c,d,e} | {a,c,d,e} | {1,a,b,c,d,e} | {a,b,c,e} |
| d | d | {1,b,c,d} | {0,1,a,b,c,d,e} | {1,a,b,d} | {1,a,b,c,d,e} | {1,b,d,e} | {1,a,b,c,d,e} |
| e | e | {1,a,b,c,d,e} | {a,c,d,e} | {0,1,a,b,c,d,e} | {a,b,c,e} | {1,a,b,c,d,e} | {1,a,c,e} |

| $HF_7^{114}$ | 0 | 1 | a | b | c | d | e |
|---|---|---|---|---|---|---|---|
| 0 | 0 | 1 | a | b | c | d | e |
| 1 | 1 | {1,a,b,d} | {1,a,b,c,d,e} | {1,a,b,d,e} | {0,1,a,b,c,d,e} | {1,b,c,d,e} | {1,a,b,c,d,e} |
| a | a | {1,a,b,c,d,e} | {a,b,c,e} | {1,a,b,c,d,e} | {1,a,b,c,e} | {0,1,a,b,c,d,e} | {1,a,c,d,e} |
| b | b | {1,a,b,d,e} | {1,a,b,c,d,e} | {1,b,c,d} | {1,a,b,c,d,e} | {1,a,b,c,d} | {0,1,a,b,c,d,e} |
| c | c | {0,1,a,b,c,d,e} | {1,a,b,c,e} | {1,a,b,c,d,e} | {a,c,d,e} | {1,a,b,c,d,e} | {a,b,c,d,e} |
| d | d | {1,b,c,d,e} | {0,1,a,b,c,d,e} | {1,a,b,c,d} | {1,a,b,c,d,e} | {1,b,d,e} | {1,a,b,c,d,e} |
| e | e | {1,a,b,c,d,e} | {1,a,c,d,e} | {0,1,a,b,c,d,e} | {a,b,c,d,e} | {1,a,b,c,d,e} | {1,a,c,e} |

| $HF_7^{115}$ | 0 | 1 | a | b | c | d | e |
|---|---|---|---|---|---|---|---|
| 0 | 0 | 1 | a | b | c | d | e |
| 1 | 1 | {1,a,b,e} | {1,a,d} | {1,b,c,e} | {0,1,a,b,c,d,e} | {1,a,c,d} | {1,c,e} |
| a | a | {1,a,d} | {1,a,b,c} | {a,b,e} | {1,a,c,d} | {0,1,a,b,c,d,e} | {a,b,d,e} |
| b | b | {1,b,c,e} | {a,b,e} | {a,b,c,d} | {1,b,c} | {a,b,d,e} | {0,1,a,b,c,d,e} |
| c | c | {0,1,a,b,c,d,e} | {1,a,c,d} | {1,b,c} | {b,c,d,e} | {a,c,d} | {1,b,c,e} |
| d | d | {1,a,c,d} | {0,1,a,b,c,d,e} | {a,b,d,e} | {a,c,d} | {1,c,d,e} | {b,d,e} |
| e | e | {1,c,e} | {a,b,d,e} | {0,1,a,b,c,d,e} | {1,b,c,e} | {b,d,e} | {1,a,d,e} |

| $HF_7^{116}$ | 0 | 1 | a | b | c | d | e |
|---|---|---|---|---|---|---|---|
| 0 | 0 | 1 | a | b | c | d | e |
| 1 | 1 | {1,a,b,e} | {1,a,d} | {1,a,b,c,e} | {0,1,a,b,c,d,e} | {1,a,c,d,e} | {1,c,e} |
| a | a | {1,a,d} | {1,a,b,c} | {a,b,e} | {1,a,b,c,d} | {0,1,a,b,c,d,e} | {1,a,b,d,e} |
| b | b | {1,a,b,c,e} | {a,b,e} | {a,b,c,d} | {1,b,c} | {a,b,c,d,e} | {0,1,a,b,c,d,e} |
| c | c | {0,1,a,b,c,d,e} | {1,a,b,c,d} | {1,b,c} | {b,c,d,e} | {a,c,d} | {1,b,c,d,e} |
| d | d | {1,a,c,d,e} | {0,1,a,b,c,d,e} | {a,b,c,d,e} | {a,c,d} | {1,c,d,e} | {b,d,e} |
| e | e | {1,c,e} | {1,a,b,d,e} | {0,1,a,b,c,d,e} | {1,b,c,d,e} | {b,d,e} | {1,a,d,e} |

| $HF_7^{117}$ | 0 | 1 | a | b | c | d | e |
|---|---|---|---|---|---|---|---|
| 0 | 0 | 1 | a | b | c | d | e |
| 1 | 1 | {1,a,b,e} | {1,a,b,d,e} | {1,b,c,d,e} | {0,1,a,b,c,d,e} | {1,a,b,c,d} | {1,a,c,d,e} |
| a | a | {1,a,b,d,e} | {1,a,b,c} | {1,a,b,c,e} | {1,a,c,d,e} | {0,1,a,b,c,d,e} | {a,b,c,d,e} |
| b | b | {1,b,c,d,e} | {1,a,b,c,e} | {a,b,c,d} | {1,a,b,c,d} | {1,a,b,d,e} | {0,1,a,b,c,d,e} |
| c | c | {0,1,a,b,c,d,e} | {1,a,c,d,e} | {1,a,b,c,d} | {b,c,d,e} | {a,b,c,d,e} | {1,a,b,c,e} |
| d | d | {1,a,b,c,d} | {0,1,a,b,c,d,e} | {1,a,b,d,e} | {a,b,c,d,e} | {1,c,d,e} | {1,b,c,d,e} |
| e | e | {1,a,c,d,e} | {a,b,c,d,e} | {0,1,a,b,c,d,e} | {1,a,b,c,e} | {1,b,c,d,e} | {1,a,d,e} |

| $HF_7^{118}$ | 0 | 1 | a | b | c | d | e |
|---|---|---|---|---|---|---|---|
| 0 | 0 | 1 | a | b | c | d | e |
| 1 | 1 | {1,a,b,e} | {1,a,b,d,e} | {1,a,b,c,d,e} | {0,1,a,b,c,d,e} | {1,a,b,c,d,e} | {1,a,c,d,e} |
| a | a | {1,a,b,d,e} | {1,a,b,c} | {1,a,b,c,e} | {1,a,b,c,d,e} | {0,1,a,b,c,d,e} | {1,a,b,c,d,e} |
| b | b | {1,a,b,c,d,e} | {1,a,b,c,e} | {a,b,c,d} | {1,a,b,c,d} | {1,a,b,c,d,e} | {0,1,a,b,c,d,e} |
| c | c | {0,1,a,b,c,d,e} | {1,a,b,c,d,e} | {1,a,b,c,d} | {b,c,d,e} | {a,b,c,d,e} | {1,a,b,c,d,e} |
| d | d | {1,a,b,c,d,e} | {0,1,a,b,c,d,e} | {1,a,b,c,d,e} | {a,b,c,d,e} | {1,c,d,e} | {1,b,c,d,e} |
| e | e | {1,a,c,d,e} | {1,a,b,c,d,e} | {0,1,a,b,c,d,e} | {1,a,b,c,d,e} | {1,b,c,d,e} | {1,a,d,e} |

| $HF_7^{119}$ | 0 | 1 | a | b | c | d | e |
|---|---|---|---|---|---|---|---|
| 0 | 0 | 1 | a | b | c | d | e |
| 1 | 1 | {1,a,c,d} | {1,a,c} | {1,b,e} | {0,1,a,b,c,d,e} | {1,c,d} | {1,b,e} |
| a | a | {1,a,c} | {a,b,d,e} | {a,b,d} | {1,a,c} | {0,1,a,b,c,d,e} | {a,d,e} |
| b | b | {1,b,e} | {a,b,d} | {1,b,c,e} | {b,c,e} | {a,b,d} | {0,1,a,b,c,d,e} |
| c | c | {0,1,a,b,c,d,e} | {1,a,c} | {b,c,e} | {1,a,c,d} | {1,c,d} | {b,c,e} |
| d | d | {1,c,d} | {0,1,a,b,c,d,e} | {a,b,d} | {1,c,d} | {a,b,d,e} | {a,d,e} |
| e | e | {1,b,e} | {a,d,e} | {0,1,a,b,c,d,e} | {b,c,e} | {a,d,e} | {1,b,c,e} |

| $HF_7^{120}$ | 0 | 1 | a | b | c | d | e |
|---|---|---|---|---|---|---|---|
| 0 | 0 | 1 | a | b | c | d | e |
| 1 | 1 | {1,a,c,d} | {1,a,c} | {1,a,b,e} | {0,1,a,b,c,d,e} | {1,c,d,e} | {1,b,e} |
| a | a | {1,a,c} | {a,b,d,e} | {a,b,d} | {1,a,b,c} | {0,1,a,b,c,d,e} | {1,a,d,e} |
| b | b | {1,a,b,e} | {a,b,d} | {1,b,c,e} | {b,c,e} | {a,b,c,d} | {0,1,a,b,c,d,e} |
| c | c | {0,1,a,b,c,d,e} | {1,a,b,c} | {b,c,e} | {1,a,c,d} | {1,c,d} | {b,c,d,e} |
| d | d | {1,c,d,e} | {0,1,a,b,c,d,e} | {a,b,c,d} | {1,c,d} | {a,b,d,e} | {a,d,e} |
| e | e | {1,b,e} | {1,a,d,e} | {0,1,a,b,c,d,e} | {b,c,d,e} | {a,d,e} | {1,b,c,e} |

| $HF_7^{121}$ | 0 | 1 | a | b | c | d | e |
|---|---|---|---|---|---|---|---|
| 0 | 0 | 1 | a | b | c | d | e |
| 1 | 1 | {1,a,c,d} | {1,a,b,c,e} | {1,b,d,e} | {0,1,a,b,c,d,e} | {1,b,c,d} | {1,a,b,d,e} |
| a | a | {1,a,b,c,e} | {a,b,d,e} | {1,a,b,c,d} | {1,a,c,e} | {0,1,a,b,c,d,e} | {a,c,d,e} |
| b | b | {1,b,d,e} | {1,a,b,c,d} | {1,b,c,e} | {a,b,c,d,e} | {1,a,b,d} | {0,1,a,b,c,d,e} |
| c | c | {0,1,a,b,c,d,e} | {1,a,c,e} | {a,b,c,d,e} | {1,a,c,d} | {1,b,c,d,e} | {a,b,c,e} |
| d | d | {1,b,c,d} | {0,1,a,b,c,d,e} | {1,a,b,d} | {1,b,c,d,e} | {a,b,d,e} | {1,a,c,d,e} |
| e | e | {1,a,b,d,e} | {a,c,d,e} | {0,1,a,b,c,d,e} | {a,b,c,e} | {1,a,c,d,e} | {1,b,c,e} |

| $HF_7^{122}$ | 0 | 1 | a | b | c | d | e |
|---|---|---|---|---|---|---|---|
| 0 | 0 | 1 | a | b | c | d | e |
| 1 | 1 | {1,a,c,d} | {1,a,b,c,e} | {1,a,b,d,e} | {0,1,a,b,c,d,e} | {1,b,c,d,e} | {1,a,b,d,e} |
| a | a | {1,a,b,c,e} | {a,b,d,e} | {1,a,b,c,d} | {1,a,b,c,e} | {0,1,a,b,c,d,e} | {1,a,c,d,e} |
| b | b | {1,a,b,d,e} | {1,a,b,c,d} | {1,b,c,e} | {a,b,c,d,e} | {1,a,b,c,d} | {0,1,a,b,c,d,e} |
| c | c | {0,1,a,b,c,d,e} | {1,a,b,c,e} | {a,b,c,d,e} | {1,a,c,d} | {1,b,c,d,e} | {a,b,c,d,e} |
| d | d | {1,b,c,d,e} | {0,1,a,b,c,d,e} | {1,a,b,c,d} | {1,b,c,d,e} | {a,b,d,e} | {1,a,c,d,e} |
| e | e | {1,a,b,d,e} | {1,a,c,d,e} | {0,1,a,b,c,d,e} | {a,b,c,d,e} | {1,a,c,d,e} | {1,b,c,e} |

| $HF_7^{123}$ | 0 | 1 | a | b | c | d | e |
|---|---|---|---|---|---|---|---|
| 0 | 0 | 1 | a | b | c | d | e |
| 1 | 1 | {1,a,c,e} | {1,a,b,e} | {1,b,c,d,e} | {0,1,a,b,c,d,e} | {1,a,b,c,d} | {1,a,d,e} |
| a | a | {1,a,b,e} | {1,a,b,d} | {1,a,b,c} | {1,a,c,d,e} | {0,1,a,b,c,d,e} | {a,b,c,d,e} |
| b | b | {1,b,c,d,e} | {1,a,b,c} | {a,b,c,e} | {a,b,c,d} | {1,a,b,d,e} | {0,1,a,b,c,d,e} |
| c | c | {0,1,a,b,c,d,e} | {1,a,c,d,e} | {a,b,c,d} | {1,b,c,d} | {b,c,d,e} | {1,a,b,c,e} |
| d | d | {1,a,b,c,d} | {0,1,a,b,c,d,e} | {1,a,b,d,e} | {b,c,d,e} | {a,c,d,e} | {1,c,d,e} |
| e | e | {1,a,d,e} | {a,b,c,d,e} | {0,1,a,b,c,d,e} | {1,a,b,c,e} | {1,c,d,e} | {1,b,d,e} |

| $HF_7^{124}$ | 0 | 1 | a | b | c | d | e |
|---|---|---|---|---|---|---|---|
| 0 | 0 | 1 | a | b | c | d | e |
| 1 | 1 | {1,a,c,e} | {1,a,b,e} | {1,a,b,c,d,e} | {0,1,a,b,c,d,e} | {1,a,b,c,d,e} | {1,a,d,e} |
| a | a | {1,a,b,e} | {1,a,b,d} | {1,a,b,c} | {1,a,b,c,d,e} | {0,1,a,b,c,d,e} | {1,a,b,c,d,e} |
| b | b | {1,a,b,c,d,e} | {1,a,b,c} | {a,b,c,e} | {a,b,c,d} | {1,a,b,c,d,e} | {0,1,a,b,c,d,e} |
| c | c | {0,1,a,b,c,d,e} | {1,a,b,c,d,e} | {a,b,c,d} | {1,b,c,d} | {b,c,d,e} | {1,a,b,c,d,e} |
| d | d | {1,a,b,c,d,e} | {0,1,a,b,c,d,e} | {1,a,b,c,d,e} | {b,c,d,e} | {a,c,d,e} | {1,c,d,e} |
| e | e | {1,a,d,e} | {1,a,b,c,d,e} | {0,1,a,b,c,d,e} | {1,a,b,c,d,e} | {1,c,d,e} | {1,b,d,e} |

| $HF_7^{125}$ | 0 | 1 | a | b | c | d | e |
|---|---|---|---|---|---|---|---|
| 0 | 0 | 1 | a | b | c | d | e |
| 1 | 1 | {1,b,c,d} | {1,a,b,c,d,e} | {1,a,b,d} | {0,1,a,b,c,d,e} | {1,b,d,e} | {1,a,b,c,d,e} |
| a | a | {1,a,b,c,d,e} | {a,c,d,e} | {1,a,b,c,d,e} | {a,b,c,e} | {0,1,a,b,c,d,e} | {1,a,c,e} |
| b | b | {1,a,b,d} | {1,a,b,c,d,e} | {1,b,d,e} | {1,a,b,c,d,e} | {1,b,c,d} | {0,1,a,b,c,d,e} |
| c | c | {0,1,a,b,c,d,e} | {a,b,c,e} | {1,a,b,c,d,e} | {1,a,c,e} | {1,a,b,c,d,e} | {a,c,d,e} |
| d | d | {1,b,d,e} | {0,1,a,b,c,d,e} | {1,b,c,d} | {1,a,b,c,d,e} | {1,a,b,d} | {1,a,b,c,d,e} |
| e | e | {1,a,b,c,d,e} | {1,a,c,e} | {0,1,a,b,c,d,e} | {a,c,d,e} | {1,a,b,c,d,e} | {a,b,c,e} |

| $HF_7^{126}$ | 0 | 1 | a | b | c | d | e |
|---|---|---|---|---|---|---|---|
| 0 | 0 | 1 | a | b | c | d | e |
| 1 | 1 | {1,a,b,c,d} | {1,a,c,d} | {1,b,e} | {0,1,a,b,c,d,e} | {1,c,d} | {1,b,c,e} |
| a | a | {1,a,c,d} | {a,b,c,d,e} | {a,b,d,e} | {1,a,c} | {0,1,a,b,c,d,e} | {a,d,e} |
| b | b | {1,b,e} | {a,b,d,e} | {1,b,c,d,e} | {1,b,c,e} | {a,b,d} | {0,1,a,b,c,d,e} |
| c | c | {0,1,a,b,c,d,e} | {1,a,c} | {1,b,c,e} | {1,a,c,d,e} | {1,a,c,d} | {b,c,e} |
| d | d | {1,c,d} | {0,1,a,b,c,d,e} | {a,b,d} | {1,a,c,d} | {1,a,b,d,e} | {a,b,d,e} |
| e | e | {1,b,c,e} | {a,d,e} | {0,1,a,b,c,d,e} | {b,c,e} | {a,b,d,e} | {1,a,b,c,e} |

| $HF_7^{127}$ | 0 | 1 | a | b | c | d | e |
|---|---|---|---|---|---|---|---|
| 0 | 0 | 1 | a | b | c | d | e |
| 1 | 1 | {1,a,b,c,d} | {1,a,c,d} | {1,a,b,e} | {0,1,a,b,c,d,e} | {1,c,d,e} | {1,b,c,e} |
| a | a | {1,a,c,d} | {a,b,c,d,e} | {a,b,d,e} | {1,a,b,c} | {0,1,a,b,c,d,e} | {1,a,d,e} |
| b | b | {1,a,b,e} | {a,b,d,e} | {1,b,c,d,e} | {1,b,c,e} | {a,b,c,d} | {0,1,a,b,c,d,e} |
| c | c | {0,1,a,b,c,d,e} | {1,a,b,c} | {1,b,c,e} | {1,a,c,d,e} | {1,a,c,d} | {b,c,d,e} |
| d | d | {1,c,d,e} | {0,1,a,b,c,d,e} | {a,b,c,d} | {1,a,c,d} | {1,a,b,d,e} | {a,b,d,e} |
| e | e | {1,b,c,e} | {1,a,d,e} | {0,1,a,b,c,d,e} | {b,c,d,e} | {a,b,d,e} | {1,a,b,c,e} |

| $HF_7^{128}$ | 0 | 1 | a | b | c | d | e |
|---|---|---|---|---|---|---|---|
| 0 | 0 | 1 | a | b | c | d | e |
| 1 | 1 | {1,a,b,c,d} | {1,a,b,c,d,e} | {1,b,d,e} | {0,1,a,b,c,d,e} | {1,b,c,d} | {1,a,b,c,d,e} |
| a | a | {1,a,b,c,d,e} | {a,b,c,d,e} | {1,a,b,c,d,e} | {1,a,c,e} | {0,1,a,b,c,d,e} | {a,c,d,e} |
| b | b | {1,b,d,e} | {1,a,b,c,d,e} | {1,b,c,d,e} | {1,a,b,c,d,e} | {1,a,b,d} | {0,1,a,b,c,d,e} |
| c | c | {0,1,a,b,c,d,e} | {1,a,c,e} | {1,a,b,c,d,e} | {1,a,c,d,e} | {1,a,b,c,d,e} | {a,b,c,e} |
| d | d | {1,b,c,d} | {0,1,a,b,c,d,e} | {1,a,b,d} | {1,a,b,c,d,e} | {1,a,b,d,e} | {1,a,b,c,d,e} |
| e | e | {1,a,b,c,d,e} | {a,c,d,e} | {0,1,a,b,c,d,e} | {a,b,c,e} | {1,a,b,c,d,e} | {1,a,b,c,e} |

| $HF_7^{129}$ | 0 | 1 | a | b | c | d | e |
|---|---|---|---|---|---|---|---|
| 0 | 0 | 1 | a | b | c | d | e |
| 1 | 1 | {1,a,b,c,d} | {1,a,b,c,d,e} | {1,a,b,d,e} | {0,1,a,b,c,d,e} | {1,b,c,d,e} | {1,a,b,c,d,e} |
| a | a | {1,a,b,c,d,e} | {a,b,c,d,e} | {1,a,b,c,d,e} | {1,a,b,c,e} | {0,1,a,b,c,d,e} | {1,a,c,d,e} |
| b | b | {1,a,b,d,e} | {1,a,b,c,d,e} | {1,b,c,d,e} | {1,a,b,c,d,e} | {1,a,b,c,d} | {0,1,a,b,c,d,e} |
| c | c | {0,1,a,b,c,d,e} | {1,a,b,c,e} | {1,a,b,c,d,e} | {1,a,c,d,e} | {1,a,b,c,d,e} | {a,b,c,d,e} |
| d | d | {1,b,c,d,e} | {0,1,a,b,c,d,e} | {1,a,b,c,d} | {1,a,b,c,d,e} | {1,a,b,d,e} | {1,a,b,c,d,e} |
| e | e | {1,a,b,c,d,e} | {1,a,c,d,e} | {0,1,a,b,c,d,e} | {a,b,c,d,e} | {1,a,b,c,d,e} | {1,a,b,c,e} |

| $HF_7^{130}$ | 0 | 1 | a | b | c | d | e |
|---|---|---|---|---|---|---|---|
| 0 | 0 | 1 | a | b | c | d | e |
| 1 | 1 | {1,a,b,c,e} | {1,a,d} | {1,b,c,e} | {0,1,a,b,c,d,e} | {1,a,c,d} | {1,c,e} |
| a | a | {1,a,d} | {1,a,b,c,d} | {a,b,e} | {1,a,c,d} | {0,1,a,b,c,d,e} | {a,b,d,e} |
| b | b | {1,b,c,e} | {a,b,e} | {a,b,c,d,e} | {1,b,c} | {a,b,d,e} | {0,1,a,b,c,d,e} |
| c | c | {0,1,a,b,c,d,e} | {1,a,c,d} | {1,b,c} | {1,b,c,d,e} | {a,c,d} | {1,b,c,e} |
| d | d | {1,a,c,d} | {0,1,a,b,c,d,e} | {a,b,d,e} | {a,c,d} | {1,a,c,d,e} | {b,d,e} |
| e | e | {1,c,e} | {a,b,d,e} | {0,1,a,b,c,d,e} | {1,b,c,e} | {b,d,e} | {1,a,b,d,e} |

| $HF_7^{131}$ | 0 | 1 | a | b | c | d | e |
|---|---|---|---|---|---|---|---|
| 0 | 0 | 1 | a | b | c | d | e |
| 1 | 1 | {1,a,b,c,e} | {1,a,d} | {1,a,b,c,e} | {0,1,a,b,c,d,e} | {1,a,c,d,e} | {1,c,e} |
| a | a | {1,a,d} | {1,a,b,c,d} | {a,b,e} | {1,a,b,c,d} | {0,1,a,b,c,d,e} | {1,a,b,d,e} |
| b | b | {1,a,b,c,e} | {a,b,e} | {a,b,c,d,e} | {1,b,c} | {a,b,c,d,e} | {0,1,a,b,c,d,e} |
| c | c | {0,1,a,b,c,d,e} | {1,a,b,c,d} | {1,b,c} | {1,b,c,d,e} | {a,c,d} | {1,b,c,d,e} |
| d | d | {1,a,c,d,e} | {0,1,a,b,c,d,e} | {a,b,c,d,e} | {a,c,d} | {1,a,c,d,e} | {b,d,e} |
| e | e | {1,c,e} | {1,a,b,d,e} | {0,1,a,b,c,d,e} | {1,b,c,d,e} | {b,d,e} | {1,a,b,d,e} |

| $HF_7^{132}$ | 0 | 1 | a | b | c | d | e |
|---|---|---|---|---|---|---|---|
| 0 | 0 | 1 | a | b | c | d | e |
| 1 | 1 | {1,a,b,c,e} | {1,a,b,d,e} | {1,b,c,d,e} | {0,1,a,b,c,d,e} | {1,a,b,c,d} | {1,a,c,d,e} |
| a | a | {1,a,b,d,e} | {1,a,b,c,d} | {1,a,b,c,e} | {1,a,c,d,e} | {0,1,a,b,c,d,e} | {a,b,c,d,e} |
| b | b | {1,b,c,d,e} | {1,a,b,c,e} | {a,b,c,d,e} | {1,a,b,c,d} | {1,a,b,d,e} | {0,1,a,b,c,d,e} |
| c | c | {0,1,a,b,c,d,e} | {1,a,c,d,e} | {1,a,b,c,d} | {1,b,c,d,e} | {a,b,c,d,e} | {1,a,b,c,e} |
| d | d | {1,a,b,c,d} | {0,1,a,b,c,d,e} | {1,a,b,d,e} | {a,b,c,d,e} | {1,a,c,d,e} | {1,b,c,d,e} |
| e | e | {1,a,c,d,e} | {a,b,c,d,e} | {0,1,a,b,c,d,e} | {1,a,b,c,e} | {1,b,c,d,e} | {1,a,b,d,e} |

| $HF_7^{133}$ | 0 | 1 | a | b | c | d | e |
|---|---|---|---|---|---|---|---|
| 0 | 0 | 1 | a | b | c | d | e |
| 1 | 1 | {1,a,b,c,e} | {1,a,b,d,e} | {1,a,b,c,d,e} | {0,1,a,b,c,d,e} | {1,a,b,c,d,e} | {1,a,c,d,e} |
| a | a | {1,a,b,d,e} | {1,a,b,c,d} | {1,a,b,c,e} | {1,a,b,c,d,e} | {0,1,a,b,c,d,e} | {1,a,b,c,d,e} |
| b | b | {1,a,b,c,d,e} | {1,a,b,c,e} | {a,b,c,d,e} | {1,a,b,c,d} | {1,a,b,c,d,e} | {0,1,a,b,c,d,e} |
| c | c | {0,1,a,b,c,d,e} | {1,a,b,c,d,e} | {1,a,b,c,d} | {1,b,c,d,e} | {a,b,c,d,e} | {1,a,b,c,d,e} |
| d | d | {1,a,b,c,d,e} | {0,1,a,b,c,d,e} | {1,a,b,c,d,e} | {a,b,c,d,e} | {1,a,c,d,e} | {1,b,c,d,e} |
| e | e | {1,a,c,d,e} | {1,a,b,c,d,e} | {0,1,a,b,c,d,e} | {1,a,b,c,d,e} | {1,b,c,d,e} | {1,a,b,d,e} |

| $HF_7^{134}$ | 0 | 1 | a | b | c | d | e |
|---|---|---|---|---|---|---|---|
| 0 | 0 | 1 | a | b | c | d | e |
| 1 | 1 | {1,a,b,d,e} | {1,a,c,d} | {1,b,c,e} | {0,1,a,b,c,d,e} | {1,a,c,d} | {1,b,c,e} |
| a | a | {1,a,c,d} | {1,a,b,c,e} | {a,b,d,e} | {1,a,c,d} | {0,1,a,b,c,d,e} | {a,b,d,e} |
| b | b | {1,b,c,e} | {a,b,d,e} | {1,a,b,c,d} | {1,b,c,e} | {a,b,d,e} | {0,1,a,b,c,d,e} |
| c | c | {0,1,a,b,c,d,e} | {1,a,c,d} | {1,b,c,e} | {a,b,c,d,e} | {1,a,c,d} | {1,b,c,e} |
| d | d | {1,a,c,d} | {0,1,a,b,c,d,e} | {a,b,d,e} | {1,a,c,d} | {1,b,c,d,e} | {a,b,d,e} |
| e | e | {1,b,c,e} | {a,b,d,e} | {0,1,a,b,c,d,e} | {1,b,c,e} | {a,b,d,e} | {1,a,c,d,e} |

| $HF_7^{135}$ | 0 | 1 | a | b | c | d | e |
|---|---|---|---|---|---|---|---|
| 0 | 0 | 1 | a | b | c | d | e |
| 1 | 1 | {1,a,b,d,e} | {1,a,c,d} | {1,a,b,c,e} | {0,1,a,b,c,d,e} | {1,a,c,d,e} | {1,b,c,e} |
| a | a | {1,a,c,d} | {1,a,b,c,e} | {a,b,d,e} | {1,a,b,c,d} | {0,1,a,b,c,d,e} | {1,a,b,d,e} |
| b | b | {1,a,b,c,e} | {a,b,d,e} | {1,a,b,c,d} | {1,b,c,e} | {a,b,c,d,e} | {0,1,a,b,c,d,e} |
| c | c | {0,1,a,b,c,d,e} | {1,a,b,c,d} | {1,b,c,e} | {a,b,c,d,e} | {1,a,c,d} | {1,b,c,d,e} |
| d | d | {1,a,c,d,e} | {0,1,a,b,c,d,e} | {a,b,c,d,e} | {1,a,c,d} | {1,b,c,d,e} | {a,b,d,e} |
| e | e | {1,b,c,e} | {1,a,b,d,e} | {0,1,a,b,c,d,e} | {1,b,c,d,e} | {a,b,d,e} | {1,a,c,d,e} |

| $HF_7^{136}$ | 0 | 1 | a | b | c | d | e |
|---|---|---|---|---|---|---|---|
| 0 | 0 | 1 | a | b | c | d | e |
| 1 | 1 | {1,a,b,d,e} | {1,a,b,c,d,e} | {1,b,c,d,e} | {0,1,a,b,c,d,e} | {1,a,b,c,d} | {1,a,b,c,d,e} |
| a | a | {1,a,b,c,d,e} | {1,a,b,c,e} | {1,a,b,c,d,e} | {1,a,c,d,e} | {0,1,a,b,c,d,e} | {a,b,c,d,e} |
| b | b | {1,b,c,d,e} | {1,a,b,c,d,e} | {1,a,b,c,d} | {1,a,b,c,d,e} | {1,a,b,d,e} | {0,1,a,b,c,d,e} |
| c | c | {0,1,a,b,c,d,e} | {1,a,c,d,e} | {1,a,b,c,d,e} | {a,b,c,d,e} | {1,a,b,c,d,e} | {1,a,b,c,e} |
| d | d | {1,a,b,c,d} | {0,1,a,b,c,d,e} | {1,a,b,d,e} | {1,a,b,c,d,e} | {1,b,c,d,e} | {1,a,b,c,d,e} |
| e | e | {1,a,b,c,d,e} | {a,b,c,d,e} | {0,1,a,b,c,d,e} | {1,a,b,c,e} | {1,a,b,c,d,e} | {1,a,c,d,e} |

| $HF_7^{137}$ | 0 | 1 | a | b | c | d | e |
|---|---|---|---|---|---|---|---|
| 0 | 0 | 1 | a | b | c | d | e |
| 1 | 1 | {1,a,b,d,e} | {1,a,b,c,d,e} | {1,a,b,c,d,e} | {0,1,a,b,c,d,e} | {1,a,b,c,d,e} | {1,a,b,c,d,e} |
| a | a | {1,a,b,c,d,e} | {1,a,b,c,e} | {1,a,b,c,d,e} | {1,a,b,c,d,e} | {0,1,a,b,c,d,e} | {1,a,b,c,d,e} |
| b | b | {1,a,b,c,d,e} | {1,a,b,c,d,e} | {1,a,b,c,d} | {1,a,b,c,d,e} | {1,a,b,c,d,e} | {0,1,a,b,c,d,e} |
| c | c | {0,1,a,b,c,d,e} | {1,a,b,c,d,e} | {1,a,b,c,d,e} | {a,b,c,d,e} | {1,a,b,c,d,e} | {1,a,b,c,d,e} |
| d | d | {1,a,b,c,d,e} | {0,1,a,b,c,d,e} | {1,a,b,c,d,e} | {1,a,b,c,d,e} | {1,b,c,d,e} | {1,a,b,c,d,e} |
| e | e | {1,a,b,c,d,e} | {1,a,b,c,d,e} | {0,1,a,b,c,d,e} | {1,a,b,c,d,e} | {1,a,b,c,d,e} | {1,a,c,d,e} |

| $HF_7^{138}$ | 0 | 1 | a | b | c | d | e |
|---|---|---|---|---|---|---|---|
| 0 | 0 | 1 | a | b | c | d | e |
| 1 | 1 | {1,a,b,c,d,e} | {1,a,c,d} | {1,b,c,e} | {0,1,a,b,c,d,e} | {1,a,c,d} | {1,b,c,e} |
| a | a | {1,a,c,d} | {1,a,b,c,d,e} | {a,b,d,e} | {1,a,c,d} | {0,1,a,b,c,d,e} | {a,b,d,e} |
| b | b | {1,b,c,e} | {a,b,d,e} | {1,a,b,c,d,e} | {1,b,c,e} | {a,b,d,e} | {0,1,a,b,c,d,e} |
| c | c | {0,1,a,b,c,d,e} | {1,a,c,d} | {1,b,c,e} | {1,a,b,c,d,e} | {1,a,c,d} | {1,b,c,e} |
| d | d | {1,a,c,d} | {0,1,a,b,c,d,e} | {a,b,d,e} | {1,a,c,d} | {1,a,b,c,d,e} | {a,b,d,e} |
| e | e | {1,b,c,e} | {a,b,d,e} | {0,1,a,b,c,d,e} | {1,b,c,e} | {a,b,d,e} | {1,a,b,c,d,e} |

| $HF_7^{139}$ | 0 | 1 | a | b | c | d | e |
|---|---|---|---|---|---|---|---|
| 0 | 0 | 1 | a | b | c | d | e |
| 1 | 1 | {1,a,b,c,d,e} | {1,a,c,d} | {1,a,b,c,e} | {0,1,a,b,c,d,e} | {1,a,c,d,e} | {1,b,c,e} |
| a | a | {1,a,c,d} | {1,a,b,c,d,e} | {a,b,d,e} | {1,a,b,c,d} | {0,1,a,b,c,d,e} | {1,a,b,d,e} |
| b | b | {1,a,b,c,e} | {a,b,d,e} | {1,a,b,c,d,e} | {1,b,c,e} | {a,b,c,d,e} | {0,1,a,b,c,d,e} |
| c | c | {0,1,a,b,c,d,e} | {1,a,b,c,d} | {1,b,c,e} | {1,a,b,c,d,e} | {1,a,c,d} | {1,b,c,d,e} |
| d | d | {1,a,c,d,e} | {0,1,a,b,c,d,e} | {a,b,c,d,e} | {1,a,c,d} | {1,a,b,c,d,e} | {a,b,d,e} |
| e | e | {1,b,c,e} | {1,a,b,d,e} | {0,1,a,b,c,d,e} | {1,b,c,d,e} | {a,b,d,e} | {1,a,b,c,d,e} |

| $HF_7^{140}$ | 0 | 1 | a | b | c | d | e |
|---|---|---|---|---|---|---|---|
| 0 | 0 | 1 | a | b | c | d | e |
| 1 | 1 | {1,a,b,c,d,e} | {1,a,b,c,d,e} | {1,b,c,d,e} | {0,1,a,b,c,d,e} | {1,a,b,c,d} | {1,a,b,c,d,e} |
| a | a | {1,a,b,c,d,e} | {1,a,b,c,d,e} | {1,a,b,c,d,e} | {1,a,c,d,e} | {0,1,a,b,c,d,e} | {a,b,c,d,e} |
| b | b | {1,b,c,d,e} | {1,a,b,c,d,e} | {1,a,b,c,d,e} | {1,a,b,c,d,e} | {1,a,b,d,e} | {0,1,a,b,c,d,e} |
| c | c | {0,1,a,b,c,d,e} | {1,a,c,d,e} | {1,a,b,c,d,e} | {1,a,b,c,d,e} | {1,a,b,c,d,e} | {1,a,b,c,e} |
| d | d | {1,a,b,c,d} | {0,1,a,b,c,d,e} | {1,a,b,d,e} | {1,a,b,c,d,e} | {1,a,b,c,d,e} | {1,a,b,c,d,e} |
| e | e | {1,a,b,c,d,e} | {a,b,c,d,e} | {0,1,a,b,c,d,e} | {1,a,b,c,e} | {1,a,b,c,d,e} | {1,a,b,c,d,e} |

| $HF_7^{141}$ | 0 | 1 | a | b | c | d | e |
|---|---|---|---|---|---|---|---|
| 0 | 0 | 1 | a | b | c | d | e |
| 1 | 1 | {1,a,b,c,d,e} | {1,a,b,c,d,e} | {1,a,b,c,d,e} | {0,1,a,b,c,d,e} | {1,a,b,c,d,e} | {1,a,b,c,d,e} |
| a | a | {1,a,b,c,d,e} | {1,a,b,c,d,e} | {1,a,b,c,d,e} | {1,a,b,c,d,e} | {0,1,a,b,c,d,e} | {1,a,b,c,d,e} |
| b | b | {1,a,b,c,d,e} | {1,a,b,c,d,e} | {1,a,b,c,d,e} | {1,a,b,c,d,e} | {1,a,b,c,d,e} | {0,1,a,b,c,d,e} |
| c | c | {0,1,a,b,c,d,e} | {1,a,b,c,d,e} | {1,a,b,c,d,e} | {1,a,b,c,d,e} | {1,a,b,c,d,e} | {1,a,b,c,d,e} |
| d | d | {1,a,b,c,d,e} | {0,1,a,b,c,d,e} | {1,a,b,c,d,e} | {1,a,b,c,d,e} | {1,a,b,c,d,e} | {1,a,b,c,d,e} |
| e | e | {1,a,b,c,d,e} | {1,a,b,c,d,e} | {0,1,a,b,c,d,e} | {1,a,b,c,e} | {1,a,b,c,d,e} | {1,a,b,c,d,e} |

## B. The hyperfields with self-opposite elements

**B1i. Hyperfields for which card(x-x) = 2 and x∈ x-x, for every non-zero element x.**

| $HF_7^{142}$ | 0 | 1 | a | b | c | d | e |
|---|---|---|---|---|---|---|---|
| 0 | 0 | 1 | a | b | c | d | e |
| 1 | 1 | {0,1} | {b,c,d,e} | {a,c,d,e} | {a,b,d,e} | {a,b,c,e} | {a,b,c,d} |
| a | a | {b,c,d,e} | {0,a} | {1,c,d,e} | {1,b,d,e} | {1,b,c,e} | {1,b,c,d} |
| b | b | {a,c,d,e} | {1,c,d,e} | {0,b} | {1,a,d,e} | {1,a,c,e} | {1,a,c,d} |
| c | c | {a,b,d,e} | {1,b,d,e} | {1,a,d,e} | {0,c} | {1,a,b,e} | {1,a,b,d} |
| d | d | {a,b,c,e} | {1,b,c,e} | {1,a,c,e} | {1,a,b,e} | {0,d} | {1,a,b,c} |
| e | e | {a,b,c,d} | {1,b,c,d} | {1,a,c,d} | {1,a,b,d} | {1,a,b,c} | {0,e} |

**B1ii. Hyperfields for which card(x-x) = 2 and x∉ x-x, for every non-zero element x.**

| $HF_7^{143}$ | 0 | 1 | a | b | c | d | e |
|---|---|---|---|---|---|---|---|
| 0 | 0 | 1 | a | b | c | d | e |
| 1 | 1 | {0,a} | {1,d} | {c,d} | {b,e} | {a,b} | {c,e} |
| a | a | {1,d} | {0,b} | {a,e} | {d,e} | {1,c} | {b,c} |
| b | b | {c,d} | {a,e} | {0,c} | {1,b} | {1,e} | {a,d} |
| c | c | {b,e} | {d,e} | {1,b} | {0,d} | {a,c} | {1,a} |
| d | d | {a,b} | {1,c} | {1,e} | {a,c} | {0,e} | {b,d} |
| e | e | {c,e} | {b,c} | {a,d} | {1,a} | {b,d} | {0,1} |

**B2i. Hyperfields for which card(x-x) = 3 and x∈ x-x, for every non-zero element x.**

| $HF_7^{144}$ | 0 | 1 | a | b | c | d | e |
|---|---|---|---|---|---|---|---|
| 0 | 0 | 1 | a | b | c | d | e |
| 1 | 1 | {0,1,a} | {1,b,c,d,e} | {a,c,d,e} | {a,b,d,e} | {a,b,c,e} | {a,b,c,d,e} |
| a | a | {1,b,c,d,e} | {0,a,b} | {1,a,c,d,e} | {1,b,d,e} | {1,b,c,e} | {1,b,c,d} |
| b | b | {a,c,d,e} | {1,a,c,d,e} | {0,b,c} | {1,a,b,d,e} | {1,a,c,e} | {1,a,c,d} |
| c | c | {a,b,d,e} | {1,b,d,e} | {1,a,b,d,e} | {0,c,d} | {1,a,b,c,e} | {1,a,b,d} |
| d | d | {a,b,c,e} | {1,b,c,e} | {1,a,c,e} | {1,a,b,c,e} | {0,d,e} | {1,a,b,c,d} |
| e | e | {a,b,c,d,e} | {1,b,c,d} | {1,a,c,d} | {1,a,b,d} | {1,a,b,c,d} | {0,1,e} |

| $HF_7^{145}$ | 0 | 1 | a | b | c | d | e |
|---|---|---|---|---|---|---|---|
| 0 | 0 | 1 | a | b | c | d | e |
| 1 | 1 | {0,1,b} | {b,c,d,e} | {1,a,c,e} | {a,b,d,e} | {a,c,d,e} | {a,b,c,d} |
| a | a | {b,c,d,e} | {0,a,c} | {1,c,d,e} | {1,a,b,d} | {1,b,c,e} | {1,b,d,e} |
| b | b | {1,a,c,e} | {1,c,d,e} | {0,b,d} | {1,a,d,e} | {a,b,c,e} | {1,a,c,d} |
| c | c | {a,b,d,e} | {1,a,b,d} | {1,a,d,e} | {0,c,e} | {1,a,b,e} | {1,b,c,d} |
| d | d | {a,c,d,e} | {1,b,c,e} | {a,b,c,e} | {1,a,b,e} | {0,1,d} | {1,a,b,c} |
| e | e | {a,b,c,d} | {1,b,d,e} | {1,a,c,d} | {1,b,c,d} | {1,a,b,c} | {0,a,e} |

| $HF_7^{146}$ | 0 | 1 | a | b | c | d | e |
|---|---|---|---|---|---|---|---|
| 0 | 0 | 1 | a | b | c | d | e |
| 1 | 1 | {0,1,b} | {b,c,d,e} | {1,a,c,d,e} | {a,b,d,e} | {a,b,c,d,e} | {a,b,c,d} |
| a | a | {b,c,d,e} | {0,a,c} | {1,c,d,e} | {1,a,b,d,e} | {1,b,c,e} | {1,b,c,d,e} |
| b | b | {1,a,c,d,e} | {1,c,d,e} | {0,b,d} | {1,a,d,e} | {1,a,b,c,e} | {1,a,c,d} |
| c | c | {a,b,d,e} | {1,a,b,d,e} | {1,a,d,e} | {0,c,e} | {1,a,b,e} | {1,a,b,c,d} |
| d | d | {a,b,c,d,e} | {1,b,c,e} | {1,a,b,c,e} | {1,a,b,e} | {0,1,d} | {1,a,b,c} |
| e | e | {a,b,c,d} | {1,b,c,d,e} | {1,a,c,d} | {1,a,b,c,d} | {1,a,b,c} | {0,a,e} |

| $HF_7^{147}$ | 0 | 1 | a | b | c | d | e |
|---|---|---|---|---|---|---|---|
| 0 | 0 | 1 | a | b | c | d | e |
| 1 | 1 | {0,1,c} | {b,c,d,e} | {a,c,d,e} | {1,a,b,c,d,e} | {a,b,c,e} | {a,b,c,d} |
| a | a | {b,c,d,e} | {0,a,d} | {1,c,d,e} | {1,b,d,e} | {1,a,b,c,d,e} | {1,b,c,d} |
| b | b | {a,c,d,e} | {1,c,d,e} | {0,b,e} | {1,a,d,e} | {1,a,c,e} | {1,a,b,c,d,e} |
| c | c | {1,a,b,c,d,e} | {1,b,d,e} | {1,a,d,e} | {0,1,c} | {1,a,b,e} | {1,a,b,d} |
| d | d | {a,b,c,e} | {1,a,b,c,d,e} | {1,a,c,e} | {1,a,b,e} | {0,a,d} | {1,a,b,c} |
| e | e | {a,b,c,d} | {1,b,c,d} | {1,a,b,c,d,e} | {1,a,b,d} | {1,a,b,c} | {0,b,e} |

**B2ii. Hyperfields for which card(x-x) = 3 and x∉ x-x, for every non-zero element x.**

| $HF_7^{148}$ | 0 | 1 | a | b | c | d | e |
|---|---|---|---|---|---|---|---|
| 0 | 0 | 1 | a | b | c | d | e |
| 1 | 1 | {0,a,b} | {1,c,d} | {1,c,d,e} | {a,b,d,e} | {a,b,c,d} | {b,c,e} |
| a | a | {1,c,d} | {0,b,c} | {a,d,e} | {1,a,d,e} | {1,b,c,e} | {b,c,d,e} |
| b | b | {1,c,d,e} | {a,d,e} | {0,c,d} | {1,b,e} | {1,a,b,e} | {1,a,c,d} |
| c | c | {a,b,d,e} | {1,a,d,e} | {1,b,e} | {0,d,e} | {1,a,c} | {1,a,b,c} |
| d | d | {a,b,c,d} | {1,b,c,e} | {1,a,b,e} | {1,a,c} | {0,1,e} | {a,b,d} |
| e | e | {b,c,e} | {b,c,d,e} | {1,a,c,d} | {1,a,b,c} | {a,b,d} | {0,1,a} |

| $HF_7^{149}$ | 0 | 1 | a | b | c | d | e |
|---|---|---|---|---|---|---|---|
| 0 | 0 | 1 | a | b | c | d | e |
| 1 | 1 | {0,a,b} | {1,b,c,e} | {1,a,c,e} | {a,b,d,e} | {a,c,d,e} | {a,b,c,d,e} |
| a | a | {1,b,c,d,e} | {0,b,c} | {1,a,c,d,e} | {1,a,b,d} | {1,b,c,e} | {1,b,d,e} |
| b | b | {1,a,c,e} | {1,a,c,d,e} | {0,c,d} | {1,a,b,d,e} | {a,b,c,e} | {1,a,c,d} |
| c | c | {a,b,d,e} | {1,a,b,d} | {1,a,b,d,e} | {0,d,e} | {1,a,b,c,e} | {1,b,c,d} |
| d | d | {a,c,d,e} | {1,b,c,e} | {a,b,c,e} | {1,a,b,c,e} | {0,1,e} | {1,a,b,c,d} |
| e | e | {a,b,c,d,e} | {1,b,d,e} | {1,a,c,d} | {1,b,c,d} | {1,a,b,c,d} | {0,1,a} |

| $HF_7^{150}$ | 0 | 1 | a | b | c | d | e |
|---|---|---|---|---|---|---|---|
| 0 | 0 | 1 | a | b | c | d | e |
| 1 | 1 | {0,a,b} | {1,b,c,d,e} | {1,a,c,d,e} | {a,b,d,e} | {a,b,c,d,e} | {a,b,c,d,e} |
| a | a | {1,b,c,d,e} | {0,b,c} | {1,a,c,d,e} | {1,a,b,d,e} | {1,b,c,e} | {1,b,c,d,e} |
| b | b | {1,a,c,d,e} | {1,a,c,d,e} | {0,c,d} | {1,a,b,d,e} | {1,a,b,c,e} | {1,a,c,d} |
| c | c | {a,b,d,e} | {1,a,b,d,e} | {1,a,b,d,e} | {0,d,e} | {1,a,b,c,e} | {1,a,b,c,d} |
| d | d | {a,b,c,d,e} | {1,b,c,e} | {1,a,b,c,e} | {1,a,b,c,e} | {0,1,e} | {1,a,b,c,d} |
| e | e | {a,b,c,d,e} | {1,b,c,d,e} | {1,a,c,d} | {1,a,b,c,d} | {1,a,b,c,d} | {0,1,a} |

| $HF_7^{151}$ | 0 | 1 | a | b | c | d | e |
|---|---|---|---|---|---|---|---|
| 0 | 0 | 1 | a | b | c | d | e |
| 1 | 1 | {0,b,d} | {b,c,d,e} | {1,a,b,c,e} | {a,b,d,e} | {1,a,c,d,e} | {a,b,c,d} |
| a | a | {b,c,d,e} | {0,c,e} | {1,c,d,e} | {1,a,b,c,d} | {1,b,c,e} | {1,a,b,d,e} |
| b | b | {1,a,b,c,e} | {1,c,d,e} | {0,1,d} | {1,a,d,e} | {a,b,c,d,e} | {1,a,c,d} |
| c | c | {a,b,d,e} | {1,a,b,c,d} | {1,a,d,e} | {0,a,e} | {1,a,b,e} | {1,b,c,d,e} |
| d | d | {1,a,c,d,e} | {1,b,c,e} | {a,b,c,d,e} | {1,a,b,e} | {0,1,b} | {1,a,b,c} |
| e | e | {a,b,c,d} | {1,a,b,d,e} | {1,a,c,d} | {1,b,c,d,e} | {1,a,b,c} | {0,a,c} |

| $HF_7^{152}$ | 0 | 1 | a | b | c | d | e |
|---|---|---|---|---|---|---|---|
| 0 | 0 | 1 | a | b | c | d | e |
| 1 | 1 | {0,b,d} | {b,c,d,e} | {1,a,b,c,d} | {a,b,d,e} | {1,a,b,c,d,e} | {a,b,c,d} |
| a | a | {b,c,d,e} | {0,c,e} | {1,c,d,e} | {1,a,b,c,d,e} | {1,b,c,e} | {1,a,b,c,d,e} |
| b | b | {1,a,b,c,d,e} | {1,c,d,e} | {0,1,d} | {1,a,d,e} | {1,a,b,c,d,e} | {1,a,c,d} |
| c | c | {a,b,d,e} | {1,a,b,c,d,e} | {1,a,d,e} | {0,a,e} | {1,a,b,e} | {1,a,b,c,d,e} |
| d | d | {1,a,b,c,d,e} | {1,b,c,e} | {1,a,b,c,d,e} | {1,a,b,e} | {0,1,b} | {1,a,b,c} |
| e | e | {a,b,c,d} | {1,a,b,c,d,e} | {1,a,c,d} | {1,a,b,c,d,e} | {1,a,b,c} | {0,a,c} |

**B3i. Hyperfields for which card(x-x) = 4 and x∈ x-x, for every non-zero element x.**

| $HF_7^{153}$ | 0 | 1 | a | b | c | d | e |
|---|---|---|---|---|---|---|---|
| 0 | 0 | 1 | a | b | c | d | e |
| 1 | 1 | {0,1,a,b} | {1,c,d} | {1,c,d,e} | {a,b,d,e} | {a,b,c,d} | {b,c,e} |
| a | a | {1,c,d} | {0,a,b,c} | {a,d,e} | {1,a,d,e} | {1,b,c,e} | {b,c,d,e} |
| b | b | {1,c,d,e} | {a,d,e} | {0,b,c,d} | {1,b,e} | {1,a,b,e} | {1,a,c,d} |
| c | c | {a,b,d,e} | {1,a,d,e} | {1,b,e} | {0,c,d,e} | {1,a,c} | {1,a,b,c} |
| d | d | {a,b,c,d} | {1,b,c,e} | {1,a,b,e} | {1,a,c} | {0,1,d,e} | {a,b,d} |
| e | e | {b,c,e} | {b,c,d,e} | {1,a,c,d} | {1,a,b,c} | {a,b,d} | {0,1,a,e} |

| $HF_7^{154}$ | 0 | 1 | a | b | c | d | e |
|---|---|---|---|---|---|---|---|
| 0 | 0 | 1 | a | b | c | d | e |
| 1 | 1 | {0,1,a,b} | {1,b,c,d,e} | {1,a,c,e} | {a,b,d,e} | {a,c,d,e} | {a,b,c,d,e} |
| a | a | {1,b,c,d,e} | {0,a,b,c} | {1,a,c,d,e} | {1,a,b,d} | {1,b,c,e} | {1,b,d,e} |
| b | b | {1,a,c,e} | {1,a,c,d,e} | {0,b,c,d} | {1,a,b,d,e} | {a,b,c,e} | {1,a,c,d} |
| c | c | {a,b,d,e} | {1,a,b,d} | {1,a,b,d,e} | {0,c,d,e} | {1,a,b,c,e} | {1,b,c,d} |
| d | d | {a,c,d,e} | {1,b,c,e} | {a,b,c,e} | {1,a,b,c,e} | {0,1,d,e} | {1,a,b,c,d} |
| e | e | {a,b,c,d,e} | {1,b,d,e} | {1,a,c,d} | {1,b,c,d} | {1,a,b,c,d} | {0,1,a,e} |

| $HF_7^{155}$ | 0 | 1 | a | b | c | d | e |
|---|---|---|---|---|---|---|---|
| 0 | 0 | 1 | a | b | c | d | e |
| 1 | 1 | {0,1,a,b} | {1,b,c,d,e} | {1,a,c,d,e} | {a,b,d,e} | {a,b,c,d,e} | {a,b,c,d,e} |
| a | a | {1,b,c,d,e} | {0,a,b,c} | {1,a,c,d,e} | {1,a,b,d,e} | {1,b,c,e} | {1,b,c,d,e} |
| b | b | {1,a,c,d,e} | {1,a,c,d,e} | {0,b,c,d} | {1,a,b,d,e} | {1,a,b,c,e} | {1,a,c,d} |
| c | c | {a,b,d,e} | {1,a,b,d,e} | {1,a,b,d,e} | {0,c,d,e} | {1,a,b,c,e} | {1,a,b,c,d} |
| d | d | {a,b,c,d,e} | {1,b,c,e} | {1,a,b,c,e} | {1,a,b,c,e} | {0,1,d,e} | {1,a,b,c,d} |
| e | e | {a,b,c,d,e} | {1,b,c,d,e} | {1,a,c,d} | {1,a,b,c,d} | {1,a,b,c,d} | {0,1,a,e} |

| $HF_7^{156}$ | 0 | 1 | a | b | c | d | e |
|---|---|---|---|---|---|---|---|
| 0 | 0 | 1 | a | b | c | d | e |
| 1 | 1 | {0,1,a,c} | {1,b,c,e} | {a,d,e} | {1,a,c,d} | {b,c,e} | {a,b,d,e} |
| a | a | {1,b,c,e} | {0,a,b,d} | {1,a,c,d} | {1,b,e} | {a,b,d,e} | {1,c,d} |
| b | b | {a,d,e} | {1,a,c,d} | {0,b,c,e} | {a,b,d,e} | {1,a,c} | {1,b,c,e} |
| c | c | {1,a,c,d} | {1,b,e} | {a,b,d,e} | {0,1,c,d} | {1,b,c,e} | {a,b,d} |
| d | d | {b,c,e} | {a,b,d,e} | {1,a,c} | {1,b,c,e} | {0,a,d,e} | {1,a,c,d} |
| e | e | {a,b,d,e} | {1,c,d} | {1,b,c,e} | {a,b,d} | {1,a,c,d} | {0,1,b,e} |

| $HF_7^{157}$ | 0 | 1 | a | b | c | d | e |
|---|---|---|---|---|---|---|---|
| 0 | 0 | 1 | a | b | c | d | e |
| 1 | 1 | {0,1,a,c} | {1,b,c,d,e} | {a,c,d,e} | {1,a,b,c,d,e} | {a,b,c,e} | {a,b,c,d,e} |
| a | a | {1,b,c,d,e} | {0,a,b,d} | {1,a,c,d,e} | {1,b,d,e} | {1,a,b,c,d,e} | {1,b,c,d} |
| b | b | {a,c,d,e} | {1,a,c,d,e} | {0,b,c,e} | {1,a,b,d,e} | {1,a,c,e} | {1,a,b,c,d,e} |
| c | c | {1,a,b,c,d,e} | {1,b,d,e} | {1,a,b,d,e} | {0,1,c,d} | {1,a,b,c,e} | {1,a,b,d} |
| d | d | {a,b,c,e} | {1,a,b,c,d,e} | {1,a,c,e} | {1,a,b,c,e} | {0,a,d,e} | {1,a,b,c,d} |
| e | e | {a,b,c,d,e} | {1,b,c,d} | {1,a,b,c,d,e} | {1,a,b,d} | {1,a,b,c,d} | {0,1,b,e} |

| $HF_7^{158}$ | 0 | 1 | a | b | c | d | e |
|---|---|---|---|---|---|---|---|
| 0 | 0 | 1 | a | b | c | d | e |
| 1 | 1 | {0,1,a,d} | {1,c,d} | {b,c,e} | {a,b,d,e} | {1,a,c} | {b,c,e} |
| a | a | {1,c,d} | {0,a,b,e} | {a,d,e} | {1,c,d} | {1,b,c,e} | {a,b,d} |
| b | b | {b,c,e} | {a,d,e} | {0,1,b,c} | {1,b,e} | {a,d,e} | {1,a,c,d} |
| c | c | {a,b,d,e} | {1,c,d} | {1,b,e} | {0,a,c,d} | {1,a,c} | {1,b,e} |
| d | d | {1,a,c} | {1,b,c,e} | {a,d,e} | {1,a,c} | {0,b,d,e} | {a,b,d} |
| e | e | {b,c,e} | {a,b,d} | {1,a,c,d} | {1,b,e} | {a,b,d} | {0,1,c,e} |

| $HF_7^{159}$ | 0 | 1 | a | b | c | d | e |
|---|---|---|---|---|---|---|---|
| 0 | 0 | 1 | a | b | c | d | e |
| 1 | 1 | {0,1,a,d} | {1,c,d} | {b,c,d,e} | {a,b,d,e} | {1,a,b,c} | {b,c,e} |
| a | a | {1,c,d} | {0,a,b,e} | {a,d,e} | {1,c,d,e} | {1,b,c,e} | {a,b,c,d} |
| b | b | {b,c,d,e} | {a,d,e} | {0,1,b,c} | {1,b,e} | {1,a,d,e} | {1,a,c,d} |
| c | c | {a,b,d,e} | {1,c,d,e} | {1,b,e} | {0,a,c,d} | {1,a,c} | {1,a,b,e} |
| d | d | {1,a,b,c} | {1,b,c,e} | {1,a,d,e} | {1,a,c} | {0,b,d,e} | {a,b,d} |
| e | e | {b,c,e} | {a,b,c,d} | {1,a,c,d} | {1,a,b,e} | {a,b,d} | {0,1,c,e} |

| $HF_7^{160}$ | 0 | 1 | a | b | c | d | e |
|---|---|---|---|---|---|---|---|
| 0 | 0 | 1 | a | b | c | d | e |
| 1 | 1 | {0,1,a,d} | {1,b,c,d,e} | {a,b,c,e} | {a,b,d,e} | {1,a,c,e} | {a,b,c,d,e} |
| a | a | {1,b,c,d,e} | {0,a,b,e} | {1,a,c,d,e} | {1,b,c,d} | {1,b,c,e} | {1,a,b,d} |
| b | b | {a,b,c,e} | {1,a,c,d,e} | {0,1,b,c} | {1,a,b,d,e} | {a,c,d,e} | {1,a,c,d} |
| c | c | {a,b,d,e} | {1,b,c,d} | {1,a,b,d,e} | {0,a,c,d} | {1,a,b,c,e} | {1,b,d,e} |
| d | d | {1,a,c,e} | {1,b,c,e} | {a,c,d,e} | {1,a,b,c,e} | {0,b,d,e} | {1,a,b,c,d} |
| e | e | {a,b,c,d,e} | {1,a,b,d} | {1,a,c,d} | {1,b,d,e} | {1,a,b,c,d} | {0,1,c,e} |

| $HF_7^{161}$ | 0 | 1 | a | b | c | d | e |
|---|---|---|---|---|---|---|---|
| 0 | 0 | 1 | a | b | c | d | e |
| 1 | 1 | {0,1,a,d} | {1,b,c,d,e} | {a,b,c,d,e} | {a,b,d,e} | {1,a,b,c,e} | {a,b,c,d,e} |
| a | a | {1,b,c,d,e} | {0,a,b,e} | {1,a,c,d,e} | {1,b,c,d,e} | {1,b,c,e} | {1,a,b,c,d} |
| b | b | {a,b,c,d,e} | {1,a,c,d,e} | {0,1,b,c} | {1,a,b,d,e} | {1,a,c,d,e} | {1,a,c,d} |
| c | c | {a,b,d,e} | {1,b,c,d,e} | {1,a,b,d,e} | {0,a,c,d} | {1,a,b,c,e} | {1,a,b,d,e} |
| d | d | {1,a,b,c,e} | {1,b,c,e} | {1,a,c,d,e} | {1,a,b,c,e} | {0,b,d,e} | {1,a,b,c,d} |
| e | e | {a,b,c,d,e} | {1,a,b,c,d} | {1,a,c,d} | {1,a,b,d,e} | {1,a,b,c,d} | {0,1,c,e} |

| $HF_7^{162}$ | 0 | 1 | a | b | c | d | e |
|---|---|---|---|---|---|---|---|
| 0 | 0 | 1 | a | b | c | d | e |
| 1 | 1 | {0,1,a,e} | {1,a,b,c,d,e} | {a,c,d,e} | {a,b,d,e} | {a,b,c,e} | {1,a,b,c,d,e} |
| a | a | {1,a,b,c,d,e} | {0,1,a,b} | {1,a,b,c,d,e} | {1,b,d,e} | {1,b,c,e} | {1,b,c,d} |
| b | b | {a,c,d,e} | {1,a,b,c,d,e} | {0,a,b,c} | {1,a,b,c,d,e} | {1,a,c,e} | {1,a,c,d} |
| c | c | {a,b,d,e} | {1,b,d,e} | {1,a,b,c,d,e} | {0,b,c,d} | {1,a,b,c,d,e} | {1,a,b,d} |
| d | d | {a,b,c,e} | {1,b,c,e} | {1,a,c,e} | {1,a,b,c,d,e} | {0,c,d,e} | {1,a,b,c,d,e} |
| e | e | {1,a,b,c,d,e} | {1,b,c,d} | {1,a,c,d} | {1,a,b,d} | {1,a,b,c,d,e} | {0,1,d,e} |

| $HF_7^{163}$ | 0 | 1 | a | b | c | d | e |
|---|---|---|---|---|---|---|---|
| 0 | 0 | 1 | a | b | c | d | e |
| 1 | 1 | {0,1,b,c} | {b,d,e} | {1,a,c,d} | {1,b,c,e} | {a,b,d,e} | {a,c,d} |
| a | a | {b,d,e} | {0,a,c,d} | {1,c,e} | {a,b,d,e} | {1,a,c,d} | {1,b,c,e} |
| b | b | {1,a,c,d} | {1,c,e} | {0,b,d,e} | {1,a,d} | {1,b,c,e} | {a,b,d,e} |
| c | c | {1,b,c,e} | {a,b,d,e} | {1,a,d} | {0,1,c,e} | {a,b,e} | {1,a,c,d} |
| d | d | {a,b,d,e} | {1,a,c,d} | {1,b,c,e} | {a,b,e} | {0,1,a,d} | {1,b,c} |
| e | e | {a,c,d} | {1,b,c,e} | {a,b,d,e} | {1,a,c,d} | {1,b,c} | {0,a,b,e} |

| $HF_7^{164}$ | 0 | 1 | a | b | c | d | e |
|---|---|---|---|---|---|---|---|
| 0 | 0 | 1 | a | b | c | d | e |
| 1 | 1 | {0,1,b,c} | {b,c,d,e} | {1,a,c,e} | {1,a,b,c,d,e} | {a,c,d,e} | {a,b,c,d} |
| a | a | {b,c,d,e} | {0,a,c,d} | {1,c,d,e} | {1,a,b,d} | {1,a,b,c,d,e} | {1,b,d,e} |
| b | b | {1,a,c,e} | {1,c,d,e} | {0,b,d,e} | {1,a,d,e} | {a,b,c,e} | {1,a,b,c,d,e} |
| c | c | {1,a,b,c,d,e} | {1,a,b,d} | {1,a,d,e} | {0,1,c,e} | {1,a,b,e} | {1,b,c,d} |
| d | d | {a,c,d,e} | {1,a,b,c,d,e} | {a,b,c,e} | {1,a,b,e} | {0,1,a,d} | {1,a,b,c} |
| e | e | {a,b,c,d} | {1,b,d,e} | {1,a,b,c,d,e} | {1,b,c,d} | {1,a,b,c} | {0,a,b,e} |

| $HF_7^{165}$ | 0 | 1 | a | b | c | d | e |
|---|---|---|---|---|---|---|---|
| 0 | 0 | 1 | a | b | c | d | e |
| 1 | 1 | {0,1,b,c} | {b,c,d,e} | {1,a,c,d,e} | {1,a,b,c,d,e} | {a,b,c,d,e} | {a,b,c,d} |
| a | a | {b,c,d,e} | {0,a,c,d} | {1,c,d,e} | {1,a,b,d,e} | {1,a,b,c,d,e} | {1,b,c,d,e} |
| b | b | {1,a,c,d,e} | {1,c,d,e} | {0,b,d,e} | {1,a,d,e} | {1,a,b,c,e} | {1,a,b,c,d,e} |
| c | c | {1,a,b,c,d,e} | {1,a,b,d,e} | {1,a,d,e} | {0,1,c,e} | {1,a,b,e} | {1,a,b,c,d} |
| d | d | {a,b,c,d,e} | {1,a,b,c,d,e} | {1,a,b,c,e} | {1,a,b,e} | {0,1,a,d} | {1,a,b,c} |
| e | e | {a,b,c,d} | {1,b,c,d,e} | {1,a,b,c,d,e} | {1,a,b,c,d} | {1,a,b,c} | {0,a,b,e} |

| $HF_7^{166}$ | 0 | 1 | a | b | c | d | e |
|---|---|---|---|---|---|---|---|
| 0 | 0 | 1 | a | b | c | d | e |
| 1 | 1 | {0,1,b,d} | {b,c,d,e} | {1,a,b,c,e} | {a,b,d,e} | {1,a,c,d,e} | {a,b,c,d} |
| a | a | {b,c,d,e} | {0,a,c,e} | {1,c,d,e} | {1,a,b,c,d} | {1,b,c,e} | {1,a,b,d,e} |
| b | b | {1,a,b,c,e} | {1,c,d,e} | {0,1,b,d} | {1,a,d,e} | {a,b,c,d,e} | {1,a,c,d} |
| c | c | {a,b,d,e} | {1,a,b,c,d} | {1,a,d,e} | {0,a,c,e} | {1,a,b,e} | {1,b,c,d,e} |
| d | d | {1,a,c,d,e} | {1,b,c,e} | {a,b,c,d,e} | {1,a,b,e} | {0,1,b,d} | {1,a,b,c} |
| e | e | {a,b,c,d} | {1,a,b,d,e} | {1,a,c,d} | {1,b,c,d,e} | {1,a,b,c} | {0,a,c,e} |

| $HF_7^{167}$ | 0 | 1 | a | b | c | d | e |
|---|---|---|---|---|---|---|---|
| 0 | 0 | 1 | a | b | c | d | e |
| 1 | 1 | {0,1,b,d} | {b,c,d,e} | {1,a,b,c,d,e} | {a,b,d,e} | {1,a,b,c,d,e} | {a,b,c,d} |
| a | a | {b,c,d,e} | {0,a,c,e} | {1,c,d,e} | {1,a,b,c,d,e} | {1,b,c,e} | {1,a,b,c,d,e} |
| b | b | {1,a,b,c,d,e} | {1,c,d,e} | {0,1,b,d} | {1,a,d,e} | {1,a,b,c,d,e} | {1,a,c,d} |
| c | c | {a,b,d,e} | {1,a,b,c,d,e} | {1,a,d,e} | {0,a,c,e} | {1,a,b,e} | {1,a,b,c,d,e} |
| d | d | {1,a,b,c,d,e} | {1,b,c,e} | {1,a,b,c,d,e} | {1,a,b,e} | {0,1,b,d} | {1,a,b,c} |
| e | e | {a,b,c,d} | {1,a,b,c,d,e} | {1,a,c,d} | {1,a,b,c,d,e} | {1,a,b,c} | {0,a,c,e} |

**B3ii. Hyperfields for which card(x-x) = 4 and x∉ x-x, for every non-zero element x.**

| $HF_7^{168}$ | 0 | 1 | a | b | c | d | e |
|---|---|---|---|---|---|---|---|
| 0 | 0 | 1 | a | b | c | d | e |
| 1 | 1 | {0,a,b,c} | {1,c,d} | {1,c,d,e} | {1,a,b,c,d,e} | {a,b,c,d} | {b,c,e} |
| a | a | {1,c,d} | {0,b,c,d} | {a,d,e} | {1,a,d,e} | {1,a,b,c,d,e} | {b,c,d,e} |
| b | b | {1,c,d,e} | {a,d,e} | {0,c,d,e} | {1,b,e} | {1,a,b,e} | {1,a,b,c,d,e} |
| c | c | {1,a,b,c,d,e} | {1,a,d,e} | {1,b,e} | {0,1,d,e} | {1,a,c} | {1,a,b,c} |
| d | d | {a,b,c,d} | {1,a,b,c,d,e} | {1,a,b,e} | {1,a,c} | {0,1,a,e} | {a,b,d} |
| e | e | {b,c,e} | {b,c,d,e} | {1,a,b,c,d,e} | {1,a,b,c} | {a,b,d} | {0,1,a,b} |

| $HF_7^{169}$ | 0 | 1 | a | b | c | d | e |
|---|---|---|---|---|---|---|---|
| 0 | 0 | 1 | a | b | c | d | e |
| 1 | 1 | {0,a,b,c} | {1,b,c,e} | {1,a,d,e} | {1,a,c,d} | {b,c,d,e} | {a,b,d,e} |
| a | a | {1,b,c,e} | {0,b,c,d} | {1,a,c,d} | {1,a,b,e} | {a,b,d,e} | {1,c,d,e} |
| b | b | {1,a,d,e} | {1,a,c,d} | {0,c,d,e} | {a,b,d,e} | {1,a,b,c} | {1,b,c,e} |
| c | c | {1,a,c,d} | {1,a,b,e} | {a,b,d,e} | {0,1,d,e} | {1,b,c,e} | {a,b,c,d} |
| d | d | {b,c,d,e} | {a,b,d,e} | {1,a,b,c} | {1,b,c,e} | {0,1,a,e} | {1,a,c,d} |
| e | e | {a,b,d,e} | {1,c,d,e} | {1,b,c,e} | {a,b,c,d} | {1,a,c,d} | {0,1,a,b} |

| $HF_7^{170}$ | 0 | 1 | a | b | c | d | e |
|---|---|---|---|---|---|---|---|
| 0 | 0 | 1 | a | b | c | d | e |
| 1 | 1 | {0,a,b,c} | {1,b,d,e} | {1,a,c} | {1,b,c,e} | {a,d,e} | {a,c,d,e} |
| a | a | {1,b,d,e} | {0,b,c,d} | {1,a,c,e} | {a,b,d} | {1,a,c,d} | {1,b,e} |
| b | b | {1,a,c} | {1,a,c,e} | {0,c,d,e} | {1,a,b,d} | {b,c,e} | {a,b,d,e} |
| c | c | {1,b,c,e} | {a,b,d} | {1,a,b,d} | {0,1,d,e} | {a,b,c,e} | {1,c,d} |
| d | d | {a,d,e} | {1,a,c,d} | {b,c,e} | {a,b,c,e} | {0,1,a,e} | {1,b,c,d} |
| e | e | {a,c,d,e} | {1,b,e} | {a,b,d,e} | {1,c,d} | {1,b,c,d} | {0,1,a,b} |

| $HF_7^{171}$ | 0 | 1 | a | b | c | d | e |
|---|---|---|---|---|---|---|---|
| 0 | 0 | 1 | a | b | c | d | e |
| 1 | 1 | {0,a,b,c} | {1,b,d,e} | {1,a,c,d} | {1,b,c,e} | {a,b,d,e} | {a,c,d,e} |
| a | a | {1,b,d,e} | {0,b,c,d} | {1,a,c,e} | {a,b,d,e} | {1,a,c,d} | {1,b,c,e} |
| b | b | {1,a,c,d} | {1,a,c,e} | {0,c,d,e} | {1,a,b,d} | {1,b,c,e} | {a,b,d,e} |
| c | c | {1,b,c,e} | {a,b,d,e} | {1,a,b,d} | {0,1,d,e} | {a,b,c,e} | {1,a,c,d} |
| d | d | {a,b,d,e} | {1,a,c,d} | {1,b,c,e} | {a,b,c,e} | {0,1,a,e} | {1,b,c,d} |
| e | e | {a,c,d,e} | {1,b,c,e} | {a,b,d,e} | {1,a,c,d} | {1,b,c,d} | {0,1,a,b} |

| $HF_7^{172}$ | 0 | 1 | a | b | c | d | e |
|---|---|---|---|---|---|---|---|
| 0 | 0 | 1 | a | b | c | d | e |
| 1 | 1 | {0,a,b,c} | {1,b,c,d,e} | {1,a,c,e} | {1,a,b,c,d,e} | {a,c,d,e} | {a,b,c,d,e} |
| a | a | {1,b,c,d,e} | {0,b,c,d} | {1,a,c,d,e} | {1,a,b,d} | {1,a,b,c,d,e} | {1,b,d,e} |
| b | b | {1,a,c,e} | {1,a,c,d,e} | {0,c,d,e} | {1,a,b,d,e} | {a,b,c,e} | {1,a,b,c,d,e} |
| c | c | {1,a,b,c,d,e} | {1,a,b,d} | {1,a,b,d,e} | {0,1,d,e} | {1,a,b,c,e} | {1,b,c,d} |
| d | d | {a,c,d,e} | {1,a,b,c,d,e} | {a,b,c,e} | {1,a,b,c,e} | {0,1,a,e} | {1,a,b,c,d} |
| e | e | {a,b,c,d,e} | {1,b,d,e} | {1,a,b,c,d,e} | {1,b,c,d} | {1,a,b,c,d} | {0,1,a,b} |

| $HF_7^{173}$ | 0 | 1 | a | b | c | d | e |
|---|---|---|---|---|---|---|---|
| 0 | 0 | 1 | a | b | c | d | e |
| 1 | 1 | {0,a,b,c} | {1,b,c,d,e} | {1,a,c,d,e} | {1,a,b,c,d,e} | {a,b,c,d,e} | {a,b,c,d,e} |
| a | a | {1,b,c,d,e} | {0,b,c,d} | {1,a,c,d,e} | {1,a,b,d,e} | {1,a,b,c,d,e} | {1,b,c,d,e} |
| b | b | {1,a,c,d,e} | {1,a,c,d,e} | {0,c,d,e} | {1,a,b,d,e} | {1,a,b,c,e} | {1,a,b,c,d,e} |
| c | c | {1,a,b,c,d,e} | {1,a,b,d,e} | {1,a,b,d,e} | {0,1,d,e} | {1,a,b,c,e} | {1,a,b,c,d} |
| d | d | {a,b,c,d,e} | {1,a,b,c,d,e} | {1,a,b,c,e} | {1,a,b,c,e} | {0,1,a,e} | {1,a,b,c,d} |
| e | e | {a,b,c,d,e} | {1,b,c,d,e} | {1,a,b,c,d,e} | {1,a,b,c,d} | {1,a,b,c,d} | {0,1,a,b} |

| $HF_7^{174}$ | 0 | 1 | a | b | c | d | e |
|---|---|---|---|---|---|---|---|
| 0 | 0 | 1 | a | b | c | d | e |
| 1 | 1 | {0,a,b,d} | {1,c,d} | {1,b,c,e} | {a,b,d,e} | {1,a,c,d} | {b,c,e} |
| a | a | {1,c,d} | {0,b,c,e} | {a,d,e} | {1,a,c,d} | {1,b,c,e} | {a,b,d,e} |
| b | b | {1,b,c,e} | {a,d,e} | {0,1,c,d} | {1,b,e} | {a,b,d,e} | {1,a,c,d} |
| c | c | {a,b,d,e} | {1,a,c,d} | {1,b,e} | {0,a,d,e} | {1,a,c} | {1,b,c,e} |
| d | d | {1,a,c,d} | {1,b,c,e} | {a,b,d,e} | {1,a,c} | {0,1,b,e} | {a,b,d} |
| e | e | {b,c,e} | {a,b,d,e} | {1,a,c,d} | {1,b,c,e} | {a,b,d} | {0,1,a,c} |

| $HF_7^{175}$ | 0 | 1 | a | b | c | d | e |
|---|---|---|---|---|---|---|---|
| 0 | 0 | 1 | a | b | c | d | e |
| 1 | 1 | {0,a,b,d} | {1,c,d} | {1,b,c,d,e} | {a,b,d,e} | {1,a,b,c,d} | {b,c,e} |
| a | a | {1,c,d} | {0,b,c,e} | {a,d,e} | {1,a,c,d,e} | {1,b,c,e} | {a,b,c,d,e} |
| b | b | {1,b,c,d,e} | {a,d,e} | {0,1,c,d} | {1,b,e} | {1,a,b,d,e} | {1,a,c,d} |
| c | c | {a,b,d,e} | {1,a,c,d,e} | {1,b,e} | {0,a,d,e} | {1,a,c} | {1,a,b,c,e} |
| d | d | {1,a,b,c,d} | {1,b,c,e} | {1,a,b,d,e} | {1,a,c} | {0,1,b,e} | {a,b,d} |
| e | e | {b,c,e} | {a,b,c,d,e} | {1,a,c,d} | {1,a,b,c,e} | {a,b,d} | {0,1,a,c} |

| $HF_7^{176}$ | 0 | 1 | a | b | c | d | e |
|---|---|---|---|---|---|---|---|
| 0 | 0 | 1 | a | b | c | d | e |
| 1 | 1 | {0,a,b,d} | {1,b,c,d,e} | {1,a,b,c,e} | {a,b,d,e} | {1,a,c,d,e} | {a,b,c,d,e} |
| a | a | {1,b,c,d,e} | {0,b,c,e} | {1,a,c,d,e} | {1,a,b,c,d} | {1,b,c,e} | {1,a,b,d,e} |
| b | b | {1,a,b,c,e} | {1,a,c,d,e} | {0,1,c,d} | {1,a,b,d,e} | {a,b,c,d,e} | {1,a,c,d} |
| c | c | {a,b,d,e} | {1,a,b,c,d} | {1,a,b,d,e} | {0,a,d,e} | {1,a,b,c,e} | {1,b,c,d,e} |
| d | d | {1,a,c,d,e} | {1,b,c,e} | {a,b,c,d,e} | {1,a,b,c,e} | {0,1,b,e} | {1,a,b,c,d} |
| e | e | {a,b,c,d,e} | {1,a,b,d,e} | {1,a,c,d} | {1,b,c,d,e} | {1,a,b,c,d} | {0,1,a,c} |

| $HF_7^{177}$ | 0 | 1 | a | b | c | d | e |
|---|---|---|---|---|---|---|---|
| 0 | 0 | 1 | a | b | c | d | e |
| 1 | 1 | {0,a,b,d} | {1,b,c,d,e} | {1,a,b,c,d,e} | {a,b,d,e} | {1,a,b,c,d,e} | {a,b,c,d,e} |
| a | a | {1,b,c,d,e} | {0,b,c,e} | {1,a,c,d,e} | {1,a,b,c,d,e} | {1,b,c,e} | {1,a,b,c,d,e} |
| b | b | {1,a,b,c,d,e} | {1,a,c,d,e} | {0,1,c,d} | {1,a,b,d,e} | {1,a,b,c,d,e} | {1,a,c,d} |
| c | c | {a,b,d,e} | {1,a,b,c,d,e} | {1,a,b,d,e} | {0,a,d,e} | {1,a,b,c,e} | {1,a,b,c,d,e} |
| d | d | {1,a,b,c,d,e} | {1,b,c,e} | {1,a,b,c,d,e} | {1,a,b,c,e} | {0,1,b,e} | {1,a,b,c,d} |
| e | e | {a,b,c,d,e} | {1,a,b,c,d,e} | {1,a,c,d} | {1,a,b,c,d,e} | {1,a,b,c,d} | {0,1,a,c} |

| $HF_7^{178}$ | 0 | 1 | a | b | c | d | e |
|---|---|---|---|---|---|---|---|
| 0 | 0 | 1 | a | b | c | d | e |
| 1 | 1 | {0,a,b,e} | {1,a,c,d} | {1,c,e} | {a,b,d,e} | {a,c,d} | {1,b,c,e} |
| a | a | {1,a,c,d} | {0,1,b,c} | {a,b,d,e} | {1,a,d} | {1,b,c,e} | {b,d,e} |
| b | b | {1,c,e} | {a,b,d,e} | {0,a,c,d} | {1,b,c,e} | {a,b,e} | {1,a,c,d} |
| c | c | {a,b,d,e} | {1,a,d} | {1,b,c,e} | {0,b,d,e} | {1,a,c,d} | {1,b,c} |
| d | d | {a,c,d} | {1,b,c,e} | {a,b,e} | {1,a,c,d} | {0,1,c,e} | {a,b,d,e} |
| e | e | {1,b,c,e} | {b,d,e} | {1,a,c,d} | {1,b,c} | {a,b,d,e} | {0,1,a,d} |

| $HF_7^{179}$ | 0 | 1 | a | b | c | d | e |
|---|---|---|---|---|---|---|---|
| 0 | 0 | 1 | a | b | c | d | e |
| 1 | 1 | {0,a,b,e} | {1,a,c,d} | {1,c,d,e} | {a,b,d,e} | {a,b,c,d} | {1,b,c,e} |
| a | a | {1,a,c,d} | {0,1,b,c} | {a,b,d,e} | {1,a,d,e} | {1,b,c,e} | {b,c,d,e} |
| b | b | {1,c,d,e} | {a,b,d,e} | {0,a,c,d} | {1,b,c,e} | {1,a,b,e} | {1,a,c,d} |
| c | c | {a,b,d,e} | {1,a,d,e} | {1,b,c,e} | {0,b,d,e} | {1,a,c,d} | {1,a,b,c} |
| d | d | {a,b,c,d} | {1,b,c,e} | {1,a,b,e} | {1,a,c,d} | {0,1,c,e} | {a,b,d,e} |
| e | e | {1,b,c,e} | {b,c,d,e} | {1,a,c,d} | {1,a,b,c} | {a,b,d,e} | {0,1,a,d} |

| $HF_7^{180}$ | 0 | 1 | a | b | c | d | e |
|---|---|---|---|---|---|---|---|
| 0 | 0 | 1 | a | b | c | d | e |
| 1 | 1 | {0,a,b,e} | {1,a,b,c,d,e} | {1,a,c,e} | {a,b,d,e} | {a,c,d,e} | {1,a,b,c,d,e} |
| a | a | {1,a,b,c,d,e} | {0,1,b,c} | {1,a,b,c,d,e} | {1,a,b,d} | {1,b,c,e} | {1,b,d,e} |
| b | b | {1,a,c,e} | {1,a,b,c,d,e} | {0,a,c,d} | {1,a,b,c,d,e} | {a,b,c,e} | {1,a,c,d} |
| c | c | {a,b,d,e} | {1,a,b,d} | {1,a,b,c,d,e} | {0,b,d,e} | {1,a,b,c,d,e} | {1,b,c,d} |
| d | d | {a,c,d,e} | {1,b,c,e} | {a,b,c,e} | {1,a,b,c,d,e} | {0,1,c,e} | {1,a,b,c,d,e} |
| e | e | {1,a,b,c,d,e} | {1,b,d,e} | {1,a,c,d} | {1,b,c,d} | {1,a,b,c,d,e} | {0,1,a,d} |

| $HF_7^{181}$ | 0 | 1 | a | b | c | d | e |
|---|---|---|---|---|---|---|---|
| 0 | 0 | 1 | a | b | c | d | e |
| 1 | 1 | {0,a,b,e} | {1,a,b,c,d,e} | {1,a,c,e} | {a,b,d,e} | {a,b,c,d,e} | {1,a,b,c,d,e} |
| a | a | {1,a,b,c,d,e} | {0,1,b,c} | {1,a,b,c,d,e} | {1,a,b,d,e} | {1,b,c,e} | {1,b,c,d,e} |
| b | b | {1,a,c,e} | {1,a,b,c,d,e} | {0,a,c,d} | {1,a,b,c,d,e} | {1,a,b,c,e} | {1,a,c,d} |
| c | c | {a,b,d,e} | {1,a,b,d} | {1,a,b,c,d,e} | {0,b,d,e} | {1,a,b,c,d,e} | {1,a,b,c,d} |
| d | d | {a,b,c,d,e} | {1,b,c,e} | {1,a,b,c,e} | {1,a,b,c,d,e} | {0,1,c,e} | {1,a,b,c,d,e} |
| e | e | {1,a,b,c,d,e} | {1,b,c,d,e} | {1,a,c,d} | {1,a,b,c,d} | {1,a,b,c,d,e} | {0,1,a,d} |

| $HF_7^{182}$ | 0 | 1 | a | b | c | d | e |
|---|---|---|---|---|---|---|---|
| 0 | 0 | 1 | a | b | c | d | e |
| 1 | 1 | {0,a,c,d} | {1,b,c,e} | {a,b,e} | {1,a,c,d} | {1,c,e} | {a,b,d,e} |
| a | a | {1,b,c,e} | {0,b,d,e} | {1,a,c,d} | {1,b,c} | {a,b,d,e} | {1,a,d} |
| b | b | {a,b,e} | {1,a,c,d} | {0,1,c,e} | {a,b,d,e} | {a,c,d} | {1,b,c,e} |
| c | c | {1,a,c,d} | {1,b,c} | {a,b,d,e} | {0,1,a,d} | {1,b,c,e} | {b,d,e} |
| d | d | {1,c,e} | {a,b,d,e} | {a,c,d} | {1,b,c,e} | {0,a,b,e} | {1,a,c,d} |
| e | e | {a,b,d,e} | {1,a,d} | {1,b,c,e} | {b,d,e} | {1,a,c,d} | {0,1,b,c} |

| $HF_7^{183}$ | 0 | 1 | a | b | c | d | e |
|---|---|---|---|---|---|---|---|
| 0 | 0 | 1 | a | b | c | d | e |
| 1 | 1 | {0,a,c,d} | {1,b,c,e} | {a,b,d,e} | {1,a,c,d} | {1,b,c,e} | {a,b,d,e} |
| a | a | {1,b,c,e} | {0,b,d,e} | {1,a,c,d} | {1,b,c,e} | {a,b,d,e} | {1,a,c,d} |
| b | b | {a,b,d,e} | {1,a,c,d} | {0,1,c,e} | {a,b,d,e} | {1,a,c,d} | {1,b,c,e} |
| c | c | {1,a,c,d} | {1,b,c,e} | {a,b,d,e} | {0,1,a,d} | {1,b,c,e} | {a,b,d,e} |
| d | d | {1,b,c,e} | {a,b,d,e} | {1,a,c,d} | {1,b,c,e} | {0,a,b,e} | {1,a,c,d} |
| e | e | {a,b,d,e} | {1,a,c,d} | {1,b,c,e} | {a,b,d,e} | {1,a,c,d} | {0,1,b,c} |

| $HF_7^{184}$ | 0 | 1 | a | b | c | d | e |
|---|---|---|---|---|---|---|---|
| 0 | 0 | 1 | a | b | c | d | e |
| 1 | 1 | {0,a,c,d} | {1,b,d,e} | {a,b,c,d} | {1,b,c,e} | {1,a,b,e} | {a,c,d,e} |
| a | a | {1,b,d,e} | {0,b,d,e} | {1,a,c,e} | {b,c,d,e} | {1,a,c,d} | {1,a,b,c} |
| b | b | {a,b,c,d} | {1,a,c,e} | {0,1,c,e} | {1,a,b,d} | {1,c,d,e} | {a,b,d,e} |
| c | c | {1,b,c,e} | {b,c,d,e} | {1,a,b,d} | {0,1,a,d} | {a,b,c,e} | {1,a,d,e} |
| d | d | {1,a,b,e} | {1,a,c,d} | {1,c,d,e} | {a,b,c,e} | {0,a,b,e} | {1,b,c,d} |
| e | e | {a,c,d,e} | {1,a,b,c} | {a,b,d,e} | {1,a,d,e} | {1,b,c,d} | {0,1,b,c} |

| $HF_7^{185}$ | 0 | 1 | a | b | c | d | e |
|---|---|---|---|---|---|---|---|
| 0 | 0 | 1 | a | b | c | d | e |
| 1 | 1 | {0,a,c,d} | {1,b,c,d,e} | {a,b,c,e} | {1,a,b,c,d,e} | {1,a,c,e} | {a,b,c,d,e} |
| a | a | {1,b,c,d,e} | {0,b,d,e} | {1,a,c,d,e} | {1,b,c,d} | {1,a,b,c,d,e} | {1,a,b,d} |
| b | b | {a,b,c,e} | {1,a,c,d,e} | {0,1,c,e} | {1,a,b,d,e} | {a,c,d,e} | {1,a,b,c,d,e} |
| c | c | {1,a,b,c,d,e} | {1,b,c,d} | {1,a,b,d,e} | {0,1,a,d} | {1,a,b,c,e} | {1,b,d,e} |
| d | d | {1,a,c,e} | {1,a,b,c,d,e} | {a,c,d,e} | {1,a,b,c,e} | {0,a,b,e} | {1,a,b,c,d} |
| e | e | {a,b,c,d,e} | {1,a,b,d} | {1,a,b,c,d,e} | {1,b,d,e} | {1,a,b,c,d} | {0,1,b,c} |

| $HF_7^{186}$ | 0 | 1 | a | b | c | d | e |
|---|---|---|---|---|---|---|---|
| 0 | 0 | 1 | a | b | c | d | e |
| 1 | 1 | {0,a,c,d} | {1,b,c,d,e} | {a,b,c,d,e} | {1,a,b,c,d,e} | {1,a,b,c,e} | {a,b,c,d,e} |
| a | a | {1,b,c,d,e} | {0,b,d,e} | {1,a,c,d,e} | {1,b,c,d,e} | {1,a,b,c,d,e} | {1,a,b,c,d} |
| b | b | {a,b,c,d,e} | {1,a,c,d,e} | {0,1,c,e} | {1,a,b,d,e} | {1,a,c,d,e} | {1,a,b,c,d,e} |
| c | c | {1,a,b,c,d,e} | {1,b,c,d,e} | {1,a,b,d,e} | {0,1,a,d} | {1,a,b,c,e} | {1,a,b,d,e} |
| d | d | {1,a,b,c,e} | {1,a,b,c,d,e} | {1,a,c,d,e} | {1,a,b,c,e} | {0,a,b,e} | {1,a,b,c,d} |
| e | e | {a,b,c,d,e} | {1,a,b,c,d} | {1,a,b,c,d,e} | {1,a,b,d,e} | {1,a,b,c,d} | {0,1,b,c} |

| $HF_7^{187}$ | 0 | 1 | a | b | c | d | e |
|---|---|---|---|---|---|---|---|
| 0 | 0 | 1 | a | b | c | d | e |
| 1 | 1 | {0,a,c,e} | {1,a,b,c,e} | {a,d,e} | {1,a,c,d} | {b,c,e} | {1,a,b,d,e} |
| a | a | {1,a,b,c,e} | {0,1,b,d} | {1,a,b,c,d} | {1,b,e} | {a,b,d,e} | {1,c,d} |
| b | b | {a,d,e} | {1,a,b,c,d} | {0,a,c,e} | {a,b,c,d,e} | {1,a,c} | {1,b,c,e} |
| c | c | {1,a,c,d} | {1,b,e} | {a,b,c,d,e} | {0,1,b,d} | {1,b,c,d,e} | {a,b,d} |
| d | d | {b,c,e} | {a,b,d,e} | {1,a,c} | {1,b,c,d,e} | {0,a,c,e} | {1,a,c,d,e} |
| e | e | {1,a,b,d,e} | {1,c,d} | {1,b,c,e} | {a,b,d} | {1,a,c,d,e} | {0,1,b,d} |

| $HF_7^{188}$ | 0 | 1 | a | b | c | d | e |
|---|---|---|---|---|---|---|---|
| 0 | 0 | 1 | a | b | c | d | e |
| 1 | 1 | {0,a,c,e} | {1,a,b,c,d,e} | {a,c,d,e} | {1,a,b,c,d,e} | {a,b,c,e} | {1,a,b,c,d,e} |
| a | a | {1,a,b,c,d,e} | {0,1,b,d} | {1,a,b,c,d,e} | {1,b,d,e} | {1,a,b,c,d,e} | {1,b,c,d} |
| b | b | {a,c,d,e} | {1,a,b,c,d,e} | {0,a,c,e} | {1,a,b,c,d,e} | {1,a,c,e} | {1,a,b,c,d,e} |
| c | c | {1,a,b,c,d,e} | {1,b,d,e} | {1,a,b,c,d,e} | {0,1,b,d} | {1,a,b,c,d,e} | {1,a,b,d} |
| d | d | {a,b,c,e} | {1,a,b,c,d,e} | {1,a,c,e} | {1,a,b,c,d,e} | {0,a,c,e} | {1,a,b,c,d,e} |
| e | e | {1,a,b,c,d,e} | {1,b,c,d} | {1,a,b,c,d,e} | {1,a,b,d} | {1,a,b,c,d,e} | {0,1,b,d} |

| $HF_7^{189}$ | 0 | 1 | a | b | c | d | e |
|---|---|---|---|---|---|---|---|
| 0 | 0 | 1 | a | b | c | d | e |
| 1 | 1 | {0,b,c,d} | {b,c,e} | {1,a,b,e} | {1,a,c,d} | {1,c,d,e} | {a,b,d} |
| a | a | {b,c,e} | {0,c,d,e} | {1,c,d} | {1,a,b,c} | {a,b,d,e} | {1,a,d,e} |
| b | b | {1,a,b,e} | {1,c,d} | {0,1,d,e} | {a,d,e} | {a,b,c,d} | {1,b,c,e} |
| c | c | {1,a,c,d} | {1,a,b,c} | {a,d,e} | {0,1,a,e} | {1,b,e} | {b,c,d,e} |
| d | d | {1,c,d,e} | {a,b,d,e} | {a,b,c,d} | {1,b,e} | {0,1,a,b} | {1,a,c} |
| e | e | {a,b,d} | {1,a,d,e} | {1,b,c,e} | {b,c,d,e} | {1,a,c} | {0,a,b,c} |

| $HF_7^{190}$ | 0 | 1 | a | b | c | d | e |
|---|---|---|---|---|---|---|---|
| 0 | 0 | 1 | a | b | c | d | e |
| 1 | 1 | {0,b,c,d} | {b,c,e} | {1,a,b,d,e} | {1,a,c,d} | {1,b,c,d,e} | {a,b,d} |
| a | a | {b,c,e} | {0,c,d,e} | {1,c,d} | {1,a,b,c,e} | {a,b,d,e} | {1,a,c,d,e} |
| b | b | {1,a,b,d,e} | {1,c,d} | {0,1,d,e} | {a,d,e} | {1,a,b,c,d} | {1,b,c,e} |
| c | c | {1,a,c,d} | {1,a,b,c,e} | {a,d,e} | {0,1,a,e} | {1,b,e} | {a,b,c,d,e} |
| d | d | {1,b,c,d,e} | {a,b,d,e} | {1,a,b,c,d} | {1,b,e} | {0,1,a,b} | {1,a,c} |
| e | e | {a,b,d} | {1,a,c,d,e} | {1,b,c,e} | {a,b,c,d,e} | {1,a,c} | {0,a,b,c} |

| $HF_7^{191}$ | 0 | 1 | a | b | c | d | e |
|---|---|---|---|---|---|---|---|
| 0 | 0 | 1 | a | b | c | d | e |
| 1 | 1 | {0,b,c,d} | {b,c,d,e} | {1,a,b,c,e} | {1,a,b,c,d,e} | {1,a,c,d,e} | {a,b,c,d} |
| a | a | {b,c,d,e} | {0,c,d,e} | {1,c,d,e} | {1,a,b,c,d} | {1,a,b,c,d,e} | {1,a,b,d,e} |
| b | b | {1,a,b,c,e} | {1,c,d,e} | {0,1,d,e} | {1,a,d,e} | {a,b,c,d,e} | {1,a,b,c,d,e} |
| c | c | {1,a,b,c,d,e} | {1,a,b,c,d} | {1,a,d,e} | {0,1,a,e} | {1,a,b,e} | {1,b,c,d,e} |
| d | d | {1,a,c,d,e} | {1,a,b,c,d,e} | {a,b,c,d,e} | {1,a,b,e} | {0,1,a,b} | {1,a,b,c} |
| e | e | {a,b,c,d} | {1,a,b,d,e} | {1,a,b,c,d,e} | {1,b,c,d,e} | {1,a,b,c} | {0,a,b,c} |

| $HF_7^{192}$ | 0 | 1 | a | b | c | d | e |
|---|---|---|---|---|---|---|---|
| 0 | 0 | 1 | a | b | c | d | e |
| 1 | 1 | {0,b,c,d} | {b,c,d,e} | {1,a,b,c,d,e} | {1,a,b,c,d,e} | {1,a,b,c,d,e} | {a,b,c,d} |
| a | a | {b,c,d,e} | {0,c,d,e} | {1,c,d,e} | {1,a,b,c,d,e} | {1,a,b,c,d,e} | {1,a,b,c,d,e} |
| b | b | {1,a,b,c,d,e} | {1,c,d,e} | {0,1,d,e} | {1,a,d,e} | {1,a,b,c,d,e} | {1,a,b,c,d,e} |
| c | c | {1,a,b,c,d,e} | {1,a,b,c,d,e} | {1,a,d,e} | {0,1,a,e} | {1,a,b,e} | {1,a,b,c,d,e} |
| d | d | {1,a,b,c,d,e} | {1,a,b,c,d,e} | {1,a,b,c,d,e} | {1,a,b,e} | {0,1,a,b} | {1,a,b,c} |
| e | e | {a,b,c,d} | {1,a,b,c,d,e} | {1,a,b,c,d,e} | {1,a,b,c,d,e} | {1,a,b,c} | {0,a,b,c} |

**B4i. Hyperfields for which card(x-x) = 5 and x∈ x-x, for every non-zero element x.**

| $HF_7^{193}$ | 0 | 1 | a | b | c | d | e |
|---|---|---|---|---|---|---|---|
| 0 | 0 | 1 | a | b | c | d | e |
| 1 | 1 | {0,1,a,b,c} | {1,c,d} | {1,c,d,e} | {1,a,b,c,d,e} | {a,b,c,d} | {b,c,e} |
| a | a | {1,c,d} | {0,a,b,c,d} | {a,d,e} | {1,a,d,e} | {1,a,b,c,d,e} | {b,c,d,e} |
| b | b | {1,c,d,e} | {a,d,e} | {0,b,c,d,e} | {1,b,e} | {1,a,b,e} | {1,a,b,c,d,e} |
| c | c | {1,a,b,c,d,e} | {1,a,d,e} | {1,b,e} | {0,1,c,d,e} | {1,a,c} | {1,a,b,c} |
| d | d | {a,b,c,d} | {1,a,b,c,d,e} | {1,a,b,e} | {1,a,c} | {0,1,a,d,e} | {a,b,d} |
| e | e | {b,c,e} | {b,c,d,e} | {1,a,b,c,d,e} | {1,a,b,c} | {a,b,d} | {0,1,a,b,e} |

| $HF_7^{194}$ | 0 | 1 | a | b | c | d | e |
|---|---|---|---|---|---|---|---|
| 0 | 0 | 1 | a | b | c | d | e |
| 1 | 1 | {0,1,a,b,c} | {1,b,c,e} | {1,a,d,e} | {1,a,c,d} | {b,c,d,e} | {a,b,d,e} |
| a | a | {1,b,c,e} | {0,a,b,c,d} | {1,a,c,d} | {1,a,b,e} | {a,b,d,e} | {1,c,d,e} |
| b | b | {1,a,d,e} | {1,a,c,d} | {0,b,c,d,e} | {a,b,d,e} | {1,a,b,c} | {1,b,c,e} |
| c | c | {1,a,c,d} | {1,a,b,e} | {a,b,d,e} | {0,1,c,d,e} | {1,b,c,e} | {a,b,c,d} |
| d | d | {b,c,d,e} | {a,b,d,e} | {1,a,b,c} | {1,b,c,e} | {0,1,a,d,e} | {1,a,c,d} |
| e | e | {a,b,d,e} | {1,c,d,e} | {1,b,c,e} | {a,b,c,d} | {1,a,c,d} | {0,1,a,b,e} |

| $HF_7^{195}$ | 0 | 1 | a | b | c | d | e |
|---|---|---|---|---|---|---|---|
| 0 | 0 | 1 | a | b | c | d | e |
| 1 | 1 | {0,1,a,b,c} | {1,b,d,e} | {1,a,c} | {1,b,c,e} | {a,d,e} | {a,c,d,e} |
| a | a | {1,b,d,e} | {0,a,b,c,d} | {1,a,c,e} | {a,b,d} | {1,a,c,d} | {1,b,e} |
| b | b | {1,a,c} | {1,a,c,e} | {0,b,c,d,e} | {1,a,b,d} | {b,c,e} | {a,b,d,e} |
| c | c | {1,b,c,e} | {a,b,d} | {1,a,b,d} | {0,1,c,d,e} | {a,b,c,e} | {1,c,d} |
| d | d | {a,d,e} | {1,a,c,d} | {b,c,e} | {a,b,c,e} | {0,1,a,d,e} | {1,b,c,d} |
| e | e | {a,c,d,e} | {1,b,e} | {a,b,d,e} | {1,c,d} | {1,b,c,d} | {0,1,a,b,e} |

| $HF_7^{196}$ | 0 | 1 | a | b | c | d | e |
|---|---|---|---|---|---|---|---|
| 0 | 0 | 1 | a | b | c | d | e |
| 1 | 1 | {0,1,a,b,c} | {1,b,d,e} | {1,a,c,d} | {1,b,c,e} | {a,b,d,e} | {a,c,d,e} |
| a | a | {1,b,d,e} | {0,a,b,c,d} | {1,a,c,e} | {a,b,d,e} | {1,a,c,d} | {1,b,c,e} |
| b | b | {1,a,c,d} | {1,a,c,e} | {0,b,c,d,e} | {1,a,b,d} | {1,b,c,e} | {a,b,d,e} |
| c | c | {1,b,c,e} | {a,b,d,e} | {1,a,b,d} | {0,1,c,d,e} | {a,b,c,e} | {1,a,c,d} |
| d | d | {a,b,d,e} | {1,a,c,d} | {1,b,c,e} | {a,b,c,e} | {0,1,a,d,e} | {1,b,c,d} |
| e | e | {a,c,d,e} | {1,b,c,e} | {a,b,d,e} | {1,a,c,d} | {1,b,c,d} | {0,1,a,b,e} |

| $HF_7^{197}$ | 0 | 1 | a | b | c | d | e |
|---|---|---|---|---|---|---|---|
| 0 | 0 | 1 | a | b | c | d | e |
| 1 | 1 | {0,1,a,b,c} | {1,b,c,d,e} | {1,a,c,e} | {1,a,b,c,d,e} | {a,c,d,e} | {a,b,c,d,e} |
| a | a | {1,b,c,d,e} | {0,a,b,c,d} | {1,a,c,d,e} | {1,a,b,d} | {1,a,b,c,d,e} | {1,b,d,e} |
| b | b | {1,a,c,e} | {1,a,c,d,e} | {0,b,c,d,e} | {1,a,b,d,e} | {a,b,c,e} | {1,a,b,c,d,e} |
| c | c | {1,a,b,c,d,e} | {1,a,b,d} | {1,a,b,d,e} | {0,1,c,d,e} | {1,a,b,c,e} | {1,b,c,d} |
| d | d | {a,c,d,e} | {1,a,b,c,d,e} | {a,b,c,e} | {1,a,b,c,e} | {0,1,a,d,e} | {1,a,b,c,d} |
| e | e | {a,b,c,d,e} | {1,b,d,e} | {1,a,b,c,d,e} | {1,b,c,d} | {1,a,b,c,d} | {0,1,a,b,e} |

| $HF_7^{198}$ | 0 | 1 | a | b | c | d | e |
|---|---|---|---|---|---|---|---|
| 0 | 0 | 1 | a | b | c | d | e |
| 1 | 1 | {0,1,a,b,c} | {1,b,c,d,e} | {1,a,c,d,e} | {1,a,b,c,d,e} | {a,b,c,d,e} | {a,b,c,d,e} |
| a | a | {1,b,c,d,e} | {0,a,b,c,d} | {1,a,c,d,e} | {1,a,b,d,e} | {1,a,b,c,d,e} | {1,b,c,d,e} |
| b | b | {1,a,c,d,e} | {1,a,c,d,e} | {0,b,c,d,e} | {1,a,b,d,e} | {1,a,b,c,e} | {1,a,b,c,d,e} |
| c | c | {1,a,b,c,d,e} | {1,a,b,d,e} | {1,a,b,d,e} | {0,1,c,d,e} | {1,a,b,c,e} | {1,a,b,c,d} |
| d | d | {a,b,c,d,e} | {1,a,b,c,d,e} | {1,a,b,c,e} | {1,a,b,c,e} | {0,1,a,d,e} | {1,a,b,c,d} |
| e | e | {a,b,c,d,e} | {1,b,c,d,e} | {1,a,b,c,d,e} | {1,a,b,c,d} | {1,a,b,c,d} | {0,1,a,b,e} |

| $HF_7^{199}$ | 0 | 1 | a | b | c | d | e |
|---|---|---|---|---|---|---|---|
| 0 | 0 | 1 | a | b | c | d | e |
| 1 | 1 | {0,1,a,b,d} | {1,c,d} | {1,b,c,e} | {a,b,d,e} | {1,a,c,d} | {b,c,e} |
| a | a | {1,c,d} | {0,a,b,c,e} | {a,d,e} | {1,a,c,d} | {1,b,c,e} | {a,b,d,e} |
| b | b | {1,b,c,e} | {a,d,e} | {0,1,b,c,d} | {1,b,e} | {a,b,d,e} | {1,a,c,d} |
| c | c | {a,b,d,e} | {1,a,c,d} | {1,b,e} | {0,a,c,d,e} | {1,a,c} | {1,b,c,e} |
| d | d | {1,a,c,d} | {1,b,c,e} | {a,b,d,e} | {1,a,c} | {0,1,b,d,e} | {a,b,d} |
| e | e | {b,c,e} | {a,b,d,e} | {1,a,c,d} | {1,b,c,e} | {a,b,d} | {0,1,a,c,e} |

| $HF_7^{200}$ | 0 | 1 | a | b | c | d | e |
|---|---|---|---|---|---|---|---|
| 0 | 0 | 1 | a | b | c | d | e |
| 1 | 1 | {0,1,a,b,d} | {1,c,d} | {1,b,c,d,e} | {a,b,d,e} | {1,a,b,c,d} | {b,c,e} |
| a | a | {1,c,d} | {0,a,b,c,e} | {a,d,e} | {1,a,c,d,e} | {1,b,c,e} | {a,b,c,d,e} |
| b | b | {1,b,c,d,e} | {a,d,e} | {0,1,b,c,d} | {1,b,e} | {1,a,b,d,e} | {1,a,c,d} |
| c | c | {a,b,d,e} | {1,a,c,d,e} | {1,b,e} | {0,a,c,d,e} | {1,a,c} | {1,a,b,c,e} |
| d | d | {1,a,b,c,d} | {1,b,c,e} | {1,a,b,d,e} | {1,a,c} | {0,1,b,d,e} | {a,b,d} |
| e | e | {b,c,e} | {a,b,c,d,e} | {1,a,c,d} | {1,a,b,c,e} | {a,b,d} | {0,1,a,c,e} |

| $HF_7^{201}$ | 0 | 1 | a | b | c | d | e |
|---|---|---|---|---|---|---|---|
| 0 | 0 | 1 | a | b | c | d | e |
| 1 | 1 | {0,1,a,b,d} | {1,b,c,d,e} | {1,a,b,c,e} | {a,b,d,e} | {1,a,c,d,e} | {a,b,c,d,e} |
| a | a | {1,b,c,d,e} | {0,a,b,c,e} | {1,a,c,d,e} | {1,a,b,c,d} | {1,b,c,e} | {1,a,b,d,e} |
| b | b | {1,a,b,c,e} | {1,a,c,d,e} | {0,1,b,c,d} | {1,a,b,d,e} | {a,b,c,d,e} | {1,a,c,d} |
| c | c | {a,b,d,e} | {1,a,b,c,d} | {1,a,b,d,e} | {0,a,c,d,e} | {1,a,b,c,e} | {1,b,c,d,e} |
| d | d | {1,a,c,d,e} | {1,b,c,e} | {a,b,c,d,e} | {1,a,b,c,e} | {0,1,b,d,e} | {1,a,b,c,d} |
| e | e | {a,b,c,d,e} | {1,a,b,d,e} | {1,a,c,d} | {1,b,c,d,e} | {1,a,b,c,d} | {0,1,a,c,e} |

| $HF_7^{202}$ | 0 | 1 | a | b | c | d | e |
|---|---|---|---|---|---|---|---|
| 0 | 0 | 1 | a | b | c | d | e |
| 1 | 1 | {0,1,a,b,d} | {1,b,c,d,e} | {1,a,b,c,d,e} | {a,b,d,e} | {1,a,b,c,d,e} | {a,b,c,d,e} |
| a | a | {1,b,c,d,e} | {0,a,b,c,e} | {1,a,c,d,e} | {1,a,b,c,d,e} | {1,b,c,e} | {1,a,b,c,d,e} |
| b | b | {1,a,b,c,d,e} | {1,a,c,d,e} | {0,1,b,c,d} | {1,a,b,d,e} | {1,a,b,c,d,e} | {1,a,c,d} |
| c | c | {a,b,d,e} | {1,a,b,c,d,e} | {1,a,b,d,e} | {0,a,c,d,e} | {1,a,b,c,e} | {1,a,b,c,d,e} |
| d | d | {1,a,b,c,d,e} | {1,b,c,e} | {1,a,b,c,d,e} | {1,a,b,c,e} | {0,1,b,d,e} | {1,a,b,c,d} |
| e | e | {a,b,c,d,e} | {1,a,b,c,d,e} | {1,a,c,d} | {1,a,b,c,d,e} | {1,a,b,c,d} | {0,1,a,c,e} |

| $HF_7^{203}$ | 0 | 1 | a | b | c | d | e |
|---|---|---|---|---|---|---|---|
| 0 | 0 | 1 | a | b | c | d | e |
| 1 | 1 | {0,1,a,b,e} | {1,a,c,d} | {1,c,e} | {a,b,d,e} | {a,c,d} | {1,b,c,e} |
| a | a | {1,a,c,d} | {0,1,a,b,c} | {a,b,d,e} | {1,a,d} | {1,b,c,e} | {b,d,e} |
| b | b | {1,c,e} | {a,b,d,e} | {0,a,b,c,d} | {1,b,c,e} | {a,b,e} | {1,a,c,d} |
| c | c | {a,b,d,e} | {1,a,d} | {1,b,c,e} | {0,b,c,d,e} | {1,a,c,d} | {1,b,c} |
| d | d | {a,c,d} | {1,b,c,e} | {a,b,e} | {1,a,c,d} | {0,1,c,d,e} | {a,b,d,e} |
| e | e | {1,b,c,e} | {b,d,e} | {1,a,c,d} | {1,b,c} | {a,b,d,e} | {0,1,a,d,e} |

| $HF_7^{204}$ | 0 | 1 | a | b | c | d | e |
|---|---|---|---|---|---|---|---|
| 0 | 0 | 1 | a | b | c | d | e |
| 1 | 1 | {0,1,a,b,e} | {1,a,c,d} | {1,c,d,e} | {a,b,d,e} | {a,b,c,d} | {1,b,c,e} |
| a | a | {1,a,c,d} | {0,1,a,b,c} | {a,b,d,e} | {1,a,d,e} | {1,b,c,e} | {b,c,d,e} |
| b | b | {1,c,d,e} | {a,b,d,e} | {0,a,b,c,d} | {1,b,c,e} | {1,a,b,e} | {1,a,c,d} |
| c | c | {a,b,d,e} | {1,a,d,e} | {1,b,c,e} | {0,b,c,d,e} | {1,a,c,d} | {1,a,b,c} |
| d | d | {a,b,c,d} | {1,b,c,e} | {1,a,b,e} | {1,a,c,d} | {0,1,c,d,e} | {a,b,d,e} |
| e | e | {1,b,c,e} | {b,c,d,e} | {1,a,c,d} | {1,a,b,c} | {a,b,d,e} | {0,1,a,d,e} |

| $HF_7^{205}$ | 0 | 1 | a | b | c | d | e |
|---|---|---|---|---|---|---|---|
| 0 | 0 | 1 | a | b | c | d | e |
| 1 | 1 | {0,1,a,b,e} | {1,a,b,c,d,e} | {1,a,c,e} | {a,b,d,e} | {a,c,d,e} | {1,a,b,c,d,e} |
| a | a | {1,a,b,c,d,e} | {0,1,a,b,c} | {1,a,b,c,d,e} | {1,a,b,d} | {1,b,c,e} | {1,b,d,e} |
| b | b | {1,a,c,e} | {1,a,b,c,d,e} | {0,a,b,c,d} | {1,a,b,c,d,e} | {a,b,c,e} | {1,a,c,d} |
| c | c | {a,b,d,e} | {1,a,b,d} | {1,a,b,c,d,e} | {0,b,c,d,e} | {1,a,b,c,d,e} | {1,b,c,d} |
| d | d | {a,c,d,e} | {1,b,c,e} | {a,b,c,e} | {1,a,b,c,d,e} | {0,1,c,d,e} | {1,a,b,c,d,e} |
| e | e | {1,a,b,c,d,e} | {1,b,d,e} | {1,a,c,d} | {1,b,c,d} | {1,a,b,c,d,e} | {0,1,a,d,e} |

| $HF_7^{206}$ | 0 | 1 | a | b | c | d | e |
|---|---|---|---|---|---|---|---|
| 0 | 0 | 1 | a | b | c | d | e |
| 1 | 1 | {0,1,a,b,e} | {1,a,b,c,d,e} | {1,a,c,d,e} | {a,b,d,e} | {a,b,c,d,e} | {1,a,b,c,d,e} |
| a | a | {1,a,b,c,d,e} | {0,1,a,b,c} | {1,a,b,c,d,e} | {1,a,b,d,e} | {1,b,c,e} | {1,b,c,d,e} |
| b | b | {1,a,c,d,e} | {1,a,b,c,d,e} | {0,a,b,c,d} | {1,a,b,c,d,e} | {1,a,b,c,e} | {1,a,c,d} |
| c | c | {a,b,d,e} | {1,a,b,d,e} | {1,a,b,c,d,e} | {0,b,c,d,e} | {1,a,b,c,d,e} | {1,a,b,c,d} |
| d | d | {a,b,c,d,e} | {1,b,c,e} | {1,a,b,c,e} | {1,a,b,c,d,e} | {0,1,c,d,e} | {1,a,b,c,d,e} |
| e | e | {1,a,b,c,d,e} | {1,b,c,d,e} | {1,a,c,d} | {1,a,b,c,d} | {1,a,b,c,d,e} | {0,1,a,d,e} |

| $HF_7^{207}$ | 0 | 1 | a | b | c | d | e |
|---|---|---|---|---|---|---|---|
| 0 | 0 | 1 | a | b | c | d | e |
| 1 | 1 | {0,1,a,c,d} | {1,c,d} | {b,c,e} | {1,a,b,c,d,e} | {1,a,c} | {b,c,e} |
| a | a | {1,c,d} | {0,a,b,d,e} | {a,d,e} | {1,c,d} | {1,a,b,c,d,e} | {a,b,d} |
| b | b | {b,c,e} | {a,d,e} | {0,1,b,c,e} | {1,b,e} | {a,d,e} | {1,a,b,c,d,e} |
| c | c | {1,a,b,c,d,e} | {1,c,d} | {1,b,e} | {0,1,a,c,d} | {1,a,c} | {1,b,e} |
| d | d | {1,a,c} | {1,a,b,c,d,e} | {a,d,e} | {1,a,c} | {0,a,b,d,e} | {a,b,d} |
| e | e | {b,c,e} | {a,b,d} | {1,a,b,c,d,e} | {1,b,e} | {a,b,d} | {0,1,b,c,e} |

| $HF_7^{208}$ | 0 | 1 | a | b | c | d | e |
|---|---|---|---|---|---|---|---|
| 0 | 0 | 1 | a | b | c | d | e |
| 1 | 1 | {0,1,a,c,d} | {1,c,d} | {b,c,d,e} | {1,a,b,c,d,e} | {1,a,b,c} | {b,c,e} |
| a | a | {1,c,d} | {0,a,b,d,e} | {a,d,e} | {1,c,d,e} | {1,a,b,c,d,e} | {a,b,c,d} |
| b | b | {b,c,d,e} | {a,d,e} | {0,1,b,c,e} | {1,b,e} | {1,a,d,e} | {1,a,b,c,d,e} |
| c | c | {1,a,b,c,d,e} | {1,c,d,e} | {1,b,e} | {0,1,a,c,d} | {1,a,c} | {1,a,b,e} |
| d | d | {1,a,b,c} | {1,a,b,c,d,e} | {1,a,d,e} | {1,a,c} | {0,a,b,d,e} | {a,b,d} |
| e | e | {b,c,e} | {a,b,c,d} | {1,a,b,c,d,e} | {1,a,b,e} | {a,b,d} | {0,1,b,c,e} |

| $HF_7^{209}$ | 0 | 1 | a | b | c | d | e |
|---|---|---|---|---|---|---|---|
| 0 | 0 | 1 | a | b | c | d | e |
| 1 | 1 | {0,1,a,c,d} | {1,b,c,e} | {a,b,e} | {1,a,c,d} | {1,c,e} | {a,b,d,e} |
| a | a | {1,b,c,e} | {0,a,b,d,e} | {1,a,c,d} | {1,b,c} | {a,b,d,e} | {1,a,d} |
| b | b | {a,b,e} | {1,a,c,d} | {0,1,b,c,e} | {a,b,d,e} | {a,c,d} | {1,b,c,e} |
| c | c | {1,a,c,d} | {1,b,c} | {a,b,d,e} | {0,1,a,c,d} | {1,b,c,e} | {b,d,e} |
| d | d | {1,c,e} | {a,b,d,e} | {a,c,d} | {1,b,c,e} | {0,a,b,d,e} | {1,a,c,d} |
| e | e | {a,b,d,e} | {1,a,d} | {1,b,c,e} | {b,d,e} | {1,a,c,d} | {0,1,b,c,e} |

| $HF_7^{210}$ | 0 | 1 | a | b | c | d | e |
|---|---|---|---|---|---|---|---|
| 0 | 0 | 1 | a | b | c | d | e |
| 1 | 1 | {0,1,a,c,d} | {1,b,c,e} | {a,b,d,e} | {1,a,c,d} | {1,b,c,e} | {a,b,d,e} |
| a | a | {1,b,c,e} | {0,a,b,d,e} | {1,a,c,d} | {1,b,c,e} | {a,b,d,e} | {1,a,c,d} |
| b | b | {a,b,d,e} | {1,a,c,d} | {0,1,b,c,e} | {a,b,d,e} | {1,a,c,d} | {1,b,c,e} |
| c | c | {1,a,c,d} | {1,b,c,e} | {a,b,d,e} | {0,1,a,c,d} | {1,b,c,e} | {a,b,d,e} |
| d | d | {1,b,c,e} | {a,b,d,e} | {1,a,c,d} | {1,b,c,e} | {0,a,b,d,e} | {1,a,c,d} |
| e | e | {a,b,d,e} | {1,a,c,d} | {1,b,c,e} | {a,b,d,e} | {1,a,c,d} | {0,1,b,c,e} |

| $HF_7^{211}$ | 0 | 1 | a | b | c | d | e |
|---|---|---|---|---|---|---|---|
| 0 | 0 | 1 | a | b | c | d | e |
| 1 | 1 | {0,1,a,c,d} | {1,b,d,e} | {a,b,c,d} | {1,b,c,e} | {1,a,b,e} | {a,c,d,e} |
| a | a | {1,b,d,e} | {0,a,b,d,e} | {1,a,c,e} | {b,c,d,e} | {1,a,c,d} | {1,a,b,c} |
| b | b | {a,b,c,d} | {1,a,c,e} | {0,1,b,c,e} | {1,a,b,d} | {1,c,d,e} | {a,b,d,e} |
| c | c | {1,b,c,e} | {b,c,d,e} | {1,a,b,d} | {0,1,a,c,d} | {a,b,c,e} | {1,a,d,e} |
| d | d | {1,a,b,e} | {1,a,c,d} | {1,c,d,e} | {a,b,c,e} | {0,a,b,d,e} | {1,b,c,d} |
| e | e | {a,c,d,e} | {1,a,b,c} | {a,b,d,e} | {1,a,d,e} | {1,b,c,d} | {0,1,b,c,e} |

| $HF_7^{212}$ | 0 | 1 | a | b | c | d | e |
|---|---|---|---|---|---|---|---|
| 0 | 0 | 1 | a | b | c | d | e |
| 1 | 1 | {0,1,a,c,d} | {1,b,c,d,e} | {a,b,c,e} | {1,a,b,c,d,e} | {1,a,c,e} | {a,b,c,d,e} |
| a | a | {1,b,c,d,e} | {0,a,b,d,e} | {1,a,c,d,e} | {1,b,c,d} | {1,a,b,c,d,e} | {1,a,b,d} |
| b | b | {a,b,c,e} | {1,a,c,d,e} | {0,1,b,c,e} | {1,a,b,d,e} | {a,c,d,e} | {1,a,b,c,d,e} |
| c | c | {1,a,b,c,d,e} | {1,b,c,d} | {1,a,b,d,e} | {0,1,a,c,d} | {1,a,b,c,e} | {1,b,d,e} |
| d | d | {1,a,c,e} | {1,a,b,c,d,e} | {a,c,d,e} | {1,a,b,c,e} | {0,a,b,d,e} | {1,a,b,c,d} |
| e | e | {a,b,c,d,e} | {1,a,b,d} | {1,a,b,c,d,e} | {1,b,d,e} | {1,a,b,c,d} | {0,1,b,c,e} |

| $HF_7^{213}$ | 0 | 1 | a | b | c | d | e |
|---|---|---|---|---|---|---|---|
| 0 | 0 | 1 | a | b | c | d | e |
| 1 | 1 | {0,1,a,c,d} | {1,b,c,d,e} | {a,b,c,d,e} | {1,a,b,c,d,e} | {1,a,b,c,e} | {a,b,c,d,e} |
| a | a | {1,b,c,d,e} | {0,a,b,d,e} | {1,a,c,d,e} | {1,b,c,d,e} | {1,a,b,c,d,e} | {1,a,b,c,d} |
| b | b | {a,b,c,d,e} | {1,a,c,d,e} | {0,1,b,c,e} | {1,a,b,d,e} | {1,a,c,d,e} | {1,a,b,c,d,e} |
| c | c | {1,a,b,c,d,e} | {1,b,c,d,e} | {1,a,b,d,e} | {0,1,a,c,d} | {1,a,b,c,e} | {1,a,b,d,e} |
| d | d | {1,a,b,c,e} | {1,a,b,c,d,e} | {1,a,c,d,e} | {1,a,b,c,e} | {0,a,b,d,e} | {1,a,b,c,d} |
| e | e | {a,b,c,d,e} | {1,a,b,c,d} | {1,a,b,c,d,e} | {1,a,b,d,e} | {1,a,b,c,d} | {0,1,b,c,e} |

| $HF_7^{214}$ | 0 | 1 | a | b | c | d | e |
|---|---|---|---|---|---|---|---|
| 0 | 0 | 1 | a | b | c | d | e |
| 1 | 1 | {0,1,a,c,e} | {1,a,b,c,e} | {a,d,e} | {1,a,c,d} | {b,c,e} | {1,a,b,d,e} |
| a | a | {1,a,b,c,e} | {0,1,a,b,d} | {1,a,b,c,d} | {1,b,e} | {a,b,d,e} | {1,c,d} |
| b | b | {a,d,e} | {1,a,b,c,d} | {0,a,b,c,e} | {a,b,c,d,e} | {1,a,c} | {1,b,c,e} |
| c | c | {1,a,c,d} | {1,b,e} | {a,b,c,d,e} | {0,1,b,c,d} | {1,b,c,d,e} | {a,b,d} |
| d | d | {b,c,e} | {a,b,d,e} | {1,a,c} | {1,b,c,d,e} | {0,a,c,d,e} | {1,a,c,d,e} |
| e | e | {1,a,b,d,e} | {1,c,d} | {1,b,c,e} | {a,b,d} | {1,a,c,d,e} | {0,1,b,d,e} |

| $HF_7^{215}$ | 0 | 1 | a | b | c | d | e |
|---|---|---|---|---|---|---|---|
| 0 | 0 | 1 | a | b | c | d | e |
| 1 | 1 | {0,1,a,c,e} | {1,a,b,c,d,e} | {a,c,d,e} | {1,a,b,c,d,e} | {a,b,c,e} | {1,a,b,c,d,e} |
| a | a | {1,a,b,c,d,e} | {0,1,a,b,d} | {1,a,b,c,d,e} | {1,b,d,e} | {1,a,b,c,d,e} | {1,b,c,d} |
| b | b | {a,c,d,e} | {1,a,b,c,d,e} | {0,a,b,c,e} | {1,a,b,c,d,e} | {1,a,c,e} | {1,a,b,c,d,e} |
| c | c | {1,a,b,c,d,e} | {1,b,d,e} | {1,a,b,c,d,e} | {0,1,b,c,d} | {1,a,b,c,d,e} | {1,a,b,d} |
| d | d | {a,b,c,e} | {1,a,b,c,d,e} | {1,a,c,e} | {1,a,b,c,d,e} | {0,a,c,d,e} | {1,a,b,c,d,e} |
| e | e | {1,a,b,c,d,e} | {1,b,c,d} | {1,a,b,c,d,e} | {1,a,b,d} | {1,a,b,c,d,e} | {0,1,b,d,e} |

| $HF_7^{216}$ | 0 | 1 | a | b | c | d | e |
|---|---|---|---|---|---|---|---|
| 0 | 0 | 1 | a | b | c | d | e |
| 1 | 1 | {0,1,b,c,d} | {b,c,e} | {1,a,b,e} | {1,a,c,d} | {1,c,d,e} | {a,b,d} |
| a | a | {b,c,e} | {0,a,c,d,e} | {1,c,d} | {1,a,b,c} | {a,b,d,e} | {1,a,d,e} |
| b | b | {1,a,b,e} | {1,c,d} | {0,1,b,d,e} | {a,d,e} | {a,b,c,d} | {1,b,c,e} |
| c | c | {1,a,c,d} | {1,a,b,c} | {a,d,e} | {0,1,a,c,e} | {1,b,e} | {b,c,d,e} |
| d | d | {1,c,d,e} | {a,b,d,e} | {a,b,c,d} | {1,b,e} | {0,1,a,b,d} | {1,a,c} |
| e | e | {a,b,d} | {1,a,d,e} | {1,b,c,e} | {b,c,d,e} | {1,a,c} | {0,a,b,c,e} |

| $HF_7^{217}$ | 0 | 1 | a | b | c | d | e |
|---|---|---|---|---|---|---|---|
| 0 | 0 | 1 | a | b | c | d | e |
| 1 | 1 | {0,1,b,c,d} | {b,c,e} | {1,a,b,d,e} | {1,a,c,d} | {1,b,c,d,e} | {a,b,d} |
| a | a | {b,c,e} | {0,a,c,d,e} | {1,c,d} | {1,a,b,c,e} | {a,b,d,e} | {1,a,c,d,e} |
| b | b | {1,a,b,d,e} | {1,c,d} | {0,1,b,d,e} | {a,d,e} | {1,a,b,c,d} | {1,b,c,e} |
| c | c | {1,a,c,d} | {1,a,b,c,e} | {a,d,e} | {0,1,a,c,e} | {1,b,e} | {a,b,c,d,e} |
| d | d | {1,b,c,d,e} | {a,b,d,e} | {1,a,b,c,d} | {1,b,e} | {0,1,a,b,d} | {1,a,c} |
| e | e | {a,b,d} | {1,a,c,d,e} | {1,b,c,e} | {a,b,c,d,e} | {1,a,c} | {0,a,b,c,e} |

| $HF_7^{218}$ | 0 | 1 | a | b | c | d | e |
|---|---|---|---|---|---|---|---|
| 0 | 0 | 1 | a | b | c | d | e |
| 1 | 1 | {0,1,b,c,d} | {b,c,d,e} | {1,a,b,c,e} | {1,a,b,c,d,e} | {1,a,c,d,e} | {a,b,c,d} |
| a | a | {b,c,d,e} | {0,a,c,d,e} | {1,c,d,e} | {1,a,b,c,d} | {1,a,b,c,d,e} | {1,a,b,d,e} |
| b | b | {1,a,b,c,e} | {1,c,d,e} | {0,1,b,d,e} | {1,a,d,e} | {a,b,c,d,e} | {1,a,b,c,d,e} |
| c | c | {1,a,b,c,d,e} | {1,a,b,c,d} | {1,a,d,e} | {0,1,a,c,e} | {1,a,b,e} | {1,b,c,d,e} |
| d | d | {1,a,c,d,e} | {1,a,b,c,d,e} | {a,b,c,d,e} | {1,a,b,e} | {0,1,a,b,d} | {1,a,b,c} |
| e | e | {a,b,c,d} | {1,a,b,d,e} | {1,a,b,c,d,e} | {1,b,c,d,e} | {1,a,b,c} | {0,a,b,c,e} |

| $HF_7^{219}$ | 0 | 1 | a | b | c | d | e |
|---|---|---|---|---|---|---|---|
| 0 | 0 | 1 | a | b | c | d | e |
| 1 | 1 | {0,1,b,c,d} | {b,c,d,e} | {1,a,b,c,d,e} | {1,a,b,c,d,e} | {1,a,b,c,d,e} | {a,b,c,d} |
| a | a | {b,c,d,e} | {0,a,c,d,e} | {1,c,d,e} | {1,a,b,c,d,e} | {1,a,b,c,d,e} | {1,a,b,c,d,e} |
| b | b | {1,a,b,c,d,e} | {1,c,d,e} | {0,1,b,d,e} | {1,a,d,e} | {1,a,b,c,d,e} | {1,a,b,c,d,e} |
| c | c | {1,a,b,c,d,e} | {1,a,b,c,d,e} | {1,a,d,e} | {0,1,a,c,e} | {1,a,b,e} | {1,a,b,c,d,e} |
| d | d | {1,a,b,c,d,e} | {1,a,b,c,d,e} | {1,a,b,c,d,e} | {1,a,b,e} | {0,1,a,b,d} | {1,a,b,c} |
| e | e | {a,b,c,d} | {1,a,b,c,d,e} | {1,a,b,c,d,e} | {1,a,b,c,d,e} | {1,a,b,c} | {0,a,b,c,e} |

**B4ii. Hyperfields for which card(x-x) = 5 and x∉ x-x, for every non-zero element x.**

| $HF_7^{220}$ | 0 | 1 | a | b | c | d | e |
|---|---|---|---|---|---|---|---|
| 0 | 0 | 1 | a | b | c | d | e |
| 1 | 1 | {0,a,b,c,d} | {1,c,d} | {1,b,c,e} | {1,a,b,c,d,e} | {1,a,c,d} | {b,c,e} |
| a | a | {1,c,d} | {0,b,c,d,e} | {a,d,e} | {1,a,c,d} | {1,a,b,c,d,e} | {a,b,d,e} |
| b | b | {1,b,c,e} | {a,d,e} | {0,1,c,d,e} | {1,b,e} | {a,b,d,e} | {1,a,b,c,d,e} |
| c | c | {1,a,b,c,d,e} | {1,a,c,d} | {1,b,e} | {0,1,a,d,e} | {1,a,c} | {1,b,c,e} |
| d | d | {1,a,c,d} | {1,a,b,c,d,e} | {a,b,d,e} | {1,a,c} | {0,1,a,b,e} | {a,b,d} |
| e | e | {b,c,e} | {a,b,d,e} | {1,a,b,c,d,e} | {1,b,c,e} | {a,b,d} | {0,1,a,b,c} |

| $HF_7^{221}$ | 0 | 1 | a | b | c | d | e |
|---|---|---|---|---|---|---|---|
| 0 | 0 | 1 | a | b | c | d | e |
| 1 | 1 | {0,a,b,c,d} | {1,c,d} | {1,b,c,d,e} | {1,a,b,c,d,e} | {1,a,b,c,d} | {b,c,e} |
| a | a | {1,c,d} | {0,b,c,d,e} | {a,d,e} | {1,a,c,d,e} | {1,a,b,c,d,e} | {a,b,c,d,e} |
| b | b | {1,b,c,d,e} | {a,d,e} | {0,1,c,d,e} | {1,b,e} | {1,a,b,d,e} | {1,a,b,c,d,e} |
| c | c | {1,a,b,c,d,e} | {1,a,c,d,e} | {1,b,e} | {0,1,a,d,e} | {1,a,c} | {1,a,b,c,e} |
| d | d | {1,a,b,c,d} | {1,a,b,c,d,e} | {1,a,b,d,e} | {1,a,c} | {0,1,a,b,e} | {a,b,d} |
| e | e | {b,c,e} | {a,b,c,d,e} | {1,a,b,c,d,e} | {1,a,b,c,e} | {a,b,d} | {0,1,a,b,c} |

| $HF_7^{222}$ | 0 | 1 | a | b | c | d | e |
|---|---|---|---|---|---|---|---|
| 0 | 0 | 1 | a | b | c | d | e |
| 1 | 1 | {0,a,b,c,d} | {1,b,c,e} | {1,a,b,e} | {1,a,c,d} | {1,c,d,e} | {a,b,d,e} |
| a | a | {1,b,c,e} | {0,b,c,d,e} | {1,a,c,d} | {1,a,b,c} | {a,b,d,e} | {1,a,d,e} |
| b | b | {1,a,b,e} | {1,a,c,d} | {0,1,c,d,e} | {a,b,d,e} | {a,b,c,d} | {1,b,c,e} |
| c | c | {1,a,c,d} | {1,a,b,c} | {a,b,d,e} | {0,1,a,d,e} | {1,b,c,e} | {b,c,d,e} |
| d | d | {1,c,d,e} | {a,b,d,e} | {a,b,c,d} | {1,b,c,e} | {0,1,a,b,e} | {1,a,c,d} |
| e | e | {a,b,d,e} | {1,a,d,e} | {1,b,c,e} | {b,c,d,e} | {1,a,c,d} | {0,1,a,b,c} |

| $HF_7^{223}$ | 0 | 1 | a | b | c | d | e |
|---|---|---|---|---|---|---|---|
| 0 | 0 | 1 | a | b | c | d | e |
| 1 | 1 | {0,a,b,c,d} | {1,b,c,e} | {1,a,b,d,e} | {1,a,c,d} | {1,b,c,d,e} | {a,b,d,e} |
| a | a | {1,b,c,e} | {0,b,c,d,e} | {1,a,c,d} | {1,a,b,c,e} | {a,b,d,e} | {1,a,c,d,e} |
| b | b | {1,a,b,d,e} | {1,a,c,d} | {0,1,c,d,e} | {a,b,d,e} | {1,a,b,c,d} | {1,b,c,e} |
| c | c | {1,a,c,d} | {1,a,b,c,e} | {a,b,d,e} | {0,1,a,d,e} | {1,b,c,e} | {a,b,c,d,e} |
| d | d | {1,b,c,d,e} | {a,b,d,e} | {1,a,b,c,d} | {1,b,c,e} | {0,1,a,b,e} | {1,a,c,d} |
| e | e | {a,b,d,e} | {1,a,c,d,e} | {1,b,c,e} | {a,b,c,d,e} | {1,a,c,d} | {0,1,a,b,c} |

| $HF_7^{224}$ | 0 | 1 | a | b | c | d | e |
|---|---|---|---|---|---|---|---|
| 0 | 0 | 1 | a | b | c | d | e |
| 1 | 1 | {0,a,b,c,d} | {1,b,d,e} | {1,a,b,c} | {1,b,c,e} | {1,a,d,e} | {a,c,d,e} |
| a | a | {1,b,d,e} | {0,b,c,d,e} | {1,a,c,e} | {a,b,c,d} | {1,a,c,d} | {1,a,b,e} |
| b | b | {1,a,b,c} | {1,a,c,e} | {0,1,c,d,e} | {1,a,b,d} | {b,c,d,e} | {a,b,d,e} |
| c | c | {1,b,c,e} | {a,b,c,d} | {1,a,b,d} | {0,1,a,d,e} | {a,b,c,e} | {1,c,d,e} |
| d | d | {1,a,d,e} | {1,a,c,d} | {b,c,d,e} | {a,b,c,e} | {0,1,a,b,e} | {1,b,c,d} |
| e | e | {a,c,d,e} | {1,a,b,e} | {a,b,d,e} | {1,c,d,e} | {1,b,c,d} | {0,1,a,b,c} |

| $HF_7^{225}$ | 0 | 1 | a | b | c | d | e |
|---|---|---|---|---|---|---|---|
| 0 | 0 | 1 | a | b | c | d | e |
| 1 | 1 | {0,a,b,c,d} | {1,b,d,e} | {1,a,b,c,d} | {1,b,c,e} | {1,a,b,d,e} | {a,c,d,e} |
| a | a | {1,b,d,e} | {0,b,c,d,e} | {1,a,c,e} | {a,b,c,d,e} | {1,a,c,d} | {1,a,b,c,e} |
| b | b | {1,a,b,c,d} | {1,a,c,e} | {0,1,c,d,e} | {1,a,b,d} | {1,b,c,d,e} | {a,b,d,e} |
| c | c | {1,b,c,e} | {a,b,c,d,e} | {1,a,b,d} | {0,1,a,d,e} | {a,b,c,e} | {1,a,c,d,e} |
| d | d | {1,a,b,d,e} | {1,a,c,d} | {1,b,c,d,e} | {a,b,c,e} | {0,1,a,b,e} | {1,b,c,d} |
| e | e | {a,c,d,e} | {1,a,b,c,e} | {a,b,d,e} | {1,a,c,d,e} | {1,b,c,d} | {0,1,a,b,c} |

| $HF_7^{226}$ | 0 | 1 | a | b | c | d | e |
|---|---|---|---|---|---|---|---|
| 0 | 0 | 1 | a | b | c | d | e |
| 1 | 1 | {0,a,b,c,d} | {1,b,c,d,e} | {1,a,b,c,e} | {1,a,b,c,d,e} | {1,a,c,d,e} | {a,b,c,d,e} |
| a | a | {1,b,c,d,e} | {0,b,c,d,e} | {1,a,c,d,e} | {1,a,b,c,d} | {1,a,b,c,d,e} | {1,a,b,d,e} |
| b | b | {1,a,b,c,e} | {1,a,c,d,e} | {0,1,c,d,e} | {1,a,b,d,e} | {a,b,c,d,e} | {1,a,b,c,d,e} |
| c | c | {1,a,b,c,d,e} | {1,a,b,c,d} | {1,a,b,d,e} | {0,1,a,d,e} | {1,a,b,c,e} | {1,b,c,d,e} |
| d | d | {1,a,c,d,e} | {1,a,b,c,d,e} | {a,b,c,d,e} | {1,a,b,c,e} | {0,1,a,b,e} | {1,a,b,c,d} |
| e | e | {a,b,c,d,e} | {1,a,b,d,e} | {1,a,b,c,d,e} | {1,b,c,d,e} | {1,a,b,c,d} | {0,1,a,b,c} |

| $HF_7^{227}$ | 0 | 1 | a | b | c | d | e |
|---|---|---|---|---|---|---|---|
| 0 | 0 | 1 | a | b | c | d | e |
| 1 | 1 | {0,a,b,c,d} | {1,b,c,d,e} | {1,a,b,c,d,e} | {1,a,b,c,d,e} | {1,a,b,c,d,e} | {a,b,c,d,e} |
| a | a | {1,b,c,d,e} | {0,b,c,d,e} | {1,a,c,d,e} | {1,a,b,c,d,e} | {1,a,b,c,d,e} | {1,a,b,c,d,e} |
| b | b | {1,a,b,c,d,e} | {1,a,c,d,e} | {0,1,c,d,e} | {1,a,b,d,e} | {1,a,b,c,d,e} | {1,a,b,c,d,e} |
| c | c | {1,a,b,c,d,e} | {1,a,b,c,d,e} | {1,a,b,d,e} | {0,1,a,d,e} | {1,a,b,c,e} | {1,a,b,c,d,e} |
| d | d | {1,a,b,c,d,e} | {1,a,b,c,d,e} | {1,a,b,c,d,e} | {1,a,b,c,e} | {0,1,a,b,e} | {1,a,b,c,d} |
| e | e | {a,b,c,d,e} | {1,a,b,c,d,e} | {1,a,b,c,d,e} | {1,a,b,c,d,e} | {1,a,b,c,d} | {0,1,a,b,c} |

| $HF_7^{228}$ | 0 | 1 | a | b | c | d | e |
|---|---|---|---|---|---|---|---|
| 0 | 0 | 1 | a | b | c | d | e |
| 1 | 1 | {0,a,b,c,e} | {1,a,c,d} | {1,c,e} | {1,a,b,c,d,e} | {a,c,d} | {1,b,c,e} |
| a | a | {1,a,c,d} | {0,1,b,c,d} | {a,b,d,e} | {1,a,d} | {1,a,b,c,d,e} | {b,d,e} |
| b | b | {1,c,e} | {a,b,d,e} | {0,a,c,d,e} | {1,b,c,e} | {a,b,e} | {1,a,b,c,d,e} |
| c | c | {1,a,b,c,d,e} | {1,a,d} | {1,b,c,e} | {0,1,b,d,e} | {1,a,c,d} | {1,b,c} |
| d | d | {a,c,d} | {1,a,b,c,d,e} | {a,b,e} | {1,a,c,d} | {0,1,a,c,e} | {a,b,d,e} |
| e | e | {1,b,c,e} | {b,d,e} | {1,a,b,c,d,e} | {1,b,c} | {a,b,d,e} | {0,1,a,b,d} |

| $HF_7^{229}$ | 0 | 1 | a | b | c | d | e |
|---|---|---|---|---|---|---|---|
| 0 | 0 | 1 | a | b | c | d | e |
| 1 | 1 | {0,a,b,c,e} | {1,a,c,d} | {1,c,d,e} | {1,a,b,c,d,e} | {a,b,c,d} | {1,b,c,e} |
| a | a | {1,a,c,d} | {0,1,b,c,d} | {a,b,d,e} | {1,a,d,e} | {1,a,b,c,d,e} | {b,c,d,e} |
| b | b | {1,c,d,e} | {a,b,d,e} | {0,a,c,d,e} | {1,b,c,e} | {1,a,b,e} | {1,a,b,c,d,e} |
| c | c | {1,a,b,c,d,e} | {1,a,d,e} | {1,b,c,e} | {0,1,b,d,e} | {1,a,c,d} | {1,a,b,c} |
| d | d | {a,b,c,d} | {1,a,b,c,d,e} | {1,a,b,e} | {1,a,c,d} | {0,1,a,c,e} | {a,b,d,e} |
| e | e | {1,b,c,e} | {b,c,d,e} | {1,a,b,c,d,e} | {1,a,b,c} | {a,b,d,e} | {0,1,a,b,d} |

| $HF_7^{230}$ | 0 | 1 | a | b | c | d | e |
|---|---|---|---|---|---|---|---|
| 0 | 0 | 1 | a | b | c | d | e |
| 1 | 1 | {0,a,b,c,e} | {1,a,b,c,e} | {1,a,d,e} | {1,a,c,d} | {b,c,d,e} | {1,a,b,d,e} |
| a | a | {1,a,b,c,e} | {0,1,b,c,d} | {1,a,b,c,d} | {1,a,b,e} | {a,b,d,e} | {1,c,d,e} |
| b | b | {1,a,d,e} | {1,a,b,c,d} | {0,a,c,d,e} | {a,b,c,d,e} | {1,a,b,c} | {1,b,c,e} |
| c | c | {1,a,c,d} | {1,a,b,e} | {a,b,c,d,e} | {0,1,b,d,e} | {1,b,c,d,e} | {a,b,c,d} |
| d | d | {b,c,d,e} | {a,b,d,e} | {1,a,b,c} | {1,b,c,d,e} | {0,1,a,c,e} | {1,a,c,d,e} |
| e | e | {1,a,b,d,e} | {1,c,d,e} | {1,b,c,e} | {a,b,c,d} | {1,a,c,d,e} | {0,1,a,b,d} |

| $HF_7^{231}$ | 0 | 1 | a | b | c | d | e |
|---|---|---|---|---|---|---|---|
| 0 | 0 | 1 | a | b | c | d | e |
| 1 | 1 | {0,a,b,c,e} | {1,a,b,d,e} | {1,a,c} | {1,b,c,e} | {a,d,e} | {1,a,c,d,e} |
| a | a | {1,a,b,d,e} | {0,1,b,c,d} | {1,a,b,c,e} | {a,b,d} | {1,a,c,d} | {1,b,e} |
| b | b | {1,a,c} | {1,a,b,c,e} | {0,a,c,d,e} | {1,a,b,c,d} | {b,c,e} | {a,b,d,e} |
| c | c | {1,b,c,e} | {a,b,d} | {1,a,b,c,d} | {0,1,b,d,e} | {a,b,c,d,e} | {1,c,d} |
| d | d | {a,d,e} | {1,a,c,d} | {b,c,e} | {a,b,c,d,e} | {0,1,a,c,e} | {1,b,c,d,e} |
| e | e | {1,a,c,d,e} | {1,b,e} | {a,b,d,e} | {1,c,d} | {1,b,c,d,e} | {0,1,a,b,d} |

| $HF_7^{232}$ | 0 | 1 | a | b | c | d | e |
|---|---|---|---|---|---|---|---|
| 0 | 0 | 1 | a | b | c | d | e |
| 1 | 1 | {0,a,b,c,e} | {1,a,b,d,e} | {1,a,c,d} | {1,b,c,e} | {a,b,d,e} | {1,a,c,d,e} |
| a | a | {1,a,b,d,e} | {0,1,b,c,d} | {1,a,b,c,e} | {a,b,d,e} | {1,a,c,d} | {1,b,c,e} |
| b | b | {1,a,c,d} | {1,a,b,c,e} | {0,a,c,d,e} | {1,a,b,c,d} | {1,b,c,e} | {a,b,d,e} |
| c | c | {1,b,c,e} | {a,b,d,e} | {1,a,b,c,d} | {0,1,b,d,e} | {a,b,c,d,e} | {1,a,c,d} |
| d | d | {a,b,d,e} | {1,a,c,d} | {1,b,c,e} | {a,b,c,d,e} | {0,1,a,c,e} | {1,b,c,d,e} |
| e | e | {1,a,c,d,e} | {1,b,c,e} | {a,b,d,e} | {1,a,c,d} | {1,b,c,d,e} | {0,1,a,b,d} |

| $HF_7^{233}$ | 0 | 1 | a | b | c | d | e |
|---|---|---|---|---|---|---|---|
| 0 | 0 | 1 | a | b | c | d | e |
| 1 | 1 | {0,a,b,c,e} | {1,a,b,c,d,e} | {1,a,c,e} | {1,a,b,c,d,e} | {a,c,d,e} | {1,a,b,c,d,e} |
| a | a | {1,a,b,c,d,e} | {0,1,b,c,d} | {1,a,b,c,d,e} | {1,a,b,d} | {1,a,b,c,d,e} | {1,b,d,e} |
| b | b | {1,a,c,e} | {1,a,b,c,d,e} | {0,a,c,d,e} | {1,a,b,c,d,e} | {a,b,c,e} | {1,a,b,c,d,e} |
| c | c | {1,a,b,c,d,e} | {1,a,b,d} | {1,a,b,c,d,e} | {0,1,b,d,e} | {1,a,b,c,d,e} | {1,b,c,d} |
| d | d | {a,c,d,e} | {1,a,b,c,d,e} | {a,b,c,e} | {1,a,b,c,d,e} | {0,1,a,c,e} | {1,a,b,c,d,e} |
| e | e | {1,a,b,c,d,e} | {1,b,d,e} | {1,a,b,c,d,e} | {1,b,c,d} | {1,a,b,c,d,e} | {0,1,a,b,d} |

| $HF_7^{234}$ | 0 | 1 | a | b | c | d | e |
|---|---|---|---|---|---|---|---|
| 0 | 0 | 1 | a | b | c | d | e |
| 1 | 1 | {0,a,b,c,e} | {1,a,b,c,d,e} | {1,a,c,d,e} | {1,a,b,c,d,e} | {a,b,c,d,e} | {1,a,b,c,d,e} |
| a | a | {1,a,b,c,d,e} | {0,1,b,c,d} | {1,a,b,c,d,e} | {1,a,b,d,e} | {1,a,b,c,d,e} | {1,b,c,d,e} |
| b | b | {1,a,c,d,e} | {1,a,b,c,d,e} | {0,a,c,d,e} | {1,a,b,c,d,e} | {1,a,b,c,e} | {1,a,b,c,d,e} |
| c | c | {1,a,b,c,d,e} | {1,a,b,d,e} | {1,a,b,c,d,e} | {0,1,b,d,e} | {1,a,b,c,d,e} | {1,a,b,c,d} |
| d | d | {a,b,c,d,e} | {1,a,b,c,d,e} | {1,a,b,c,e} | {1,a,b,c,d,e} | {0,1,a,c,e} | {1,a,b,c,d,e} |
| e | e | {1,a,b,c,d,e} | {1,b,c,d,e} | {1,a,b,c,d,e} | {1,a,b,c,d} | {1,a,b,c,d,e} | {0,1,a,b,d} |

| $HF_7^{235}$ | 0 | 1 | a | b | c | d | e |
|---|---|---|---|---|---|---|---|
| 0 | 0 | 1 | a | b | c | d | e |
| 1 | 1 | {0,a,b,d,e} | {1,a,c,d} | {1,b,c,e} | {a,b,d,e} | {1,a,c,d} | {1,b,c,e} |
| a | a | {1,a,c,d} | {0,1,b,c,e} | {a,b,d,e} | {1,a,c,d} | {1,b,c,e} | {a,b,d,e} |
| b | b | {1,b,c,e} | {a,b,d,e} | {0,1,a,c,d} | {1,b,c,e} | {a,b,d,e} | {1,a,c,d} |
| c | c | {a,b,d,e} | {1,a,c,d} | {1,b,c,e} | {0,a,b,d,e} | {1,a,c,d} | {1,b,c,e} |
| d | d | {1,a,c,d} | {1,b,c,e} | {a,b,d,e} | {1,a,c,d} | {0,1,b,c,e} | {a,b,d,e} |
| e | e | {1,b,c,e} | {a,b,d,e} | {1,a,c,d} | {1,b,c,e} | {a,b,d,e} | {0,1,a,c,d} |

| $HF_7^{236}$ | 0 | 1 | a | b | c | d | e |
|---|---|---|---|---|---|---|---|
| 0 | 0 | 1 | a | b | c | d | e |
| 1 | 1 | {0,a,b,d,e} | {1,a,c,d} | {1,b,c,d,e} | {a,b,d,e} | {1,a,b,c,d} | {1,b,c,e} |
| a | a | {1,a,c,d} | {0,1,b,c,e} | {a,b,d,e} | {1,a,c,d,e} | {1,b,c,e} | {a,b,c,d,e} |
| b | b | {1,b,c,d,e} | {a,b,d,e} | {0,1,a,c,d} | {1,b,c,e} | {1,a,b,d,e} | {1,a,c,d} |
| c | c | {a,b,d,e} | {1,a,c,d,e} | {1,b,c,e} | {0,a,b,d,e} | {1,a,c,d} | {1,a,b,c,e} |
| d | d | {1,a,b,c,d} | {1,b,c,e} | {1,a,b,d,e} | {1,a,c,d} | {0,1,b,c,e} | {a,b,d,e} |
| e | e | {1,b,c,e} | {a,b,c,d,e} | {1,a,c,d} | {1,a,b,c,e} | {a,b,d,e} | {0,1,a,c,d} |

| $HF_7^{237}$ | 0 | 1 | a | b | c | d | e |
|---|---|---|---|---|---|---|---|
| 0 | 0 | 1 | a | b | c | d | e |
| 1 | 1 | {0,a,b,d,e} | {1,a,b,c,d,e} | {1,a,b,c,e} | {a,b,d,e} | {1,a,c,d,e} | {1,a,b,c,d,e} |
| a | a | {1,a,b,c,d,e} | {0,1,b,c,e} | {1,a,b,c,d,e} | {1,a,b,c,d} | {1,b,c,e} | {1,a,b,d,e} |
| b | b | {1,a,b,c,e} | {1,a,b,c,d,e} | {0,1,a,c,d} | {1,a,b,c,d,e} | {a,b,c,d,e} | {1,a,c,d} |
| c | c | {a,b,d,e} | {1,a,b,c,d} | {1,a,b,c,d,e} | {0,a,b,d,e} | {1,a,b,c,d,e} | {1,b,c,d,e} |
| d | d | {1,a,c,d,e} | {1,b,c,e} | {a,b,c,d,e} | {1,a,b,c,d,e} | {0,1,b,c,e} | {1,a,b,c,d,e} |
| e | e | {1,a,b,c,d,e} | {1,a,b,d,e} | {1,a,c,d} | {1,b,c,d,e} | {1,a,b,c,d,e} | {0,1,a,c,d} |

| $HF_7^{238}$ | 0 | 1 | a | b | c | d | e |
|---|---|---|---|---|---|---|---|
| 0 | 0 | 1 | a | b | c | d | e |
| 1 | 1 | {0,a,b,d,e} | {1,a,b,c,d,e} | {1,a,b,c,d,e} | {a,b,d,e} | {1,a,b,c,d,e} | {1,a,b,c,d,e} |
| a | a | {1,a,b,c,d,e} | {0,1,b,c,e} | {1,a,b,c,d,e} | {1,a,b,c,d,e} | {1,b,c,e} | {1,a,b,c,d,e} |
| b | b | {1,a,b,c,d,e} | {1,a,b,c,d,e} | {0,1,a,c,d} | {1,a,b,c,d,e} | {1,a,b,c,d,e} | {1,a,c,d} |
| c | c | {a,b,d,e} | {1,a,b,c,d,e} | {1,a,b,c,d,e} | {0,a,b,d,e} | {1,a,b,c,d,e} | {1,a,b,c,d,e} |
| d | d | {1,a,b,c,d,e} | {1,b,c,e} | {1,a,b,c,d,e} | {1,a,b,c,d,e} | {0,1,b,c,e} | {1,a,b,c,d,e} |
| e | e | {1,a,b,c,d,e} | {1,a,b,c,d,e} | {1,a,c,d} | {1,a,b,c,d,e} | {1,a,b,c,d,e} | {0,1,a,c,d} |

**B5i. Hyperfields for which card(x-x) = 6 and x∈ x-x, for every non-zero element x.**

| $HF_7^{239}$ | 0 | 1 | a | b | c | d | e |
|---|---|---|---|---|---|---|---|
| 0 | 0 | 1 | a | b | c | d | e |
| 1 | 1 | {0,1,a,b,c,d} | {1,c,d} | {1,b,c,e} | {1,a,b,c,d,e} | {1,a,c,d} | {b,c,e} |
| a | a | {1,c,d} | {0,a,b,c,d,e} | {a,d,e} | {1,a,c,d} | {1,a,b,c,d,e} | {a,b,d,e} |
| b | b | {1,b,c,e} | {a,d,e} | {0,1,b,c,d,e} | {1,b,e} | {a,b,d,e} | {1,a,b,c,d,e} |
| c | c | {1,a,b,c,d,e} | {1,a,c,d} | {1,b,e} | {0,1,a,c,d,e} | {1,a,c} | {1,b,c,e} |
| d | d | {1,a,c,d} | {1,a,b,c,d,e} | {a,b,d,e} | {1,a,c} | {0,1,a,b,d,e} | {a,b,d} |
| e | e | {b,c,e} | {a,b,d,e} | {1,a,b,c,d,e} | {1,b,c,e} | {a,b,d} | {0,1,a,b,c,e} |

| $HF_7^{240}$ | 0 | 1 | a | b | c | d | e |
|---|---|---|---|---|---|---|---|
| 0 | 0 | 1 | a | b | c | d | e |
| 1 | 1 | {0,1,a,b,c,d} | {1,c,d} | {1,b,c,d,e} | {1,a,b,c,d,e} | {1,a,b,c,d} | {b,c,e} |
| a | a | {1,c,d} | {0,a,b,c,d,e} | {a,d,e} | {1,a,c,d,e} | {1,a,b,c,d,e} | {a,b,c,d,e} |
| b | b | {1,b,c,d,e} | {a,d,e} | {0,1,b,c,d,e} | {1,b,e} | {1,a,b,d,e} | {1,a,b,c,d,e} |
| c | c | {1,a,b,c,d,e} | {1,a,c,d,e} | {1,b,e} | {0,1,a,c,d,e} | {1,a,c} | {1,a,b,c,e} |
| d | d | {1,a,b,c,d} | {1,a,b,c,d,e} | {1,a,b,d,e} | {1,a,c} | {0,1,a,b,d,e} | {a,b,d} |
| e | e | {b,c,e} | {a,b,c,d,e} | {1,a,b,c,d,e} | {1,a,b,c,e} | {a,b,d} | {0,1,a,b,c,e} |

| $HF_7^{241}$ | 0 | 1 | a | b | c | d | e |
|---|---|---|---|---|---|---|---|
| 0 | 0 | 1 | a | b | c | d | e |
| 1 | 1 | {0,1,a,b,c,d} | {1,b,c,e} | {1,a,b,e} | {1,a,c,d} | {1,c,d,e} | {a,b,d,e} |
| a | a | {1,b,c,e} | {0,a,b,c,d,e} | {1,a,c,d} | {1,a,b,c} | {a,b,d,e} | {1,a,d,e} |
| b | b | {1,a,b,e} | {1,a,c,d} | {0,1,b,c,d,e} | {a,b,d,e} | {a,b,c,d} | {1,b,c,e} |
| c | c | {1,a,c,d} | {1,a,b,c} | {a,b,d,e} | {0,1,a,c,d,e} | {1,b,c,e} | {b,c,d,e} |
| d | d | {1,c,d,e} | {a,b,d,e} | {a,b,c,d} | {1,b,c,e} | {0,1,a,b,d,e} | {1,a,c,d} |
| e | e | {a,b,d,e} | {1,a,d,e} | {1,b,c,e} | {b,c,d,e} | {1,a,c,d} | {0,1,a,b,c,e} |

| $HF_7^{242}$ | 0 | 1 | a | b | c | d | e |
|---|---|---|---|---|---|---|---|
| 0 | 0 | 1 | a | b | c | d | e |
| 1 | 1 | {0,1,a,b,c,d} | {1,b,c,e} | {1,a,b,d,e} | {1,a,c,d} | {1,b,c,d,e} | {a,b,d,e} |
| a | a | {1,b,c,e} | {0,a,b,c,d,e} | {1,a,c,d} | {1,a,b,c,e} | {a,b,d,e} | {1,a,c,d,e} |
| b | b | {1,a,b,d,e} | {1,a,c,d} | {0,1,b,c,d,e} | {a,b,d,e} | {1,a,b,c,d} | {1,b,c,e} |
| c | c | {1,a,c,d} | {1,a,b,c,e} | {a,b,d,e} | {0,1,a,c,d,e} | {1,b,c,e} | {a,b,c,d,e} |
| d | d | {1,b,c,d,e} | {a,b,d,e} | {1,a,b,c,d} | {1,b,c,e} | {0,1,a,b,d,e} | {1,a,c,d} |
| e | e | {a,b,d,e} | {1,a,c,d,e} | {1,b,c,e} | {a,b,c,d,e} | {1,a,c,d} | {0,1,a,b,c,e} |

| $HF_7^{243}$ | 0 | 1 | a | b | c | d | e |
|---|---|---|---|---|---|---|---|
| 0 | 0 | 1 | a | b | c | d | e |
| 1 | 1 | {0,1,a,b,c,d} | {1,b,d,e} | {1,a,b,c} | {1,b,c,e} | {1,a,d,e} | {a,c,d,e} |
| a | a | {1,b,d,e} | {0,a,b,c,d,e} | {1,a,c,e} | {a,b,c,d} | {1,a,c,d} | {1,a,b,e} |
| b | b | {1,a,b,c} | {1,a,c,e} | {0,1,b,c,d,e} | {1,a,b,d} | {b,c,d,e} | {a,b,d,e} |
| c | c | {1,b,c,e} | {a,b,c,d} | {1,a,b,d} | {0,1,a,c,d,e} | {a,b,c,e} | {1,c,d,e} |
| d | d | {1,a,d,e} | {1,a,c,d} | {b,c,d,e} | {a,b,c,e} | {0,1,a,b,d,e} | {1,b,c,d} |
| e | e | {a,c,d,e} | {1,a,b,e} | {a,b,d,e} | {1,c,d,e} | {1,b,c,d} | {0,1,a,b,c,e} |

| $HF_7^{244}$ | 0 | 1 | a | b | c | d | e |
|---|---|---|---|---|---|---|---|
| 0 | 0 | 1 | a | b | c | d | e |
| 1 | 1 | {0,1,a,b,c,d} | {1,b,d,e} | {1,a,b,c,d} | {1,b,c,e} | {1,a,b,d,e} | {a,c,d,e} |
| a | a | {1,b,d,e} | {0,a,b,c,d,e} | {1,a,c,e} | {a,b,c,d,e} | {1,a,c,d} | {1,a,b,c,e} |
| b | b | {1,a,b,c,d} | {1,a,c,e} | {0,1,b,c,d,e} | {1,a,b,d} | {1,b,c,d,e} | {a,b,d,e} |
| c | c | {1,b,c,e} | {a,b,c,d,e} | {1,a,b,d} | {0,1,a,c,d,e} | {a,b,c,e} | {1,a,c,d,e} |
| d | d | {1,a,b,d,e} | {1,a,c,d} | {1,b,c,d,e} | {a,b,c,e} | {0,1,a,b,d,e} | {1,b,c,d} |
| e | e | {a,c,d,e} | {1,a,b,c,e} | {a,b,d,e} | {1,a,c,d,e} | {1,b,c,d} | {0,1,a,b,c,e} |

| $HF_7^{245}$ | 0 | 1 | a | b | c | d | e |
|---|---|---|---|---|---|---|---|
| 0 | 0 | 1 | a | b | c | d | e |
| 1 | 1 | {0,1,a,b,c,d} | {1,b,c,d,e} | {1,a,b,c,e} | {1,a,b,c,d,e} | {1,a,c,d,e} | {a,b,c,d,e} |
| a | a | {1,b,c,d,e} | {0,a,b,c,d,e} | {1,a,c,d,e} | {1,a,b,c,d} | {1,a,b,c,d,e} | {1,a,b,d,e} |
| b | b | {1,a,b,c,e} | {1,a,c,d,e} | {0,1,b,c,d,e} | {1,a,b,d,e} | {a,b,c,d,e} | {1,a,b,c,d,e} |
| c | c | {1,a,b,c,d,e} | {1,a,b,c,d} | {1,a,b,d,e} | {0,1,a,c,d,e} | {1,a,b,c,e} | {1,b,c,d,e} |
| d | d | {1,a,c,d,e} | {1,a,b,c,d,e} | {a,b,c,d,e} | {1,a,b,c,e} | {0,1,a,b,d,e} | {1,a,b,c,d} |
| e | e | {a,b,c,d,e} | {1,a,b,d,e} | {1,a,b,c,d,e} | {1,b,c,d,e} | {1,a,b,c,d} | {0,1,a,b,c,e} |

| $HF_7^{246}$ | 0 | 1 | a | b | c | d | e |
|---|---|---|---|---|---|---|---|
| 0 | 0 | 1 | a | b | c | d | e |
| 1 | 1 | {0,1,a,b,c,d} | {1,b,c,d,e} | {1,a,b,c,d,e} | {1,a,b,c,d,e} | {1,a,b,c,d,e} | {a,b,c,d,e} |
| a | a | {1,b,c,d,e} | {0,a,b,c,d,e} | {1,a,c,d,e} | {1,a,b,c,d,e} | {1,a,b,c,d,e} | {1,a,b,c,d,e} |
| b | b | {1,a,b,c,d,e} | {1,a,c,d,e} | {0,1,b,c,d,e} | {1,a,b,d,e} | {1,a,b,c,d,e} | {1,a,b,c,d,e} |
| c | c | {1,a,b,c,d,e} | {1,a,b,c,d,e} | {1,a,b,d,e} | {0,1,a,c,d,e} | {1,a,b,c,e} | {1,a,b,c,d,e} |
| d | d | {1,a,b,c,d,e} | {1,a,b,c,d,e} | {1,a,b,c,d,e} | {1,a,b,c,e} | {0,1,a,b,d,e} | {1,a,b,c,d} |
| e | e | {a,b,c,d,e} | {1,a,b,c,d,e} | {1,a,b,c,d,e} | {1,a,b,c,d,e} | {1,a,b,c,d} | {0,1,a,b,c,e} |

| $HF_7^{247}$ | 0 | 1 | a | b | c | d | e |
|---|---|---|---|---|---|---|---|
| 0 | 0 | 1 | a | b | c | d | e |
| 1 | 1 | {0,1,a,b,c,e} | {1,a,c,d} | {1,c,e} | {1,a,b,c,d,e} | {a,c,d} | {1,b,c,e} |
| a | a | {1,a,c,d} | {0,1,a,b,c,d} | {a,b,d,e} | {1,a,d} | {1,a,b,c,d,e} | {b,d,e} |
| b | b | {1,c,e} | {a,b,d,e} | {0,a,b,c,d,e} | {1,b,c,e} | {a,b,e} | {1,a,b,c,d,e} |
| c | c | {1,a,b,c,d,e} | {1,a,d} | {1,b,c,e} | {0,1,b,c,d,e} | {1,a,c,d} | {1,b,c} |
| d | d | {a,c,d} | {1,a,b,c,d,e} | {a,b,e} | {1,a,c,d} | {0,1,a,c,d,e} | {a,b,d,e} |
| e | e | {1,b,c,e} | {b,d,e} | {1,a,b,c,d,e} | {1,b,c} | {a,b,d,e} | {0,1,a,b,d,e} |

| $HF_7^{248}$ | 0 | 1 | a | b | c | d | e |
|---|---|---|---|---|---|---|---|
| 0 | 0 | 1 | a | b | c | d | e |
| 1 | 1 | {0,1,a,b,c,e} | {1,a,c,d} | {1,c,d,e} | {1,a,b,c,d,e} | {a,b,c,d} | {1,b,c,e} |
| a | a | {1,a,c,d} | {0,1,a,b,c,d} | {a,b,d,e} | {1,a,d,e} | {1,a,b,c,d,e} | {b,c,d,e} |
| b | b | {1,c,d,e} | {a,b,d,e} | {0,a,b,c,d,e} | {1,b,c,e} | {1,a,b,e} | {1,a,b,c,d,e} |
| c | c | {1,a,b,c,d,e} | {1,a,d,e} | {1,b,c,e} | {0,1,b,c,d,e} | {1,a,c,d} | {1,a,b,c} |
| d | d | {a,b,c,d} | {1,a,b,c,d,e} | {1,a,b,e} | {1,a,c,d} | {0,1,a,c,d,e} | {a,b,d,e} |
| e | e | {1,b,c,e} | {b,c,d,e} | {1,a,b,c,d,e} | {1,a,b,c} | {a,b,d,e} | {0,1,a,b,d,e} |

| $HF_7^{249}$ | 0 | 1 | a | b | c | d | e |
|---|---|---|---|---|---|---|---|
| 0 | 0 | 1 | a | b | c | d | e |
| 1 | 1 | {0,1,a,b,c,e} | {1,a,b,c,e} | {1,a,d,e} | {1,a,c,d} | {b,c,d,e} | {1,a,b,d,e} |
| a | a | {1,a,b,c,e} | {0,1,a,b,c,d} | {1,a,b,c,d} | {1,a,b,e} | {a,b,d,e} | {1,c,d,e} |
| b | b | {1,a,d,e} | {1,a,b,c,d} | {0,a,b,c,d,e} | {a,b,c,d,e} | {1,a,b,c} | {1,b,c,e} |
| c | c | {1,a,c,d} | {1,a,b,e} | {a,b,c,d,e} | {0,1,b,c,d,e} | {1,b,c,d,e} | {a,b,c,d} |
| d | d | {b,c,d,e} | {a,b,d,e} | {1,a,b,c} | {1,b,c,d,e} | {0,1,a,c,d,e} | {1,a,c,d,e} |
| e | e | {1,a,b,d,e} | {1,c,d,e} | {1,b,c,e} | {a,b,c,d} | {1,a,c,d,e} | {0,1,a,b,d,e} |

| $HF_7^{250}$ | 0 | 1 | a | b | c | d | e |
|---|---|---|---|---|---|---|---|
| 0 | 0 | 1 | a | b | c | d | e |
| 1 | 1 | {0,1,a,b,c,e} | {1,a,b,d,e} | {1,a,c} | {1,b,c,e} | {a,d,e} | {1,a,c,d,e} |
| a | a | {1,a,b,d,e} | {0,1,a,b,c,d} | {1,a,b,c,e} | {a,b,d} | {1,a,c,d} | {1,b,e} |
| b | b | {1,a,c} | {1,a,b,c,e} | {0,a,b,c,d,e} | {1,a,b,c,d} | {b,c,e} | {a,b,d,e} |
| c | c | {1,b,c,e} | {a,b,d} | {1,a,b,c,d} | {0,1,b,c,d,e} | {a,b,c,d,e} | {1,c,d} |
| d | d | {a,d,e} | {1,a,c,d} | {b,c,e} | {a,b,c,d,e} | {0,1,a,c,d,e} | {1,b,c,d,e} |
| e | e | {1,a,c,d,e} | {1,b,e} | {a,b,d,e} | {1,c,d} | {1,b,c,d,e} | {0,1,a,b,d,e} |

| $HF_7^{251}$ | 0 | 1 | a | b | c | d | e |
|---|---|---|---|---|---|---|---|
| 0 | 0 | 1 | a | b | c | d | e |
| 1 | 1 | {0,1,a,b,c,e} | {1,a,b,d,e} | {1,a,c,d} | {1,b,c,e} | {a,b,d,e} | {1,a,c,d,e} |
| a | a | {1,a,b,d,e} | {0,1,a,b,c,d} | {1,a,b,c,e} | {a,b,d,e} | {1,a,c,d} | {1,b,c,e} |
| b | b | {1,a,c,d} | {1,a,b,c,e} | {0,a,b,c,d,e} | {1,a,b,c,d} | {1,b,c,e} | {a,b,d,e} |
| c | c | {1,b,c,e} | {a,b,d,e} | {1,a,b,c,d} | {0,1,b,c,d,e} | {a,b,c,d,e} | {1,a,c,d} |
| d | d | {a,b,d,e} | {1,a,c,d} | {1,b,c,e} | {a,b,c,d,e} | {0,1,a,c,d,e} | {1,b,c,d,e} |
| e | e | {1,a,c,d,e} | {1,b,c,e} | {a,b,d,e} | {1,a,c,d} | {1,b,c,d,e} | {0,1,a,b,d,e} |

| $HF_7^{252}$ | 0 | 1 | a | b | c | d | e |
|---|---|---|---|---|---|---|---|
| 0 | 0 | 1 | a | b | c | d | e |
| 1 | 1 | {0,1,a,b,c,e} | {1,a,b,c,d,e} | {1,a,c,e} | {1,a,b,c,d,e} | {a,c,d,e} | {1,a,b,c,d,e} |
| a | a | {1,a,b,c,d,e} | {0,1,a,b,c,d} | {1,a,b,c,d,e} | {1,a,b,d} | {1,a,b,c,d,e} | {1,b,d,e} |
| b | b | {1,a,c,e} | {1,a,b,c,d,e} | {0,a,b,c,d,e} | {1,a,b,c,d,e} | {a,b,c,e} | {1,a,b,c,d,e} |
| c | c | {1,a,b,c,d,e} | {1,a,b,d} | {1,a,b,c,d,e} | {0,1,b,c,d,e} | {1,a,b,c,d,e} | {1,b,c,d} |
| d | d | {a,c,d,e} | {1,a,b,c,d,e} | {a,b,c,e} | {1,a,b,c,d,e} | {0,1,a,c,d,e} | {1,a,b,c,d,e} |
| e | e | {1,a,b,c,d,e} | {1,b,d,e} | {1,a,b,c,d,e} | {1,b,c,d} | {1,a,b,c,d,e} | {0,1,a,b,d,e} |

| $HF_7^{253}$ | 0 | 1 | a | b | c | d | e |
|---|---|---|---|---|---|---|---|
| 0 | 0 | 1 | a | b | c | d | e |
| 1 | 1 | {0,1,a,b,c,e} | {1,a,b,c,d,e} | {1,a,c,d,e} | {1,a,b,c,d,e} | {a,b,c,d,e} | {1,a,b,c,d,e} |
| a | a | {1,a,b,c,d,e} | {0,1,a,b,c,d} | {1,a,b,c,d,e} | {1,a,b,d,e} | {1,a,b,c,d,e} | {1,b,c,d,e} |
| b | b | {1,a,c,d,e} | {1,a,b,c,d,e} | {0,a,b,c,d,e} | {1,a,b,c,d,e} | {1,a,b,c,e} | {1,a,b,c,d,e} |
| c | c | {1,a,b,c,d,e} | {1,a,b,d,e} | {1,a,b,c,d,e} | {0,1,b,c,d,e} | {1,a,b,c,d,e} | {1,a,b,c,d} |
| d | d | {a,b,c,d,e} | {1,a,b,c,d,e} | {1,a,b,c,e} | {1,a,b,c,d,e} | {0,1,a,c,d,e} | {1,a,b,c,d,e} |
| e | e | {1,a,b,c,d,e} | {1,b,c,d,e} | {1,a,b,c,d,e} | {1,a,b,c,d} | {1,a,b,c,d,e} | {0,1,a,b,d,e} |

| $HF_7^{254}$ | 0 | 1 | a | b | c | d | e |
|---|---|---|---|---|---|---|---|
| 0 | 0 | 1 | a | b | c | d | e |
| 1 | 1 | {0,1,a,b,d,e} | {1,a,c,d} | {1,b,c,e} | {a,b,d,e} | {1,a,c,d} | {1,b,c,e} |
| a | a | {1,a,c,d} | {0,1,a,b,c,e} | {a,b,d,e} | {1,a,c,d} | {1,b,c,e} | {a,b,d,e} |
| b | b | {1,b,c,e} | {a,b,d,e} | {0,1,a,b,c,d} | {1,b,c,e} | {a,b,d,e} | {1,a,c,d} |
| c | c | {a,b,d,e} | {1,a,c,d} | {1,b,c,e} | {0,a,b,c,d,e} | {1,a,c,d} | {1,b,c,e} |
| d | d | {1,a,c,d} | {1,b,c,e} | {a,b,d,e} | {1,a,c,d} | {0,1,b,c,d,e} | {a,b,d,e} |
| e | e | {1,b,c,e} | {a,b,d,e} | {1,a,c,d} | {1,b,c,e} | {a,b,d,e} | {0,1,a,c,d,e} |

| $HF_7^{255}$ | 0 | 1 | a | b | c | d | e |
|---|---|---|---|---|---|---|---|
| 0 | 0 | 1 | a | b | c | d | e |
| 1 | 1 | {0,1,a,b,d,e} | {1,a,c,d} | {1,b,c,d,e} | {a,b,d,e} | {1,a,b,c,d} | {1,b,c,e} |
| a | a | {1,a,c,d} | {0,1,a,b,c,e} | {a,b,d,e} | {1,a,c,d,e} | {1,b,c,e} | {a,b,c,d,e} |
| b | b | {1,b,c,d,e} | {a,b,d,e} | {0,1,a,b,c,d} | {1,b,c,e} | {1,a,b,d,e} | {1,a,c,d} |
| c | c | {a,b,d,e} | {1,a,c,d,e} | {1,b,c,e} | {0,a,b,c,d,e} | {1,a,c,d} | {1,a,b,c,e} |
| d | d | {1,a,b,c,d} | {1,b,c,e} | {1,a,b,d,e} | {1,a,c,d} | {0,1,b,c,d,e} | {a,b,d,e} |
| e | e | {1,b,c,e} | {a,b,c,d,e} | {1,a,c,d} | {1,a,b,c,e} | {a,b,d,e} | {0,1,a,c,d,e} |

| $HF_7^{256}$ | 0 | 1 | a | b | c | d | e |
|---|---|---|---|---|---|---|---|
| 0 | 0 | 1 | a | b | c | d | e |
| 1 | 1 | {0,1,a,b,d,e} | {1,a,b,c,d,e} | {1,a,b,c,e} | {a,b,d,e} | {1,a,c,d,e} | {1,a,b,c,d,e} |
| a | a | {1,a,b,c,d,e} | {0,1,a,b,c,e} | {1,a,b,c,d,e} | {1,a,b,c,d} | {1,b,c,e} | {1,a,b,d,e} |
| b | b | {1,a,b,c,e} | {1,a,b,c,d,e} | {0,1,a,b,c,d} | {1,a,b,c,d,e} | {a,b,c,d,e} | {1,a,c,d} |
| c | c | {a,b,d,e} | {1,a,b,c,d} | {1,a,b,c,d,e} | {0,a,b,c,d,e} | {1,a,b,c,d,e} | {1,b,c,d,e} |
| d | d | {1,a,c,d,e} | {1,b,c,e} | {a,b,c,d,e} | {1,a,b,c,d,e} | {0,1,b,c,d,e} | {1,a,b,c,d,e} |
| e | e | {1,a,b,c,d,e} | {1,a,b,d,e} | {1,a,c,d} | {1,b,c,d,e} | {1,a,b,c,d,e} | {0,1,a,c,d,e} |

| $HF_7^{257}$ | 0 | 1 | a | b | c | d | e |
|---|---|---|---|---|---|---|---|
| 0 | 0 | 1 | a | b | c | d | e |
| 1 | 1 | {0,1,a,b,d,e} | {1,a,b,c,d,e} | {1,a,b,c,d,e} | {a,b,d,e} | {1,a,b,c,d,e} | {1,a,b,c,d,e} |
| a | a | {1,a,b,c,d,e} | {0,1,a,b,c,e} | {1,a,b,c,d,e} | {1,a,b,c,d,e} | {1,b,c,e} | {1,a,b,c,d,e} |
| b | b | {1,a,b,c,d,e} | {1,a,b,c,d,e} | {0,1,a,b,c,d} | {1,a,b,c,d,e} | {1,a,b,c,d,e} | {1,a,c,d} |
| c | c | {a,b,d,e} | {1,a,b,c,d,e} | {1,a,b,c,d,e} | {0,a,b,c,d,e} | {1,a,b,c,d,e} | {1,a,b,c,d,e} |
| d | d | {1,a,b,c,d,e} | {1,b,c,e} | {1,a,b,c,d,e} | {1,a,b,c,d,e} | {0,1,b,c,d,e} | {1,a,b,c,d,e} |
| e | e | {1,a,b,c,d,e} | {1,a,b,c,d,e} | {1,a,c,d} | {1,a,b,c,d,e} | {1,a,b,c,d,e} | {0,1,a,c,d,e} |

**B5ii. Hyperfields for which card(x-x) = 6 and x∉ x-x, for every non-zero element x.**

| $HF_7^{258}$ | 0 | 1 | a | b | c | d | e |
|---|---|---|---|---|---|---|---|
| 0 | 0 | 1 | a | b | c | d | e |
| 1 | 1 | {0,a,b,c,d,e} | {1,a} | {1,b} | {1,c} | {1,d} | {1,e} |
| a | a | {1,a} | {0,1,b,c,d,e} | {a,b} | {a,c} | {a,d} | {a,e} |
| b | b | {1,b} | {a,b} | {0,1,a,c,d,e} | {b,c} | {b,d} | {b,e} |
| c | c | {1,c} | {a,c} | {b,c} | {0,1,a,b,d,e} | {c,d} | {c,e} |
| d | d | {1,d} | {a,d} | {b,d} | {c,d} | {0,1,a,b,c,e} | {d,e} |
| e | e | {1,e} | {a,e} | {b,e} | {c,e} | {d,e} | {0,1,a,b,c,d} |

| $HF_7^{259}$ | 0 | 1 | a | b | c | d | e |
|---|---|---|---|---|---|---|---|
| 0 | 0 | 1 | a | b | c | d | e |
| 1 | 1 | {0,a,b,c,d,e} | {1,a} | {1,b,d} | {1,c} | {1,b,d} | {1,e} |
| a | a | {1,a} | {0,1,b,c,d,e} | {a,b} | {a,c,e} | {a,d} | {a,c,e} |
| b | b | {1,b,d} | {a,b} | {0,1,a,c,d,e} | {b,c} | {1,b,d} | {b,e} |
| c | c | {1,c} | {a,c,e} | {b,c} | {0,1,a,b,d,e} | {c,d} | {a,c,e} |
| d | d | {1,b,d} | {a,d} | {1,b,d} | {c,d} | {0,1,a,b,c,e} | {d,e} |
| e | e | {1,e} | {a,c,e} | {b,e} | {a,c,e} | {d,e} | {0,1,a,b,c,d} |

| $HF_7^{260}$ | 0 | 1 | a | b | c | d | e |
|---|---|---|---|---|---|---|---|
| 0 | 0 | 1 | a | b | c | d | e |
| 1 | 1 | {0,a,b,c,d,e} | {1,a,c} | {1,b,e} | {1,a,c,d} | {1,c,d} | {1,b,e} |
| a | a | {1,a,c} | {0,1,b,c,d,e} | {a,b,d} | {1,a,c} | {a,b,d,e} | {a,d,e} |
| b | b | {1,b,e} | {a,b,d} | {0,1,a,c,d,e} | {b,c,e} | {a,b,d} | {1,b,c,e} |
| c | c | {1,a,c,d} | {1,a,c} | {b,c,e} | {0,1,a,b,d,e} | {1,c,d} | {b,c,e} |
| d | d | {1,c,d} | {a,b,d,e} | {a,b,d} | {1,c,d} | {0,1,a,b,c,e} | {a,d,e} |
| e | e | {1,b,e} | {a,d,e} | {1,b,c,e} | {b,c,e} | {a,d,e} | {0,1,a,b,c,d} |

| $HF_7^{261}$ | 0 | 1 | a | b | c | d | e |
|---|---|---|---|---|---|---|---|
| 0 | 0 | 1 | a | b | c | d | e |
| 1 | 1 | {0,a,b,c,d,e} | {1,a,c} | {1,b,d,e} | {1,a,c,d} | {1,b,c,d} | {1,b,e} |
| a | a | {1,a,c} | {0,1,b,c,d,e} | {a,b,d} | {1,a,c,e} | {a,b,d,e} | {a,c,d,e} |
| b | b | {1,b,d,e} | {a,b,d} | {0,1,a,c,d,e} | {b,c,e} | {1,a,b,d} | {1,b,c,e} |
| c | c | {1,a,c,d} | {1,a,c,e} | {b,c,e} | {0,1,a,b,d,e} | {1,c,d} | {a,b,c,e} |
| d | d | {1,b,c,d} | {a,b,d,e} | {1,a,b,d} | {1,c,d} | {0,1,a,b,c,e} | {a,d,e} |
| e | e | {1,b,e} | {a,c,d,e} | {1,b,c,e} | {a,b,c,e} | {a,d,e} | {0,1,a,b,c,d} |

| $HF_7^{262}$ | 0 | 1 | a | b | c | d | e |
|---|---|---|---|---|---|---|---|
| 0 | 0 | 1 | a | b | c | d | e |
| 1 | 1 | {0,a,b,c,d,e} | {1,a,c,d} | {1,b,c,e} | {1,a,b,c,d,e} | {1,a,c,d} | {1,b,c,e} |
| a | a | {1,a,c,d} | {0,1,b,c,d,e} | {a,b,d,e} | {1,a,c,d} | {1,a,b,c,d,e} | {a,b,d,e} |
| b | b | {1,b,c,e} | {a,b,d,e} | {0,1,a,c,d,e} | {1,b,c,e} | {a,b,d,e} | {1,a,b,c,d,e} |
| c | c | {1,a,b,c,d,e} | {1,a,c,d} | {1,b,c,e} | {0,1,a,b,d,e} | {1,a,c,d} | {1,b,c,e} |
| d | d | {1,a,c,d} | {1,a,b,c,d,e} | {a,b,d,e} | {1,a,c,d} | {0,1,a,b,c,e} | {a,b,d,e} |
| e | e | {1,b,c,e} | {a,b,d,e} | {1,a,b,c,d,e} | {1,b,c,e} | {a,b,d,e} | {0,1,a,b,c,d} |

| $HF_7^{263}$ | 0 | 1 | a | b | c | d | e |
|---|---|---|---|---|---|---|---|
| 0 | 0 | 1 | a | b | c | d | e |
| 1 | 1 | {0,a,b,c,d,e} | {1,a,c,d} | {1,b,c,d,e} | {1,a,b,c,d,e} | {1,a,b,c,d} | {1,b,c,e} |
| a | a | {1,a,c,d} | {0,1,b,c,d,e} | {a,b,d,e} | {1,a,c,d,e} | {1,a,b,c,d,e} | {a,b,c,d,e} |
| b | b | {1,b,c,d,e} | {a,b,d,e} | {0,1,a,c,d,e} | {1,b,c,e} | {1,a,b,d,e} | {1,a,b,c,d,e} |
| c | c | {1,a,b,c,d,e} | {1,a,c,d,e} | {1,b,c,e} | {0,1,a,b,d,e} | {1,a,c,d} | {1,a,b,c,e} |
| d | d | {1,a,b,c,d} | {1,a,b,c,d,e} | {1,a,b,d,e} | {1,a,c,d} | {0,1,a,b,c,e} | {a,b,d,e} |
| e | e | {1,b,c,e} | {a,b,c,d,e} | {1,a,b,c,d,e} | {1,a,b,c,e} | {a,b,d,e} | {0,1,a,b,c,d} |

| $HF_7^{264}$ | 0 | 1 | a | b | c | d | e |
|---|---|---|---|---|---|---|---|
| 0 | 0 | 1 | a | b | c | d | e |
| 1 | 1 | {0,a,b,c,d,e} | {1,a,b,c,e} | {1,a,b,e} | {1,a,c,d} | {1,c,d,e} | {1,a,b,d,e} |
| a | a | {1,a,b,c,e} | {0,1,b,c,d,e} | {1,a,b,c,d} | {1,a,b,c} | {a,b,d,e} | {1,a,d,e} |
| b | b | {1,a,b,e} | {1,a,b,c,d} | {0,1,a,c,d,e} | {a,b,c,d,e} | {a,b,c,d} | {1,b,c,e} |
| c | c | {1,a,c,d} | {1,a,b,c} | {a,b,c,d,e} | {0,1,a,b,d,e} | {1,b,c,d,e} | {b,c,d,e} |
| d | d | {1,c,d,e} | {a,b,d,e} | {a,b,c,d} | {1,b,c,d,e} | {0,1,a,b,c,e} | {1,a,c,d,e} |
| e | e | {1,a,b,d,e} | {1,a,d,e} | {1,b,c,e} | {b,c,d,e} | {1,a,c,d,e} | {0,1,a,b,c,d} |

| $HF_7^{265}$ | 0 | 1 | a | b | c | d | e |
|---|---|---|---|---|---|---|---|
| 0 | 0 | 1 | a | b | c | d | e |
| 1 | 1 | {0,a,b,c,d,e} | {1,a,b,c,e} | {1,a,b,d,e} | {1,a,c,d} | {1,b,c,d,e} | {1,a,b,d,e} |
| a | a | {1,a,b,c,e} | {0,1,b,c,d,e} | {1,a,b,c,d} | {1,a,b,c,e} | {a,b,d,e} | {1,a,c,d,e} |
| b | b | {1,a,b,d,e} | {1,a,b,c,d} | {0,1,a,c,d,e} | {a,b,c,d,e} | {1,a,b,c,d} | {1,b,c,e} |
| c | c | {1,a,c,d} | {1,a,b,c,e} | {a,b,c,d,e} | {0,1,a,b,d,e} | {1,b,c,d,e} | {a,b,c,d,e} |
| d | d | {1,b,c,d,e} | {a,b,d,e} | {1,a,b,c,d} | {1,b,c,d,e} | {0,1,a,b,c,e} | {1,a,c,d,e} |
| e | e | {1,a,b,d,e} | {1,a,c,d,e} | {1,b,c,e} | {a,b,c,d,e} | {1,a,c,d,e} | {0,1,a,b,c,d} |

| $HF_7^{266}$ | 0 | 1 | a | b | c | d | e |
|---|---|---|---|---|---|---|---|
| 0 | 0 | 1 | a | b | c | d | e |
| 1 | 1 | {0,a,b,c,d,e} | {1,a,b,c,d,e} | {1,a,b,c,e} | {1,a,b,c,d,e} | {1,a,c,d,e} | {1,a,b,c,d,e} |
| a | a | {1,a,b,c,d,e} | {0,1,b,c,d,e} | {1,a,b,c,d,e} | {1,a,b,c,d} | {1,a,b,c,d,e} | {1,a,b,d,e} |
| b | b | {1,a,b,c,e} | {1,a,b,c,d,e} | {0,1,a,c,d,e} | {1,a,b,c,d,e} | {a,b,c,d,e} | {1,a,b,c,d,e} |
| c | c | {1,a,b,c,d,e} | {1,a,b,c,d} | {1,a,b,c,d,e} | {0,1,a,b,d,e} | {1,a,b,c,d,e} | {1,b,c,d,e} |
| d | d | {1,a,c,d,e} | {1,a,b,c,d,e} | {a,b,c,d,e} | {1,a,b,c,d,e} | {0,1,a,b,c,e} | {1,a,b,c,d,e} |
| e | e | {1,a,b,c,d,e} | {1,a,b,d,e} | {1,a,b,c,d,e} | {1,b,c,d,e} | {1,a,b,c,d,e} | {0,1,a,b,c,d} |

| $HF_7^{267}$ | 0 | 1 | a | b | c | d | e |
|---|---|---|---|---|---|---|---|
| 0 | 0 | 1 | a | b | c | d | e |
| 1 | 1 | {0,a,b,c,d,e} | {1,a,b,c,d,e} | {1,a,b,c,d,e} | {1,a,b,c,d,e} | {1,a,b,c,d,e} | {1,a,b,c,d,e} |
| a | a | {1,a,b,c,d,e} | {0,1,b,c,d,e} | {1,a,b,c,d,e} | {1,a,b,c,d,e} | {1,a,b,c,d,e} | {1,a,b,c,d,e} |
| b | b | {1,a,b,c,d,e} | {1,a,b,c,d,e} | {0,1,a,c,d,e} | {1,a,b,c,d,e} | {1,a,b,c,d,e} | {1,a,b,c,d,e} |
| c | c | {1,a,b,c,d,e} | {1,a,b,c,d,e} | {1,a,b,c,d,e} | {0,1,a,b,d,e} | {1,a,b,c,d,e} | {1,a,b,c,d,e} |
| d | d | {1,a,b,c,d,e} | {1,a,b,c,d,e} | {1,a,b,c,d,e} | {1,a,b,c,d,e} | {0,1,a,b,c,e} | {1,a,b,c,d,e} |
| e | e | {1,a,b,c,d,e} | {1,a,b,c,d,e} | {1,a,b,c,d,e} | {1,a,b,c,d,e} | {1,a,b,c,d,e} | {0,1,a,b,c,d} |

**B6. Hyperfields for which card(x-x) = 7, for every non-zero element x.**

| $HF_7^{268}$ | 0 | 1 | a | b | c | d | e |
|---|---|---|---|---|---|---|---|
| 0 | 0 | 1 | a | b | c | d | e |
| 1 | 1 | {0,1,a,b,c,d,e} | {1,a} | {1,b} | {1,c} | {1,d} | {1,e} |
| a | a | {1,a} | {0,1,a,b,c,d,e} | {a,b} | {a,c} | {a,d} | {a,e} |
| b | b | {1,b} | {a,b} | {0,1,a,b,c,d,e} | {b,c} | {b,d} | {b,e} |
| c | c | {1,c} | {a,c} | {b,c} | {0,1,a,b,c,d,e} | {c,d} | {c,e} |
| d | d | {1,d} | {a,d} | {b,d} | {c,d} | {0,1,a,b,c,d,e} | {d,e} |
| e | e | {1,e} | {a,e} | {b,e} | {c,e} | {d,e} | {0,1,a,b,c,d,e} |

| $HF_7^{269}$ | 0 | 1 | a | b | c | d | e |
|---|---|---|---|---|---|---|---|
| 0 | 0 | 1 | a | b | c | d | e |
| 1 | 1 | {0,1,a,b,c,d,e} | {1,a} | {1,b,d} | {1,c} | {1,b,d} | {1,e} |
| a | a | {1,a} | {0,1,a,b,c,d,e} | {a,b} | {a,c,e} | {a,d} | {a,c,e} |
| b | b | {1,b,d} | {a,b} | {0,1,a,b,c,d,e} | {b,c} | {1,b,d} | {b,e} |
| c | c | {1,c} | {a,c,e} | {b,c} | {0,1,a,b,c,d,e} | {c,d} | {a,c,e} |
| d | d | {1,b,d} | {a,d} | {1,b,d} | {c,d} | {0,1,a,b,c,d,e} | {d,e} |
| e | e | {1,e} | {a,c,e} | {b,e} | {a,c,e} | {d,e} | {0,1,a,b,c,d,e} |

| $HF_7^{270}$ | 0 | 1 | a | b | c | d | e |
|---|---|---|---|---|---|---|---|
| 0 | 0 | 1 | a | b | c | d | e |
| 1 | 1 | {0,1,a,b,c,d,e} | {1,a,c} | {1,b,e} | {1,a,c,d} | {1,c,d} | {1,b,e} |
| a | a | {1,a,c} | {0,1,a,b,c,d,e} | {a,b,d} | {1,a,c} | {a,b,d,e} | {a,d,e} |
| b | b | {1,b,e} | {a,b,d} | {0,1,a,b,c,d,e} | {b,c,e} | {a,b,d} | {1,b,c,e} |
| c | c | {1,a,c,d} | {1,a,c} | {b,c,e} | {0,1,a,b,c,d,e} | {1,c,d} | {b,c,e} |
| d | d | {1,c,d} | {a,b,d,e} | {a,b,d} | {1,c,d} | {0,1,a,b,c,d,e} | {a,d,e} |
| e | e | {1,b,e} | {a,d,e} | {1,b,c,e} | {b,c,e} | {a,d,e} | {0,1,a,b,c,d,e} |

| $HF_7^{271}$ | 0 | 1 | a | b | c | d | e |
|---|---|---|---|---|---|---|---|
| 0 | 0 | 1 | a | b | c | d | e |
| 1 | 1 | {0,1,a,b,c,d,e} | {1,a,c} | {1,b,d,e} | {1,a,c,d} | {1,b,c,d} | {1,b,e} |
| a | a | {1,a,c} | {0,1,a,b,c,d,e} | {a,b,d} | {1,a,c,e} | {a,b,d,e} | {a,c,d,e} |
| b | b | {1,b,d,e} | {a,b,d} | {0,1,a,b,c,d,e} | {b,c,e} | {1,a,b,d} | {1,b,c,e} |
| c | c | {1,a,c,d} | {1,a,c,e} | {b,c,e} | {0,1,a,b,c,d,e} | {1,c,d} | {a,b,c,e} |
| d | d | {1,b,c,d} | {a,b,d,e} | {1,a,b,d} | {1,c,d} | {0,1,a,b,c,d,e} | {a,d,e} |
| e | e | {1,b,e} | {a,c,d,e} | {1,b,c,e} | {a,b,c,e} | {a,d,e} | {0,1,a,b,c,d,e} |

| $HF_7^{272}$ | 0 | 1 | a | b | c | d | e |
|---|---|---|---|---|---|---|---|
| 0 | 0 | 1 | a | b | c | d | e |
| 1 | 1 | {0,1,a,b,c,d,e} | {1,a,c,d} | {1,b,c,e} | {1,a,b,c,d,e} | {1,a,c,d} | {1,b,c,e} |
| a | a | {1,a,c,d} | {0,1,a,b,c,d,e} | {a,b,d,e} | {1,a,c,d} | {1,a,b,c,d,e} | {a,b,d,e} |
| b | b | {1,b,c,e} | {a,b,d,e} | {0,1,a,b,c,d,e} | {1,b,c,e} | {a,b,d,e} | {1,a,b,c,d,e} |
| c | c | {1,a,b,c,d,e} | {1,a,c,d} | {1,b,c,e} | {0,1,a,b,c,d,e} | {1,a,c,d} | {1,b,c,e} |
| d | d | {1,a,c,d} | {1,a,b,c,d,e} | {a,b,d,e} | {1,a,c,d} | {0,1,a,b,c,d,e} | {a,b,d,e} |
| e | e | {1,b,c,e} | {a,b,d,e} | {1,a,b,c,d,e} | {1,b,c,e} | {a,b,d,e} | {0,1,a,b,c,d,e} |

| $HF_7^{273}$ | 0 | 1 | a | b | c | d | e |
|---|---|---|---|---|---|---|---|
| 0 | 0 | 1 | a | b | c | d | e |
| 1 | 1 | {0,1,a,b,c,d,e} | {1,a,c,d} | {1,b,c,d,e} | {1,a,b,c,d,e} | {1,a,b,c,d} | {1,b,c,e} |
| a | a | {1,a,c,d} | {0,1,a,b,c,d,e} | {a,b,d,e} | {1,a,c,d,e} | {1,a,b,c,d,e} | {a,b,c,d,e} |
| b | b | {1,b,c,d,e} | {a,b,d,e} | {0,1,a,b,c,d,e} | {1,b,c,e} | {1,a,b,d,e} | {1,a,b,c,d,e} |
| c | c | {1,a,b,c,d,e} | {1,a,c,d,e} | {1,b,c,e} | {0,1,a,b,c,d,e} | {1,a,c,d} | {1,a,b,c,e} |
| d | d | {1,a,b,c,d} | {1,a,b,c,d,e} | {1,a,b,d,e} | {1,a,c,d} | {0,1,a,b,c,d,e} | {a,b,d,e} |
| e | e | {1,b,c,e} | {a,b,c,d,e} | {1,a,b,c,d,e} | {1,a,b,c,e} | {a,b,d,e} | {0,1,a,b,c,d,e} |

| $HF_7^{274}$ | 0 | 1 | a | b | c | d | e |
|---|---|---|---|---|---|---|---|
| 0 | 0 | 1 | a | b | c | d | e |
| 1 | 1 | {0,1,a,b,c,d,e} | {1,a,b,c,e} | {1,a,b,e} | {1,a,c,d} | {1,c,d,e} | {1,a,b,d,e} |
| a | a | {1,a,b,c,e} | {0,1,a,b,c,d,e} | {1,a,b,c,d} | {1,a,b,c} | {a,b,d,e} | {1,a,d,e} |
| b | b | {1,a,b,e} | {1,a,b,c,d} | {0,1,a,b,c,d,e} | {a,b,c,d,e} | {a,b,c,d} | {1,b,c,e} |
| c | c | {1,a,c,d} | {1,a,b,c} | {a,b,c,d,e} | {0,1,a,b,c,d,e} | {1,b,c,d,e} | {b,c,d,e} |
| d | d | {1,c,d,e} | {a,b,d,e} | {a,b,c,d} | {1,b,c,d,e} | {0,1,a,b,c,d,e} | {1,a,c,d,e} |
| e | e | {1,a,b,d,e} | {1,a,d,e} | {1,b,c,e} | {b,c,d,e} | {1,a,c,d,e} | {0,1,a,b,c,d,e} |

| $HF_7^{275}$ | 0 | 1 | a | b | c | d | e |
|---|---|---|---|---|---|---|---|
| 0 | 0 | 1 | a | b | c | d | e |
| 1 | 1 | {0,1,a,b,c,d,e} | {1,a,b,c,e} | {1,a,b,d,e} | {1,a,c,d} | {1,b,c,d,e} | {1,a,b,d,e} |
| a | a | {1,a,b,c,e} | {0,1,a,b,c,d,e} | {1,a,b,c,d} | {1,a,b,c,e} | {a,b,d,e} | {1,a,c,d,e} |
| b | b | {1,a,b,d,e} | {1,a,b,c,d} | {0,1,a,b,c,d,e} | {a,b,c,d,e} | {1,a,b,c,d} | {1,b,c,e} |
| c | c | {1,a,c,d} | {1,a,b,c,e} | {a,b,c,d,e} | {0,1,a,b,c,d,e} | {1,b,c,d,e} | {a,b,c,d,e} |
| d | d | {1,b,c,d,e} | {a,b,d,e} | {1,a,b,c,d} | {1,b,c,d,e} | {0,1,a,b,c,d,e} | {1,a,c,d,e} |
| e | e | {1,a,b,d,e} | {1,a,c,d,e} | {1,b,c,e} | {a,b,c,d,e} | {1,a,c,d,e} | {0,1,a,b,c,d,e} |

| $HF_7^{276}$ | 0 | 1 | a | b | c | d | e |
|---|---|---|---|---|---|---|---|
| 0 | 0 | 1 | a | b | c | d | e |
| 1 | 1 | {0,1,a,b,c,d,e} | {1,a,b,c,d,e} | {1,a,b,c,e} | {1,a,b,c,d,e} | {1,a,c,d,e} | {1,a,b,c,d,e} |
| a | a | {1,a,b,c,d,e} | {0,1,a,b,c,d,e} | {1,a,b,c,d,e} | {1,a,b,c,d} | {1,a,b,c,d,e} | {1,a,b,d,e} |
| b | b | {1,a,b,c,e} | {1,a,b,c,d,e} | {0,1,a,b,c,d,e} | {1,a,b,c,d,e} | {a,b,c,d,e} | {1,a,b,c,d,e} |
| c | c | {1,a,b,c,d,e} | {1,a,b,c,d} | {1,a,b,c,d,e} | {0,1,a,b,c,d,e} | {1,a,b,c,d,e} | {1,b,c,d,e} |
| d | d | {1,a,c,d,e} | {1,a,b,c,d,e} | {a,b,c,d,e} | {1,a,b,c,d,e} | {0,1,a,b,c,d,e} | {1,a,b,c,d,e} |
| e | e | {1,a,b,c,d,e} | {1,a,b,d,e} | {1,a,b,c,d,e} | {1,b,c,d,e} | {1,a,b,c,d,e} | {0,1,a,b,c,d,e} |

| $HF_7^{277}$ | 0 | 1 | a | b | c | d | e |
|---|---|---|---|---|---|---|---|
| 0 | 0 | 1 | a | b | c | d | e |
| 1 | 1 | {0,1,a,b,c,d,e} | {1,a,b,c,d,e} | {1,a,b,c,d,e} | {1,a,b,c,d,e} | {1,a,b,c,d,e} | {1,a,b,c,d,e} |
| a | a | {1,a,b,c,d,e} | {0,1,a,b,c,d,e} | {1,a,b,c,d,e} | {1,a,b,c,d,e} | {1,a,b,c,d,e} | {1,a,b,c,d,e} |
| b | b | {1,a,b,c,d,e} | {1,a,b,c,d,e} | {0,1,a,b,c,d,e} | {1,a,b,c,d,e} | {1,a,b,c,d,e} | {1,a,b,c,d,e} |
| c | c | {1,a,b,c,d,e} | {1,a,b,c,d,e} | {1,a,b,c,d,e} | {0,1,a,b,c,d,e} | {1,a,b,c,d,e} | {1,a,b,c,d,e} |
| d | d | {1,a,b,c,d,e} | {1,a,b,c,d,e} | {1,a,b,c,d,e} | {1,a,b,c,d,e} | {0,1,a,b,c,d,e} | {1,a,b,c,d,e} |
| e | e | {1,a,b,c,d,e} | {1,a,b,c,d,e} | {1,a,b,c,d,e} | {1,a,b,c,d,e} | {1,a,b,c,d,e} | {0,1,a,b,c,d,e} |

CORE DEPARTMENT, EURIPUS CAMPUS, NATIONAL AND KAPODISTRIAN UNIVERSITY OF ATHENS,

*Email address:* chrmas@uoa.gr or ch.massouros@gmail.com

SCHOOL OF SOCIAL SCIENCES, HELLENIC OPEN UNIVERSITY,

*Email address:* germasouros@gmail.com